\documentclass[a4paper, 11pt]{article} 

\usepackage{times}
\usepackage{graphicx}
\usepackage{amsfonts,amssymb,amsmath,stmaryrd}
\usepackage{mathrsfs} 
\usepackage{upref}
\usepackage[sf,SF]{subfigure} 
\usepackage{caption}
\usepackage{amsthm,euscript}
\usepackage[a4paper, portrait, margin=1in]{geometry}
\usepackage{xcolor}
\usepackage[bookmarks=true,hidelinks]{hyperref}
\usepackage{units}
\usepackage{bbm}
\usepackage{enumitem}
\usepackage{caption}
\usepackage[title]{appendix}

\usepackage{dsfont} 

\usepackage{mathtools}
\mathtoolsset{showonlyrefs}

\newcommand{\R}{\mathbb{R}}

\newcommand{\N}{\mathbb{N}}
\newcommand{\Z}{\mathbb{Z}}
\renewcommand{\leq}{\leqslant}
\renewcommand{\geq}{\geqslant}

\newcommand{\hf}{{\unitfrac{1}{2}}}

\renewcommand{\epsilon}{\varepsilon}

\newcommand{\bk}{\mathbf{k}}

\newcommand{\diff}{\mathrm{d}}
\newcommand{\rhostat}{\rho_\infty}

\newcommand{\bb}{\frac{\beta^2}{2}}
\newcommand{\HH}{\displaystyle H\left( \omega \right)}

\newcommand{\norme}[1]{\left\lVert#1\right\rVert}
\newcommand{\pscal}[2]{\left\langle ~#1~\middle\vert~#2~\right\rangle}

\definecolor{navy}{rgb}{0,0,0.5}   
\definecolor{firebrick}{rgb}{0.698,0.133,0.133}  

\DeclareMathOperator{\supp}{supp}
\DeclareMathOperator\erf{erf}

\newtheorem{theorem}{Theorem}

\newtheorem{remark}[theorem]{Remark}
\newtheorem{definition}[theorem]{Definition}
\newtheorem{lemma}[theorem]{Lemma}
\newtheorem{corollary}[theorem]{Corollary}

\newtheorem{hyp}[theorem]{Assumption}

\numberwithin{equation}{section}
\numberwithin{theorem}{section}

\makeatletter
\def\@cite#1#2{\textup{[{#1\if@tempswa , #2\fi}]}}
\makeatother

\allowdisplaybreaks

\let\OLDthebibliography\thebibliography
\renewcommand\thebibliography[1]{
  \OLDthebibliography{#1}
  \setlength{\parskip}{1.5pt}
  \setlength{\itemsep}{0.5pt}
}

\title{Nonlinear partial differential equations in neuroscience: from modelling to mathematical theory}

\author{José A. Carrillo\thanks{Mathematical Institute, University of Oxford, Oxford OX2 6GG, UK (carrillo@maths.ox.ac.uk)} \and Pierre Roux\thanks{Institut Camille Jordan, \'Ecole Centrale de Lyon, Ecully 69134, France (pierre.roux@ec-lyon.fr)}}

\begin{document}

\maketitle

\begin{abstract}
    Many systems of partial differential equations have been proposed as simplified representations of complex collective behaviours in large networks of neurons. In this survey, we briefly discuss their derivations and then review the mathematical methods developed to handle the unique features of these models, which are often nonlinear and non-local. The first part focuses on parabolic Fokker-Planck equations: the Nonlinear Noisy Leaky Integrate and Fire neuron model, stochastic neural fields in PDE form with applications to grid cells, and rate-based models for decision-making. The second part concerns hyperbolic transport equations, namely the model of the Time Elapsed since the last discharge and the jump-based Leaky Integrate and Fire model. The last part covers some kinetic mesoscopic models, with particular attention to the kinetic Voltage-Conductance model and FitzHugh-Nagumo kinetic Fokker-Planck systems.
\end{abstract}

\vspace{12cm}

\noindent\textit{Keywords:} Neuroscience, Partial Differential Equations, Neural Networks, Mathematical Biology.

\noindent\textit{Mathematics Subject Classification:}  Primary: 35-02, 
92-02, 
35Q80; 
Secondary: 
35B40,  
35B10,  
65C30.  

\newpage

\setcounter{tocdepth}{2}

\tableofcontents

\section*{Introduction}

Since the first mathematical models for individual neurons \cite{L,HH,F,NAY}, a humongous research effort has been devoted to understanding complex phenomena in the brain and the modelling of large neural assemblies. Contrary to several older fields like theoretical physics, which can for example rely upon quantum field theory and general relativity, there is no abstract paradigm or first physics principles allowing to construct knowledge in neuroscience. Furthermore, animal brains exhibit an astonishing level of complexity and, despite recent advances, precise \textit{in vivo} measurements of entire neural networks are still challenging.

For these reasons, researchers seek to develop simple mathematical models that represent the qualitative properties of neural tissues at the mesoscopic and macroscopic scales. They need to be mathematically and numerically tractable while still encompassing as many key features as possible. A popular approach is to start with a simple model for how a state variable $x$ evolves for individual neurons, for example with stochastic differential equations (SDE), to choose interaction rules and to derive a mean-field equation for a network containing an infinite number of neurons \cite{BH,B3,CTRM,OKS,OKKS,OBH,BFFT,FTC,FT}. In many cases, the probability density $\rho(x,t)$ of finding in the network, at time $t$, a neuron with state variable $x$, solves a nonlinear partial differential equation (PDE).

These nonlinear equations often present unique structures and non-local terms, making it impossible to straightforwardly apply the standard theories of PDEs to them. In the past two decades, a lot of progress has been made in the mathematical understanding of nonlinear partial differential equations modelling neuronal assemblies. Specific tools have been developed and many methods from mathematical physics have been adapted to their specific difficulties.

Although there are already excellent survey articles \cite{B,C,CPWT,VL} and books \cite{GK,C3} about mathematical modelling in neuroscience, they are either biology-oriented or focused on older types of equations like neural fields, Wilson-Cowan models \cite{WC1,WC2} and dynamical systems. We propose here a review article focused on the mathematical theory for neuroscience models in the form of nonlinear PDEs, often rigorously obtained from particle systems, with an emphasis on rigorous proofs of the properties of the solutions. We have no intention, nor do we pretend, to be exhaustive. Our attempt is rather to accurately describe the process of constructing methods for a few key types of models and to show how they are related to each other through asymptotic reductions and approximations, while keeping in mind the application purpose: getting back from proving theorems to challenging our understanding of the brain.

The structure of this survey is organized around three large families of PDEs. In the first part, we discuss systems of parabolic type, with Fokker-Planck differential operators. This entails smoothness of solutions and allows to use tools such as Green kernels and entropy methods. We discuss the Nonlinear Noisy Leaky Integrate and Fire (NNLIF) model for neural networks, a Fokker-Planck system representing grid cells, and rate-based equations used to model decision making. In the second part, we focus on first-order hyperbolic type PDEs, whose solutions are not always smooth nor unique. It allows us to introduce other key methods like characteristics, the measure solutions framework, and Doeblin-Harris theorems. We mainly discuss the Time Elapsed neuron model and a jump-based Leaky Integrate and Fire model. In the last part, we focus on kinetic models with degenerate diffusion, for which the regularity theory is challenging and the long-time behaviour still poorly understood. We propose an overview of Voltage-Conductance and FitzHug-Nagumo models, both in the form of kinetic Fokker-Planck equations.

\paragraph{Notation:} Most conventions will be explained when required, but we want to fix for clarity some choices that are used throughout this document. We denote $\R_+ = [0,+\infty)$ (note that $0\in\R_+$). We write $C([a,b])$ or $C[a,b]$ for the space of continuous functions from the interval $[a,b]$ to $\R$; when the norm is not ambiguous, $C(X,Y)$ is the space of continuous functions from space $X$ to space $Y$, or just $C(X)$ if $Y=\R$. Similarly, $C^k$ will be used to denote spaces with higher regularity than just continuous. The Lebesgue spaces on a set $\Omega$ are denoted $L^p(\Omega)$ or just $L^p$, $p\in[1,+\infty]$, with the classical norms $\norme{\cdot}_{L^p}$. Those norms are also written $\norme{\cdot}_p$ when there is no ambiguity. Given a weight function $\omega:\Omega\to\R_+$, we denote $L^p(\Omega,\omega)$ or $L^p(\omega)$ the space of measurable functions $f$ with finite weighted $L^p$ norm $\norme{f\omega^{1/p}}_{L^p}$. The Sobolev spaces on $\Omega$ are $W^{k,p}(\Omega)$; for $p=2$, it is the standard Hilbert space $H^k$. Similar to Lebesgue spaces, we consider weighted versions $W^{k,p}(\Omega,\omega)$, $H^k(\Omega,\omega)$. When the notation $W$ is already taken by a connectivity kernel, we use $\mathcal W^{k,p}(\Omega)$ for Sobolev spaces. Given an interval $I\subset \R$ and a Banach space $X$ we use the Bochner spaces $L^p(I,X)$, also denoted by $L^p(a,b\,;X)$ for $I=(a,b)$. The space of probability measures on a domain $\Omega\subset \R^d$ will be denoted by $\mathcal P(\Omega)$ typically in dimension $d\in\{1,2,3\}$. We will write $f\equiv a$ to mean "the function $f$ is constant and equal to $a\in\R$".

\part{Parabolic equations of Fokker-Planck type}

\section{The Nonlinear Noisy Leaky Integrate \& Fire model for neural networks}
\label{sec:NNLIF}

We focus here on the so-called Nonlinear Noisy Leaky Integrate \& Fire model for neural networks, or in short form NNLIF\footnote{In some references, the model has no name, or the first N stands for "Network". The probabilistic community refers to it as a nonlinear McKean-Vlasov equation. The name NNLIF was first proposed in \cite{CCP}.}. It has been introduced about twenty-five years ago, in the seminal works of Brunel and Hakim \cite{BH,B}. Their main motivation was to explain a counter-intuitive phenomenon: self-sustained oscillations in networks of neurons with strong inhibitory connections and low individual firing rate.

Self-sustained oscillations are ubiquitous in neuroscience and play pivotal roles in many core functions like breathing and memory. They are found in many brain areas, like the visual cortex and the olfactory cortex. Previous studies had been focusing on networks where individual neurons are themselves oscillators, leading to total or partial synchronisation in a few clusters under the effect of excitatory connections. However, the role of inhibition in the emergence of periodic activity was still poorly understood. Even more striking was the fact that in networks of neurons that sporadically fire, oscillations at the population scale can be very fast: gamma frequency range ($\sim$ 30Hz) or even faster.

We start by presenting the historical derivation of the NNLIF model through physics arguments. We then set the framework for mathematical analysis of this type of PDE systems and, in order: discuss the existence and multiplicity of stationary states; explain how to prove existence and regularity of solutions through changes of variables and representation formulae; discuss the occurrence and meaning of finite-time blow-up in the non-delayed case; present the way entropy methods can be used to study long-time asymptotics in the weakly connected case; explain the specific methods and heuristic arguments allowing us to understand the NNLIF model in the delayed case. After all these considerations on the standard theory around this type of equations, we propose a brief overview of the rigourous derivation through a mean-field limit, and we present recent methods developed on the one hand to extend solutions after blow-up \textit{via} time dilations and on the other hand to study the asymptotic behaviour beyond the weakly-connected case \textit{via} semigroup methods and a sequence of pseudo-equilibria. Last, we present the numerous variants of the standard model and discuss some numerical methods tailored to this specific problem.

\subsection{Derivation of the model \textit{à la physicienne}}

Let us propose first a physical derivation of the model as can be found in early references \cite{BH,CCP}. The rigorous mean-field limit from a particle system was established later in 2015, in \cite{DIRT2}. The main idea is to consider a large network where a representative neuron of electric potential $V(t)$ is modelled by the Lapicque Integrate \& Fire model \cite{L} given by
\begin{equation}\label{voltage1}
	C_m\dfrac{\diff V}{\diff t} = - g_L(V-V_L) + I(t),
\end{equation}
where $C_m$ is the capacitance of the membrane, $g_L$ the leak conductance and $V_L$ the leak potential. The characteristic time of relaxation $\tau_m = \frac{g_L}{C_m}$ is of the order of $2ms$ and $V_L \simeq -70 mV$. The incoming synaptic current is a stochastic process of the form
\begin{equation}\label{current1}
	I(t) = J_E\sum_{i=1}^{C_E}\sum_{j\in\N}\delta(t-t_{E_j}^i-d) + J_I\sum_{i=1}^{C_I}\sum_{j\in\N}\delta(t-t_{I_j}^i-d),
\end{equation}
where $\delta$ is the Dirac mass measure at zero, $J_E$ and $J_I$ are the strength of excitatory and inhibitory synapses, $C_E$ and $C_I$ are the number of excitatory and inhibitory connections received on average, and $t_{E_j}^i$, $t_{I_j}^i$ are the times of the $j^{th}$ action potential emitted by the $i^{th}$ pre-synaptic excitatory or inhibitory neuron; those times are subject to random fluctuations. The parameter $d$ represents the average time it takes for an action potential to go from a neuron to an other and be integrated; as we will see, this parameter can have a crucial effect on the solutions and the model is more singular when $d=0$.

When a neuron reaches the firing potential $V_F\simeq -50 mV$, it emits an action potential (it is also said to \textit{fire}), and then comes back to the reset potential $V_R\simeq -60 mV$. Assuming that each neuron in the network fires independently following a Poisson law of intensity $\nu(t)$, the mean and variance of $I(t)$ are $\mu_C(t-d) = b\nu(t-d)$, where
\[ b = C_EJ_E - C_IJ_I,  \]
is the average connectivity of the network, and
$  {\sigma_C}^2(t-d) = (C_E{J_E}^2 + C_I{J_I}^2)\nu(t-d). $
When $b<0$ we say that the network is average inhibitory, when $b>0$ we describe it as average excitatory. When $b=0$, the network is disconnected and only mathematicians care about it, although it is useful as linearized equation for small connectivity.

The stochastic process \eqref{current1} being difficult to study, a classical approach \cite{B3,BH,RBW,mg,sirovich,OKS} is to replace it by an Ornstein-Uhlenbeck process with same mean and variance
\[  I(t)\diff t \simeq \mu_c(t-d)\, \diff t + \sigma_C(t-d)\, \diff B_t,  \]
where $B_t$ is a standard Brownian motion.

If we rescale the voltage and time variables so to have $C_m=g_L=1$, we can write the following stochastic differential equation for a neuron of this network:
\begin{equation}\label{voltage2}
	\diff V = (-V + V_L + \mu_C(t-d))\,\diff t + \sigma_C(t-d) \,\diff B_t, \qquad V<V_F
\end{equation}
with the jump process
$  \limsup_{t\to t_0^-} V(t) = V_F \implies \liminf_{t\to t_0^+} V(t) = V_R, $ 
and the discharge intensity
$ \nu(t) = \nu_{ext} + N(t), $
where $N(t)$ is the flux of neurons in the network that cross the firing potential $V_F$ (and emit action potentials) at time $t$. The term $\nu_{\mathrm{ext}}\in\R$ represents an external input from outside the network. We consider it constant here, but it may vary in time in more realistic settings. Although it plays an important role in the modelling, as an interface between a group of neurons, other networks and sensory inputs, we will often rescale the voltage $V$ in order to absorb it, which simplifies the mathematical analysis.

The stochastic differential equation \eqref{voltage2} is associated to a Fokker-Planck partial differential equation, describing the probability density $(v,t)\mapsto \rho(v,t) = \mathbb P(V(t)\in dv)$ of finding in the network, at time $t\geqslant 0$, a neuron with membrane potential $v\in\,(-\infty,V_F]$, given by 
\begin{equation}\label{NNLIFphy}
	\dfrac{\partial \rho}{\partial t}(v,t) + \dfrac{\partial}{\partial v}\Big[(V_L  -v + b_0 + bN(t-d))\rho(v,t)\Big]-a\big(N(t-d)\big)\dfrac{\partial^2 \rho}{\partial v^2}(v,t)= \ \delta_{v=V_R}  N(t),
\end{equation}
with $b_0=b\nu_{\mathrm{ext}}$ and $a\big(N(t-d)\big) = \dfrac{\sigma_C^2(t)}{2}= a_0 + a_1 N(t-d)$, $\rho(V_F,t) = \rho(-\infty,t)=0$ and
\begin{equation}\label{NNLIFphy2} N(t) = - a\big(N(t-d)\big) \dfrac{\partial \rho}{\partial v}(V_F,t). \end{equation}
Since $\rho$ must remain a probability density , integrating \eqref{NNLIFphy} in $v$ yields the relation \eqref{NNLIFphy2}. The singular resetting of neurons is embodied by the Dirac mass at potential $V_R$.

\subsection{Mathematical setting for the NNLIF type models}\label{sec:NNLIF_math}

Up to appropriate variables re-scaling, the Fokker-Planck equation \eqref{NNLIFphy}--\eqref{NNLIFphy2} can be rewritten as the Cauchy problem
\begin{equation}\label{NNLIF}\tag{NNLIF}
	\left\{\begin{array}{ll}
		\displaystyle \dfrac{\partial \rho}{\partial t} + \dfrac{\partial}{\partial v}\Big[(-v + bN(t-d))\rho\Big]-a\big(N(t-d)\big)\dfrac{\partial^2 \rho}{\partial v^2}= \ \delta_{v=V_R}  N(t), & v\in(-\infty,V_F],\\ &t > 0,\\
		\displaystyle N(t)= -a\big(N(t-d)\big) \dfrac{\partial \rho}{\partial v}(V_F,t),  \quad \rho(V_F,t)=0, \quad \rho(-\infty,t)=0,& t > 0,\\
		\displaystyle \rho(v,0)=\rho^0(v) \geqslant 0, \quad \int_{-\infty}^{V_F} \rho^0(v)\diff v = 1,&v\in(-\infty,V_F],\\
		N(t)=N^0(t),&t\in[-d,0].
	\end{array}\right.
\end{equation}
We call this model the standard NNLIF model, with ($d > 0$) or without ($d=0$) synaptic delay. When $d=0$, the initial condition $N^0$ for the firing rate is not required.

Note that the external input $\nu_{\mathrm{ext}}$ has been absorbed in the values $V_F$ and $V_R$ through an affine voltage rescaling in order to reduce the number of parameters and simplify the mathematical description. Every time a result is stated in terms of $V_R,V_F,b$ and the function $a$, the role of external inputs $\nu_{ext}$ is hidden in the values of these parameters. The impact of $\nu_{\mathrm{ext}}$ is especially visible when taking negative values for $V_F$, which is equivalent to a highly excitatory external input. Note also that we allow for more general choices for the function $a$ than the form $a(N)=a_0+a_1N$ obtained in the previous section. The rigorous derivation from a particle system yields $a_1 = 0$, in contrast to the physical one. \\

To better understand the role of the Dirac mass $\delta_{v=V_R}$ in the right-hand side, there are two viewpoints: 
\begin{itemize}
	\item \eqref{NNLIF} is satisfied in the classical sense on $(-\infty,V_R)\,\cup\,(V_R,V_F]$ and in the sense of distributions on $(-\infty,V_F]$;
	\item \eqref{NNLIF} is satisfied in the classical sense on $(-\infty,V_R)\,\cup\,(V_R,V_F]$ and the function $\partial_v \rho(\cdot,t)$ has left and right limits at $V_R$ with a jump of size $\frac{N(t)}{a(N(t-d))}$, \textit{i.e.}
	\[  \lim_{v \to {V_R}^-} \left(\dfrac{\partial \rho}{\partial v}(v,t) \right) - \lim_{v \to {V_R}^+} \left(\dfrac{\partial \rho}{\partial v}(v,t) \right)  = \dfrac{N(t)}{a\big(N(t-d)\big)}.   \]
\end{itemize}

In order to avoid technical difficulties, we will consider fast-decaying (at $-\infty$) initial conditions and solutions, since in neuroscience applications voltage distributions are not fat-tailed.

\begin{definition}\label{def:classical}
	We say that $(\rho,N)$ is a \emph{strong} fast-decaying \emph{solution}, of system \eqref{NNLIF} on $[0,T^*)\,$, $T^*\in(0,+\infty]$ if
	\begin{itemize}
		\item[$\bullet$] $\rho\in C^0(\,(-\infty,V_F]\times [0,T^*)\,)\cap C^{2,1}\Big( \big(\,(-\infty,V_R)\,\cup\,(V_R,V_F]\big)\times (0,T^*) \Big)$ and $N\in C^0([-d,T^*))$;
		\item[$\bullet$] Functions $\rho$ and $N$ are non-negative, satisfy \eqref{NNLIF} in the classical sense on $(-\infty,V_R)\,\cup\,(V_R,V_F]$ and in the sense of distributions in $(-\infty,V_F]$;
		\item[$\bullet$] $\forall t\in[0,T^*)\,,\forall Q\in\R[X]$, $\lim_{v\to-\infty}Q(v)\rho(v,t)=0$ and $\lim_{v\to-\infty}Q(v)\dfrac{\partial \rho}{\partial v}(v,t)=0$.
	\end{itemize}
\end{definition}

We can also define weak solutions.
\begin{definition}
	Let $T^*\in(0,+\infty]$, let $\rho\in L^\infty([0,T^*)\,, L^1_+(\,(-\infty,V_F)\,))$ and $N\in L^1_{\mathrm{loc},+}([0,T^*)\,)$. The pair $(p,N)$ is said to be a fast-decaying \emph{weak solution} of \eqref{NNLIF} if for every test function $\phi\in C^2(\,(-\infty,V_F])$ with polynomial growth\footnote{There exists $Q\in\R[X]$ such that $\forall v\in (-\infty,V_F], \ \phi(v)\leqslant Q(v)$.}
	we have
	\begin{multline*}
		\dfrac{d}{d t} \int_{-\infty}^{V_F} \rho(v,t)\phi(v)d v = \int_{-\infty}^{V_F}\Big[ \big(-v+bN(t-d)\big)\dfrac{\partial \phi}{\partial v} + a\big(N(t-d)\big)\dfrac{\partial^2 \phi}{\partial v^2} \Big]\rho(v,t)d v \\ + N(t) \Big( \phi(V_R) \ -\ \phi(V_F) \Big),
	\end{multline*}
	where all the quantities are finite, the time derivative is taken in the sense of distributions, $N(t-d)$ is given by $N^0$ when $t\leqslant d$ and $\lim_{t\to 0} \norme{ \rho(\cdot,t) - \rho^0(\cdot) }_{ L^1} = 0 $.\\
\end{definition}
All strong solutions are weak solutions, which can be checked by integration by parts.

\begin{remark}
	If we apply the definition of weak solution with the test function $\phi=1$, we obtain for $t > 0$,
	\[ \int_{-\infty}^{V_F}\rho(v,t)d v=\int_{-\infty}^{V_F}\rho^0(v)d v=1. \]
\end{remark}

We will refer to the following regularity assumptions on initial data.
\begin{definition}[Initial datum for classical solutions]\label{Assumption1} We say that $\rho^0$ is a compatible initial condition for \eqref{NNLIF} if
	$\rho^0\in C^0(\,(-\infty,V_F])\cap C^1(\,(-\infty,V_R)\,\cup\,(V_R,V_F])\cap \mathcal P((-\infty,V_F])$ is fast-decaying at $-\infty$, $\rho^0(V_F)=0$, $\frac{ d  \rho^0}{ d  v} $ admits finite left and right limits at $V_R$ and is fast decreasing at $-\infty$. In the case with delay $d > 0$, we require that $N^0\in C([-d,0])$ and $N^0(0) = -a(N(-d))\frac{ d  \rho^0}{ d  v} (V_F) $.
\end{definition}

\subsection{Characterisation of the stationary states}

We call a stationary state, or steady state, of \eqref{NNLIF} a solution $(\rho_\infty,N_\infty)$ of the system
\begin{equation}\label{thm:NNLIF_stat}
	\left\{\begin{array}{ll}
		\displaystyle \dfrac{\partial}{\partial v}\Big[(-v + bN_\infty)\rho_\infty\Big]-a(N_\infty)\dfrac{\partial^2 \rho_\infty}{\partial v^2}= \ \delta_{v=V_R}  N_\infty,\qquad&v\in(-\infty,V_F],\\ 
		\displaystyle N_\infty= -a(N_\infty) \dfrac{\partial \rho_\infty}{\partial v}(V_F),  \qquad&\\
		\displaystyle \rho_\infty(V_F)=\rho_\infty(-\infty)=0, \quad \rho_\infty \geqslant 0, \quad \int_{-\infty}^{V_F} \rho_\infty(v)\diff v = 1.&
	\end{array}\right.
\end{equation}
They are of course the same with and without delay. Direct integration and rescalings allow to characterise the stationary states in a semi-explicit way.
\begin{theorem}[C\'aceres, Carrillo, Perthame \cite{CCP}]\label{StaticValue}
	Let $(\rho_\infty,N_\infty)$ be a stationary state of \eqref{NNLIF}. It must be of the form
	\begin{equation}\label{eq:NNLIF_rho_stat} \rho_\infty (v) = \dfrac{N_\infty}{a(N_\infty)} e^{-\frac{(v-bN_\infty)^2}{2a(N_\infty)}} \int_{\max(v,V_R)}^{V_F} e^{\frac{(w-bN_\infty)^2}{2a(N_\infty)}} \diff w. \end{equation}
	where $N_\infty$ satisfies the fixed-point equation
	\begin{equation}\label{eq:Fix_N_Stat} \dfrac{1}{N_\infty} = I(N_\infty), \qquad I(N) = \int_{0}^{+\infty}  \dfrac{ e^{-\frac{s^2}{2}} }{s} \left( e^{s\frac{V_F-bN}{\sqrt{a(N)}} } - e^{s\frac{V_R-bN}{\sqrt{a(N)}} } \right)  \diff s. \end{equation}
\end{theorem}
Notice that given the stationary firing rate $N_\infty$, the stationary density is explicitly given. The existence and number of classical stationary states thus depends solely on the solutions of the fixed-point equation \eqref{eq:Fix_N_Stat}. A couple pages of real analysis lead to a good picture of what can happen.
\begin{theorem}[C\'aceres, Carrillo, Perthame, \cite{CCP}]\label{thm:ClassifStatic}
	Assume $a(N) \equiv a>0$; then, for \eqref{NNLIF},
	\begin{itemize}
		\item[$\bullet$] there is a constant $c(a)>0$ such that if $b\in (-\infty,c(a)]$, there is a unique stationary state;
		\item[$\bullet$] if  $0 < b < V_F - V_R$  or $0 < 2a b < (V_F-V_R)^2V_R$, there is at least one stationary state;
		\item[$\bullet$] if $b > V_F - V_R$ and $0 < 2a b < (V_F-V_R)^2V_R$, there are at least two stationary states;
		\item[$\bullet$] if $b> \max\big(2(V_F-V_R)\,,\,2V_F I(0)\big)$, there is no stationary state.
	\end{itemize}
\end{theorem}
When the diffusion coefficient linearly depends on $N$, the classification remains similar, but without the uniqueness for inhibitory or small connectivity.
\begin{theorem}[C\'aceres, Carrillo, Perthame \cite{CCP}]
	Assume $a(N)=a_0 + a_1N$, with $a_0,a_1>0$; then,
	\begin{itemize}
		\item[$\bullet$] If $ b < V_F - V_R$   or $\big( 2a_0 b + 2 a_1 V_R < (V_F-V_R)^2V_R$ and $b>0\big)$, 
		there is at least one steady state;
		\item[$\bullet$] If
		$ b > V_F - V_R$ and $2a_0 b + 2 a_1 V_R < (V_F-V_R)^2V_R$,
		there are at least two stationary states;
		\item[$\bullet$] if $b> \max\big(2(V_F-V_R)\,,\,2V_F I(0)\big)$, there is no stationary state.
	\end{itemize}
\end{theorem}

\subsection{Construction of solutions}

\subsubsection{Local-in-time existence \textit{via} a Stefan-like free boundary problem}

Because it was the basis for the construction of solutions in other parameter settings or variants of the NNLIF model, we will present here the main ideas of the first deterministic existence proof, as can be found in \cite{CGGS}. It concerns only the case $d=0$, $a(N)\equiv a > 0$ constant, but it was generalised later to the case with delay $d >0$ in \cite{CRSS}, to the case with internal noise $a(N)=a_0+a_1N$ in \cite{DZ} and to other PDEs from neuroscience, see, for example \eqref{eq:4PDE} below. Proofs \textit{via} stochastic methods were independently achieved in \cite{DIRT1, DIRT3} for classical solutions and in \cite{DIRT2} for generalised solutions. 

Here we will present the main ideas of the deterministic proof for the case $d=0$. The case with delay proceeds similarly. The idea is to transform the original problem on a fixed domain $(-\infty,V_F]$ into a free-boundary Stefan-like problem, then isolate the firing rate in a closed equation and apply the Banach fixed point theorem.\\

First, let us rescale the variables, in order to have $a=1$ and $V_F=0$, which comes at the price of changing $b$ and putting back an external input $\nu_{\mathrm{ext}}$. It is enough to apply the scaling $\rho_{\mathrm{new}}(v,t)=\sqrt a \rho_{\mathrm{old}}(\sqrt a v + V_F,t)$, $N_{\mathrm{new}}(t)= -\partial_v \rho_{\mathrm{new}}(0,t)$. If we define $\nu_{\mathrm{ext}} = \frac{V_F}{\sqrt a }$ and the new parameters $\bar{b} = \frac{b}{\sqrt a}$,  $\bar{V_R} = \frac{V_R - V_F}{\sqrt a}$, then dropping the bars, the "new" pair $(\rho,N)$ satisfies, 
\[\dfrac{\partial \rho}{\partial t} + \dfrac{\partial}{\partial v}\Big[(-v + \nu_{\mathrm{ext}} + bN(t))\rho\Big]-\dfrac{\partial^2 \rho}{\partial v^2}= \ \delta_{v=V_R}  N(t), \qquad v\in(-\infty,0],\ t > 0.\]
Note that this preserves the sign of $b$ and compensates the change $V_F\to 0$ of the firing potential by an equivalent external force $\nu_{\mathrm{ext}}$. This external force is not the same as the one derived in the first Subsection, as $V_F\neq 0$ in real life; it is here purely mathematical.

Secondly, we apply a classical change of variables that transforms the linear Fokker-Planck equation into a heat equation \cite{CT},
\[  y=e^{t}v,      \qquad  \tau = \dfrac{1}{2}(e^{2t}-1)   ,   \qquad w(y,\tau) = \alpha(\tau) \rho\big(y\alpha(\tau), -\log\alpha(\tau)\big),\]
with the notation
$ \alpha(\tau)=\dfrac{1}{\sqrt{2\tau+1}}=e^{-t}$.  
The firing rate in the new coordinates is 
$ M(\tau) = \alpha(\tau)^2 N(t)$.

Now, we apply the second change of variables
\begin{equation*}
	x\, =\, y- \int_0^\tau\left( \nu_{\mathrm{ext}}\alpha(s) + b\frac{ M(s)}{\alpha(s)} \right)\diff s\ =\, y - \nu_{ext} (  \sqrt{1+2\tau}-1  ) - b\int_0^\tau \dfrac{M(s)}{\alpha(s)}\diff s, 
\end{equation*}
and we denote $u(x,\tau) = w(y,\tau)$. The pair $(u,M)$ is a solution to the following Stefan-like free boundary problem with a moving distribution in the right-hand side given by
\begin{equation}\label{EqnEqu}
	\left\{\begin{array}{ll}
		\displaystyle \dfrac{\partial u}{\partial \tau}(x,\tau) = \dfrac{\partial^2 u}{\partial y^2}(x,\tau) +M(\tau)\delta_{s(\tau)+\frac{V_R}{\alpha(\tau)}}( x) &\qquad x\in(-\infty,s(\tau)],\ \tau > 0 \\
		\displaystyle M(\tau)=-\dfrac{\partial u}{\partial x}(s(\tau),\tau)& \qquad \tau > 0 \\
		\displaystyle s(\tau) = s(0) - \nu_{\mathrm{ext}} (  \sqrt{1+2\tau}-1  ) - b\int_0^\tau \dfrac{M(s)}{\alpha(s)}\diff s & \qquad \tau > 0\\
		\displaystyle u(-\infty,\tau)=u(s(\tau),\tau)=0, \qquad u(x,0) = u^0(x)  & \qquad x\in(-\infty,s(0)],  \tau > 0 .
	\end{array}\right.  
\end{equation}
Notice that
\begin{enumerate}
	\item the mass is conserved: for all $\tau$, $\int_{-\infty}^{s(\tau)} u(x,\tau)\diff x=\int_{-\infty}^{0} u^0(x)\diff x\,$; 
	
	\item the flux through the point $s_1(\tau):=s(\tau)+\frac{V_R}{\alpha(\tau)}$ is exactly
	\[
	M(\tau):=-\partial_xu(s(\tau),\tau)=\partial_xu(s_1(\tau)^-,\tau)-\partial_xu(s_1(\tau)^+,\tau).
	\]
\end{enumerate}

Away from the singular points $s(\tau)$ and $s_1(\tau)=s(\tau)+\frac{V_R}{\alpha(\tau)}$, the problem is a linear heat equation. Hence, given the heat kernel and its space derivative
\[  G(x,\tau,\xi,\eta) = \dfrac{1}{\sqrt{4\pi (\tau-\eta)}}e^{-\frac{(x-\xi)^2}{4(\tau-\eta)}},\qquad  \dfrac{\partial G}{\partial \xi}(x,\tau,\xi,\eta) = \dfrac{x-\xi}{4\sqrt{\pi}(\tau-\eta)^{\unitfrac{3}{2}} }e^{-\frac{(x-\xi)^2}{4(\tau-\eta)}}  , \]
we have the Green identity
\[ \dfrac{\partial}{\partial \xi} \left( G\dfrac{\partial u}{\partial \xi} - u\dfrac{\partial G}{\partial \xi} \right) - \dfrac{\partial }{\partial \eta}\big( Gu \big)  = 0.\]
Fix any $\tau\in\R_+$ and any $x\in(-\infty,s_1(\tau))\cup(s_1(\tau),s(\tau))$. Integrating the identity and using the regularity and decay of both the solution and the heat kernel $G$ yields
\begin{multline}\label{eq:Duhamel}
	u(x,\tau)= \int_{-\infty}^{0}G(x,\tau,\xi,0)u^0(\xi)\diff \xi - \int_{0}^{\tau}M(\eta)G(x,\tau,s(\eta),\eta)\diff \eta\\ + \int_{0}^{\tau}M(\eta)G(x,\tau,s_1(\eta),\eta)\diff \eta.
\end{multline}
The first term is the solution to the linear heat equation on the whole real line with initial condition $u^0$; the other two terms represent the loss of neurons at the free boundary and their return at a moving reset point $s_1(\tau)$. Now that we have $u$ in terms of only $u^0$ and the firing rate $M$, it is possible to isolate $M$ by taking the partial derivative along $x$ and then the limit $x\to s(\tau)$. It leads to the close formula
\begin{multline}\label{PointFix}
	M(\tau) = - 2 \int_{-\infty}^{0} G(s(\tau),\tau,\xi,0)\dfrac{\diff u^0}{\diff x}(\xi)\diff \xi  + 2\int_{0}^{\tau} M(\eta)\dfrac{\partial}{\partial x} G(s(\tau),\tau,s(\eta),\eta)\diff \eta \\  - 2\int_{0}^{\tau} M(\eta)\dfrac{\partial}{\partial x} G(s(\tau),\tau,s(\eta)+\frac{V_R}{\alpha(\eta)},\eta)\diff \eta.
\end{multline}

From there, it is possible to apply the Banach fixed point theorem in $\{ M\in C([0,\tau_0])\ | \ \norme{M}_\infty \leqslant m \}$ with $\tau_0$ small enough and $m$ large enough. For
\[  \tau_0 < \norme{ \frac{\diff u_0}{\diff x} }_{L^\infty}^{-1}, \]
there exists a unique continuous local-in-time solution $M$ on $[0,\tau_0]$ and the maximal time of existence is $\tau^*=\sup\{ \tau\in\R_+\ | \ M(\tau) <  +\infty \}$. The existence of the firing rate implies, through the Duhamel formula \eqref{eq:Duhamel}, the existence of a solution $(u,M)$ to the Stefan-like problem \eqref{EqnEqu}. Undoing all changes of variables yields the following result.

\begin{theorem}[Carrillo, Gonzalez, Gualdani, Schonbeck \cite{CGGS}]\label{thm:NNLIFexistence}
    Assume $d=0$, $a(N) = a > 0$ is constant, and $\rho^0$ is a compatible initial condition. Then there exists a unique strong solution to \eqref{NNLIF} whose maximal time of existence $T^*$ satisfies
    \begin{equation}\label{eq:maxTNNLIF}
        T^* = \sup\{ \tau\in\R_+\ | \ N(t) <  +\infty \}. 
    \end{equation}
\end{theorem}

\subsubsection{Global-in-time existence and uniform bounds in the inhibitory case}

Consider again the case $d=0$,  $a(N) \equiv a > 0$ constant. Then, in the case $b < 0$ of an inhibitory network, strong solutions are in fact global in time. Two different methods were used to prove it.

First, in \cite{CGGS}, global-in-time existence for inhibitory networks is proved using the Stefan problem reformulation. The main idea is that the nonlinear term in the free boundary motion is monotonically increasing, yielding the key estimate
\[  s(\tau) - s_1(\eta) \geqslant |V_R| - |\nu_{\mathrm{ext}}| (\tau-\eta), \qquad \eta < \tau. \]
Then, all terms in the fixed-point formulation \eqref{PointFix} can be controlled uniformly on a short enough time and a contradiction argument with \eqref{eq:maxTNNLIF} proves global-in-time existence.

Then, in \cite{CPSS}, the same result was achieved and improved by using universal super-solutions. The idea is to use the fact that in the inhibitory case, non-increasing super-solutions of the linear problem are super-solutions for the nonlinear problem and can control the firing rate.

\begin{definition}[universal super-solutions] Assume $d=0$,  $a(N) \equiv a > 0$ and $b < 0$. Let $V_0\in(-\infty,V_F)$ and $T > 0$. $(\bar \rho, \bar N)$ is a universal super-solution to \eqref{NNLIF} on $[V_0,V_F]\times[0,T]$ if
\[  \dfrac{\partial \bar \rho}{\partial t} - \dfrac{\partial }{\partial v}[ v\bar\rho] - a \dfrac{\partial^2 \bar \rho}{\partial v^2} \geqslant \delta_{v=V_R} \bar N(t) \]
in the distributional sense on $[V_0,V_F]\times[0,T]$, in the classical sense away from $V_R$, if $\bar N(t) = -a\partial_v \bar\rho (V_F,t)$ and if for all $t$, $\bar \rho(\cdot,t)$ is non-increasing on $[V_0,V_F]$.
\end{definition}
This definition allows the following comparison principle.
\begin{lemma}
    Assume $d=0$,  $a(N) = a > 0$ and $b < 0$. Let $V_0\in(-\infty,V_F)$ and $T > 0$. Let $(\rho,N)$ be a classical solution of \eqref{NNLIF} defined up to time $T$ and $(\bar \rho,\bar N)$ be a universal super-solution on $[V_0,V_F]\times[0,T]$. Assume that for all $(v,t)\in[V_0,V_F]\times[0,T]$,
    \[
    \bar \rho(v,0) \geqslant \rho^0(v), \qquad \bar\rho(V_0,t)  \geqslant \rho(V_0,t)
    \]
    Then $\bar \rho \geqslant \rho$ on $[V_0,V_F]\times[0,T]$ and if $\bar\rho(\cdot,0)-\rho^0$ is not identically zero, then $\bar \rho > \rho$ on $(V_0,V_F)\times(0,T]$.
\end{lemma}
Using two specific families of super-solutions constructed at hand, and taking advantage of the Stefan reformulation and classical heat equation estimates, \cite{CPSS} proves

\begin{theorem}[Carrillo, Perthame, Salort, Smets \cite{CPSS}; generalising \cite{CGGS}]\label{thm:NNLIFglobalinhib}
    Assume $d=0$,  $a(N) = a > 0$ is constant and $b < 0$. Let $\rho^0$ be a compatible initial condition. Then there exists a unique global-in-time solution $(\rho,N)$ to \eqref{NNLIF} and
    \begin{equation*}  \rho\in L^\infty((-\infty,V_F] 
    \times 
    [0,+\infty)),\qquad N\in L^\infty(0,+\infty). \end{equation*}
\end{theorem}

\noindent The result can be improved to $\rho^0\in L^1_+\cap L^\infty(-\infty,V_F)$ with
\[  \limsup_{v\to V_F} \dfrac{\rho^0}{V_F-v} < +\infty \]
if there exists a strong solution. The influence of less stringent decays near the boundary has been studied in the context of finance models in \cite{HLS,DNS}.

Note the crucial assumption $a(N) = a > 0$ constant; as we shall discuss, in the case with linear internal noise $a(N) = a_0 + a_1 N$, $a_0,a_1 > 0$, finite time blow-up occurs even in the inhibitory case. 

\subsection{Finite time blow-up in non-delayed networks (d=0)}

When there is no synaptic delay ($d=0$), the NNLIF model is prone to finite time blow-up. As we can see in Theorem \ref{thm:NNLIFexistence}, finite-time blow-up is caused by the divergence of the firing rate $N(t)$ at the blow-up time $T^*$:
\[  \limsup_{t\to T^*} N(t) = +\infty. \]
From a biology perspective, it is a rough depiction of a macroscopic part of the neural network synchronising and firing at the same time.

A first proof of this phenomenon was provided in \cite{CCP}, using a contradiction argument on the exponential moment 
\[ M_\mu(t) = \int_{-\infty}^{V_F}e^{\mu v} \rho(v,t)dv.\]
Indeed, for $\mu$ large enough, it satisfies the differential inequality
\[ \dfrac{\diff }{\diff t} M_\mu (t) \geqslant b\mu N(t)\left[ M_\mu(t) -  \dfrac{e^{\mu V_F}-e^{\mu V_R}}{b\mu}\right], \]
which allows to prove that if the right-hand side is initially positive, then $M_\mu$ is increasing and, using another lower-bound, that $M_\mu$ is unbounded, contradicting the natural bound $M_\mu(t)\leqslant e^{\mu V_F}$.
\begin{theorem}[C\'aceres, Carrillo, Perthame \cite{CCP}]\label{ExploPert}
	Assume $d=0$, $a(N)\geqslant a_0>0$ and $b>0$. Let $\mu>\max(\frac{V_F}{a_0},\frac{1}{b})$. If
	\begin{equation}\label{eq:blowup}b\mu\int_{-\infty}^{V_F}e^{\mu v}\rho^0(v)\diff v > e^{\mu V_F} - e^{\mu V_R} , \end{equation}
	then there is no global-in-time weak solution of \eqref{NNLIF}.
\end{theorem}

The condition \eqref{eq:blowup} is satisfied, for example, whenever $\rho^0$ is concentrated close enough to the firing threshold: if $\rho^0$ is fully supported in $[V_F-\varepsilon,V_F]$ with $\varepsilon$ small enough, then $b\mu\geqslant 1$ and
\begin{align*}\int_{-\infty}^{V_F}e^{\mu v}\rho^0(v)\diff v\,  &=\, \int_{V_F-\varepsilon}^{V_F}e^{\mu v}\rho^0(v)\diff v\\ &\geqslant \,e^{-\mu \varepsilon}e^{\mu V_F}\int_{V_F-\varepsilon}^{V_F} \rho^0(v)\diff v  \,=\, e^{-\mu \varepsilon}e^{\mu V_F}\int_{-\infty}^{V_F}\rho^0(v)\diff v =\, e^{-\mu \varepsilon}e^{\mu V_F}\\ &>\, e^{\mu V_F} - e^{\mu V_R}.
\end{align*}
The smaller $b>0$, the higher $\mu$, and then the smaller $\varepsilon$ we need. For a low excitatory connectivity $b>0$, blow-up occurs when the initial density is very concentrated near the firing threshold. In contrast, the stronger the connectivity, the easier it is for the network to synchronise.

In fact, it was proved in \cite{RS} that for a strong enough excitatory connectivity, all initial conditions lead to finite-time blow-up.

\begin{theorem}[Roux, Salort \cite{RS}]{\color{white},}
    \begin{itemize}
        \item If $V_F\leqslant 0$, and $b\geqslant V_F-V_R$, no fast-decaying weak solution of \eqref{NNLIF} can be global-in-time.
        \item If $V_F>0$, there exists a value $b_c>0$ such that if $b\geqslant b_c$\,, no fast-decaying weak solution of \eqref{NNLIF} can be global-in-time.
    \end{itemize}
\end{theorem}

The main idea is to study the linear moment
\[ M_1(t) = \int_{-\infty}^{V_F}(V_F-v)\rho(v,t) dv \]
and to prove that its evolution in time leads to enough concentration of the density to use the contradiction method on the exponential moment. Note that the case $V_F\leqslant 0$ implies the existence of an excitatory external input $\nu_{\mathrm{ext}}$ which is strong enough to cause neurons to fire even in the absence of noise. It explains why finite-time blow-up, an avatar of synchronisation, occurs more easily in this parameter range.

As seen in Theorem \ref{thm:NNLIFglobalinhib}, constant diffusion $a(N)\equiv a>0$ and an inhibitory connectivity $b<0$ prevent blow-up from occurring. However, when there is linear internal noise in the form $a(N)=a_0+a_1N$, $a_0, a_1>0$, there can be finite time blow-up even in the inhibitory case.

\begin{theorem}[Carrillo, Perthame, Salort, Smets \cite{CPSS}] \label{thm:blowup_internal_noise}
    Assume $d=0$, $a(N)=a_0+a_1N$, $a_0, a_1>0$. Let $\mu> 0$ be such that $a_1\mu^2+\mu b>1$ and $a_0\mu^2-V_F>0$. If
    \begin{equation}\label{eq:blowup2}\int_{-\infty}^{V_F}e^{\mu v}\rho^0(v)\diff v > e^{\mu V_F} - e^{\mu V_R} , \end{equation}
	then there is no global-in-time weak solution of \eqref{NNLIF}.
\end{theorem}
The reason for this discrepancy between the cases $a_1=0$ and $a_1\neq 0$ is that in the latter $N(t)$ is defined by the self-consistency equation
\[  N(t) = -\big(a_0+a_1 N(t)\big) \dfrac{\partial \rho}{\partial v}(V_F,t), \]
which can be rewritten as
\[ N(t) = \dfrac{-a_0\partial_v \rho(V_F,t)}{1+a_1 \partial_v \rho(V_F,t)}. \]
Because $N(t)$ must remain finite and non-negative for classical solutions, internal noise adds the well-posedness constraint $a_1 \partial_v \rho(V_F,t) < 1 $.\\

As per the biological interpretation of finite time blow-up as a partial or total synchronisation event, the solutions should not stop at the first blow-up time. In principle, it should be possible to compute which part of the mass blows-up in the density and to reset it at $V_R$ as a proportionate delta Dirac mass. From the mathematical analysis perspective, this problem is challenging and, although some progress have been made, a lot of questions are still open. In what follows, we will discuss some existing probabilistic \cite{DIRT2,ST,TW} and deterministic \cite{DZ,DPSZ2,CDRZ} methods. Apart from conjectures \cite{HLS} and results for related models (see Theorem \ref{thm:NNLIF_dPFM} below), the blow-up asymptotics, namely the divergence speed of $N(t)$ close to $T^*$, are poorly understood (see however \cite{DNS,HV} for the Stefan formulation).

\subsection{Desynchronisation in the weakly connected case}

We have seen that in some cases the NNLIF model synchronises through finite-time blow-up. An opposite phenomenon can occur: convergence of the solution towards a stationary state. A stationary density $\rho_\infty(v)$ does not mean that nothing is happening in the network. At the microscopic level, many neurons fire continuously, but at the global scale, the distribution is not evolving. This is akin to the convergence of a gas towards a state of higher entropy. As we shall see, the relative entropy method that works well for the linear Fokker-Planck equation on the whole space (see the classical review \cite{MV}) can be adapted, with enough work, to the weakly nonlinear regime of the NNLIF model, as was the case for many biological models in PDE form \cite{MMP2,MMP1,P}. In the case of the NNLIF model, specific difficulties due to the non-local structure require appropriate control of the firing rate $N(t)$.

In all this subsection, we choose $d=0$ and $a(N)\equiv a > 0$ constant.

\subsubsection{\textit{A priori} convergence in relative entropy}

\paragraph{Linear case $b=0$.}

Let us first showcase the method on the linear case $b=0$ where neurons are disconnected from each other. Theorem \ref{thm:ClassifStatic} ensures that there exists a unique stationary state $(\rho_\infty,N_\infty)$. For any convex and $C^2$ function $G:\R\to\R_+$, we have the entropy dissipation
\[
\begin{array}{rcl}
	\displaystyle{ \dfrac{\diff}{\diff t}\int_{-\infty}^{V_F}   G\left(\dfrac{ \rho}{ \rho_\infty}\right) \rho_\infty \diff v } \ &=&  \
	\displaystyle{ -a\int_{-\infty}^{V_F}  
		\left[\dfrac{\partial}{\partial v}\left(\dfrac{ \rho}{ \rho_\infty}\right)
		\right]^2G''\left(\dfrac{ \rho}{ \rho_\infty}\right) \rho_\infty \diff v   } \\
	&&	 \displaystyle{  -N_{\infty}
		\Big[G\left(\dfrac{N}{N_\infty}\right) - G\left(\dfrac{ \rho}{ \rho_\infty}\right) - \Big( \dfrac{N}{N_\infty} - \dfrac{ \rho}{ \rho_\infty} \Big)
		G'\left(\dfrac{ \rho}{ \rho_\infty}\right)   \Big]\bigg\rvert_{v=V_R}  }.
\end{array}
\]
Note that by convexity of $G$, the second part of the right-hand side is always non-positive.
We can choose $G(x)=(x-1)^2$ and apply a Poincaré-like inequality to the first term, which yields
\[ \dfrac{\diff}{\diff t} \int_{-\infty}^{V_F}  \left(  \dfrac{\rho(v,t) - \rho_\infty(v)}{\rho_\infty(v)} \right)^2 \rho_\infty(v) \diff v \leqslant -  2 a \nu \int_{-\infty}^{V_F}  \left(  \dfrac{\rho(v,t) - \rho_\infty(v)}{\rho_\infty(v)} \right)^2 \rho_\infty(v) \diff v, \]
where $\nu$ is the Poincaré constant. By Grönwall's lemma,
\[  \int_{-\infty}^{V_F}  \left(  \dfrac{\rho(v,t) - \rho_\infty(v)}{\rho_\infty(v)} \right)^2 \rho_\infty(v) \diff v \leqslant e^{- 2 a \nu t} \int_{-\infty}^{V_F}  \left(  \dfrac{\rho^0(v) - \rho_\infty(v)}{\rho_\infty(v)} \right)^2 \rho_\infty(v) \diff v,   \]
which is the convergence in relative entropy of $\rho(\cdot,t)$ towards $\rho_\infty$ at exponential speed.

\paragraph{Weakly nonlinear case.}

Proceeding like in the linear case, for any convex and $C^2$ function $G$ and a solution with appropriate decay, we compute the entropy dissipation
\[
\begin{array}{rcl}
	\displaystyle{ \dfrac{\diff}{\diff t}\int_{-\infty}^{V_F} 
		G\left(\dfrac{ \rho}{ \rho_\infty}\right)  \rho_\infty  \diff v } \ &=&  \
	\displaystyle{ -a\int_{-\infty}^{V_F} 
		\left[\dfrac{\partial}{\partial v}\left(\dfrac{ \rho}{ \rho_\infty}\right)
		\right]^2G''\left(\dfrac{ \rho}{ \rho_\infty}\right)\rho_\infty  \diff v   } \\
	&&	 \displaystyle{  -N_{\infty}
		\Big[G\left(\dfrac{N}{N_\infty}\right) - G\left(\dfrac{ \rho}{ \rho_\infty}\right) - \Big( \dfrac{N}{N_\infty} - \dfrac{ \rho}{ \rho_\infty} \Big)
		G'\left(\dfrac{ \rho}{ \rho_\infty}\right)   \Big]\bigg\rvert_{v=V_R}  }\\
	&&	 \displaystyle{    + b(N(t)-N_\infty)\int_{-\infty}^{V_F}\dfrac{\partial  \rho_\infty}{\partial v}\Big[G\left(\dfrac{ \rho}{ \rho_\infty}\right) -\dfrac{ \rho}{ \rho_\infty}G'\left(\frac{ \rho}{ \rho_\infty}\right)\Big]\diff v }.
\end{array}
\]

For $|b|$ small enough, using natural upper bounds, Sobolev injections, a Poincaré-like inequality of constant $\nu$ and Grönwall's lemma, we can get to
\begin{equation}\label{eq:almost}
	\int_{-\infty}^{V_F}   
	G\left(\dfrac{ \rho}{ \rho_\infty}\right) \rho_\infty \diff v  \leqslant e^{-\nu at + \frac{8b^2}{a} \int_0^t |(N(s)-N_\infty)|^2 \diff s} \int_{-\infty}^{V_F} 
	G\left(\dfrac{ \rho^0}{ \rho_\infty}\right)  \rho_\infty  \diff v.  
\end{equation}
In order to conclude, we need a uniform estimate of the form
\begin{equation}\label{eq:L2N}
    \int_0^t N(s)^2 \diff s \leqslant C(1+t),
\end{equation}
where $C$ does not depend on $b$, so we can take $|b|$ small enough with respect to $\nu$, $a$ and $C$.

We can apply again a relative entropy method, but in a neighbourhood of $V_F$, in order to isolate the firing rate. However, it would not yield a $b$-uniform estimate. The trick is to compare the solution $\rho$ not with the stationary state $\rho_\infty$ but rather with $\rho_\infty^1$ obtained with another parameter $b_1$.

Fix $b_1  >  0$ such that there is at least one stationary state $(\rho_\infty^1,N_\infty^1)$. We work near $V_F$ by fixing some $V_M\in(V_R,V_F)$ and using the cutoff function
$\gamma(v) = \mathds 1_{ v > V_M }e^{\frac{-1}{(V_F-V_M)^2-(V_F-v)^2}}$. Define
\[  W(v,t)= \left(\dfrac{p(v,t)}{p_\infty^1(v)}\right)^2p_\infty^1(v) =\dfrac{p(v,t)^2}{p_\infty^1(v)}  \qquad \mathrm{and} \qquad I(t)=\int_{-\infty}^{V_F}W(v,t)\gamma(v) \diff v. \]
The key reason for this quantity $W$ is that applying the boundary conditions for the equations \eqref{NNLIF} associated to the parameters $b$ and $b_1$, we can define by continuity
\[  - a \dfrac{\partial W}{\partial v}(V_F,t) = \dfrac{N(t)^2}{N_\infty^1}, \qquad W(V_F,t) = 0.  \]
Then, it all comes down to finding differential inequalities for the modified relative entropy $I(t)$.
\begin{itemize}
	\item[$\bullet$] \textit{Inhibitory case}: Assume $b\leqslant0$. There exist $C_1,C_2,C_3 > 0$ independent of $b$ such that
	\begin{equation} \label{integramos}
		\frac{d}{dt}I(t) \leq -C_1N(t)^2+C_2 +(bN(t)-C_3)I(t).
	\end{equation}
	\item[$\bullet$]  \textit{Excitatory case}: Assume $b  > 0$. There exist $C_1,C_2,C_3,C_4 > 0$ independent of $b$ such that
	\begin{equation}\frac{\diff}{\diff t}I (t)\leqslant N(t)^2\Big( C_1b^2 - C_2 + C_3 b^2 I(t) \Big)  + C_3  - C_4 I(t). \label{Ayda} \end{equation}
\end{itemize}
These differential inequalities allow to get uniform $L^2$ control over $N$. However, because of the multiplicative term $I(t)$ in \eqref{Ayda}, we need $b$ small enough with respect to $I(0)$ in the excitatory case $b > 0$. In both cases, it yields an estimate of the form \eqref{eq:L2N} (see Theorem \ref{thm:LqN} for a generalised statement).
 Applying these $L^2$ estimates to \eqref{eq:almost}, and choosing also $b$ small enough yields:
\begin{theorem}[Carrillo, Perthame, Salort, Smets \cite{CPSS}]\label{thm:NNLIFentropy}
	Assume $d=0$, $a(N)\equiv a >  0$ constant. Let us fix $b_1  >  0$ such that there exists at least one corresponding stationary state $(\rho_\infty^1,N_\infty^1)$ and $V_M\in(V_R,V_F)$ . Let $(\rho,N)$ a classical solution of \eqref{NNLIF} such that,
	\[  \int_{V_M}^{V_F}  \dfrac{p^0(v)^2}{p_\infty^1(v)} \diff v < +\infty,\quad \int_{-\infty}^{V_F}  \left(\dfrac{\rho^0(v)-\rho_\infty(v)}{\rho_\infty(v)}\right)^2 \rho_\infty(v) \diff v < +\infty. \] 
	There exist positive constants $b^*,\mu>0$ such that for all $b\in[-b^*,b^*]$,
	\[  \int_{-\infty}^{V_F}  \left(  \dfrac{\rho(v,t) - \rho_\infty(v)}{\rho_\infty(v)} \right)^2 \rho_\infty(v) \diff v \leqslant e^{- \mu t} \int_{-\infty}^{V_F}  \left(  \dfrac{\rho^0(v) - \rho_\infty(v)}{\rho_\infty(v)} \right)^2 \rho_\infty(v) \diff v.   \]
\end{theorem}

When $b\leqslant 0$, $b^*$ and $\mu$ can be chosen independently of $\rho^0$ if we start at a larger time $T=T(\rho^0)$,
\[\int_{-\infty}^{V_F}  \left(  \dfrac{\rho(v,t) - \rho_\infty(v)}{\rho_\infty(v)} \right)^2 \rho_\infty(v)\diff v \leqslant e^{- \mu (t-T)} \int_{-\infty}^{V_F}  \left(  \dfrac{\rho(v,T) - \rho_\infty(v)}{\rho_\infty(v)} \right)^2 \rho_\infty(v)\diff v.\] In the excitatory case $b > 0$, $\mu$ and $b^*$ heavily depend on $\rho^0$ and there is nothing we can do about it, since there is blow-up for some initial conditions however small we choose $b > 0$ (Theorem \ref{ExploPert}). The more concentrated $\rho^0$ is around $V_F$, the smaller $b$ has to be. In the inhibitory case, this \textit{a priori} estimate is true for any small enough $b$ independently of $\rho^0$. Although the $L^2$ estimate on $N$ holds for any inhibitory connectivity parameter $b\in(-\infty,0]$, we still need $b$ small for $\mu$ to be positive.

\subsubsection{$L^q$ estimates on the firing rate and global existence for weak interaction}

Theorem \ref{thm:NNLIFentropy} is an \textit{a priori} convergence result, valid as long as the solution exists; it does not enforce global-in-time existence by itself. Indeed, convergence happens in the relative entropy space $L^2(\rho_\infty^{-1})$ in which the boundary derivative operator is not bounded (see Subsection \ref{sec:NNLIF_general_b} below for another approach involving a new stronger space). Global-in-time existence for the weakly nonlinear case was achieved first in \cite{DIRT1} using probability methods on the associated stochastic differential equation (see Subsection \ref{sec:NNLIF_particles} below for an introduction of this SDE).

\begin{theorem}[Delarue, Inglis, Rubenthaler, Tanré \cite{DIRT1}]
    For any initial condition $V_0 < V_F$, there exists $b^*=b^*(V_0)$ such that for all $b\in [0,b^*]$, there is a unique global-in-time strong solution to the stochastic differential equation \eqref{eq:NNLIF_MF}. 
\end{theorem}

Later, the method of \cite{CPSS} was extended in \cite{RS} to prove a similar result in a deterministic PDE framework. The idea is to use an iteration argument to get from the $L^2$ estimates of \cite{CPSS} to $b$-uniform $L^q$ estimates for any integer $q\geqslant 2$.
\begin{theorem}[Roux, Salort \cite{RS}, generalising \cite{CPSS}]\label{thm:LqN}
	Assume $d=0$, $a(N) \equiv a  >  0$ is constant. Fix $q\geqslant 2$ integer and $b_1  >  0$ such that there is at least one stationary state $(\rho_\infty^1,N_\infty^1)$ and $V_M\in(V_R,V_F)$ . For any compatible initial condition $\rho^0$ such that
	$ \int_{V_M}^{V_F}  \frac{p^0(v)^q}{p_\infty^{q-1}(v)} \diff v < +\infty$, 
	the solution of \eqref{NNLIF} satisfies
    \begin{itemize}
    \item There exists $C_q$ depending only on $b_1,q,a,V_F$ and $T=T(\rho^0)$ such that for all $b \leqslant 0$,
    \[  \forall t\geqslant T,\qquad  \int_{T}^t N(s)^q \diff s \leqslant C_q \,(1+t-T). \]
    \item There exist $C_q$ depending only on $b_1,q,a,V_F,\rho^0$ and $b^*=b^*(\rho^0)$ such that for all $0 < b \leqslant b^*$,
    \[  \forall t\geqslant 0,\qquad  \int_{0}^t N(s)^q \diff s \leqslant C_q \,(1+t). \]
    \end{itemize}
\end{theorem}
Then, coming back to the Stefan-like free-boundary problem reformulation exposed above, it is possible to use $L^3$ estimates on the firing rate in the representation formula \eqref{PointFix} and to repeat the contradiction argument used by \cite{CGGS} for the inhibitory case.
\begin{theorem}[Roux, Salort \cite{RS}]
    Assume $d=0$, $a(N) \equiv a  >  0$ is constant. Let $\rho^0$ be a compatible initial condition. There exists $b^*=b^*(\rho^0)$ such that for all $b\in(-\infty,b^*]$, the unique strong solution of \eqref{NNLIF} is global-in-time.
\end{theorem}
The reason $L^2$ estimates are not enough but $L^3$ estimates on $N$ yield global-in-time existence can be seen in the Hölder inequality
\[ \int_0^\tau M(\eta) \dfrac{\partial G}{\partial x}(s(\tau),\tau,s(\eta),\eta) \diff \eta \leqslant \left(\int_0^\tau M(\eta)^3\diff \eta\right)^\frac13\left( \int_0^\tau  \left|\dfrac{\partial G}{\partial x}(s(\tau),\tau,s(\eta),\eta)^\frac23 \right|\diff \eta\right)^\frac23. \]
Because the free boundary position $s(\tau)$ is Lipschitz, the second integral is convergent. However, if the exponents $\frac12+\frac12=1$ are used for the inequality, then the second part is a diverging integral.

\subsection{The delayed NNLIF model and self-sustained oscillations}
\label{sec:NNLIF_delay}

As we said above, the goal of Brunel and Hakim in \cite{BH} was to explain fast periodic activity in the strongly inhibitory case. In order to observe it in the model, we need to consider the delayed case $d > 0$.

In this subsection, we consider the case $a(N)\equiv a > 0$ constant. We talk about the first wave of results on the delayed NNLIF model. Recent advances are described below in Subsection \ref{sec:NNLIF_advances}

\subsubsection{Global-in-time existence and relaxation for weak interaction}

In \cite{CRSS}, the aforementioned techniques for the non-delayed case $d=0$ have been extended to the case $d > 0$. It is first still possible to use the Stefan free-boundary reformulation to prove local-in-time existence, with similar representation formulae for $\rho$ and $N$ in modified variables. The delayed structure even allows for shortcuts. The maximal time of existence is still characterised by
\[  T^* = \sup\{ \tau\in\R_+\ | \ N(t) <  +\infty \}. \]
A striking difference with the case $d=0$ is that there is unconditional global-in-time existence of strong solutions. Because of the delayed impact of the firing rate, the notion of universal super-solution can be extended to include the nonlinear part of the problem. On time intervals of length $d$, nonlinear super-solutions can be constructed using the firing rate at previous times in $[t-d,t)$. These super-solutions prove that $T^*=+\infty$ in all cases.

\begin{theorem}[C\'aceres, Roux, Salort, Schneider \cite{CRSS}]
    Assume $d > 0$ and $a(N)\equiv a > 0$. Let $\rho^0$ be a compatible initial condition. Then there exists a unique global-in-time strong solution to \eqref{NNLIF}.
\end{theorem}

Then, the work on $b$-uniform $L^2$ estimates and generalised entropy can be extended to the case with delay. Unfortunately, during the process, the iteration argument on time intervals of length $d$ forces the appearance of additional terms depending on $d$, the consequence being that smallness of $d$ needs also to be assumed and not only that of $b$.

\begin{theorem}[C\'aceres, Roux, Salort, Schneider \cite{CRSS}]\label{thm:NNLIF_d_stat}
   Under the assumptions of Theorem \ref{thm:NNLIFentropy} but with $d>0$, there exist $b_1=b_1(d)>0$, $b_2=b_2(d,\rho^0)>0$, $\mu=\mu(d,\rho^0)>0$ such that for all $b\in[-b_1,b_2]$, there exists $S_0(\rho^0,N^0,b,d)$ satisfying 
	\[  \int_{-\infty}^{V_F}  \left(  \dfrac{\rho(v,t) - \rho_\infty(v)}{\rho_\infty(v)} \right)^2 \rho_\infty(v) \diff v \leqslant e^{- \mu t} S_0( \rho^0,N^0,b,d).  \]
\end{theorem}
As we shall see in Theorem \ref{thm:NNLIF_CCR-L_small_b}, the fact that the smallness assumption depends on $d$ is a technical hypothesis that other methods can remove \cite{CCR-L, CCR-L2}.

\subsubsection{Numerical insights, heuristics and partial results}

\paragraph{The excitatory delayed NNLIF: convergence \textit{versus} infinite-time blow-up}

We assume for now that we are in the excitatory and delayed case: $b > 0$ and $d > 0$. As we said, all solutions are global-in-time. Two very different regimes were identified (See Subsection \ref{subsec:pseudo} below for another perspective on this dichotomy using the sequence of pseudo-equilibria).\\

$\bullet$\ \textit{When $b$ is very large ($b > b^*$):} Since the stationary states are the same for $d=0$ and $d > 0$, Theorem \ref{thm:ClassifStatic} tells us that there are no stationary states. In this regime, there are also no periodic solutions.
\begin{theorem}[C\'aceres, Roux, Salort, Schneider \cite{CRSS}]
	If $b>V_F-V_R$ and $V_F\leq 0$, then for any $d\geq 0$ there are no classical periodic solutions to \eqref{NNLIF}.
\end{theorem}
What then happens? This question was answered with detailed numerical simulations in \cite{CR-L}. For a large excitatory connection, there is infinite-time blow-up of the firing rate $N(t)$:
\[ \lim_{t\to+\infty} N(t) = +\infty \qquad (\text{numerical result}).   \]
The density then converges towards a plateau state
\begin{equation}\label{eq:plateau} \rho(\cdot,t)\underset{t\to+\infty}{\longrightarrow}  \dfrac{1}{V_F-V_R}\mathds{1}_{[V_R,V_F]}(\cdot) \qquad (\text{numerical result}).   \end{equation}
A natural interpretation is that neurons fire, return to $V_R$ and are then sent increasingly fast to $V_F$. The plateau is made up of a trail of neurons that are firing and resetting. Asymptotically, all of them join the train, and the density becomes 0 outside of $[V_R,V_F]$ and uniform inside. This feature seems meaningless from the neuroscience point of view. A simple fix is to add a refractory period, which will change the structure of stationary states in the highly excitatory regime \cite{CP,CS2} (see the system \eqref{NNLIFR} below).\\

$\bullet$\ \textit{When $b$ is not too large ($b < b^*$):} In this case, there is at least one stationary state. Periodic solutions have never been observed numerically in this regime either. Numerical solutions exhibit different behaviours depending on the value of $b$, $a,V_R,V_F$ and the initial condition:
\begin{itemize}
	\item[i)] if $b<V_F-V_R$, there is a unique equilibrium; solutions either converge directly towards this stationary state or the firing rate increases sharply, all or part of the density is absorbed by the boundary $V_F$ and after resetting at $V_R$ of a portion of the mass, there is convergence towards a stationary state.
    \item[ii)] if $V_F-V_R\leqslant b<b*$ (note that $b^*\geqslant V_F-V_R$), there are often two equilibria; either the solutions converge directly towards a stationary state or the density converges towards a plateau state as in \eqref{eq:plateau}.
\end{itemize} 

When the delay $d$ tends to 0, this dichotomy sharpens into convergence \textit{versus} finite-time blow-up. The authors in \cite{CR-L} argue that the global-in-time blow-up friendly \textit{physical solutions} of \cite{DIRT2}, see Definition \ref{def:physsol},  will converge towards the stationary state after blow-up in the non-delayed case $d=0$. Note that these physical solutions are not defined in the case of large excitatory $b$, consistent with the fact that the plateau state is meaningless from the probabilistic perspective.

\paragraph{The inhibitory delayed NNLIF and the ongoing quest for periodic solutions}

Despite more than twenty years of numerical evidence \cite{BH,CS2,HLXZ}, no analytical existence result was ever achieved for periodic solutions of \eqref{NNLIF}. In \cite{IRSS}, a heuristic approach was proposed to understand the appearance of self-sustained oscillations in the strongly inhibitory case.\\

\begin{figure}[ht!]
	\begin{center}
		\includegraphics[width=400pt]{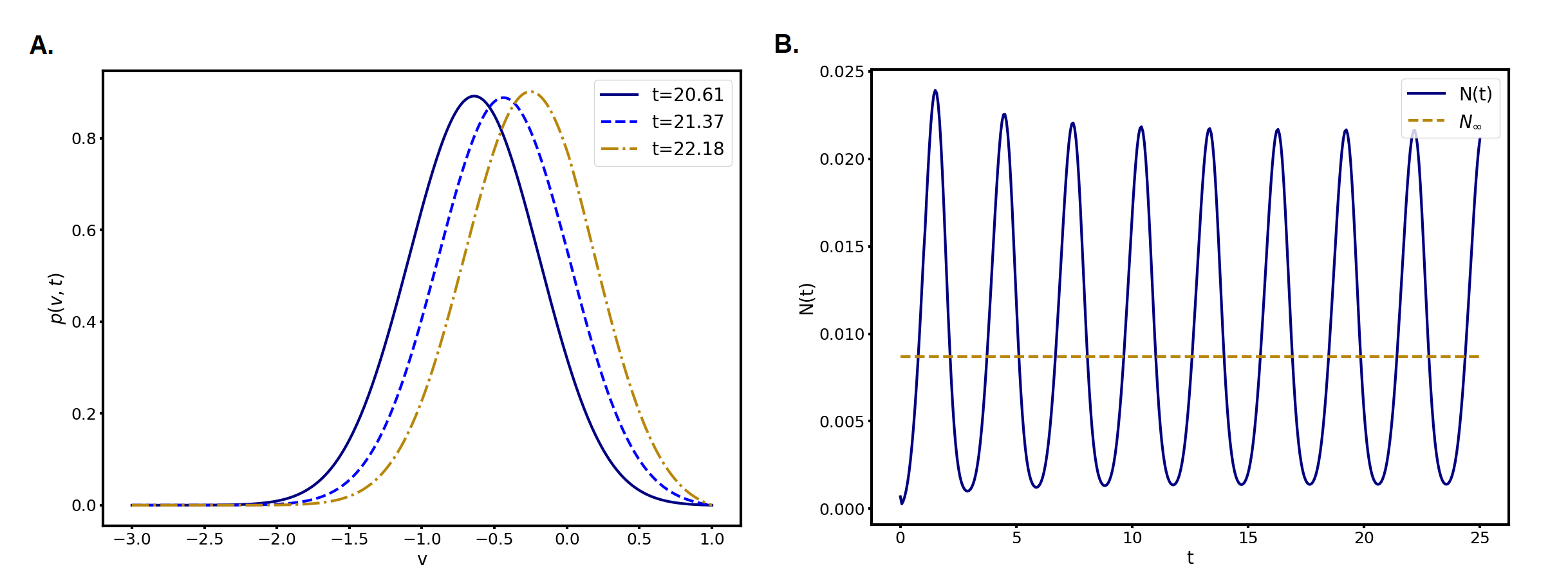}
	\end{center}
	\caption{Simulation of the PDE \eqref{NNLIF} with $V_F=1,V_R=0,a=0.2,d=1$ and $b=-45$; (A.) solution at three different times; (B.) evolution in time of the firing rate $N(t)$. Ikeda, Roux, Salort, Smets \cite{IRSS}.  \label{fig:1}}	
\end{figure}

We can look for a solution composed of a periodic wave of unit mass $\varphi(v,t)=\phi\big(v - c(t)\big)$ defined on the whole real line $\R$ plus a corrective term $R(v,t)$ needed to account for the boundary and reset conditions:
\begin{equation}\label{first}  p(v,t) = \phi\big(v - c(t)\big) + R(v,t),   \end{equation}
The wave $\varphi :(v,t)\mapsto \phi\big(v - c(t)\big)$ solves
\begin{equation}\label{eqn:wave_type}
	\left\{\begin{array}{l}
		\dfrac{\partial \varphi }{\partial t} + \dfrac{\partial }{\partial v} \Big[ (  -v + bN(t-d) ) \varphi  \Big] -a \dfrac{\partial^2 \varphi }{\partial v^2} = 0, \qquad v\in\R, t>0,\\
		\varphi(-\infty,t) = \varphi(+\infty,t)=0, \qquad \int_{-\infty}^{+\infty}\varphi(v,t)dv = 1.
	\end{array}\right.
\end{equation}
The remainder term $R:(v,t)\mapsto R(v,t)$ must be a solution of
\begin{equation}\label{eqn:remainder_term}
	\left\{\begin{array}{l}
		\dfrac{\partial R }{\partial t} + \dfrac{\partial }{\partial v} \Big[ (  -v + bN(t-d) ) R  \Big] -a \dfrac{\partial^2 R }{\partial v^2} = (\delta_{V_R} - \delta_{V_F})N(t), \qquad v\in\R,\\
		R(-\infty,t) = R(+\infty,t) = 0, \qquad R(v,0) = p^0(v)-\varphi(v,0),\qquad v\in\R.
	\end{array}\right.
\end{equation}

\begin{figure}[ht!]
	\begin{center}
		\includegraphics[width=400pt]{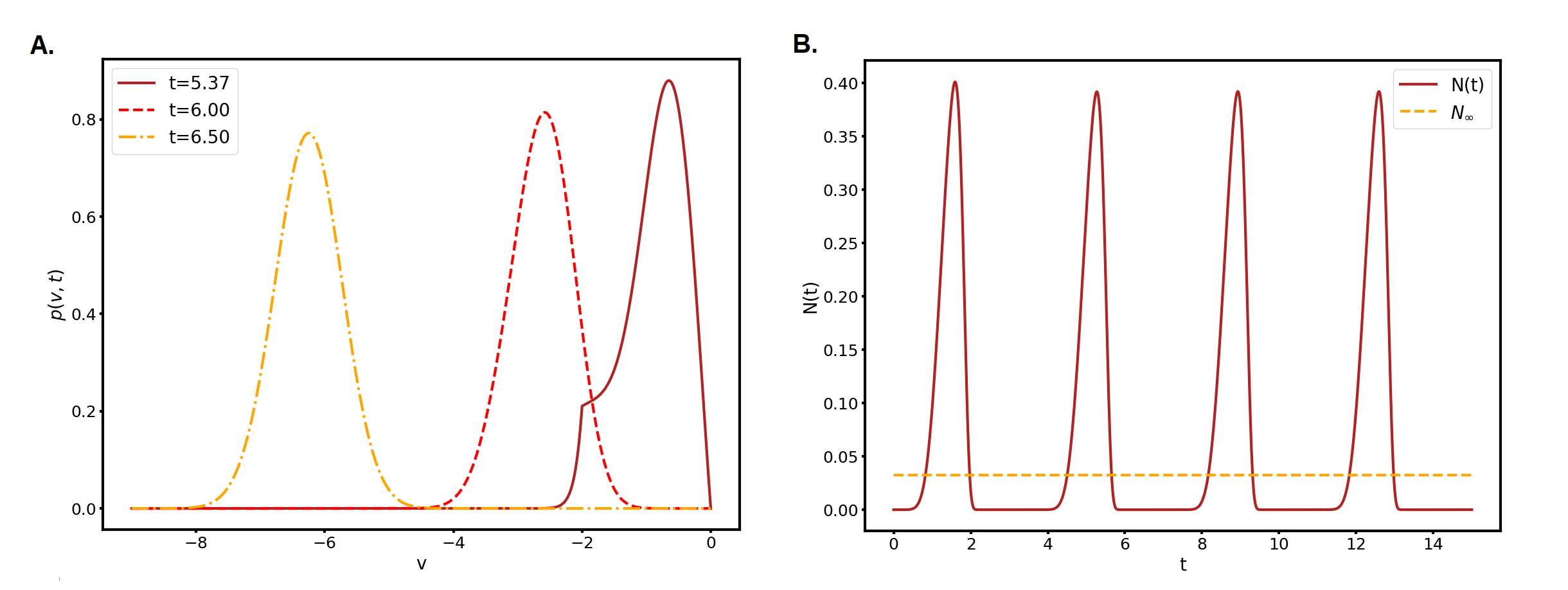}
	\end{center}
	\caption{Simulation of the PDE \eqref{NNLIF} with $V_F=0,V_R=-2,a=0.2,d=1$ and $b=-35$; (A.) solution at three different times; (B.) evolution in time of the firing rate $N(t)$. Ikeda, Roux, Salort, Smets \cite{IRSS}. \label{fig:3}}	
\end{figure}

\begin{theorem}[Ikeda, Roux, Salort, Smets \cite{IRSS}]\label{thm:1}
	Let $\phi:\R\to\R$ be the function defined by
	$  \phi(v)=\frac{1}{\sqrt{2\pi a}}e^{-\frac{v^2}{2a}}.$
	Let $c$ be a solution of
	\begin{equation}  \label{eqn:coupled}  c'(t) + c(t) = b N(t-d).   \end{equation}
	Then the function $\varphi$ defined by  $\varphi(v,t) = \phi(v-c(t))$ is a solution of \eqref{eqn:wave_type}.
\end{theorem}

Unfortunately, systems \eqref{eqn:wave_type} and \eqref{eqn:remainder_term} are strongly coupled through the firing rate:
\[N(t)=-a\partial_v \varphi(V_F,t) -a \partial_v R(V_F,t).\]
In order to make the problem autonomous and thus theoretically tractable, we assume that $b\partial_v R(V_F,t)\simeq 0$ in an appropriate sense when $b\to-\infty$. Hence, we replace $N(t)$ by the simpler firing rate
\begin{equation}\label{eqn:simpler_N}
	\mathcal N(c(t))=-a\dfrac{\partial \varphi}{\partial v}(V_F,t).
\end{equation}
in equation \ref{eqn:coupled}. It reduces the problem to the autonomous delayed differential equation (DDE)
\begin{equation} \label{eqn:harnak}
	c'(t)+c(t) = b \mathcal N(c(t-d)), \qquad \mathcal N(c) = \dfrac{1}{\sqrt{2\pi a}}\big(V_F-c\big) e^{-\frac{(V_F-c)^2}{2a}}.
\end{equation}
Exploiting the classical methods for DDEs \cite{HL2,S2,HT}, it is possible to prove existence of a periodic solution $c$, and thus a periodic solution $\varphi(v,t)=\frac{1}{\sqrt{2\pi a}}e^{-\frac{(v-c(t))^2}{2a}}$ of the simplified nonlinear PDE
\begin{equation}\label{eqn:wave_type2}
		\dfrac{\partial \varphi }{\partial t} + \dfrac{\partial }{\partial v} \Big[ (  -v + b\mathcal N(t-d) ) \varphi  \Big] -a \dfrac{\partial^2 \varphi }{\partial v^2} = 0, \quad
		\mathcal N(t) = -\dfrac{\partial \varphi}{\partial v}(V_F,t)\qquad v\in\R, t>0.
\end{equation}

\begin{theorem}[Ikeda, Roux, Salort, Smets \cite{IRSS}]\label{thm:existence_periodic_sol_0}
	Assume ${V_F}^2 \ge a$, then there exist $b_* < 0$ and $d_* > 0$ such that 
	\eqref{eqn:harnak} has a non-constant periodic solution $c (t)$ for any $b < b_*$ and $d > d_*$.
\end{theorem}

The way to prove this result is to use the Browder fixed point theorem.

\begin{definition}[Ejective fixed-point \cite{B5}]    \label{def:ejective}
	Let $E$ a Banach space, $D\subset E$ a closed subset, $\mathcal F : D \to D$ a continuous map. A fixed-point $\bar x$ of $\mathcal F$ is \emph{ejective} if there is a neighbourhood $U\subset D$ of $\bar x$ such that
	\[  \forall x\in U\setminus \{ \bar x \}, \ \exists n\in\N, \ \mathcal F^n(x) \notin U.  \]         
\end{definition}                                                                                                   
\begin{theorem}[Browder's fixed-point theorem \cite{B5}]\label{thm:browder}
	Let $D$ be a closed, bounded and convex subset of an infinite dimensional Banach space and let $\mathcal F : D\to D$ a continuous and compact map. Then $\mathcal F$ admits at least one fixed-point which is not \emph{ejective}.
\end{theorem}
Spectral study of \eqref{eqn:harnak} allows to prove that for large enough delay and strong enough inhibition, the stationary state $c_\infty$ is an ejective fixed-point of a carefully chosen functional $\mathcal F$. Then, there is a non-ejective fixed-point, which is associated to a periodic solution of \eqref{eqn:harnak}.\\

In \cite{IRSS}, asymptotic information is provided in the limit $b\to -\infty$. Assume for the sake of clarity that $V_F = 0$, $d=1$ and make the change of variables
\begin{equation}\label{eqn:rescale}
	\beta = \dfrac{1}{\log(-b)}, \qquad u(s) = \beta \log\left( - \dfrac{1}{\sqrt a} c\left(\dfrac{s}{\beta}\right)   \right),
\end{equation}
then we come to the following equivalent system:
\begin{equation}\label{eq1bis}
	u_\beta '(t)+1= \exp\left( \frac{1}{\beta}(1-u_\beta(t)) + f\big(u_\beta(t-\beta)\big)\right), \quad  1 \geq u_\beta([-\beta,0]) >0
\end{equation}
with
$	f(x)= \frac{1}{\beta}x - \frac{1}{2}e^{\frac{2}{\beta} x}.
$
Equation \eqref{eq1bis} has a unique positive stationary state $	\overline{u_\beta}= \frac{\beta}{2} \log\left(\frac{2}{\beta}\right)$.
Note that
$ \lim_{\beta \to 0}\overline{u_\beta} = 0$.
\begin{theorem}[Ikeda, Roux, Salort, Smets \cite{IRSS}]\label{thm:profileP0_0}
	Assume that for all $\beta>0$, $u_\beta([-\beta,0]) \equiv1$ and let $T>2$. Then, there exists a constant $C_M\in(0,1]$ such that
	$\lim_{\beta \to 0} u_{\beta} = P$, in $L^1(0,T)$,
	where $P$ has the following form:\\
	 {\color{white}.}\qquad $P(t)= 1-t$ on  $(0,1);$ \quad $P(t)= C_M+1-t$ on $(1,1+C_M);$ \quad $P$ is $C_M$ periodic on $(1, T)$. 
\end{theorem}

The assumption $u_\beta([-\beta,0]) \equiv1$ is purely technical, as can be checked by numerical simulations. In original variables, the result can be interpreted as follows
\begin{itemize}
	\item The period $T_b$ of solutions of \eqref{eqn:harnak} evolves in $\mathcal O\big(\log(-b)\big)$ when $b\to -\infty$.
	\item More precisely, we have in original variables
	\[  c(t) \simeq -\sqrt{a} \exp\left( \log(-b) P\left(  \dfrac{t}{\log(-b)}  \right)  \right).  \]
	In decay phases,  $c$ behaves like $ c(t)\simeq  -\sqrt{a} e^{   C_M \log(-b) - t}$,
	which constitutes an exponential decay over a spatial area of length $S_b\simeq \sqrt a (-b)^{C_M}$ during a time $T_b \simeq C_M \log(-b)$.
	\item between decays phases, $c$ has very rapid growth which is discontinuous in the re-scaled limit.
\end{itemize} 
For another asymptotic approach using the limit $d\to+\infty$ with $b$ fixed instead of $b\to-\infty$ with $d$ fixed, see Subsection \ref{subsec:pseudo} below.

\subsection{Rigorous derivation of the NNLIF model}

\subsubsection{Physical solutions from the limit of a particle system}
\label{sec:NNLIF_particles}
Although the NNLIF model was introduced in the late 90', the first rigorous construction from a mean-field limit of interacting particles was published way later, in 2015, in the article \cite{DIRT2}. Note that an overview on technical questions related to stochastic processes depending on singular hitting times was provided by the authors in the companion document \cite{DIRT3}.

For the sake of simplicity, the authors fix $V_R=0,V_F=1$ and make the assumption of a fully connected network with mean-field scaling. They consider the system of  $K$ particles satisfying the stochastic differential equation
\begin{equation}\label{eq:NNLIF_particles}
    V_t^{i,K} = V_0^{i,K} - \int_0^t V_s^{i,K}\diff s + \dfrac{b}{K}\sum_{j=0}^{K} M_t^{j,K}  + B_t^{i}- (V_F-V_R)M_t^{i,K}, 
\end{equation}
where $V_t^i$ represents the electrical potential of the $i^{th}$ neuron, all $V_0^{i,K}$ are independent and identically distributed with law $X_0$, and $(B_t^{i})_{t\geqslant 0}$ are independent Brownian motions on a probability space $(\Omega,\mathcal A,\mathbb P)$ with respect to some adequate filtration. The processes $(M_t^{i,K})$ count the number of times the $i^{th}$ neuron has spiked during the interval of time $[0,t]$, which can be precisely defined as càdlàg processes:
\begin{equation*}\label{eq:NNLIF_M_t} M_t^{i,K} = \sum_{l\geqslant 1} \mathds{1}_{[0,t]}\left( \tau_l^{i,K} \right), \qquad \tau_l^{i,K} = \inf\left\{ t > \tau_{l-1}^{i,K}\ | \ V_{t^-}^{i,K} + \frac{\alpha}{K} \sum_{j=1}^{K} \left( M_t^{j,K} - M_{t^-}^{j,K}  \right) \geqslant V_F  \right\}, \end{equation*}
and $\tau_0^{i,K}=0$. This system was previously introduced and studied in a physics framework in \cite{OBH}. In order to avoid neurons spiking multiple times at once and ensure that spikes happen in a well-ordered physical way, \cite{DIRT2} imposes $b < V_F-V_R$ (see also the problem of eternal blow-up in \cite{DZ}, described later in Section \ref{sec:gen_NNLIF}) and defines rigorously the notion of cascade: a first set $\Gamma_0$ of neurons spike, then a second set of neurons $\Gamma_1$ spike given that neurons in $\Gamma_0$ spike, etc. A cautious definition yields uniqueness for the solution of the particle system \eqref{eq:NNLIF_particles}.

Since the counting processes $M_t^{i,K}$ are naturally monotone, \cite{DIRT2} uses the M1 Skorokhod topology to perform the mean-field limit $K\to +\infty$. The limiting system writes
\begin{equation}\label{eq:NNLIF_MF}
    V_t = V_0 - \int_0^t V_s \diff s + b\mathbb E [M_t] + B_t - (V_F-V_R)M_t.
\end{equation}
Because $\mathbb E [M_t]$ might be a discontinuous function, defining a notion of solution for \eqref{eq:NNLIF_MF} is challenging. Recall, for any càdlàg function $F$, the notation
 $\Delta F(t) = F(t)-F(t^-)$.
\begin{definition}[Physical solution \cite{DIRT2}]\label{def:physsol}
    We say that the càdlàg adapted process $(V_t,M_t)_{t\geqslant 0}$ is a physical solution of \eqref{eq:NNLIF_MF} if
    \begin{itemize}
        \item $(M_t)_{t\geqslant0}$ has integrable marginal distributions and for all $t\geqslant 0$, $\mathbb P(\Delta M_t \leqslant 1)=1$;
        \item \eqref{eq:NNLIF_MF} holds true $\mathbb P$-almost surely with
        \[ M_t =\sum_{l\geqslant 1} \mathds{1}_{[0,t]}\left( \tau_l \right),\qquad \tau_l = \inf\left\{ t > \tau_{l-1} \ | \ V_{t^-} + b \Delta \mathbb E[M_t] \geqslant V_F \right\}, \qquad \tau_0 = 0  ;\]
        \item The discontinuity points of $t\mapsto \mathbb E[M_t]$ satisfy
        \[ \Delta \mathbb E[M_t] = \inf\left\{ \eta\geqslant 0 \ | \ \mathbb P(V_{t^-} + b \eta \geqslant V_F ) < \eta \right\}.\]
    \end{itemize}
\end{definition}
Through the mean-field limit in M1 Skorokhod topology, \cite{DIRT2} proves the existence of a physical solution for any initial condition of law in $L^p$, $p\geqslant 1$, and supported away from $V_F$. Note that for a fixed end-time $T>0$, almost surely
\[\lim_{K\to+\infty}  \dfrac1K \sum_{i=1}^{K} M_t^{i,K} = \mathbb E[M_t]\]
for all $t\in[0,T)$, except at any discontinuity points of $t\mapsto \mathbb E[M_t]$, of which there are countably many.

Physical solutions can also be obtained in the vanishing delay limit $d\to0$, also through convergence in the M1 Skorokhod topology \cite{DIRT2}. In \cite{DIRT1}, with more details in the notes \cite{DIRT3}, it is also proved that, at least up to the first blow-up time, which corresponds to the first discontinuity point of $\mathbb E[M_t]$, the law of a physical solution of \eqref{eq:NNLIF_MF} is a solution in the weak sense of the PDE \eqref{NNLIF} with $d=0$ and $a\equiv a>0$.

\subsubsection{Smooth deterministic solutions from an iteration approach in the linear case}

\label{sec:NNLIF_iteration}

The link between the stochastic differential equation \eqref{eq:NNLIF_MF} and the PDE \eqref{NNLIF} is investigated in depth in \cite{LWZZ} for the simpler linear case $b=0$, $d=0$, $a(N) \equiv a > 0$ constant. They choose for simplicity $V_R=0, V_F=1,a=1$.

The idea of \cite{LWZZ} is to rely on an iteration perspective, considering a collection of processes representing the laws of the successive firing events for a representative neuron. They consider a sequence of independent Ornstein-Uhlenbeck processes $(Y^{n}_t)_{n\in\N}$, each representing the trajectory from $V_R$ to $V_F$ after a spike event and before the next. They introduce a sequence of hitting times $(\tau_n)_{n\in\N}$ at $V_F=1$ for the processes $Y^{n}$ and the jump times for the process $V_t$ are defined as
\[ T_n = \sum_{i=0}^{n}\tau_i .  \]
The original process can then be recovered with the formula 
$V_t = Y^{k}_{t- T_{k-1}}$
on the time interval $[T_{k-1}, T_k)$. The process $n_t=k-1$ keeps track of which process $Y^k$ is being used, resulting in a pair $(V_t,n_t) = (Y^{k}_{t- T_{k-1}},k-1)$. Denoting
\begin{equation} \label{eq:NNLIF_Fk}F_k(v,t) = \mathbb{P}(V_t \leqslant v, n_t = k) , \end{equation}
the cumulative distribution function after $k$ jumps, it is possible to construct an infinite series of PDE systems such that $F_k$ is associated to the solution $\rho_k$ of each of them. The method relies on the key iteration formula
\[ 
F_n (x,t) = \int_{0}^t F_{n-1} (x,t-s) \diff G_1(s),
\] 
where $G_1$ is the cumulative distribution for the first spike time given by
$G_1(t) = \mathbb{P}(\tau_1 \leqslant t)$.
The collection of systems write for $n\geqslant 1$
\begin{equation}\label{eq:NNLIF_iteration}
	\left\{\begin{array}{ll}
		\displaystyle \dfrac{\partial \rho_n}{\partial t} - \dfrac{\partial}{\partial v}\big(v\rho_n\big)-\dfrac{\partial^2 \rho_n}{\partial v^2}= \ \delta_{v=V_R}  N_{n-1}(t), & v\in(-\infty,1],\ t > 0,\\
		\displaystyle N_{n}(t)= -\dfrac{\partial \rho_{n}}{\partial v}(1,t),  \quad \rho_n(1,t) = \rho_n(-\infty,t)=0,& t > 0,\\
		\displaystyle \rho_n(v,0)=0 \geqslant 0,&v\in(-\infty,1].
	\end{array}\right.
\end{equation}
and for the initial density $\rho_0$,
\begin{equation*}
    \dfrac{\partial \rho_0}{\partial t} - \dfrac{\partial}{\partial v}\big(v\rho_n\big)-\dfrac{\partial^2 \rho_0}{\partial v^2} = 0,\quad  N_0(t) = -\dfrac{\partial \rho_{0}}{\partial v}(1,t),\quad \rho_0(-\infty,t)=\rho_0(1,t)=0, 
\end{equation*}
and  $\rho_0(v,0) = \rho^0(v)\geqslant 0$, $\int_{-\infty}^{V_F} \rho^0(v)\diff v = 1$.
The function $\rho^0$ is the probability density of finding a neuron which has never spiked in the network. The density $\rho_n$ is the probability density of neurons which have spiked $n$ times but not $n+1$ times. Each firing rate $N_k$ is a source term in the next equation for $\rho_{k+1}$ whose initial condition is uniformly 0. It is then possible to prove rigorously the decomposition
\begin{equation}
    \rho(v,t) = \sum_{k=0}^{+\infty} \rho_k(v,t).
\end{equation}
Each of these densities $\rho_k$ is associated to the cumulative distribution function $F_k$ in \eqref{eq:NNLIF_Fk}. It can be used to construct strong solutions to the original PDE \eqref{NNLIF} from the stochastic process \eqref{eq:NNLIF_MF}.
\begin{theorem}[Liu, Wang, Zhang, Zhou \cite{LWZZ}]
    Assume $d=0$, $b=0$, $a(N)\equiv 1$, $V_F=1$ and $V_R=0$. If $\rho^0$ is a compatible initial condition compactly supported in $(-\infty,V_F)$ or if $\rho^0 = \delta_{v=V_0}$ is a Dirac mass at a point $V_0\in(-\infty,V_F)$, then the unique solution to \eqref{eq:NNLIF_MF} has a continuously evolved probability density $\rho(\cdot,t)$ which is for positive times a global-in-time strong solution to \eqref{NNLIF}.
\end{theorem}

\subsection{Generalised solutions and continuation after blow-up}
\label{sec:gen_NNLIF}

\subsubsection{Case with internal noise}

Most of the results on the NNLIF model were performed in the case of a constant noise $a(N)\equiv a > 0$. The mean-field limit in \cite{DIRT2} leads to a constant noise coefficient. However, earlier derivations by physicists often propose a noise strength $a(N)=a_0 + a_1 N$, $a_0,a_1 > 0$.
As we have seen in Theorem \ref{thm:blowup_internal_noise}, this internal noise case is more ill-posed as even in the inhibitory case $b < 0$, finite time blow-up can occur. Paradoxically, this parameter range is more amenable for the construction of generalised solutions with PDE tools.

In \cite{DZ}, an approach is proposed that allows both the construction of smooth classical solutions and the extension beyond blow-up points if some conditions are satisfied. The idea is to take advantage of the fact that, although in the constant noise case blow-up is equivalent to $\limsup_{t\to T^*} [-\partial_v \rho(V_F,t)] = +\infty$, the internal noise case allows blow-up to happen while $\partial_v \rho(V_F,t)$ remains bounded. It can be seen by solving for $N$ in $N = - (a_0+a_1N)\partial_v\rho(V_F,t)$:
\begin{equation}
	N(t) = \dfrac{-a_0\partial_v \rho(V_F,t)}{1+a_1 \partial_v \rho(V_F,t)}.
\end{equation}
The main idea of \cite{DZ} is to rescale time in function of the firing rate, in order to turn the divergent speed in the drift term into a locally-bounded one while enjoying uniform parabolicity in the second order term. They propose the change of variable, for some $c > 0$,
\[  \diff \tau = \big(c+N(t)\big)\diff t. \]  
Before applying the change of variable, they substitute $N(t)$ for $ - (a_0+a_1N)\partial_v\rho(V_F,t)$ in the right-hand side, which leads to the $\tau$ time-scale evolution equation
\begin{equation}\label{eq:dilated_NNLIF}
	\dfrac{\partial \tilde\rho}{\partial \tau} + \dfrac{\partial }{\partial v}\left[ (-v\tilde N -cb\tilde N + b )\tilde \rho \right] - \left(a_c\tilde N + a_1\right)\dfrac{\partial^2 \tilde\rho }{\partial v^2} =  - \left(a_c\tilde N + a_1\right) \dfrac{\partial \tilde \rho}{\partial v}(V_F,\tau) \delta_{v = V_R},
\end{equation}
with $a_c = a_0 - c a_1$, $\tilde\rho(-\infty,t) = \tilde \rho(V_F,t) = 0$, and the new firing rate
\begin{equation}\label{eq:dilated_NNLIF_N} \tilde N (\tau) = \dfrac{1}{N+c} = \dfrac{[1+a_1\partial_v \tilde\rho(V_F,\tau)]_+}{-a_0 \partial \tilde\rho(V_F,\tau) + c[1+a_1\partial_v\tilde\rho(V_F,\tau)]_+}.\end{equation}
This modified firing rate is uniformly bounded:
$  0 \leqslant \tilde N(\tau) \leqslant \frac1c$, which implies uniform parabolicity because
\[ \min\left(a_1,\frac{a_0}{c}\right) \leqslant a_c \tilde N(\tau) + a_1 \leqslant \max\left(a_1,\frac{a_0}{c}\right). \]

The original time can then be recovered \textit{via} the formula
$ t(\tau) = \int_0^\tau \tilde N(\eta) \diff\eta$. Adapting the method of \cite{CGGS}, \cite{DZ} constructs global-in-time solutions to the time-dilated problem.
\begin{theorem}[Dou, Zhou \cite{DZ}]
    For any compatible initial condition $\rho^0$, any fixed $c$, there exists a unique global-in-time strong solution to \eqref{eq:dilated_NNLIF}--\eqref{eq:dilated_NNLIF_N}.
\end{theorem}
It is then possible to recover a $c$-independent unique generalised solution for the original time problem \eqref{NNLIF} which is defined up to time
\begin{equation}
    T^* = \int_0^{+\infty} \tilde N(\tau) \diff \tau.
\end{equation}
The generalised firing rate in original time $t$ is then
\begin{equation}
    N(t) = \dfrac{1}{\tilde N(\tau(t))} - \dfrac1c \,\in[0,+\infty], \qquad \tau(t) = \sup\left\{ \tau\in(0,+\infty)\ |\ \int_0^\tau \tilde N(\eta)\diff \eta = t  \right\}.
\end{equation}

\begin{theorem}[Dou, Zhou \cite{DZ}]
    Assume $d=0$, $a(N)=a_0+a_1 N$, $a_0,a_1  > 0$. For any compatible initial condition $\rho^0$, there exists a unique generalised solution to \eqref{NNLIF} defined on the maximal interval $[0,T^*)$.
    \begin{itemize}
        \item If $b < 0$ and $V_F\geqslant 0$, then $T^*=+\infty$.
        \item If $0 < b < V_F - V_R$, $\rho^0(v) e^{-\frac{b}{a_1}v}$ is bounded and $V_F\geqslant -\frac{a_0}{a_1}\, b$, then $T^* = +\infty$.
        \item If $ V_F-V_R < b $, there exist initial conditions such that $T^* < +\infty$.
    \end{itemize}
\end{theorem}
The case where generalised solutions cannot be global-in-time in the high connectivity case is called an eternal blowup. In the time-dilated problem \eqref{eq:dilated_NNLIF}, the solution never exits blow-up regime and $\tilde N$ stays zero indefinitely. Note that this nonphysical case was excluded by design in the mean-field limit of \cite{DIRT2}, as they assume $b < V_F - V_R$. In the other cases with weaker or inhibitory interaction, the idea of the proof is to use an entropy method in time-dilated timescale to prove that if the solution enters blow-up, it must exit in finite time. A crucial ingredient is to use, the fact that if $T^*$ is assumed by contradiction to be finite, then $\tilde N$ must be small for large times, enabling weakly-nonlinear analysis.

Those generalised solutions come from specific choices: how to rewrite the $N(t)\delta_{v=V_R}$ source term and how to do the dynamics inside blow-up. Other choices could be made, leading to other types of generalised solutions. The approach in \cite{DZ} enjoys many good mathematical properties, but comes with the caveat of allowing neurons to spike multiple times during a blow-up event.

\subsubsection{Case without internal noise}
\label{sec:NNLIF_gen_nonoise}

When the diffusion is constant $a(N)\equiv a$, the method presented above no longer works because the equation in a dilated timescale is no longer parabolic. Moreover, in the case without internal noise, the Dirichlet boundary condition at $V_F$ is not satisfied during blow-up. To overcome these issues, an other approach that involves time dilation was recently developed in \cite{DPSZ2}.
Before blow-up, the problem \eqref{NNLIF} can be recast on the whole real line by putting a $-N(t)\delta_{v=V_F}$ in the right-hand side. In order to generalise this idea after blow-up, a more general measure $\mathcal S$ has to be considered. More precisely, introduce
\[ \diff \tau = N(t)\diff t, \qquad n(v,\tau) = \rho(v,t),\qquad Q(\tau) = \frac1{N(t)}.  \]
Then, $(n,Q,\mathcal S)$ satisfies
\begin{equation}\label{eq:NNLIF_gen_1}
    \dfrac{\partial n}{\partial \tau} + \dfrac{\partial}{\partial v}[ (-Q(\tau)+b)n  ] - aQ(\tau)\dfrac{\partial^2 n}{\partial v^2} = \delta_{v=V_R} - \mathcal S(v,\tau), \qquad \tau\geqslant 0,\quad v\in\R.
\end{equation}
In the classical regime when $\int_{V_F}^{+\infty} n(v,\tau)\diff v = n(V_F,\tau) = 0$, then $\mathcal S(v,\tau) = \delta_{v=V_F}$ and $$-a\partial_v n(V_F,\tau)Q(\tau)=1.$$ In the blow-up regime when $\int_{V_F}^{+\infty} n(v,\tau)\diff v > 0$, $Q(\tau)=0$. The other conditions for the generalised boundary condition write
\begin{equation}\label{eq:NNLIF_gen_2}
    \mathcal S(\cdot,\tau) \geqslant 0,\qquad \int_\R \mathcal S(\cdot,\tau)(\diff v) = 1,\qquad \mathcal S(\cdot,\tau)\equiv 0 \ \mathrm{on}\ (-\infty,V_F).
\end{equation}
Whenever $Q(\tau)>0$, a solution to \eqref{eq:NNLIF_gen_1}--\eqref{eq:NNLIF_gen_2} and related conditions corresponds to a solution to \eqref{NNLIF}. In order to obtain a generalised solution in this context, \cite{DPSZ2} uses the NNLIF model with random firing \eqref{NNLIFRand} presented below in Subsection \ref{sec:NNLIF_RF}, with the choice $\phi_\varepsilon = \frac1\varepsilon \mathds 1_{v\geqslant V_F}$. Working in the dilated timescale $\tau$ allows to obtain $\varepsilon$-uniform estimates and pass to the limit $\varepsilon\to 0$.

Given an interval $I\subset\R$, consider $\mathcal M(I)$ the space of bounded measures on $I$, $\mathcal M_+(I)\subset \mathcal M(I)$ the space of positive bounded measures on $I$ and $\mathcal P_2(I)\subset \mathcal M_+(I)$ the space of probability measures with second moment.
\begin{theorem}[Dou, Perthame, Salort, Zhou \cite{DPSZ2}]
    Assume the initial condition $n^0\in L^2(\R)\cap \mathcal P_2(\R)$. Then, there is a global-in-time weak solution to \eqref{eq:NNLIF_gen_1}--\eqref{eq:NNLIF_gen_2} such that $n\in L^\infty(0,+\infty;L^1(\R))\cap L^\infty_{\mathrm{loc}}(0,+\infty;L^2(\R))$ and for all $\tau_0 > 0$, in $(0,\tau_0)$ we have $Q\in \mathcal M_+(0,\tau_0)$, $\int_\R v^2 n(v,\tau)\diff v \leqslant C(\tau_0)$. Moreover, for all $\psi\in L^2(\R)\cap L^\infty(\R)\cap C(\R)$, $\tau\mapsto \int_\R \psi(v)n(v,\tau)\diff v\in C(0,+\infty)$ and
    \[  Q(\tau)\int_{V_F}^{+\infty} n(v,\tau)\diff v = 0. \]
\end{theorem}
Similar to the generalised solution in the internal noise case described above, the generalised solution $\rho$ defined on $[0,T^*)$ in original timescale can be global-in-time ($T^*=+\infty$) or enter an eternal blow-up ($T^*<+\infty$). Again we have
\[  T^* = \int_0^{+\infty} Q(\tau)\diff \tau.  \]
When $b<V_F-V_R$, \cite{DPSZ2} proves that $T^*=+\infty$; when $b>V_F-V_R$, examples of eternal blow-up are found, which correspond to the plateau states identified in \cite{CR-L} that we have described above.

\subsection{Asymptotic behaviour for general connectivity strength and delay}
\label{sec:NNLIF_advances}

In most previous study, analytical results on the asymptotic behaviour of \eqref{NNLIF} required smallness assumption on both the connectivity parameter $b$ and the delay parameter $d$. In two companion papers \cite{CCR-L,CCR-L2}, new methods were proposed in order to study the long-time behaviour for large $b$ and/or large $d$. The authors fix $a(N)\equiv 1$ for convenience, but claim their strategy could be extended to any constant $a(N)\equiv a>0$ with only technical changes. For simplicity of notations in their approach, they also consider initial conditions of the form
\begin{equation}
    \rho^0(v,t), \qquad N^0(t) = \dfrac{\partial \rho^0}{\partial t}(V_F,t),\qquad v\in(-\infty,V_F],\ t\in[-d,0].
\end{equation}
from which initial conditions in the sense of Definition \ref{Assumption1} can easily be retrieved provided enough regularity and decay at infinity.

\subsubsection{Improved convergence estimates with a spectral gap approach}
\label{sec:NNLIF_general_b}

In order to better understand the nonlinear system \eqref{NNLIF} in the case of strong nonlinearities, a possibility is to first look at a linearised version of the firing rate impact on the network. Fix an external firing rate $N\geqslant 0$, a connectivity parameter $b\in\R$ and define the linear operator
\begin{equation*}
    L_N \rho = \dfrac{\partial^2 \rho}{\partial v^2} + \dfrac{\partial}{\partial v}[(v-bN)\rho] + \delta_{v=V_R} N_{\rho}, \qquad N_{\rho} = -\dfrac{\partial \rho}{\partial v}(V_F,t),
\end{equation*}
naturally associated to the linear PDE    $\partial_t \rho = L_N\rho$, $\rho(V_F,t) = 0$, $\rho(\cdot,0) = \rho^0(\cdot)$.
Now write $\rho = \rho_\infty + u$ as a perturbation $u$ from a stationary state $(\rho_\infty, N_\infty)$ and denote $L_\infty = L_{N_\infty}$. Then $u$ solves the equation
\begin{equation*}\label{eq:NNLIF_u}
    \dfrac{\partial u}{\partial t} = L_\infty u - b N_u(t-d)\dfrac{\partial \rho_\infty}{\partial v} - bN_u(t-d) \dfrac{\partial u}{\partial v},
\end{equation*}
which yields, ditching the nonlinear term for small $u$, the linearised problem
\begin{equation}\label{eq:NNLIF_lin}
    \dfrac{\partial u}{\partial t} = Tu, \qquad Tu =  L_\infty u - b N_u(t-d)\dfrac{\partial \rho_\infty}{\partial v}, \qquad u(-\infty,t) = u(V_F,t) = 0.
\end{equation}
Because of the structure of the stationary states, the need to control the derivative at the boundary and the framework for classical well-posedness, a natural space for the study of strong solutions is
\begin{equation}
    X = \left\{ \rho \in C(-\infty,V_F]\cap C^1(-\infty,V_R]\cap C^1[V_R,V_F] \ \ | \ \ \rho(V_F)=0, \ \norme{\rho}_X < +\infty  \right\},
\end{equation}
where
\begin{equation}
    \norme{\rho}_X = \norme{\rho}_\infty + \norme{\partial_v \rho}_\infty + \norme{\rho}_{L^2(\phi)} + \norme{\partial_v\rho}_{L^2(\phi)}, \qquad \phi(v) = e^{-\frac{(v-bN)^2}{2}}.
\end{equation}
It can be readily checked that this space embeds into the natural space for relative entropy estimates: $\norme{\rho}_{L^2(\rho_\infty^{-1})}\leqslant C \norme{\rho}_{X}$ and $X\subset L^2(\rho_\infty^{-1})$. Note that the spectrum of the linear problem in $L^2(\rho_\infty^{-1})$ had been studied before \cite{CCR-L}, although less thoroughly, in \cite{CGGS}.

\begin{definition}
    We say that a stationary state $(\rho_\infty,N_\infty)$ is linearly stable if there exists $C \geqslant 1$ and $\lambda > 0$ such that all solutions $u$ to \eqref{eq:NNLIF_lin} from an initial condition $u_0\in C([-d,0],X)$ such that $\int_{-\infty}^{V_F}u_0(v,\cdot)\diff v \equiv 0$ on $[-d,0]$, satisfy
    \begin{equation}
        \norme{u(\cdot,t)}_X \leqslant C e^{-\lambda t} \sup_{\tau\in[-d,0]} \norme{u_0(\cdot,\tau)}_X.
    \end{equation}
    It is linearly unstable if for any $C\geqslant 1$, there exists $u_0\in C([-d,0],X)$ and $T\geqslant 0$ such that
    \begin{equation}
        \norme{u(\cdot,T)}_X \geqslant C \sup_{\tau\in[-d,0]} \norme{u_0(\cdot,\tau)}_X.
    \end{equation}
\end{definition}

Given this definition, the next step is to study the evolution along the linear flow of the firing rate of the small perturbation $u$, which can be done by introducing the functions
\begin{equation}
    q(v,t) = e^{ t L_\infty} \dfrac{\partial \rho_\infty}{\partial v},\qquad N_q(t) = - \dfrac{\partial q}{\partial v}(V_F,t),
\end{equation}
and the Laplace transform
\begin{equation}\label{eq:NNLIF_laplace}  \hat N_q(\xi) = \int_{0}^{+\infty} e^{-\xi t} N_q(t)\diff t, \end{equation}
which \cite{CCR-L} proves to be well-defined for a large enough real part $\mathcal R(\xi) > -\lambda$, as there exist $C\geqslant 1,\lambda >0$ such that $\hat N_q(t) \leqslant C e^{-\lambda t}$.

The cornerstone of the use of this Laplace transform is the following result.
\begin{theorem}[C\'aceres, Cañizo, Ramos-Lora \cite{CCR-L}]
    Let $(\rho_\infty,N_\infty)$ be a stationary state and associated $\hat N_q$ as in \eqref{eq:NNLIF_laplace}; it is linearly stable if and only if all zeros of the analytic function
    \begin{equation}
        \Phi_d(\xi) = 1 + b \hat N_q(\xi) e^{-\xi d}
    \end{equation}
    are located on the real-negative half plane $\{\xi\in\mathbb C \ | \ \mathcal R(\xi) < 0\}$.
\end{theorem}
From this characterisation, a surprising consequence is that the linear stability of stationary states of \eqref{NNLIF} can be related to the real function $I(N)$ used for the study of stationary states.

\begin{theorem}[C\'aceres, Cañizo, Ramos-Lora \cite{CCR-L}]
    Let $(\rho_\infty,N_\infty)$ be a stationary state and associated $\hat N_q$ as in \eqref{eq:NNLIF_laplace};
    \begin{itemize}
        \item $\rho_\infty$ is linearly unstable for all delay $d\geqslant 0$ under the sufficient condition
        \[\dfrac{\diff}{\diff N}\left(\dfrac{1}{I(N)}\right)_{|N=N\infty} > 1;\]
        \item $\rho_\infty$ is linearly stable for all delay $d\geqslant 0$ under the sufficient condition
        \[  |b| \int_0^{+\infty} |N_q(t)|\diff t < 1 ;\]
        \item $\rho_\infty$ is linearly unstable if the delay $d$ is large enough  under the sufficient condition
        \[\dfrac{\diff}{\diff N}\left(\dfrac{1}{I(N)}\right)_{|N=N\infty} < -1;\]
    \end{itemize}
    where the function $N\mapsto I(N)$ is the one defined in \eqref{eq:Fix_N_Stat}.
\end{theorem}
This classification is not exhaustive theoretically, but numerical results show that it is almost comprehensive in practice, because $N_q$ has a fixed sign, leading to 
$$
|b| \int_0^{+\infty} |N_q(t)|\diff t = \left|\frac{\diff}{\diff N}\left(\frac{1}{I(N)}\right)(N_\infty)\right|
$$
and the limit case $=\pm 1$ is marginal. From these linear stability results, \cite{CCR-L} draws rigorous consequences for the nonlinear system \eqref{NNLIF}.

First, the weakly-nonlinear stability performed with relative entropy in previous studies \cite{CCP,CPSS,CRSS,CS,CS2} can be performed in the smaller space $X$ with improved assumptions on $d$.
\begin{theorem}[C\'aceres, Cañizo, Ramos-Lora \cite{CCR-L}]\label{thm:NNLIF_CCR-L_small_b}
    Assume $a(N)\equiv 1$ and $|b|$ is small enough to have a unique stationary state $(\rho_\infty,N_\infty)$. There exists $C\geqslant 1,\lambda>0$ depending only on $V_R,V_F$ such that for all initial condition $\rho^0\in C([-d,0],X)$, there exists $b^*>0$ depending only on $V_R,V_F$ and 
    \begin{equation}\sup_{t\in[-d,0]} \norme{\rho^0(\cdot,t)-\rho_\infty(\cdot)}_X,
    \end{equation} 
    such that for all $b\in(-b^*,b^*)$, the solution $\rho$ to \eqref{NNLIF} satisfies
    \begin{equation}
        \forall t\geqslant 0,\qquad \norme{\rho(\cdot,t)-\rho_\infty(\cdot)}_X \leqslant C e^{-\frac{\lambda}{2+d}t} \norme{\rho^0(\cdot,0)-\rho_\infty(\cdot)}_X.
    \end{equation}
\end{theorem}
The small connectivity limit $b^*$ depends on the initial condition $\rho^0$, but is independent on the delay value $d$, significantly improving the earlier findings of \cite{CRSS} (see Theorem \ref{thm:NNLIF_d_stat} above).

Second, it can be proved that linear stability entails nonlinear stability beyond the weakly-nonlinear case when there is no delay $d=0$.
\begin{definition}
    We say that the flow $(e^{tT})_{t\geqslant 0}$ of the linearised problem \eqref{eq:NNLIF_lin} has a spectral gap $\lambda$ with constant $C\geqslant 1$ if for all $u_0\in X$ such that $\int_{-\infty}^{V_F} u_0(v)\diff v$, the solution $e^{tT}u_0$ to \eqref{eq:NNLIF_lin} satisfies
    \[  \norme{e^{tT}u_0}_X \leqslant C e^{-\lambda t} \norme{u_0}_X.  \]
\end{definition}

\begin{theorem}[C\'aceres, Cañizo, Ramos-Lora \cite{CCR-L}]
    Assume $a(N)\equiv 1$ and $d=0$. Let $(\rho_\infty,N_\infty)$ be a stationary state.
    Assume  the flow of the linearised problem \eqref{eq:NNLIF_lin} has a spectral gap of size $\lambda>0$ with constant $C\geqslant 1$.
    There exists $\varepsilon>0$ depending on $b,V_F,V_R,\rho_\infty,\lambda$ such that for any initial condition $\rho^0\in X$, if $\norme{\rho^0-\rho_\infty}_X<\varepsilon$, then the solution $\rho$ to \eqref{NNLIF} satisfies
        \begin{equation}
        \forall t\geqslant 0,\qquad \norme{\rho(\cdot,t)-\rho_\infty(\cdot)}_X \leqslant 2C e^{-\frac{\lambda}{2}t} \norme{\rho^0-\rho_\infty}_X.
    \end{equation}
\end{theorem}
The authors of \cite{CCR-L} claim that the result should hold for non-zero delay $d > 0$ and that an extension of the result is mainly a technical challenge.

\subsubsection{The discrete sequence of pseudo-equilibria}\label{subsec:pseudo}

In a companion paper \cite{CCR-L2} to \cite{CCR-L}, the troubling link between the function $I(N)$ defined in \eqref{eq:Fix_N_Stat} and the stability of stationary states is studied further with the introduction of a sequence of pseudo-equilibria.

Consider the simple case $a\equiv 1$. Their main idea is to remark that when the delay $d$ is large enough, the system behaves as a linear equation with a given external drift on each time interval of size $d$. Given a fixed firing rate $N\geqslant 0$ acting on the network as an eternal source, the solution of \eqref{NNLIF} would converge exponentially fast to the probability density
\begin{equation}\label{eq:NNLIF_pseudo}
    \rho_{\mathrm{pseudo}}(v) = \tilde N e^{-\frac{(v-bN)^2}{2}} \int_{\max(v,V_R)}^{V_F}e^{\frac{(w-bN)^2}{2}}\diff w,\qquad \tilde N = \dfrac{1}{I(N)},
\end{equation}
which is a stationary state of the full nonlinear system \eqref{NNLIF} if and only if $\tilde N = N$ or equivalently if $N I(N) = 1$. If the firing rate was $N_{0,\infty}$ on $[-d,0]$, then on $[0,d]$ with $d$ large enough, the firing rate $N(t)$ of a solution to \eqref{NNLIF} has the time to approach a new equilibrium value $N_{1,\infty} \simeq -\partial_v \rho(V_F,d)$ given by \eqref{eq:NNLIF_pseudo}. One can thus define a pseudo-equilibria sequence in order to approach the behaviour of a solution on a sequence of time intervals $[kd,(k+1)d]$.

\begin{definition}
    We call the \emph{firing rate sequence} associated to the initial firing rate $N_{0,\infty}$ the sequence $(N_{k,\infty})_{k\geqslant0}$ recursively defined by
    \begin{equation}
        N_{k+1,\infty} = \dfrac{1}{I(N_{k,\infty})}.
    \end{equation}
    The associated \emph{pseudo-equilibria sequence} $(\rho_{k,\infty})_{k\geqslant 0}$ is defined by
    \begin{equation}
        \rho_{k,\infty} = N_{k,\infty} e^{-\frac{(v-bN_{k-1,\infty})^2}{2}} \int_{\max(v,V_R)}^{V_F}e^{\frac{(w-bN_{k-1,\infty})^2}{2}}\diff w.
    \end{equation}
\end{definition}

A careful analysis of the function $\frac{1}{I(N)}$ allows to do an almost comprehensive classification of the asymptotic behaviour of the firing rate sequence.
\begin{theorem}[C\'aceres, Cañizo, Ramos-Lora \cite{CCR-L2}]
    Assume $a(N)\equiv 1$.
    \begin{itemize}
        \item In the excitatory case $b > 0$:
            \begin{itemize}
        \item if there is a unique stationary firing rate $N_\infty$ and $N\mapsto \frac1{I(N)}$ crosses the diagonal $N\mapsto N$, $(N_{k,\infty})_{k\geqslant0}$ converges monotonically to $N_\infty$;
        \item if there is a unique stationary firing rate $N_\infty$ and $N\mapsto \frac1{I(N)}$ touches tangentially the diagonal, $(N_{k,\infty})_{k\geqslant0}$ converges to $N_\infty$ when $N_{0,\infty} < N_\infty$ and diverges when $N_\infty < N_{0,\infty}$;
        \item if there are two stationary firing rates $N_\infty^1 <N_\infty^2$, $(N_{k,\infty})_{k\geqslant0}$ monotonically converges to $N_\infty^1$ when $N_{0,\infty} < N_\infty^2$ and diverges when $N_\infty^2 < N_{0,\infty}$;
        \item if there is no stationary state, $(N_{k,\infty})_{k\geqslant0}$ diverges.
    \end{itemize}
        \item In the inhibitory case $b<0$: there exists a critical connectivity $b^* < 0$ such that
        \begin{itemize}
            \item if $b^* < b < 0 $, $(N_{k,\infty})_{k\geqslant0}$ converges monotonically to the unique stationary firing rate $N_\infty$;
            \item if $ b < b^* $, there exists $N^-,N^+$ such that $N_\infty\in (N^-,N^+)$ and $(N_{k,\infty})_{k\geqslant0}$ converges to the 2-cycle $\{N^-,N^+\}$.
        \end{itemize}
    \end{itemize}
\end{theorem}

Across a wider parameter range, numerical simulations allow to check that the behaviour of the nonlinear system \eqref{NNLIF} asymptotically follows the sequence of pseudo-equilibria when the delay is large enough, $d=10$ being enough in many cases.

In the excitatory case $b > 0$, when there is a unique stationary state $(\rho_\infty,N_\infty)$, in most cases the sequence of firing rates converges to $N_\infty$. Starting from a pseudo-equilibrium as an initial condition, the nonlinear system for positive delays (hence, no finite-time blowup) will numerically converge to the stationary state. When $b$ is large enough to have two stationary states $(\rho_\infty^1,N_\infty^1)$ and $(\rho_\infty^2,N_\infty^2)$ with $N_\infty^1 < N_\infty^2$, the assymptotic behaviour of \eqref{NNLIF} also agrees with the sequence of pseudo-equilibria:
\begin{itemize}
    \item if the initial firing rate is above $N_\infty^2$, there is divergence in infinite time $\lim_{t\to+\infty}N(t) = +\infty$ and the density converges to a plateau state, \textit{i.e.} a uniform distribution between $V_R$ and $V_F$;
    \item if the initial firing rate is below $N_\infty^2$, there is convergence to the stationary state with lower firing rate $(\rho_\infty^1,N_\infty^1)$.
\end{itemize}
In the case where there is no stationary state, both the sequence of pseudo-equilibria and the solutions of \eqref{NNLIF} converge towards the plateau state.
Note that in the excitatory case, the delay does not matter unless it is $d=0$ (then the technical problem of finite-time blow-up arises) and only the initial condition and the multiplicity of stationary states matter. The behaviour of the sequence is a good predictor of the nonlinear asymptotics.

In the inhibitory case $b < 0$, the initial condition is irrelevant, but the value of the delay matters. Fix, for example, $V_F=1, V_R=0, a\equiv 1$. For small enough $b$ and any delay $d\geqslant 0$, there is convergence of the firing rate sequence to the unique stationary firing rate $N_\infty$ and convergence of the nonlinear solution of \eqref{NNLIF} towards the unique stationary state. There is a critical value $b^*\simeq -9.4$ such that when $b < b^*$, the sequence of firing rates converges to a 2-cycle $\{N^-,N^+\}$ and the sequence of pseudo-equilibria to a 2-cycle $\{\rho^-,\rho^+\}$. Then, the solution of \eqref{NNLIF} follows the same behaviour only when the delay $d$ is large enough. For example, with the parameters specified, $(\rho,N)$ converges to the unique stationary state $(\rho_\infty,N_\infty)$ when $d=2$ but to a periodic solution when $d=10$. The only impact of the initial condition is which state in the two-cycle is reached first.

In the case of a small enough connectivity, the link between the sequence of pseudo-equilibria and the stationary state is rigorously shown in \cite{CCR-L2}. The idea is to introduce, given a solution $(\rho,N)$ to \eqref{NNLIF}, 
\[  \rho_{k+1}(v,t) = \rho(v, t + kd),\qquad N_k(t) = N(t+kd),  \]
which allows to rewrite \eqref{NNLIF} in the time interval $(kd,(k+1)d)$ as
\begin{equation}\label{eq:delayres}
    \dfrac{\partial \rho_k}{\partial t} + \dfrac{\partial}{\partial v}\left[(-v + bN_{k-1}(t-d))\rho_k\right] -\dfrac{\partial^2 \rho_k}{\partial v^2} = \delta_{v=V_R} N_k(t).
\end{equation}
The stationary states to \eqref{eq:delayres} correspond to the sequence of pseudo-equilibria. Assuming a spectral gap and a regularisation property on the linearised system, \cite{CCR-L2} proves under the functional analysis artillery of \cite{CCR-L} that there exists a connectivity $b^*>0$, a delay $d^*>0$ and constants $Q>0,\lambda>0$, all depending on $V_R,V_F$ and 
$ \sup_{t\in[0,d]} \norme{\rho(\cdot,t) }_X  $
such that for all $ -b^* < b < b^*$ and $d > d^*$,
\begin{equation}
    \norme{\rho(\cdot,t)}_X \leqslant Q e^{-\mu t} \norme{\rho^0-\rho_\infty}_X. 
\end{equation}

Note that this result was already proved in previous literature \cite{CPSS,CRSS,CCR-L}. However, the relationship with the sequence of pseudo-equilibria sheds new light onto the problem. Moreover, in \cite{CRSS}, the smallness condition for $b$ was given by the constraint $N_\infty b^2 e^{\lambda d} \leqslant C$, where $C,\lambda$ are functional analysis constants, but with the pseudo-equilibria method it is replaced by the $d$-independent condition $b^4\leqslant C \lambda$, where $C,\lambda$ are related to the spectral gap and other functional analysis constraints.

\subsection{Variants of the NNLIF model}

The system \eqref{NNLIF} is a simplistic representation of a large network of neurons. Some behaviours of the solutions, such as convergence to a plateau state for a large connectivity $b> 0$ and positive delay $d>0$, are not realistic features in neuroscience. Many studies have proposed and investigated variants of the original model displaying additional features or modified equations. It allows in particular to find the minimal structures necessary to observe more realistic neuroscience patterns.

\subsubsection{Adding a refractory state}
\label{sec:NNLIF_R}

In \cite{CP}, the addition to \eqref{NNLIF} of a refractory state was rigorously studied. This models an important feature of neurons: after emitting an action potential, they cannot easily spike shortly after and have to enter a state of recovery. Mathematically, we can represent it by replacing the reset term $\delta_{v=V_R}N(t)$ by a source from an external pool of neurons $R(t)$ which are in a refractory state. Note that previous physics-oriented studies like \cite{BH,B3} had already proposed this idea and had done preliminary numerical studies. The NNLIF model with refractory period reads
\begin{equation}\label{NNLIFR}
	\left\{\begin{array}{l}
		\displaystyle \dfrac{\partial \rho}{\partial t} + \dfrac{\partial}{\partial v}\Big[(-v+bN(t))\rho\Big]-a\big(N(t)\big)\dfrac{\partial^2 \rho}{\partial v^2}= \ \delta_{v=V_R}  M(t),\\ 
		\displaystyle N(t)= -a(N(t)) \dfrac{\partial p}{\partial v}(V_F,t),\\
		\dfrac{\diff R}{\diff t}(t) = N(t) - M(t), \quad R(0)=R^0\geqslant 0,\\
		\displaystyle \rho(V_F,t)= \rho(-\infty,t)=0,\quad
		\displaystyle \int_{-\infty}^{V_F} \rho^0(v)\diff v + R^0 = 1, \quad  \rho(v,0)=\rho^0(v) \geqslant 0.
	\end{array}\right.
\end{equation}
The choice in \cite{B3} was for example $M(t) = N(t-\tau_0)$, meaning that neurons firing at time $t$ are ready to fire again after a fixed time $\tau_0$. In \cite{CP}, the authors opt for the smoother
\[  M(t) = \dfrac{R(t)}{\tau},\qquad a(N)\equiv a>0.  \]
Weak and strong solutions can be defined similar to Section \ref{sec:NNLIF_math} for \eqref{NNLIF}. The system enjoys the mass conservation property
\begin{equation}
		 \int_{-\infty}^{V_F} \rho(v,t)\diff v + R(t) = \int_{-\infty}^{V_F} \rho^0(v)\diff v + R^0 = 1.
\end{equation}
If $b>0$,  Theorem \ref{ExploPert} still holds with a similar proof, which indicates that a refractory period, unlike a time delay, is not enough to prevent finite-time blow-up. There is finite time blow-up when, assuming $a(N) > a_m > 0$,
\begin{itemize}
    \item for a fixed $\rho^0$, when $b>0$ is large enough;
    \item for a fixed $b>0$, when
    \begin{equation}
        b\mu \int_{-\infty}^{V_F} \rho^0(t)e^{\mu v}\diff v > e^{\mu V_F} , \qquad \mu > 2\max\left(\frac{V_F}{a_m},\frac1b\right)
    \end{equation}
\end{itemize}
Stationary states of \eqref{NNLIFR} are similar to those of \eqref{NNLIF} but with $R_\infty = \tau N_\infty$. However, the classification of their existence is strikingly different.
\begin{theorem}[C\'aceres, Perthame \cite{CP}]
	Assume $M(t) = \frac{R(t)}{\tau}$, $a(N)\equiv a>0$; then
	\begin{itemize}
		\item[$\bullet$] if $b\leqslant0$ there is a unique stationary state of \eqref{NNLIFR} ;
		\item[$\bullet$] if $b>0$ there is an odd\footnote{Up to multiplicity. There are a finite number of values $b$ where $1/I(N)$ touches tangentially the diagonal, with only two distinct stationary states, often during transitions from 1 to 3 stationary states. } number of stationary states of \eqref{NNLIFR}.
	\end{itemize}
\end{theorem}

In contrast with the standard NNLIF model, there is always at least one stationary state, even for very high connectivities. The refractory state removes the possibility of the plateau state, which was a form of singular stationary state for \eqref{NNLIF}. For system \eqref{NNLIFR}, there are in practice one or three stationary states (two in special cases; see footnote). The stationary density is still described by \eqref{eq:NNLIF_rho_stat}, but the stationary firing rate is now solution of $N(I(N)+\tau) = 1$ with the same function $I$ as in \eqref{eq:Fix_N_Stat}.

The natural relative entropy for \eqref{NNLIFR} given a stationary state $(\rho_\infty,N_\infty, R_\infty)$ is
\begin{equation}
    \mathcal H(\rho|\rho_\infty)(t) = \int_{-\infty}^{V_F} \left(\dfrac{ \rho(v,t)-\rho_\infty(v)}{R_\infty}\right)^2 \rho_\infty(v) \diff v + \left(\dfrac{R(t)-R_\infty}{R_\infty}\right)^2 R_\infty.
\end{equation}
It is used in \cite{CP} to prove exponential convergence to the unique stationary state in the linear case.
\begin{theorem}[C\'aceres, Perthame \cite{CP}]
    Assume $b=0$, $a(N)\equiv a>0$. There exists $\mu>0$ such that if $\mathcal H(\rho^0|\rho_\infty)<+\infty$, the solution to \eqref{NNLIFR} satisfies
    \[ \mathcal H(\rho|\rho_\infty)(t) \leqslant e^{-\mu t}\mathcal H(\rho^0|\rho_\infty).  \]
\end{theorem}
The result was extended to the weakly-nonlinear case in \cite{CS2} for solutions starting close enough to the stationary state in relative entropy and constant diffusion $a(N)\equiv a>0$.
\begin{theorem}[C\'aceres, Schneider \cite{CS2}]\label{thm:NNLIFr_CV}
    Assume $|b|$ is small enough and $a(N)\equiv a>0$. There exists $\mu>0$ such that if $\mathcal H(\rho^0|\rho_\infty)\leqslant \frac{1}{2|b|}$, the solution to \eqref{NNLIFR} satisfies
    \[ \mathcal H(\rho|\rho_\infty)(t) \leqslant e^{-\mu t}\mathcal H(\rho^0|\rho_\infty).  \]
\end{theorem}

\subsubsection{NNLIF model with random firing}
\label{sec:NNLIF_RF}

As we have seen above, a positive delay coefficient $d>0$, representing travel of action potential along axons and synaptic integration, prevents finite-time blow-up and replaces it either by a sharp peak of the firing rate followed by convergence to a stationary state or by an unphysical plateau state. A smoother way to prevent finite-time blow-up is to reduce the deterministic firing of neurons when they reach the firing potential $V_F$ by randomizing firing events. The model is then extended to the whole real line $\R=(-\infty,+\infty)$ and the discharge of neurons is represented by the decay term $\phi_\varepsilon(v)\rho(v,t)$, where the function $\phi_\varepsilon(v)$ is uniformly equal to zero on $(-\infty,V_F)$ and takes positive values on $(V_F,+\infty)$. To our knowledge, the first rigorous study was performed in \cite{CP}, where the authors propose the model
\begin{equation}\label{NNLIFRand}
	\left\{\begin{array}{l}
		\displaystyle \dfrac{\partial \rho}{\partial t} + \dfrac{\partial}{\partial v}\Big[(-v+bN(t))\rho\Big]-a\big(N(t)\big)\dfrac{\partial^2 \rho}{\partial v^2}  + \phi_\varepsilon(v)\rho  = N(t) \delta_{v=V_R}  ,\qquad v\in\R, \\ 
		\displaystyle N(t)= \int_{-\infty}^{+\infty} \phi_\varepsilon(v)\rho(v,t)\diff v,\\
		\displaystyle \rho(+\infty,t)= \rho(-\infty,t)=0,\quad
		\displaystyle \int_{-\infty}^{+\infty} \rho^0(v)\diff v = 1, \quad  \rho(v,0)=\rho^0(v) \geqslant 0,
	\end{array}\right.
\end{equation}
along with a second version with refractory period:
\begin{equation}\label{NNLIFRandR}
	\left\{\begin{array}{l}
		\displaystyle \dfrac{\partial \rho}{\partial t} + \dfrac{\partial}{\partial v}\Big[(-v+bN(t))\rho\Big]-a\big(N(t)\big)\dfrac{\partial^2 \rho}{\partial v^2}  + \phi_\varepsilon(v)\rho  = \dfrac{R(t)}{\tau} \delta_{v=V_R}  ,\qquad v\in\R, \\[0.3cm] 
		\displaystyle N(t)= \int_{-\infty}^{+\infty} \phi_\varepsilon(v)\rho(v,t)\diff v,\quad
		\dfrac{\diff R}{\diff t}(t) = N(t) - \dfrac{R(t)}{\tau}, \quad R(0)=R^0\geqslant 0,\\
		\displaystyle \rho(+\infty,t)=\rho(-\infty,t)=0,\quad
		\displaystyle \int_{-\infty}^{+\infty} \rho^0(v)\diff v + R^0 = 1, \quad  \rho(v,0)=\rho^0(v) \geqslant 0.
	\end{array}\right.
\end{equation}
Since the obstruction for global-in-time existence of solutions to \eqref{NNLIF} is finite-time divergence of the firing rate $N(t)$, the article \cite{CP} provides arguments for the absence of blow-up in \eqref{NNLIFRand} by proving locally uniform bounds for the firing rate.
\begin{theorem}[C\'aceres, Perthame \cite{CP}]
	Assume $a(N)\leqslant a_0 + a_2N^2$ for some $a_0,a_2\geqslant 0$, $\phi_\varepsilon(v) = \frac{1}{\varepsilon} (v-V_F)_+$ and
	\[   \int_{-\infty}^{+\infty} (1+|v|^3)\rho^0(v)\diff v < +\infty.   \]
	Then, there exists $C > 0$ such that any weak solution of \eqref{NNLIFRand} or of \eqref{NNLIFRandR} satisfies
	\[  N(t)\leqslant \dfrac{1}{\varepsilon} \max\left(C,\int_{-\infty}^{+\infty} (1+|v|^3)\rho^0(v)\diff v \right)^\frac13 e^{\frac{C}{3\varepsilon^2}t}.  \]
\end{theorem}
 The result is proved with the choice of a linear by part firing function
  \[\phi_\varepsilon(v) = \dfrac{1}{\varepsilon} (v-V_F)_+.\]
For the choice
$   \phi_\varepsilon(v) = \frac{1}{\varepsilon} \mathds 1_{(V_F,+\infty)}(v)$, \cite{CP} notices the simpler direct control
\[ N(t)\leqslant \dfrac{1}{\varepsilon}.  \]
Numerical simulations in \cite{CP} show that solutions of \eqref{NNLIFRandR} can display convergence towards a time-periodic solution, which was never the case for system \eqref{NNLIF} without delay ($d=0$). We could think that the random firing mechanism effectively plays the role of a delay as some randomly chosen neurons will have a delayed impact on the network upon crossing the threshold $V_F$, but a notable difference is that those periodic oscillations are in the excitatory case, which the delayed NNLIF model with deterministic firing cannot reproduce. This type of periodic activity was rigorously proved to arise through a Hopf bifurcation for a hyperbolic variant of the model \cite{CTV,CTV2}, where there is no diffusion $a(N)\equiv 0$ (See \eqref{LIFRand} In Subsection \ref{sec:LIFrand} below).

When the parameter $\varepsilon$ tends to 0, the NNLIF model with random firing formally converges towards the standard model \eqref{NNLIF}. This limit was rigorously proved in the linear case $b=0$, $d=0$, $a(N)\equiv a > 0$ in \cite{LWXZZ}, using the iteration scheme method exposed above in Section \ref{sec:NNLIF_iteration}. The authors also show with numerical simulations that the convergence should hold in the nonlinear case $b\neq 0$. Rigorous convergence in the nonlinear case was then proved rigorously in \cite{DPSZ2}, allowing to construct generalised solutions to \eqref{NNLIF} (See Subsection \ref{sec:NNLIF_gen_nonoise} above).

\subsubsection{A system of excitatory-inhibitory populations}

In \cite{B3}, a model was proposed for two interacting populations of neurons, one excitatory and the other inhibitory. In many areas of the brain, distinct interlaced populations play excitatory and inhibitory roles, like inhibitory interneurons among pyramidal cells \cite{DC}. The excitatory-inhibitory NNLIF system involves two probability densities $\rho_I(v,t)$ and $\rho_E(v,t)$ solving 
\begin{equation}\label{EI}
	\left\{\begin{array}{l}
		\dfrac{\partial \rho_I}{\partial t} + \dfrac{\partial}{\partial v}\Big[ h^I\big(v,N_E,N_I\big) \rho_I  \Big]- a_I(N_E,N_I)\dfrac{\partial^2 \rho_I}{\partial v^2}=N_I(t)\delta_{V_R},\\[0.3cm]
		\dfrac{\partial \rho_E}{\partial t} + \dfrac{\partial}{\partial v}\Big[ h^E\big(v,N_E,N_I\big) \rho_E  \Big]- a_E(N_E,N_I)\dfrac{\partial^2 \rho_E}{\partial v^2}=N_E(t)\delta_{V_R},\\[0.3cm]
		N_\alpha(t) = -a_\alpha(N_E(t),N_I(t))\dfrac{\partial \rho_\alpha}{\partial v}(V_F,t), \quad  \alpha=E,I\\
		\rho_\alpha(V_F,t)= \rho_\alpha(-\infty,t)=0,\quad \rho_\alpha(v,0)=\rho_\alpha^0(v),\quad 
		\displaystyle{\int_{-\infty}^{V_F}\rho_\alpha^0(v)\diff v=1, \quad \ \alpha=E,I}
	\end{array}\right. 
\end{equation}
with
\[
\begin{array}{rcll}
	h^\alpha\big(v,N_E,N_I\big)&=&-v + b_E^\alpha N_E - b_I^\alpha N_I,\\
	a_{\alpha}(N_E,N_I)&=& a_0^\alpha + a_E^\alpha N_E + a_I^\alpha N_I,& \quad \alpha=E,I
\end{array}
\]
where $b_E^E$, $b_E^I$, $b_I^E$, $b_I^I$, $a_E^E$, $a_I^E$, $a_E^I$, $a_I^I$, $a_0^E$ and $a_0^I$ are all non-negative parameters. The first rigorous study of such a model was proposed in 
\cite{CS}.

The action of the inhibitory population, even when the inhibitory coefficients are large is not necessarily enough to prevent entirely finite time blow-up. Under an \textit{a priori} assumption of the form $\int_0^t N_I(s)\diff s \leqslant C t$, \cite{CS} proves that as soon as $b_E^E>0$, finite-time blowup can occur for initial conditions where $\rho^0_E$ is concentrated around $V_F$. In unpublished results one can find in \cite{R2}, this technical assumption was removed and replaced by assumptions on coefficients only: $a_E$ and $a_I$ positive constant functions, $b_E^E>0$ and $b_E^E>\frac1{V_F-V_R}b_E^I b_I^E$. However, these results do not offer a complete understanding on finite-time blow-up in excitatory-inhibitory NNLIF systems as many simple at hand counter-example prove that the conditions are not sharp. The interplay between the parameters and the two initial conditions makes any attempt at a general classification challenging.

Regarding stationary states, for $\alpha = E, I$, the stationary densities are of the form
\[ \rho_{\alpha,\infty} (v) = \dfrac{N_{\alpha,\infty}}{a_\alpha(N_{E,\infty},N_{I,\infty})} e^{-\frac{(v-b_E^\alpha N_E + b_I^\alpha N_I)^2}{2a_\alpha(N_{E,\infty},N_{I,\infty})}} \int_{\max(v,V_R)}^{V_F} e^{-\frac{(w-b_E^\alpha N_E + b_I^\alpha N_I)^2}{2a_\alpha(N_{E,\infty},N_{I,\infty})}} \diff w. \]
where $N_{\alpha,\infty}$ satisfies
\[ \dfrac{1}{N_{\alpha,\infty}} =  \int_{0}^{+\infty}  \dfrac{ e^{-\frac{s^2}{2}} }{s} \left( e^{s\frac{V_F-b_E^\alpha N_E + b_I^\alpha N_I}{\sqrt{a_\alpha(N_{E,\infty},N_{I,\infty})}} } - e^{s\frac{V_R-b_E^\alpha N_E + b_I^\alpha N_I}{\sqrt{a_\alpha(N_{E,\infty},N_{I,\infty})}} } \right)  \diff s. \]
It is possible to provide a partial classification.
\begin{theorem}[C\'aceres, Schneider \cite{CS}]
	Assume $b_I^E, b_E^I>0$, and $a_E$, $a_I$ are constant positive functions. Then for $\alpha=E,I$,
	\begin{itemize}
		\item[$\bullet$] if
		$  (V_F - V_R)^2 < (V_F-V_R)(b_E^E - b_I^I) + b_E^E b_I^I - b_I^E b_E^I$,
		either there is no stationary state of \eqref{EI} or an even number of stationary states; if $b_E^E$ is large enough with respect to other parameters, there is no stationary state;
		\item[$\bullet$] if
		  $(V_F-V_R)(b_E^E - b_I^I) + b_E^E b_I^I - b_I^E b_E^I < (V_F-V_R)^2$,  
		there is an odd number of stationary states of \eqref{EI} ; if $b_E^E$ is small enough, there is a unique stationary state.
	\end{itemize}
\end{theorem}
As this result and numerical examples in \cite{CS} show, there can be one, two, three or no stationary states; the two population system exhibits a more complicated picture than the standard model \eqref{NNLIF}.

Given a stationary state $(\rho_{E,\infty},\rho_{I,\infty},N_{E,\infty},N_{I,\infty})$, a natural relative entropy for asymptotic convergence is
\[  \mathcal H[\rho_E,\rho_I](t) =   \int_{-\infty}^{V_F} \left(  \dfrac{\rho_E(v,t) - \rho_{E,\infty}(v)}{\rho_{E,\infty}(v)} \right)^2 \rho_{E,\infty}(v)\diff v + \int_{-\infty}^{V_F} \left(  \dfrac{\rho_I(v,t) - \rho_{I,\infty}(v)}{\rho_{I,\infty}(v)} \right)^2 \rho_{I,\infty}(v)\diff v, \]

\begin{theorem}[C\'aceres, Schneider \cite{CS}]\label{GoingPlaces3}
	Assume $b_E^E,b_E^I,b_I^E,b_I^I$ are all small enough, $a_E$, $a_I$ are constant positive functions. Then, if
	\[\mathcal H[\rho_E,\rho_I](0) \leqslant \dfrac{1}{2\max(b_E^E+b_I^E,b_E^I + b_I^I)},\]
	 there exists $\mu > 0$ such that the solution of \eqref{EI} satisfies,
	\[ \forall t\geqslant 0, \qquad \mathcal H[\rho_E,\rho_I](t) \leqslant e^{-\mu t}\mathcal H[\rho_E,\rho_I](0). \]
\end{theorem}

\subsubsection{Combination of excitatory-inhibitory populations, refractory state and delay}\label{subsec:twopopul}

In \cite{CS2}, a more realistic NNLIF-type model is studied theoretically and numerically. It combines separate excitatory and inhibitory populations, a refractory state and transmission delays. The model writes, for $\alpha = E,I$,
\begin{equation}\label{EIrd}
	\left\{\begin{array}{l}
		\dfrac{\partial \rho_\alpha}{\partial t} + \dfrac{\partial}{\partial v}\Big[ h^\alpha \rho_\alpha  \Big]- a_\alpha\dfrac{\partial^2 \rho_I}{\partial v^2}=\dfrac{R_\alpha(t)}{\tau_\alpha}\delta_{V_R},\\
		N_\alpha(t) = -a_\alpha \dfrac{\partial \rho_\alpha}{\partial v}(V_F,t),\\
        \dfrac{\diff R_\alpha}{\diff t}(t) = N_\alpha(t) - \dfrac{R_\alpha(t)}{\tau_\alpha}, \qquad R_\alpha(0) = R_\alpha^0 \\
		\rho_\alpha(V_F,t)= \rho_\alpha(-\infty,t)=0,\quad \rho_\alpha(v,0)=\rho_\alpha^0(v),\quad 
		\displaystyle{\int_{-\infty}^{V_F}\rho_\alpha^0(v)\diff v + R_\alpha^0 =1, }
	\end{array}\right. 
\end{equation}
with
\[
\begin{array}{rcll}
	h^\alpha(v,t)&=&-v + b_E^\alpha N_E(t-d_E^\alpha) - b_I^\alpha N_I(t-d_I^\alpha),\\
	a_{\alpha}(t)&=& a_0^\alpha + a_E^\alpha N_E(t-d_E^\alpha) + a_I^\alpha N_I(t-d_I^\alpha),& \quad \alpha=E,I
\end{array}
\]
where $b_E^E$, $b_E^I$, $b_I^E$, $b_I^I$, $a_E^E$, $a_I^E$, $a_E^I$, $a_I^I$, $a_0^E$, $a_0^I$, $d_E^E$, $d_E^I$, $d_I^E$ and $d_I^I$ are all non-negative parameters, and $\tau_E,\tau_I>0$ are positive parameters. It is of course also possible to replace the smooth refractory input $\frac{R_\alpha(t)}{\tau_\alpha}$ by $N_\alpha(t-\tau_\alpha)$. The mass of the density plus the refractory population is conserved for $\alpha=E,I$.

The study of stationary states proceeds similarly to the other NNLIF-type models. Adding a refractory state to the excitatory-inhibitory population prevents the case where no stationary state can exist.

\begin{theorem}[C\'aceres, Schneider \cite{CS2}] Assume $b_E^I,b_I^E>0$ and $a_E,a_I$ are positive constant functions. Then, there is always an odd number\footnote{Like mentioned before, up to a finite number of values of the parameters where it is up to multiplicity only.} of stationary states to \eqref{EIrd}. If $b_E^E$ is small enough, or $\tau_E$ is large enough, there is a unique stationary state.
\end{theorem}
After proving Theorem \ref{thm:NNLIFr_CV}, \cite{CS2} proves an equivalent result for the general system \eqref{EIrd} in the non-delayed case $d_E^E=d_I^E=d_E^I=d_I^I=0$, using an appropriate relative entropy functional and adequate smallness assumptions on coefficients and initial conditions. The finite time blow-up result proved for system \eqref{EI} is also extended to \eqref{EIrd} in the case $d_E^E=0$.

Note that in the case where all delay coefficients are positive, it is possible to rewrite \eqref{EIrd} in a form that allows to extend the existence proof of \cite{CRSS}, thus proving global-in-time existence of smooth solutions. More generally, the construction of local-in-time solutions of \eqref{EIrd}, although not done yet, would be a (very) technical extension of \cite{CGGS} whose main benefit would be the procurement of representation formulae. The general study of the properties of \eqref{EIrd} is still mostly open.

On the numerical side, \cite{CS2} confirms that this model also exhibits time-periodic solutions in the delayed case with strong inhibition, like for \eqref{NNLIF}, but also time-periodic solutions in the average excitatory case. A complete picture of the influence of each of the numerous parameters in \eqref{EIrd} on the occurrence of self-sustained oscillations is currently lacking.

\subsubsection{NNLIF with learning rule}

In all NNLIF-type models we have described so far, the connectivity parameter is fixed throughout the time evolution. Following works on models for learning among a finite number of neurons \cite{GK,GK2,GW1,GW2}, attempts were made in \cite{PSW} at describing learning in a large (mathematically infinite) network of neurons. The goal was to answer some questions inspired from \cite{H}:
\begin{itemize}
    \item can any pattern of neural activity be generated by a heterogeneous synaptic weights distribution?
    \item what is the equilibrium synaptic weights distribution for given learning rule and external inputs?
    \item can the network remember an external signal presented through a stationary weight distribution?
\end{itemize}
In order to investigate those questions, \cite{PSW} proposes two NNLIF-like models. First, a model structured by voltage $v\in(-\infty,V_F]$ and synaptic weights $w\in\R$,
\begin{equation}\label{eq:NNLIF_learning_1}
	\left\{\begin{array}{l}
		\displaystyle \dfrac{\partial \rho}{\partial t} + \dfrac{\partial}{\partial v}\Big[(-v + J(w) + w \sigma\big(\bar N(t)\big) )\rho\Big]-a\dfrac{\partial^2 \rho}{\partial v^2}= \ \delta_{v=V_R}  N(w,t), \\
		\displaystyle N(w,t)= -a \dfrac{\partial \rho}{\partial v}(V_F,w,t), \quad \bar N(t) = \int_{-\infty}^{+\infty} N(w,t)\diff w,  \\
		\displaystyle \rho(V_F,w,t)=\rho(-\infty,w,t)=0,\quad \rho(v,w,0)=\rho^0(v,w) \geqslant 0,  \\\displaystyle \int_{-\infty}^{V_F}\int_{-\infty}^{+\infty} \rho^0(v,w)\diff v \diff w = 1,
	\end{array}\right.
\end{equation}
where $a>0$ is a constant, $\sigma\in C^2(\R_+,\R_+)$, $\sigma'\geqslant 0$, $\sup_{N\geqslant 0}\sigma(N)=\sigma_m < +\infty$. The neurons interact through the mean-firing rate $\bar N$ and the distribution along the connectivity parameter $w$. For a same average firing rate, neurons with large positive $w$ will receive strong excitation and neurons with a large negative $w$ will receive strong inhibition. The function $J$ is an external input. In this model, the synaptic distribution cannot change. The synaptic distribution can be represented by the function
\begin{equation}
    H(w) = \int_{-\infty}^{V_F} \rho(v,w,t)\diff v = \int_{-\infty}^{V_F} \rho^0(v,w)\diff v .
\end{equation}
Note that \[\int_{-\infty}^{+\infty} H(w)\diff w=1.\]

Second, they propose a kinetic model of the form
\begin{equation}\label{eq:NNLIF_learning_2}
	\left\{\begin{array}{l}
		\displaystyle \dfrac{\partial \rho}{\partial t} + \dfrac{\partial}{\partial v}\Big[(-v + J(w) + w \sigma\big(\bar N(t)\big) )\rho\Big] + \varepsilon\dfrac{\partial}{\partial w}\Big[(\Phi-w)\rho\Big] -a\dfrac{\partial^2 \rho}{\partial v^2}= \ \delta_{v=V_R}  N(w,t),\\
		\displaystyle N(w,t)= -a \dfrac{\partial \rho}{\partial v}(V_F,w,t), \quad \bar N(t) = \int_{-\infty}^{+\infty} N(w,t)\diff w, \\
		\displaystyle  \rho(V_F,w,t)= \rho(-\infty,w,t)=0,\quad \rho(v,w,0)=\rho^0(v,w) \geqslant 0, \\\displaystyle \int_{-\infty}^{V_F}\int_{-\infty}^{+\infty} \rho^0(v,w)\diff v \diff w = 1.
	\end{array}\right.
\end{equation}
This model allows for a change in $H$ over time. The scaling parameter $\varepsilon$ allows for a slow change in the weights distribution compared to the evolution of voltage in neurons. The learning rule is defined by the function $\Phi$. One possible choice is the celebrated Hebbian learning \cite{H2,GK2}: the strength of the connectivity between two neurons increases when they have high activity simultaneously, which can be represented in the context of \eqref{eq:NNLIF_learning_2} by
\begin{equation}
    \Phi(w,t) = \bar N(t) N(w,t) K(w),
\end{equation}
for some learning kernel $K(\cdot)$. The idea is that if the specific activity $N(w,t)$ of neurons with connectivity $w$ is high and at the same time the average activity represented by $\bar N(t)$ is high, then $\Phi$ is high, driving neurons towards stronger connectivities \textit{via} a flux towards positive $w$.

Another, more general choice could be
\begin{equation*}
    \Phi(w,t) = N(w,t)\int_{-\infty}^{+\infty}K(w,\omega)N(\omega,t)\diff \omega,
\end{equation*}
where $K(\cdot,\cdot)$ is a learning kernel, or other choices inspired from spike timing dependent plasticity (STDP) rules \cite{DP,CD,MGS,GW2}.

\begin{theorem}[Perthame, Salort, Wainrib \cite{PSW}]
    Assume $J\in L^\infty(\R)$ and
    \[ \int_{0}^{+\infty} w^2 H(w) \diff w < +\infty,\]
    then there exists at least one stationary state of \eqref{eq:NNLIF_learning_1}. If $H$ is supported on $(-\infty,0]$, there is a unique stationary state.
\end{theorem}
From a stationary state $\rho_\infty$ with firing rate $N_\infty$, we can define the normalised output signal
\begin{equation*}
    S(w) = \dfrac{N_\infty(w)}{\bar N_\infty}.
\end{equation*}

Given a reasonable output signal $S(w)$, \cite{PSW} proves that it is always possible to find a synaptic distribution $H(w)$ which allows the network to produce the output $S(w)$.

\begin{theorem}[Perthame, Salort, Wainrib \cite{PSW}]
    Assume $J\in L^\infty(\R)\cap C^1(\R)$ and $S\in L^1(\R)$ satisfies $\int_{-\infty}^{+\infty}S(w)\diff w=1$ and $\int_{-\infty}^{0}S(w)e^{-\gamma w}\diff w< +\infty$ for some $\gamma>\frac{\sigma_m}{a}(V_F-V_R) $.
    Then there exists a weigth distribution $H$ and a stationary state of \eqref{eq:NNLIF_learning_1} with normalised output signal $S(w)$. If $S$ is supported in $(-\infty,0]$, then $H$ is unique.
\end{theorem}
A precise notion of discrimination property is also given and studied theoretically: from two input signals $J_1(w)$ and $J_2(w)$, we have under reasonable assumptions
\begin{equation*}
    \int_{-\infty}^{+\infty}|N_1(w)-N_2(w)|\diff w \geqslant C \int_{-\infty}^{+\infty}L(w) N(w)|J_1(w)-J_2(w)|\diff w,
\end{equation*}
with $L>0$ a well-chosen constant and $C>0$ a constant depending on $H,J_1,J_2$ and other parameters.

For the kinetic variant \eqref{eq:NNLIF_learning_2}, \cite{PSW} proves that a given input $J$ can lead to many different stationary weight distributions $H$ after learning.
\begin{theorem}[Perthame, Salort, Wainrib \cite{PSW}]
    Assume $J\in L^\infty(\R)\cap C^1(\R)$, $wK(w)>0$ for $w\neq 0$ and $0 < k_m \leqslant K(w) \leqslant K_M$ almost everywhere. Then there are infinitely many different stationary states of \eqref{eq:NNLIF_learning_2}.
\end{theorem}
It is however possible to describe some aspects of these stationary states. In specific cases, \cite{PSW} gives a way to organise and select them, by adding a Gaussian noise $\varepsilon\nu \frac{\partial^2 \rho}{\partial w^2}$ in the system in the fully inhibitory context. Then, they select a stationary state of \eqref{eq:NNLIF_learning_2} in the small noise limit $\nu\to0$.

Numerical tests are then proposed by letting the system learn from an input $J_1$ in \eqref{eq:NNLIF_learning_2}, then using the corresponding stationary weigth distribution $H$ as an initial condition in \eqref{eq:NNLIF_learning_1} subject to a different input $J_2$. The idea is to let a system learn the input $J_1$ on a long timescale and then, on a shorter timescale, to let it react to the other input $J_2$. Further numerics with a carefully crafted scheme have been proposed in \cite{HHZ}.

\subsubsection{Poissonian mean-field variant}

In all the variants of the NNLIF model we have described above, the idea was to complexify the model in order to describe more realistic properties of neural networks. In two companion papers \cite{TW,ST}, a simpler variant is proposed with the idea of doing more explicit computations and achieving a better characterisation of finite time blow-up. They assume that neurons are driven by Poisson-like noisy inputs and enter a refractory period after spiking, leading to what they name a delayed Poissonian mean-field dynamics (dPMF).
\begin{equation}\label{eq:NNLIF_dPMF}\tag{dPMF}
    \left\{\begin{array}{l}
         \dfrac{\partial \rho}{\partial t} = (\nu + bN(t))\left( \dfrac{\partial \rho}{\partial v} + \dfrac12\dfrac{\partial^2 \rho}{\partial v^2} \right) + N(t-\varepsilon_0) \delta_{v=V_R} \\
         N(t) = \dfrac12 \dfrac{\partial \rho}{\partial v}(0,t),\\
         \rho(0,t)=\rho(+\infty,t) = 0,\quad \rho(v,0) = \rho^0(v)\geqslant 0 ,\quad\int_{0}^{+\infty}\rho^0(v)\diff v = 1,
    \end{array}\right.
\end{equation}
with $b\in\R$ the connectivity parameter, $\nu>0$ a constant drift and $\varepsilon_0>0$ the length of the refractory period. Note that they choose $0 = V_F < V_R $ and study the model in the non-negative half space $v\in[0,+\infty)$ for simplicity. For this system, the firing rate can be written as
\begin{equation}
    N(t) = \dfrac{\nu\partial_v \rho(0,t-\varepsilon_0)}{2-b \partial_v \rho(0,t)}
\end{equation}
which shows that as \eqref{NNLIF} with internal noise $a(N) = a_0 +a_1 N$, finite-time blow-up happens with a finite boundary derivative $\frac12 \partial_v \rho(V_F,T^*) = \frac1b$.
Their method is to introduce the time change
\begin{equation}
    \tau = \Phi(t) = \nu + b F(t), \qquad F(t) = \int_0^t N(s) \diff s,
\end{equation}
which allows to dilate the blow-up and continue the dynamics within. This idea is quite similar to the strategy of \cite{DZ,DPSZ2} although developed with completely different methods and, up to our best knowledge, independently and simultaneously. 
After the time dilation, \cite{DZ,DPSZ2} use more deterministic PDE tools while \cite{TW,ST} use more probabilistic methods. This idea was later used on a hyperbolic PDE in \cite{CDRZ}.

In the case of the model \eqref{eq:NNLIF_dPMF}, the time change allows to get precise information about the blow-up dynamics and to perform a vanishing refractory period limit $\varepsilon_0\to 0$. A notable byproduct is a characterisation of the divergence speed for $N(t)$.

\begin{theorem}[Taillefumier, Whitman \cite{TW}]\label{thm:NNLIF_dPFM}
    Under a full blow-up condition, blowing-up solutions to \eqref{eq:NNLIF_dPMF} satisfy at the blow-up time $T^*$,
    \begin{equation}
        N(t) \sim \dfrac{b}{\sqrt{2C(T^*-t)}},
    \end{equation}   
    for some constant $C>0$ which can be characterised.
\end{theorem}
The method also allows to characterise the mass of neurons participating into a full blow-up event. In the companion article \cite{ST}, it is proved that there are global generalised solutions with infinitely many full blow-up events.

\subsubsection{NNLIF-like models in finance modelling}

Many researchers in financial mathematics \cite{HS, LS, NS,NS2,NS3,HL2, DNS}
 consider probabilistic NNLIF-type models in the form of McKean-Vlasov equations depending on singular hitting times; they have developed their own terminology and set of techniques. Although a careful description of their methods and results is beyond the scope of this review, we want to give a broad idea of the kind of models they work with, provide some references, and briefly mention some results they have obtained that can be extended to the NNLIF model with enough technical work.

Default contagion in large financial networks, bank networks, or in a large portfolio can be modelled by the distance to default $X_t\in [0,+\infty)$, $0$ being the default threshold. In the mean-field limit, the most simple probabilistic formulation \cite{HLS} is
\begin{equation}
    \left\{
    \begin{array}{l}
         X_t = X_0 +B_t - b L_t\\
         \tau = \inf\{ t\geqslant 0\ |\ X_t \leqslant 0 \},\\
         L_t = \mathbb P(\tau \leqslant t)
    \end{array}
    \right.
\end{equation}
whose law is a solution to the PDE
\begin{equation}
    \left\{
    \begin{array}{l}
         \dfrac{\partial \rho}{\partial t} - b N(t)  \dfrac{\partial \rho}{\partial x} - \dfrac12 \dfrac{\partial^2\rho}{\partial x^2} = 0\\[0.3cm]
         N(t) = \dfrac12 \dfrac{\partial \rho}{\partial x}(0,t),\\[0.3cm]
         \rho(0,t) = 0,\qquad \rho(x,0) = \delta_{x=X_0}.
    \end{array}
    \right.
\end{equation}
Many more realistic financial models are built on the same basis and mean-field limits have been proved. Contrarily to neurons, the financial agents, banks or portfolio entities are not resurrected to a reset value after defaulting. However, like the NNLIF model without delay, these finance models are prone to finite-time blow-up and a notion of \textit{physical solution} can be constructed like in \cite{DIRT2}. 

The size of the jumps in physical solutions is of a particular importance for finance applications as it is indicative of the severity of, for example, a financial crisis. It can be proved that the \textit{physical solutions} have minimal jumps after blow-up \cite[Theorem 1.2]{HLS} and the article makes progresses towards a proof of uniqueness for the \textit{physical solutions} \cite[Theorem 1.8, discussion therein]{HLS}. Recently, more methods were developed to extend solutions of similar stochastic equations in spite of the presence of blow-up \cite{DNS,LS,NS3}. This community also performed a careful study of how the decay of the solution near the boundary $0$ has an impact on the regularity after blow-up events.

\subsection{Numerical methods for NNLIF type systems}\label{subsec:numnnlif}

Several numerical methods have been proposed for approximating NNLIF equations and their variants. The first attempts by physicists and applied mathematicians have been applying Monte-Carlo methods based on the stochastic processes underlying these PDE models, see for instance \cite{ST01,RC07,CTRM,RKC08,CTSM,CTM-PNAS04,RCM-PNAS08,NKKRC10,NKKZRC10}.
The first numerical results by deterministic schemes of the NNLIF model \eqref{NNLIFphy}-\eqref{NNLIFphy2} and their kinetic conductance-voltage counterpart in \cite{CCP,CCT} used explicit in time finite difference methods in voltage and conductance based on \textit{Weighted Essentially Non-Oscillatory} (WENO) schemes \cite{S3} originating in conservation laws \cite{SO88,JS96}. They approximate the advection terms to high order in order to simulate efficiently and accurately large gradient regions in voltage or conductance that show up due to the blow-up of solutions for large connectivity. Second-order terms are discretized by centered finite differences while time discretization is performed via 3rd-order \textit{Total Variation Diminishing} Runge-Kutta schemes. The resulting schemes have been very successful in capturing blow-up phenomena, kink formation and convergence to stationary states for both the NNLIF \eqref{NNLIFphy}-\eqref{NNLIFphy2} without delay $d=0$ \cite{CCP} and the kinetic conductance-voltage model \cite{CCT}. Moreover, they have been widely used to analyse the NNLIF \eqref{NNLIFphy}-\eqref{NNLIFphy2} with delay in \cite{CS}, the two excitatory-inhibitory population model \eqref{EIrd} in \cite{CS2}, see Subsection \ref{subsec:twopopul}, and the sequence of pseudo-equilibria describing the long-time behaviour of the NNLIF model \eqref{NNLIFphy}-\eqref{NNLIFphy2} with large delay in \cite{CCR-L2}, see Subsection \ref{subsec:pseudo}. We refer to the corresponding previous subsections for details on the interpretation of the numerical results related to stability of steady states, periodic solutions and synchronization of the network.

A second family of deterministic numerical methods proposed in the literature are based on structure preserving methods, meaning schemes that preserve some properties of the evolution of the problem. In the case of the NNLIF model, the basic property to preserve at least in the linear case, $b=0$, it is the decay of the relative entropy to equilibrium. This family of methods based on finite difference/finite volume schemes for Fokker-Planck type equations can be found with the name Chang-Cooper method, see for instance \cite{BCD}, or Scharfetter-Gummel method, see for instance \cite{JY}.  Both methods use the following rewriting of the Fokker-Planck equation \eqref{NNLIFphy} in terms of the Maxwellian
$$
M(v,t)=e^{-\frac{(v-bN(t))^2}{2a(N(t))}}
$$ 
as follows,
$$
\frac{\partial p}{\partial t}(v,t) - a(N(t)) \frac{\partial}{\partial v}
\left[
M(v,t)\frac{\partial}{\partial v}\left(\frac{p(v,t)}{M(v,t)} \right)
\right]=N(t)\, \delta(v-V_R).
$$
Then, the Chang-Cooper method performs a kind of $\theta$-finite difference approximation of $p/M$, see \cite{BCD} for details, while the Scharfetter-Gummel method does a harmonic average, see \cite{HLXZ} for details. Both methods have been used in \cite{CCP,CCT,HLXZ,HHZ} for a variety of the computational neuroscience problems discussed in the previous sections, showing some of the conservation and dissipation properties of these schemes. In comparison with the WENO schemes, they are lower order, however they are able to keep known dissipation properties in the linear case $b=0$ without network connectivity. Other numerical approaches have been introduced such as spectral methods \cite{ZWZ}, multiscale methods \cite{DXZ} mixing stochastic microscopic schemes with deterministic macroscopic schemes to analyse synchronization in NNLIF, and discontinuous-Galerkin methods \cite{SPAK,SS}.


\section{Fokker-Planck representation of a stochastic neural field and the modelling of grid cells}
\label{sec:GC}

We present here a PDE framework recently developed for the study of the asymptotic behaviour and phase transitions in stochastic neural fields, with particular interest to its applications to patterning in noisy grid cells.

We first explain the main application goal of this family of models in the context of previous modelling attempts on grid cells. Then we describe the rigorous mean-field limit which gave rise to these systems. Finally, we introduce the mathematical setting and propose an overview of mathematical results on the solutions: stationary states and their bifurcations, existence and regularity of solutions, asymptotic behaviour.

\subsection{Modelling the hexagonal firing pattern of grid cells}

Grid cells, discovered in 2005 \cite{HFMMM}, have been described as the GPS in the navigation system for mammals. In fact, the hexagonal patterns they create play an important role in the integration of the moving paths of animals (including humans \cite{DBB}) as demonstrated for mammalian brains in \cite{MMM}. The understanding of their precise behaviour, see \cite{RRMM,MMM} and the references therein, is still a matter of modern research. One of the main challenges highlighted in \cite{RRMM}, was the effect of noise on grid cells and the robustness of these hexagonal patterns. The grid cell network is commonly described by deterministic continuous attractor network dynamics through a system of neural field models \cite{Ermentrout2010}. In order to better understand the formation of the hexagonal firing pattern and the role of grid cells in path integration, modelling studies have tried to represent in a simple framework how patterned activity arises through autonomous neuronal interactions and external inputs of speed and orientation. 

Around 2009, mathematical models were designed to explain the patterns of grid cells, among which the works of Burak and Fiete, \cite{BF3,FBB}. In order to stay consistent with assumptions and notation, we do not present exactly their model but a slightly modified version. Consider neurons arranged in a two-dimensional sheet at locations $(x_1,\dots,x_N)$ on a square $Q=[0,L]\times[0,L]$ of length $L$. Each space location contains four neurons with orientation preference $\beta\in\{1,2,3,4\}$, where these numbers represent North, West, South, East (all preference angles in $[0,2\pi]$ are possible in reality; the model uses only four directions for simplicity). The electrical activity (how much it fires per unit time) $s_i^\beta$ of a neuron at location $x_i$ with orientation preference $\beta$ evolves according to the differential equation 
\begin{equation}\label{eq:BF}
	\tau \dfrac{\diff s_i^\beta}{\diff t} = \Phi\left( \sum_ {\beta' = 1}^{4} \sum_ {j = 1}^{N} W_{i,j}^{\beta'} s_j^{\beta'} + B^\beta(t) \right)
\end{equation}
where the connectivity parameter between locations $x_i$ and $x_j$ given the orientation preference $\beta$ is
$W_{i,j}^\beta = W(x_i-x_j - r^\beta)$, computed according to the connectivity kernel $W$ (a function), periodically extended outside of the square domain $Q$, and the shift $r^\beta$ (a two-dimensional vector). The modulation function $\Phi$ can be thought of as a filter, like $\Phi(z)=\max(z,0)$ (the so-called ReLU function), or a logistic / sigmoid function like $\phi(z) = (1 + e^{- x} )^{-1}$. The characteristic response time for neurons is $\tau\simeq 10 ms$.

The input from the external world or other neuronal cells related to navigation, such as place or border cells, is contained in the terms $B^\beta$, which in turn depend on the velocity $v(t)$ and the orientation $\theta(t)$ of the animal. We consider here the simple form
\begin{equation*}
	B^\beta(t) = B + \alpha v(t)\cos( \theta - \theta^\beta),
\end{equation*}
used in \cite{BF3}, where $B,\alpha\in [0,\infty)$ and $(\theta^1,\theta^2,\theta^3,\theta^4) = (\tfrac\pi2, \pi,\tfrac{3\pi}{2},2\pi)$. This setting allowed to reproduce the firing pattern of grid cells and to explain how they can be used by the brain for path integration; see \cite{BF3,CetAl}.

\subsection{Rigorous mean-field limit towards a new model}

 Noise plays a role in many emerging brain phenomena and is more a feature than an accidental element of the way neural networks perform tasks; see \cite{RD} for an overview.  In general, studies on understanding noisy neural field models have been lacking \cite[Sec. 6]{B} until very recently \cite{KE,K14,KS,KS2,MB,BF2,THF,T3,B19,BAC,FI}. However, theoretical foundations are still sparse. In \cite{CCS, CHS}, Brownian fluctuations were added to the attractor model \eqref{eq:BF}, in order to investigate how robust the grid cell pattern is to noise and how it appears as a strong noise decreases to lower levels.

Consider $N$ space locations $x_1,\dots,x_N\in \mathbb T^d:= \R^d/ (L\Z)^d$ of the neural cortex. The model is stated for a general dimension $d$ because it does not impact the results, but $d=2$ is taken for the application to grid cells. Note that $d=1$ yields interesting insights for the celebrated ring model \cite{BBS}. The choice of a torus for $\mathbb T^d$ for the (cortical) spatial domain is motivated by recent insights from topological data analysis \cite{GHPBBDMM}, which indicate that the electrical activity of a grid cell module can be represented by a toroidal topology. At each location, we now consider not just one neuron but a cortical column of $4M$ neurons ($M$ per orientation preference). The $k^{th}$ neuron at location $x_i$ with orientation preference $\beta\in\{1,2,3,4\}$ has an activity level $s_{ik}^\beta(t)$ at time $t$, which evolves according to the coupled system
	\begin{align}
		\tau s_{ik}^\beta(t)=&\nonumber \,\tau s_{ik}(0)+\int_0^t \left[ -s_{ik}^\beta(\tau) + \Phi\left( \dfrac{1}{4NM}\sum_ {\beta' = 1}^{4} \sum_ {j = 1}^{N}\sum_ {l = 1}^{M} W_{i,j}^{\beta'} s_{jl}^{\beta'}(\tau) + B^\beta(\tau) \right) \right]\,\diff\tau\\&\ +\sqrt{2\sigma} \mathcal B_{ik}^\beta(t)-\ell_{ik}^\beta(t), \label{eq:SDE1}
		\\[4mm]
		\ell_{ik}^{\beta}(t)=&\,-|\ell_{ik}^{\beta}|(t),\quad |\ell_{ik}^{\beta}|(t)=\int_0^t \mathds 1_{\{s^\beta_{ik}(\tau)=0\}}\diff|\ell_{ik}^{\beta}|(\tau). \label{eq:SDE2}
	\end{align}
where the Browian motions $(\mathcal B_{ik}^\beta)_{\beta=1,\dots,4}$ can be correlated (see below).

For each $i$, $k$ and $\beta$, the reflection term $\ell_{ik}^\beta$ prevents the activity level $s_{ik}^\beta$ from taking negative values. 
It is equal for all time $t$ to the opposite of its total variation $\ell_{ik}^{\beta}(t)=-\big|\ell_{ik}^{\beta}\big|(t)$ which stays constant when $s_{ik}^\beta>0$ and increases as $\int_0^t \mathds1_{\{s_{ik}^\beta(\tau)=0\}}\diff\big|\ell_{ik}^{\beta}\big|(\tau)$ when $s_{ik}^\beta=0$; the result is that  $s_{ik}^\beta$ is pushed away from the boundary value $s=0$. Proving existence and uniqueness of this kind of term is a \emph{Skorokhod problem}. Such reflection terms are rigorously constructed for the model \eqref{eq:SDE1}--\eqref{eq:SDE2} in \cite{CCS}, applying the strategy of Lions and Sznitman \cite{LS2,S}.

Denote $\rho_{N,M}^\beta$ the empirical measure associated to the $4MN$ particles, that is
\begin{equation}\label{formula/ single empirical measure}
	\rho_{N,M}(t,dy,ds)=\frac{1}{NM}\sum_{j=1}^{N}\sum_{m=1}^M\delta_{(x_j,s_{jm}(t))}\quad\text{regarded as a measure on $\mathbb T^d\times\R^{4}$.}
\end{equation}
where $s_{jm}(t)=(s_{jm}^1,s_{jm}^2,s_{jm}^3,s_{jm}^4)\in\R^4$.

We are then interested in the mean-field limit $M\to+\infty,\ N\to+\infty$. However, this is not a classical mean field problem with $4MN$ particles in $\R$ or $MN$ particles in $\R^4$. As noted in \cite{CCS,C2} Its structure is more a collection of $N$ mean-field systems (the cortical columns) that interact with each other, each containing $M$ particles. For this modelling reason and for mathematical reasons related to measurability in $x$, the four-dimensional Brownian motions $\mathcal B_{ik}$ and the initial conditions $s_{ik}^\beta(0)$ are not independent. There will be independence only between the $M$ neurons at the same location $x_j$ but there are correlations between initial conditions and noise in two different space locations $x_i$ and $x_j$. The authors of \cite{CCS} show that it results in the convergence rate
\[\frac{1}{\sqrt{M}} + \frac{1}{N^{\tfrac{\alpha}{d}}},\]
for some positive $\alpha$ in the mean-field limit $M,N\to+\infty$ instead of the decay rate $\tfrac{1}{\sqrt{MN}}$ classicaly obtained with $MN$ particles.

Assuming the points $(x_1,\dots,x_N)$ are adequately chosen in the torus $\mathbb T^d$, and under other technical assumptions, it is shown in \cite{CCS} that, up to technical modifications due to the measurability problems and the correlations, the system at each location $x\in\mathbb T^d$ converges to independent copies evolving according to the stochastic differential equation
\begin{align}
		\tau s^\beta(x,t)=& \,s^\beta(x,0)+\int_0^t \left[ -s^\beta(x,\tau) + \Phi^\beta(x,\tau) \right]\,\diff\tau +\sqrt{2\sigma} \mathcal B(x,t)^\beta(t)-\ell^\beta(x,t), \label{eq:SDE1-MF}
		\\
		\Phi^\beta(x,t) =&\, \Phi\left( \dfrac{1}{4}\sum_ {\beta' = 1}^{4} \int_{\mathbb T^d} W^{\beta'}(x-y) \int_{\R^4} s^{\beta'}(x,y,t) \rho(t, y,\diff s )\diff y + B^\beta(t) \right),\\
		\ell^{\beta}(x,t)=&\,-|\ell^{\beta}(x,\cdot)|(t),\quad |\ell^{\beta}(x,\cdot)|(t)=\int_0^t \mathds 1_{\{s^\beta(x,\tau)(\tau)=0\}}\diff|\ell^{\beta}(x,\cdot)|(\tau). \label{eq:SDE2-MF}
\end{align}
where we denote $\rho(t,y,\diff s) = \mathrm{Law}_{\R^4}(s(y,t))$, which induces the probability density $\rho(t,\diff x,\diff s)$ on $\mathbb T^d\times\R^4$ defined by
\[  \forall \phi\in C_b(\mathbb T^d\times\R^4), \quad \int_{\mathbb T^d \times \R^4}\phi(x,s) \rho(t,\diff x, \diff s) = \int_{\mathbb T^d}\int_{\R^4}\phi(x,s) \rho(t, x, \diff s)\diff x.  \] 
The marginal distributions $\rho^\beta (x,s,t) = \mathrm{Law}_{\R}(s^\beta(x,t))$ are solution to the $4-$population system of partial differential equations, for $\beta\in\{1,2,3,4\}$,
\begin{align}\label{eq:4PDE}\tag{GC}
	\tau \frac{\partial \rho^\beta}{\partial t} =
	-\frac{\partial}{\partial s}\Bigg(
	\Big[\Phi^\beta(x,t) -s\Big] \rho^\beta
	\Bigg) + \sigma \frac{\partial^2 \rho^\beta}{\partial s^2},
\end{align}
where $\Phi^\beta(x,t)$ is given by
\begin{align}\label{eq:phi}
	\Phi^\beta(x,t) = \Phi \left(\frac{1}{4}\sum_{\beta'=1}^4 \int_{\mathbb T^d} W^{\beta'}(x-y) \int_{0}^\infty s \rho^{\beta'} (y,s,t)\, \diff s \diff y + B^\beta(t) \right),
\end{align}
and the system is supplemented with the nonlinear boundary condition
\begin{align}\label{eq:BC_4pop}
		\left( \Phi^\beta(x,t)\rho^\beta(x,s,t) - \sigma \dfrac{\partial \rho^\beta}{\partial s}(x,s,t) \right)\bigg|_{s=0} = 0.
\end{align}
The next order in the expansion in terms of $M$ and $N$ in the mean field limit can be analysed by understanding the behaviour of a stochastic PDE for the fluctuations of the system; see \cite{C2}. See also \cite{OR} for a neural coding perspective on grid cell modelling.

\subsection{Mathematical setting and definitions}

Consider again the $d-$dimensional torus of length $L$, $\mathbb T^d := \R^d/( L \Z)^d, L>0$. Writing proofs about the system \eqref{eq:4PDE} being quite cumbersome, an efficient route is to first study a simpler model, without the orientation preference $\beta$, and then see what results can be extended to the general case. Thus, we consider the simpler system for one population given by
\begin{equation}\label{eqn:2}\tag{SNF}
	\tau_c\dfrac{\partial \rho}{\partial t} = -\dfrac{\partial }{\partial s}\left[  \left(  \Phi\left( \int_{\mathbb T^d} W(x-y)\int_0^{+\infty}s\rho(y,s,t)\diff s \diff y + B\right) -s  \right)\rho(x,s,t)   \right] + \sigma \dfrac{\partial^2 \rho}{\partial s^2},
\end{equation}
with the nonlinear and non-local no-flux boundary condition
\begin{equation*}
	\left( \Phi\left( \int_{\mathbb T^d} W(x-y)\int_0^{+\infty}s\rho(y,s,t)\diff s \diff y + B\right)\rho(x,s,t) - \sigma \dfrac{\partial \rho}{\partial s}(x,s,t) \right)\bigg|_{s=0} = 0,
\end{equation*}
and the mass normalisation condition
\begin{equation}\label{eqn:mass_x}
	\forall x\in\mathbb T^d, \qquad \int_{0}^{+\infty} \rho^0(x,s)\diff s = \dfrac1{L^d}.
\end{equation}
The system \eqref{eqn:2} can be thought of as a Stochastic Neural Field in the form of a partial differential equation. Other types of stochastic neural fields have been proposed and investigated with probabilistic techniques \cite{KE,MB,BK,B2,KS,KS2,CEK,KT}; see also books and reviews \cite{B,C,C3,CPWT} for discussions about the extensive literature on neural fields in general, and \cite{BAC,BOFRC,C7,TAD,SBOT,COTA} for efforts to build more realistic second-generation neural fields.

We will use the classical notation
$u\ast v : x\mapsto \int_{\mathbb T^d} u(x-y)v(y)\diff y$
for the convolution of $u$ and $v$ on the torus $\mathbb T^d$. We also denote
\[ \Phi_{\bar\rho}(x,t) = \Phi(W \ast \bar\rho(x,t) + B),\quad  \Phi'_{\bar \rho}(x,t) = \Phi'(W \ast \bar\rho(x,t) + B), \quad  \Phi''_{\bar \rho}(x,t) = \Phi''(W \ast \bar\rho(x,t) + B),  \]
with the mean activity level defined by
$$
\bar\rho(x,t):=\int_0^{+\infty}s\rho(x,s,t)\diff s.
$$
The main equation \eqref{eqn:2} can then be rewritten as
\[  	\tau_c\dfrac{\partial \rho}{\partial t} = -\dfrac{\partial }{\partial s}\left[  \left(  \Phi_{\bar\rho} -s  \right)\rho   \right] + \sigma \dfrac{\partial^2 \rho}{\partial s^2}, \qquad \sigma\dfrac{\partial \rho}{\partial s}(x,0,t) = \Phi_{\bar\rho}\,\rho(x,0,t). \]
The nonlinear non-local boundary condition is the principal theoretical difficulty: it destroys the gradient flow structure of the Fokker-Planck equation, complicates the search for stationary states and the construction of smooth solutions. However, note that weak solutions are already provided by the mean-field limit in \cite{CCS}.

For the spatial component $\mathbb T^d$ of the system, we will use the Hilbert space $L^2(\mathbb T^d)$ or the subspace $L^2_S(\mathbb T^d)$ of coordinate-wise even $L^2$-functions:
\[ L^2_S(\mathbb T^d) = \left\{ \ u\in L^2(\mathbb T^d) \ \big| \ \forall i\in\llbracket 1,d \rrbracket,\ u(x_1,\dots,-x_i,\dots,x_d) = u(x_1,\dots,x_i,\dots,x_d)\ \ \text{a.e. in } \mathbb T^d \ \right\}.  \]
Similarly, for any $m\in \N$, the coordinate-wise even Sobolev space is $H^m_S(\mathbb T^d) = H^m(\mathbb T^d)\cap L^2_S(\mathbb T^d)$.
\begin{hyp}[On $\Phi$, $W$ and $B$]\label{as:1} There exists $m\in\N$, such that $2m > d$ and $W\in H^m_S(\mathbb T^d)$. The function $\Phi$ is continuous on $\R$, $C^3$ on $\R^*$, increasing on $\R_+$ and $\Phi\equiv 0$ on $\R_-$. Moreover, $B > 0$ and  $W_0 = \int_{\mathbb T^d} W(y)\diff y < 0$.
\end{hyp}

There are biological and mathematical reasons for these choices:
\begin{itemize}
	\item[$\bullet$] Because $W\in H^m_S(\mathbb T^d)$, with $2m > d$, the smoothing property of the convolution and Sobolev embeddings imply that $W\ast h$ will be continuous for all $h\in L^2(\mathbb T^d)$. It also avoids pathological connectivity kernels with singularities.
	\item[$\bullet$] Consistent with the modelling purpose, the external signal $B$ is excitatory, the connectivity kernel $W$ is inhibitory on average and the modulation function $\Phi$ acts as a filter, transmitting only excitatory signals to the network. 
\end{itemize}

For the evolution problem, we will use the following definitions.

\begin{definition}(Fast decay at $s=+\infty$)
	Let $a\in\R$. We say that a function $\rho : \mathbb T^d \times [a,+\infty) \to \R$, is fast-decaying if for all polynomial functions $s\mapsto P(s)$, for all $x\in\mathbb T^d$,
	  $\lim_{s\to +\infty} \rho(x,s) P(s) = 0$.
\end{definition}

\begin{definition}[Compatible initial condition]\label{def:compatible}
	We say that the function $\rho^0\in C^{0,1}(\mathbb T^d \times \R_+)$ is a compatible initial condition for \eqref{eqn:2} if $\rho^0$ is a non-negative probability density, if $\rho^0$, $\partial_s\rho^0$ are fast-decaying,
	and if
	\[ \forall x\in\mathbb T^d, \quad \int_0^{+\infty} \rho^0(x,s)ds = \dfrac{1}{L^d}.  \]
\end{definition}

\begin{definition}[Classical solution]\label{def:classicalGridCells}
	We say that a non-negative function $\rho\in C^{0}(\mathbb T^d\times\R_+\times [0,T_{\mathrm{max}}))$, $T_{\mathrm{max}}\in (0,+\infty]$ is a classical solution of \eqref{eqn:2} if    
	\begin{itemize}
		\item it has the regularity 
		$
			\rho\in \, C^{0,2,1}(\mathbb T^d\times(0,+\infty)\times (0,T_{\mathrm{max}}))\cap C^{0,1,1}(\mathbb T^d\times\R_+\times (0,T_{\mathrm{max}})),
		$
		and it satisfies the equation in the classical sense;
		\item for all $t\in[0,T_{\mathrm{max}})$, $\rho(\cdot,\cdot,t)$, $\partial_s \rho(\cdot,\cdot,t)$ are fast-decaying;
		\item  for all $(x,t)\in\mathbb T^d\times [0,T_{\mathrm{max}})$, 
		$
		\int_{0}^{+\infty} \rho(x,s,t)ds = \frac{1}{L^d}$.
	\end{itemize}
\end{definition}

\subsection{Analytical study of the stationary problem}

\begin{figure}[ht!]
	\centering
	\includegraphics[width=350pt]{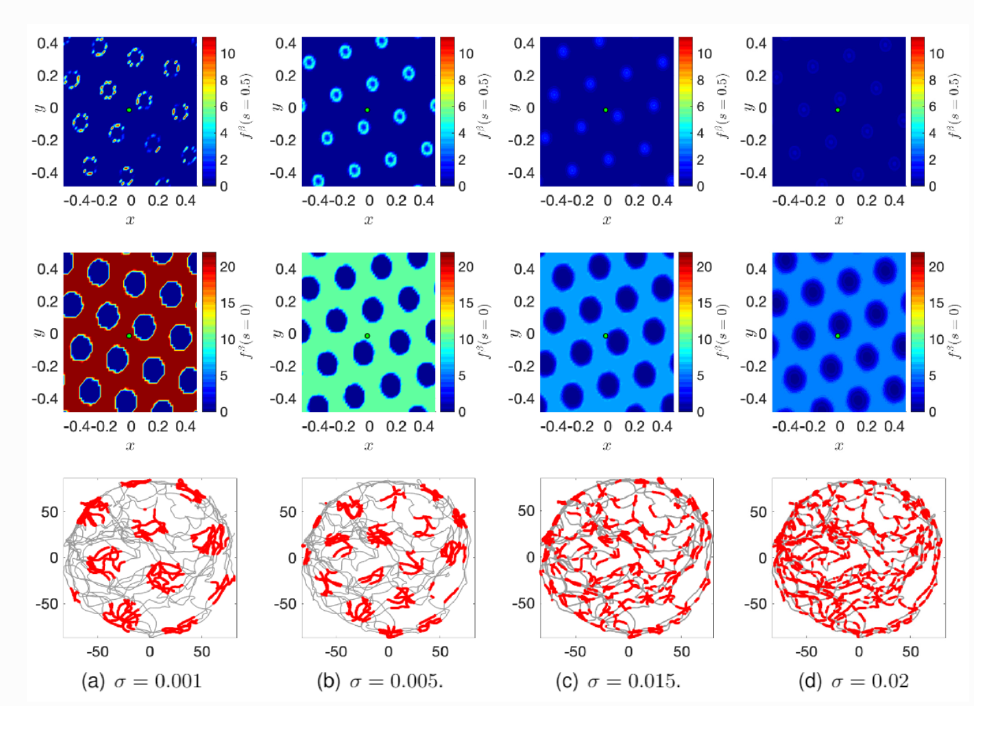}
	\caption{\textit{Numerical simulation of the four population model for grid cells \eqref{eq:4PDE} for different values of $\sigma$. Top Row is the repartition of the density for $s=0.5$; middle row is the same for $s=0$; bottom row represents in grey the path of a real rat and in red the prediction of the model of where in physical space one specific grid cell fires -- the one at $x=(0,0)$. Carrillo, Holden, Solem \cite{CHS} }}\label{fig:hex}
\end{figure}

Numerical simulations in \cite{CHS} show that the model \eqref{eq:4PDE} reproduces the hexagonal firing pattern of grid cells for appropriate choice of parameters as seen in Figure \ref{fig:hex}. Although there is no rigorous proof of the existence of hexagonal patterns yet, partial answers were obtained rigorously \textit{via} a bifurcation analysis in \cite{CRS}. In the following, we present the theory for the system \eqref{eqn:2} for simplicity, but the technical extension for \eqref{eq:4PDE} has been done in \cite{CHS,CRS}.

\subsubsection{Characterisation of stationary states}

Integrating with respect to $s$ and using the boundary condition, the stationary states satisfy
\[   \sigma \partial_s \rho(x,s) =   - \left(s - \Phi\left( W \ast \bar \rho (x) + B \right) \right) \rho(x,s),  \]
which yields the form
\begin{equation}\label{eq:stationarystateequ}
	\rho(x,s) = \dfrac{1}{Z_\rho} \exp\left( - \dfrac{ \big( s - \Phi\left(  W\ast \bar \rho (x) + B \big) \right)^2}{2\sigma} \right), 
\end{equation}
with, owing to \eqref{eqn:mass_x},
$  Z_\rho = L^d \int_{0}^{+\infty} \exp\left( - \frac{ ( s - \Phi\left(  W\ast \bar \rho (x) + B ) \right)^2}{2\sigma} \right) \diff s$. 
Stationary states must be zeros of
\begin{equation*}
	\mathcal G(\rho ,\sigma)  = \rho -\dfrac{1}{Z_\rho}e^{ - \frac{ \left( s - \Phi_{\bar\rho} \right)^2}{2\sigma}} ,\qquad Z_\rho = L^d \int_{0}^{+\infty} e^{ - \frac{ \left( s - \Phi_{\bar\rho} \right)^2}{2\sigma}} \diff s.
\end{equation*}
Consequently, the average $x\mapsto\bar \rho(x)$ is a zero of the functional
\begin{equation}\label{eq:onecompfunctional}
	\bar{\mathcal G}(\bar \rho ,\sigma)  = \bar\rho -\dfrac{1}{Z_\rho} \int_0^{+\infty} s e^{ - \frac{ \left( s - \Phi_{\bar\rho} \right)^2}{2\sigma}}\diff s,\qquad  Z_\rho =  L^d \int_{0}^{+\infty} e^{ - \frac{ \left( s - \Phi_{\bar\rho} \right)^2}{2\sigma}} \diff s.
\end{equation}
We can fully characterise the stationary states of \eqref{eqn:2} considering the zeros of the functional $\bar{\mathcal G}(\bar \rho ,\kappa)$: for any zero of the averaged functional, a stationary state can be defined through \eqref{eq:stationarystateequ} and for any stationary state, its activity average is a solution to \eqref{eq:onecompfunctional}.

A magic ingredient for the study of this problem is the short calculation
\begin{align*}  \rho(x,0) &= - \int_{0}^{+\infty} \dfrac{\partial \rho}{\partial s}(x,s) \diff s = \frac1\sigma \int_{0}^{+\infty} (s-\Phi_{\bar \rho}(x) )\rho(x,s) \diff s \\
	& = \frac1\sigma \Big( \underbrace{\int_{0}^{+\infty} s \rho(x,s)\diff s}_{ = \bar\rho(x)} - \Phi_{\bar\rho}(x) \underbrace{\int_{0}^{+\infty} \rho(x,s)\diff s }_{=\frac1{L^d}}  \Big),
\end{align*}
yielding a relation between the value at $s=0$ of a stationary state, its firing rate and its mean value,
\begin{equation}\label{eqn:rho_0}
	 \rho(x,0) = \frac1\sigma\left(\bar \rho(x) - \dfrac{\Phi_{\bar\rho}(x)}{L^d}  \right).
\end{equation}

For space-homogeneous stationary states $\rho_\infty(s)$, the characterisation \eqref{eq:onecompfunctional} simplifies to an easier scalar problem.

\begin{theorem}[Carrillo, Holden, Solem \cite{CHS}]
	Grant Assumption \ref{as:1}. Then, for all $\sigma>0$, there exists a unique space-homogeneous stationary state $s\mapsto\rho_\infty^\sigma(s)$ of \eqref{eqn:2}.
\end{theorem}
The actual result is stated under slightly more general assumptions in \cite{CHS}. Denoting
\[ \Phi_0 = \Phi(W \ast \bar\rho_\infty + B),\qquad  \Phi'_0 = \Phi'(W \ast \bar\rho_\infty + B), \qquad  \Phi''_0 = \Phi''(W \ast \bar\rho_\infty + B),  \]
the equation \eqref{eqn:rho_0} becomes the scalar identity
\begin{equation}\label{eqn:rho_0_0}
	 \rho_\infty(0) = \frac1\sigma\left(\bar \rho_\infty - \dfrac{\Phi_0}{L^d}  \right).
\end{equation}

\subsubsection{Noise-driven bifurcations}

As we mentioned, noise is ubiquitous in the brain and plays an important role, which was the main reason for the development of this new PDE model. An interesting question is: How robust are the grid cell patterns to noise? If increasing the noise can flatten the activity and kill all stationary states but the space-homogeneous one, then while decreasing the noise bifurcations should occur.

The first study \cite{CHS} on this PDE showed through numerical simulations that noise-driven bifurcations from the space-homogeneous stationary state do occur for the general four populations model \eqref{eq:4PDE}, and that the hexagonal grid cell pattern appears at some point via secondary bifurcations. They study a linearised version of the PDE and obtain a rigorous criterion for when the space-homogeneous stationary state can be guaranteed to be stable in the linear problem, which we will describe in more detail in Section \ref{sec:GC_lin}. The bifurcation analysis was then made rigorous in \cite{CRS} as follows, building upon methods in \cite{CGPS} for the McKean-Vlazov equation on the torus.

\paragraph{Functional analysis setup}

For simplicity in computation, the bifurcation parameter will be 
 \[\kappa = \dfrac{1}{\sigma}.\]
Using the equivalence between the zeros of $\mathcal G$ and $\bar{\mathcal G}$, we only need a functional setting for the space component. While working in the space $L^2_S(\mathbb T^d)$, we use the Hilbert basis
$
\left(\omega_k\right)_{k\in \N^d}
$
defined by
\begin{equation}\label{eqn:HilbertElements}  \omega_k(x) = \dfrac{\Theta(k)}{L^{\frac d2}} \prod_{i=1}^{d} \omega_{k_i}(x_i),  \qquad \omega_{k_i}(x) = \cos\left(\dfrac{2\pi k_i}{L} x_i\right),\qquad \Theta(k) = \prod_{i=1}^{d} \sqrt{2 - \delta_{k_i,0}},
\end{equation}
where $\delta_{i,j}$ is 0 if $i\neq j$ and 1 if $i=j$. The Fourier modes of $W$ are denoted by
\[\tilde  W(k) = \pscal{ W }{ \omega_k} = \int_{\mathbb T^d} W(x)\omega_k(x)\diff x , \quad k\in\N^d,  \]
If $W\in L_S^2(\mathbb T^d)$, for all $g\in L_S^2(\mathbb T^d)$,
\begin{equation*}
	W\ast g (x) = \int_{\mathbb T^d} W(x-y) g(y) \diff y = \dfrac{L^\frac d2}{\Theta(k)} \sum_{k\in \N^d} \tilde W(k) \tilde g(k) \omega_k(x). 
\end{equation*}
In particular, for all $k\in\N^d$,
\begin{equation}\label{trigo}
	W\ast \omega_k (x) =  \dfrac{L^{\frac d2}\tilde W(k)}{\Theta(k)} \omega_k(x) .
\end{equation}

In the larger space $L^2(\mathbb T^d)$, we use the Hilbert basis $	\left(\omega_k\right)_{k\in \Z^d} $ defined by \eqref{eqn:HilbertElements},
\[
\omega_{k_i}(x) = \left\{\begin{array}{ll}
	\cos\left(\dfrac{2\pi k_i}{L} x_i\right) & \mathrm{if}\ k_i > 0, \\
	1 & \mathrm{if}\ k_i=0,\\
	\sin\left(\dfrac{2\pi k_i}{L} x_i\right) & \mathrm{if}\ k_i < 0. \\
\end{array}\right.
\qquad \mathrm{and} \qquad \Theta(k) = \prod_{i=1}^{d} \sqrt{2 - \delta_{k_i,0}}.
\]

\paragraph{Applying the Crandall-Rabinowitz theorem}

The simplest way to prove the existence of noise-driven bifurcations is to apply the Crandall-Rabinowitz theorem to a modified version of the $\kappa-$rewritten functional
\begin{equation}\label{eq:onecompfunctionalKappa}
	\bar{\mathcal G}(\bar \rho ,\kappa)  = \bar\rho -\dfrac{1}{Z_\rho} \int_0^{+\infty} s e^{ - \kappa\frac{ \left( s - \Phi_{\bar\rho} \right)^2}{2}}\diff s,\qquad  Z_\rho =  L^d \int_{0}^{+\infty} e^{ - \kappa\frac{ \left( s - \Phi_{\bar\rho} \right)^2}{2}} \diff s,
\end{equation}
Denote $\rho_\infty^\kappa$ the space-homogeneous stationary state associated to the parameter $\kappa$, wich has firing rate $\Phi_0^\kappa$, and define the modified functional
\begin{equation}\label{eq:onecompfunctionalH}
	\begin{array}{rccl}
		\mathcal H : & L_S^2(\mathbb T^d)\times (0,+\infty) & \to & L_S^2(\mathbb T^d)  \\
		& (\bar \rho,\kappa) &\mapsto & \bar {\mathcal G} (\bar \rho_\infty^\kappa + \bar\rho, \kappa).
	\end{array}
\end{equation}
Note that for all $\kappa>0$, $\mathcal H(0,\kappa)=0$. The Fréchet derivatives of $\mathcal H$ at $(0,\kappa)$ can be computed in terms of the ones of $\bar{\mathcal G}$:
\begin{align}\label{eqn:der_H_0}
	D_{\bar\rho} \mathcal H (0, \kappa)[h] & =  D_{\bar\rho} \bar{\mathcal G}(\bar \rho_\infty^\kappa,\kappa)[h], \nonumber \\
	D_{\kappa} \mathcal H (0, \kappa) & =  D_{\kappa} \bar{\mathcal G}(\bar \rho_\infty^\kappa,\kappa) + D_{\bar\rho} \bar{\mathcal G}(\bar \rho_\infty^\kappa,\kappa)\left[\dfrac{\diff \bar\rho_\infty}{\diff \kappa}(\kappa)\right],\nonumber \\
	D^2_{\bar\rho \kappa} \mathcal H (0, \kappa)[h] & =  D^2_{\bar\rho\kappa} \bar{\mathcal G}(\bar \rho_\infty^\kappa,\kappa)[h]+ D^2_{\bar\rho\bar\rho} \bar{\mathcal G}(\bar \rho_\infty^\kappa,\kappa) \left[\dfrac{\diff \bar\rho_\infty}{\diff \kappa}(\kappa), h\right],
\end{align}
for all $h\in L_S^2(\mathbb T^d)$.
The road map for proving a bifurcation at parameter $\kappa^*$ consists in the following steps:
\begin{enumerate}
	\item[i)] Check smoothness assumptions for $\mathcal H$ and that $D_x \mathcal H(0,\kappa^*)$ is a Fredholm operator of index 0.
	\item[ii)] Check that $\dim\big(\ker(D_x \mathcal H(0,\kappa^*))\big)=1$.
	\item[iii)] Check a second order condition.\\
\end{enumerate}

\begin{theorem}[Carrillo, Roux, Solem \cite{CRS}]\label{thm:mainmain} Grant Assumption \ref{as:1}, and let $\kappa^*\in\left(\tfrac{2|W_0|^2}{L^{2d}\pi B^2},+\infty\right)$. Assume that there exists a unique $k^*\in\N^d$ such that
	\begin{equation}\label{eqn:W_tilde_equal_0}
		\dfrac{\tilde W (k^*)}{\Theta(k^*)} = \dfrac{1}{L^{\frac d2} (\Phi_0^{\kappa^*})' \left( \dfrac1{L^d} - L^d\bar \rho_\infty^{\kappa^*}  \left(\bar \rho_\infty^{\kappa^*} - \dfrac{ \Phi_0^{\kappa^*}}{L^d}\right)\kappa^* \right)}.
	\end{equation}
	Assume also that $ (\Phi_0^{\kappa^*})''\geq 0$. Then $(\bar\rho_\infty^{\kappa^*},\kappa^*)$ is a bifurcation point of $\bar{\mathcal G}(\bar\rho,\kappa)=0$. 
	
	Moreover, in a neighbourhood $U\times (\kappa^*-\varepsilon,\kappa^*+\varepsilon)$ of $(\bar\rho_\infty^{\kappa^*},\kappa^*)$ in $L^2_S(\mathbb T^d)\times \R_+^*$, the average of the stationary state is either of the form $(\bar \rho_\infty^\kappa,\kappa)$ or on the non-homogeneous solution curve
	\begin{equation*}
		\{\ (\bar \rho_{\kappa(z)},\kappa(z)) \ | \ z\in(-\delta,\delta), \ (\bar \rho_{\kappa(0)},\kappa(0))=(\bar \rho_\infty^{\kappa^*},\kappa^*),\ \delta>0 \ \},
	\end{equation*}
	with
	\begin{equation*}
		\bar \rho_{\kappa(z)}(x) = \bar \rho_\infty^{\kappa(z)} + z \omega_{k^*}(x) + o(z), \qquad x\in \mathbb T^d.
	\end{equation*}
\end{theorem}
Notice that the new branch cannot be made up of space-homogeneous stationary states due to the presence of the vector $\omega_{k^*}$ in the asymptotic expression for small $z$. Moreover, the vector $\omega_{k^*}$ that spans the kernel of the functional derivative gives information about the shape of the new pattern close to the bifurcation point.

\paragraph{Equivariant bifurcations}

When there are symmetries in the connectivity kernel $W$, the uniqueness condition on the Fourier coefficient $W(k^*)$ may not be satisfied, but a bifurcation may still be found. Based on equivariant bifurcation theory \cite{CL,DG2,FSV}, symmetries can be taken into account to give a more general classification.

\begin{theorem}[Carrillo, Roux, Solem \cite{CRS}]\label{thm:ebl}
	Under the conditions of Assumption \ref{as:1}, consider $\mathcal L$ the trivial $d-$dimensional cubic lattice of the flat torus generated by the vectors $(L,0,\dots,0), \dots, (0,\dots,0,L)$. Let $\Gamma$ be a subgroup of $E_d$. Assume $\mathcal H$ is $\Gamma-$equivariant on $L^{2}(\mathbb{T}^d)$ and that there exists a multi-index $k^*\in \mathbb N^d$ and $\kappa>0$ such that \eqref{eqn:W_tilde_equal_0} holds
	with $ (\Phi_0^\kappa)''\geq 0.$ Denote $\mathcal V = \mathrm{Ker}\big( D_{\bar\rho} \mathcal H (0,\kappa) \big) \subset L^{2}(\mathbb{T}^d)$. Then, for all isotropy subgroups $\Sigma$ of $\Gamma$, such that $ \mathrm{dim}\big(\mathrm{Fix}_{\mathcal V}( \Sigma ) \big) = 1$, there exists a unique branch of stationary states of \eqref{eqn:2} bifurcating from the trivial curve and whose average $\bar\rho$ in $s$ has the symmetry of $\Sigma$.
\end{theorem}

The application framework for grid cells is in dimension $d=2$. The Euclidean group $E_2$ associated with the cubic lattice on the torus is the compact semidirect sum of the dihedral group $D_4$ and the compact group of translations on $\mathbb T^2$, $D_4 \overset{.}{+} \mathbb T^2$. Note that the semidirect sum $\overset{.}{+}$ is in fact a semidirect product, but noted differently for notational clarity in a geometry context. Assume $W$ is radially symmetric. Then $\mathcal H$ is $\Gamma-$equivariant with $\Gamma=D_4 \overset{.}{+} \mathbb T^2$. The patterns obtained from Theorem \ref{thm:ebl} can be classified using the procedure in \cite{DG2,D}:
\begin{itemize}
    \item write the kernel $\mathcal V \subset L^2(\mathbb T^d)$ as a direct sum of $\Gamma-$irreducible subspaces $\mathcal V = \mathcal V_1 \oplus \cdots \oplus \mathcal V_n$, which in turn implies, for any isotropy subgroup $\Sigma$ of $\Gamma$,
\[  
\mathrm{Fix}_{\mathcal V}( \Sigma ) = \mathrm{Fix}_{\mathcal V_1}( \Sigma )\oplus \cdots \oplus \mathrm{Fix}_{\mathcal V_n}( \Sigma ). 
\]
    \item classify $\Gamma$-irreducible representations of $\mathcal V$ depending on its dimension.
\end{itemize}
In dimension 2, this work has already been done, for example in a general context \cite{DG2,D} or in the context of deterministic neural field models in \cite{FSV,VCF}. If $W$ is radially symmetric and the terms $\tilde W(k)/\Theta(k)$ are unique up to permutations of $k$, then $\mathrm{dim}(\mathcal V) = 4$ or $\mathrm{dim}(\mathcal V) = 8$. The method allows to prove the existence of branches bifurcating from the spatially homogeneous stationary state with the following symmetries: stripes, also called rolls \cite[Figure 4]{CHS}; simple squares \cite[Figure 1]{D}; squares \cite[Figure 2]{D}; and anti-squares \cite[Figure 3]{D}.

\begin{figure}[ht!]
	\centering
	\includegraphics[width=400pt]{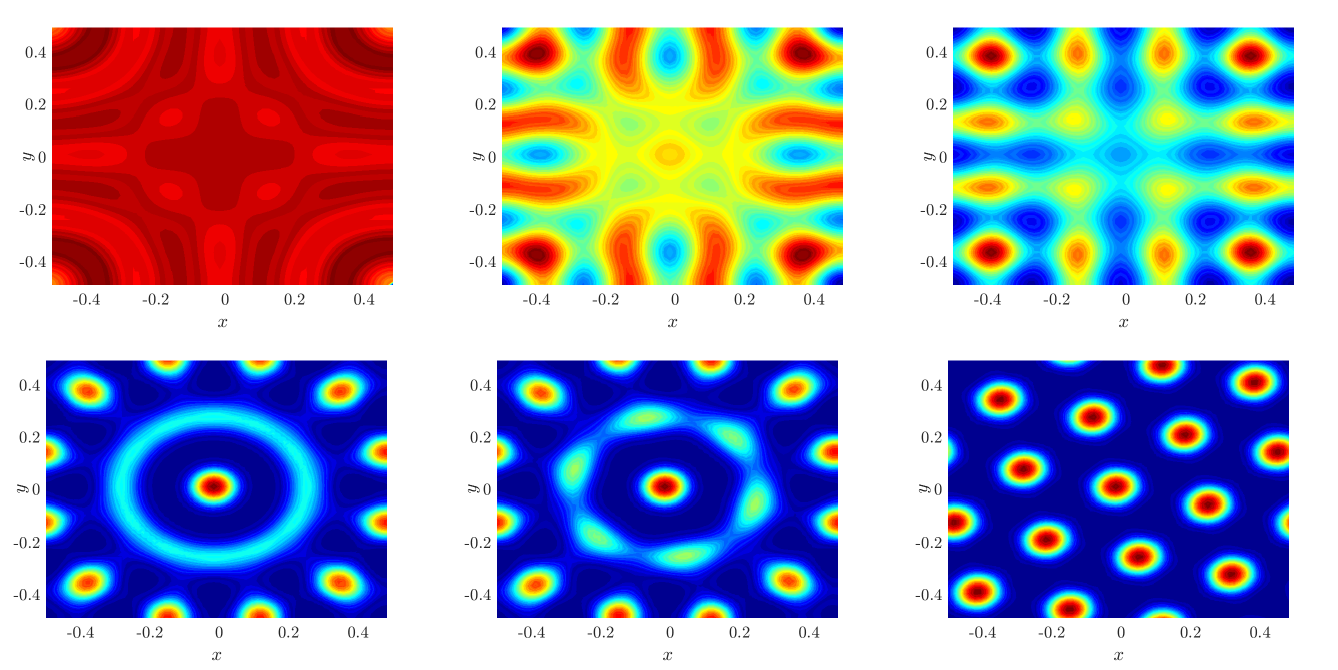}
	\caption{\textit{Numerical simulation of \eqref{eqn:2} starting from a small perturbation of a constant function, $t=40, 220, 1500, 1810, 2190,$ and $2400$ m$s$. Carrillo, Roux, Solem \cite{CRS}.}\label{fig:GC_time_evol}}
\end{figure}

Numerical results indicate that the first bifurcation points are associated to the multi-indices $k=(4,0),(0,4)$, $k=(4,1),(1,4)$, and $k=(3,3)$. Although equivariant bifurcation theory does not allow us to prove the existence of stationary states with hexagonal symmetry, it is likely that the hexagonal firing pattern appears as a secondary bifurcation, since the indexes $(4,1),(1,4),(3,3)$ are related to hexagonal patterns in many mathematical models exhibiting bifurcations \cite{M}. It can also be noted that the simple addition of these modes,
\begin{align*}
    \omega(x_1,x_2) = & \, \cos(6\pi x_1)\cos(6\pi x_2) + \cos(8\pi x_1)\cos(2\pi x_2) + \cos(2\pi x_1)\cos(8\pi x_2),
\end{align*}
looks very close visually to the transient pattern in Figure \ref{fig:GC_time_evol} at $t=1810ms$, and that
\begin{align*}
    \omega^{\mathrm{hex}}(x_1,x_2) = & \, \cos(6\pi x_1)\cos(6\pi x_2) + \cos(8\pi x_1)\cos(2\pi x_2) + \cos(2\pi x_1)\cos(8\pi x_2)\\ & 
    + \sin(6\pi x_1)\sin(6\pi x_2) - \sin(8\pi x_1)\sin(2\pi x_2) - \sin(2\pi x_1)\sin(8\pi x_2),
\end{align*}
looks very close to the final hexagonal pattern in Figure \ref{fig:GC_time_evol} at $t=2400ms$.

\subsection{Global-in-time existence of classical solutions}

In \cite{CRS2} global-in-time existence of smooth solutions is proved for both the one population system \eqref{eqn:2} and the general system \eqref{eq:4PDE} under fairly general assumptions. 

The proof consists in applying the same changes of variables than in the NNLIF existence proof. First, an affine change of variables to fix $\sigma=\tau_c=1$. We use again the variables 
$
y=e^{t}s$ and $\tau = \tfrac{1}{2}(e^{2t}-1),
$
and the function $\alpha(\tau)=(2\tau+1)^{-\tfrac12} = e^{-t}$, to define the new unknown
\[ 
q(x,y,\tau) = \alpha(\tau) \rho\Big(x,y\alpha(\tau), -\log(\alpha(\tau))\Big) \quad
\mbox{or equivalently} \quad
\rho(x,s,t) = e^t q(x,e^t s , \tfrac12 (e^{2t}-1)). 
\]
and we define $\beta(\tau) = B(-\log(\alpha(\tau)))$ and $\Psi(x,\tau)=\Phi\big( \alpha(\tau) W \ast \bar q(x,\tau) + \beta(\tau)  \big)\alpha(\tau)$.
After the additional change of variable
$
z = y - \int_{0}^{\tau} \Psi(x,\eta) \diff \eta,
$
defining $u(x,z,\tau)=q(x,y,\tau)$ and
\[ 
\gamma(x,\tau) = - \int_{0}^{\tau} \Psi(x,\eta) \diff \eta, \qquad \bar u = \int_{\gamma(x,\tau)}^{+\infty} zu(z,\tau)\diff z,
\]
we come to another nonlocal free-boundary Stefan-like problem,
\begin{equation*}\label{eq:StefanGridCell}
	\left\{\begin{array}{rcll}
		\displaystyle \dfrac{\partial u}{\partial \tau}(x,z,\tau) &=& \dfrac{\partial^2 u}{\partial z^2}(x,z,\tau), \qquad\qquad\qquad\  z\in(\gamma(x,\tau), +\infty),&  \tau\in\R_+, \\[3mm]
		\displaystyle \gamma(x,\tau) &=& - \displaystyle\int_{0}^{\tau} \Psi(x,\eta) \diff \eta, & \tau\in\R_+,\\
		\dfrac{\partial u}{\partial y}(x,\gamma(x,\tau),\tau) &=& \Psi(x,\tau) u(x,\gamma(x,\tau),\tau), &\tau\in\R_+,\\
		\Psi(x,\tau) &=& \Phi\big( \alpha(\tau) W \ast [ \bar u (x,\tau) - \gamma(x,\tau)] + \beta(\tau)\big)\alpha(\tau),\ &\tau\in\R_+,
	\end{array}\right.
\end{equation*}
with $u(x,z,0)=u^0(x,z)$ on $\mathbb T^d\times\R_+$ defined from $\rho^0$. Using again the heat kernel and Green identity
\[  G(z,\tau,\xi,\eta) = \frac{1}{\sqrt{4\pi (\tau-\eta) }}e^{-\frac{(z-\xi)^2}{4(\tau-\eta)}}, \qquad \dfrac{\partial}{\partial \xi} \left( G\dfrac{\partial u}{\partial \xi} - u\dfrac{\partial G}{\partial \xi} \right) - \dfrac{\partial }{\partial \eta}\big( Gu \big)  = 0, \]
we can integrate for each $x\in \mathbb T^d$ the identity in the domain $(\gamma(x,t),+\infty)\times (0,\tau)$; we thereby obtain
\begin{align}\label{eq:DuhamelGridCells}
	u(x,z,\tau) =\, \int_{0}^{+\infty} G(z,\tau, \xi ,0) u^0(x,\xi) \diff \xi + \int_{0}^{\tau} \dfrac{\partial G}{\partial \xi}(z,\tau,\gamma(x,\eta),\eta) u(x,\gamma(x,\eta),\eta) \diff\eta.
\end{align}

Defining the boundary value term $v(x,\tau)= u(x,\gamma(x,\tau),\tau)$ and doing some work on \eqref{eq:DuhamelGridCells}, it is possible to obtain a closed integral system 
\begin{equation*}
	\left\{\begin{array}{rcll}
		v(x,\tau)      &=&  F_v[v,\gamma, \bar u](x,\tau)\quad & \mathrm{value\ on\ the\ free\ boundary},  \\
		\gamma(x,\tau) &=&   F_\gamma[v,\gamma,\bar u](x,\tau)& \mathrm{position\ of\ the\ free\ boundary},\\
		\bar u(x,\tau) &=&  F_{\bar u}[v,\gamma,\bar u](x,\tau)& \mathrm{average\ in\ activity},
	\end{array}\right.
\end{equation*}
where the activity $s$ has been eliminated. It allows to prove local-in-time existence of a solution to the original problem by application of the Banach fixed point theorem. \textit{A priori} bounds and moment analysis lead to the conclusion that solutions are global-in-time for any Lipschitz continuous $\Phi$.
\begin{theorem}[Carrillo, Roux, Solem \cite{CRS2}]\label{thm:mainCRS2}
	Assume $\Phi$ is locally Lipschitz, $B\in C^0(\R_+)$ and $W\in L^1(\mathbb T^d)$. For any compatible initial condition $\rho^0$,
	equation \eqref{eqn:2} admits a unique maximal solution in the sense of Definition \ref{def:classicalGridCells}.
	Moreover, if $\Phi$ is globally Lipschitz, the solution is global-in-time and
	\begin{itemize}
		\item[$\bullet$] if $\Phi$ non-negative, then
		\begin{equation*}
			\forall t\in\R_+,\ \forall x\in\mathbb T^d, \quad \rho(x,0,t) \leqslant \dfrac{1}{L^d}\sqrt{ \dfrac{2}{ \pi \sigma } } \dfrac1{\sqrt{1-e^{-2\tau_c t}}};
		\end{equation*}
		\item[$\bullet$] if $\Phi$ is non-positive, then
		\[ \forall t\in \R_+,\ \forall x\in \mathbb T^d, \quad  \bar\rho(x,t) \leqslant  \dfrac{1}{L^d}(1+\sigma - \sigma e^{-2\tau_c t}) +  e^{-2\tau_c t} \sup_{x\in\mathbb T^d}\int_0^{+\infty} s^2\rho^0(x,s)\diff s . \]
	\end{itemize}
\end{theorem}
The result is then extended to the four population Grid cells model \eqref{eq:4PDE} with similar assumptions and definitions.

\subsection{Asymptotic stability of stationary states}

The relative entropy method \cite{MV} can be used to study the nonlinear stability of the stationary states of \eqref{eq:4PDE} and \eqref{eqn:2}. As above, the idea is to study it first for the simpler \eqref{eqn:2} and then to generalise to \eqref{eq:4PDE}, with the exception of linear stability, which was performed directly on \eqref{eq:4PDE} in \cite{CHS}.

An important quantity for the study of asymptotic stability for these models is the variance associated to any stationary state $\rho_\infty$:
\begin{equation*}
    M_\infty(x) = \int_0^{+\infty} \big(s-L^d\bar\rho_\infty(x)\big)^2 \rho_\infty(x,s) \diff s.
\end{equation*}
As proved in \cite{CRS}, we can obtain the following identity
\begin{equation*}
    \dfrac{M_\infty}{\sigma} = \frac{1}{L^d} g\left( \frac{\Phi_0}{\sqrt{2\sigma}}\right), \qquad g(\eta) = 1-\frac{2}{\sqrt{\pi}}\frac{\exp (-\eta^2)}{1+\erf (\eta)}\left[\frac{1}{\sqrt{\pi}}\frac{\exp (-\eta^2)}{1+\erf (\eta)}+\eta\right]
\end{equation*}
where $g$ is a differentiable positive bounded function whose image lies in $(1-\frac{2}{\pi},1)$.

\subsubsection{Naive entropy method}

The easiest way to get nonlinear stability is to do, like for the NNLIF model, a weakly-nonlinear analysis \textit{via} a Lyapunov control on the dissipation of the relative entropy. The key step is to control the non-local and nonlinear effect of the firing rate.
\begin{theorem}[Carrillo, Roux, Solem \cite{CRS2}]\label{prop:stab1}
	Let $\Phi$ be Lipschitz, $W \in L^2(\mathbb T^d)$, and $\rho_\infty$ be a stationary state of \eqref{eqn:2}. Assume that $\sigma$ satisfies 
\begin{align}\label{eq:poincareestimate}
		\|\Phi' \|_\infty \frac{\norme{W}_{L^2({\mathbb T^d})}}{L^\frac d2} \norme{ M_\infty^\frac12}_\infty < \frac{\sigma}{2} \tilde\gamma (\rho_\infty )^\frac12,
	\end{align}
	where $\tilde\gamma (\rho_\infty ) = \inf_{x \in \mathbb T^d} \gamma(\rho_\infty (x))$, and $\gamma(\rho_\infty(x))$ is a Poincar\'e constant associated to $\rho_\infty (x)$. Then, the $L^1$ norm of the relative entropy 
	\begin{align}\label{l2relentropy}
		\int_{\mathbb T^d} \int_{0}^{+\infty} \left( \dfrac{ \rho(x,s,t) - \rho_\infty(x,s) }{ \rho_\infty(x,s) } \right)^2 \rho_\infty(x,s) \diff s \diff x,
	\end{align} 
	decays exponentially fast whenever the compatible initial datum $\rho_0$ is close enough to the stationary state $\rho_\infty$ in relative entropy.
\end{theorem}

The idea of the proof is to write $ h =\frac{\rho}{\rho_\infty}$ and choose $G(h) = \frac12 (h-1)^2$. Then the entropy dissipation can be bounded by
	\begin{align*}
		  \frac{d}{dt}\int_0^{+\infty} \frac{1}{2}(h-1)^2 \rho_\infty \diff s  
		\leqslant & \ \frac{1}{2\sigma}\left|\Phi_{0}-\Phi_{\bar{\rho}}\right|^2 \int_0^{+\infty} (h-1)^2 \rho_\infty  \diff s \\&\ -\frac{\sigma}{2} \int_0^{+\infty} \left( \dfrac{\partial h}{\partial s} \right)^2\rho_\infty  \diff s -  \left(\Phi_{0}-\Phi_{\bar{\rho}}\right)\int_0^{+\infty}  \dfrac{\partial h}{\partial s}\,  \rho_\infty  \diff s.
	\end{align*}
    The deviation of the firing rate from the stationary firing rate can be controlled by
	\begin{align*}
		\left|\Phi_{0}-\Phi_{\bar{\rho}}\right| &\leq  \underbrace{\norme{\Phi'}_\infty \norme{W}_{L^2({\mathbb T^d})} \sup_{x \in {\mathbb T^d}} M_\infty^\frac12(x)}_{ C(\Phi,W, \rho_\infty)} \left(\int_{\mathbb T^d} \int_0^{+\infty} (h-1)^2 \rho_\infty \diff s \diff y \right)^\frac12.
	\end{align*}
By applying a Poincar\'e inequality, we arrive at 
\begin{multline*}
\frac{d}{dt}\int_{\mathbb T^d}\int_0^{+\infty} \frac{1}{2}(h-1)^2 \rhostat \diff s \diff y 
\leq \,\frac{1}{2\sigma} C^2 \left(\int_{\mathbb T^d} \int_0^{+\infty} (h-1)^2 \rhostat \diff s \diff y \right)^2 \\
 + \left(C-\frac{\sigma}{2} \tilde\gamma (\rhostat)^\hf\right) \left(\int_{\mathbb T^d}\int_0^{+\infty} | \partial_s h |^2\rhostat  \diff s \diff x\right)^\hf \left(\int_{\mathbb T^d} \int_0^{+\infty} (h-1)^2 \rhostat \diff s \diff y \right)^\hf.
\end{multline*}
	By assumption \eqref{eq:poincareestimate}, $\tilde C:=C-\frac{\sigma}{2} \tilde\gamma (\rho_\infty)^\frac12 < 0$. Thus, applying again the same Poincar\'e inequality,
	\begin{align*}
		\frac{d}{dt}\int_{\mathbb T^d} & \int_0^{+\infty} \frac{1}{2}(h-1)^2 \rho_\infty \diff s \diff y \\
		& \leq \left[\frac{C^2}{2\sigma}  \int_{\mathbb T^d} \int_0^{+\infty} (h-1)^2 \rho_\infty \diff s \diff y + \tilde\gamma (\rho_\infty)^\frac12\tilde C \right] \int_{\mathbb T^d} \int_0^{+\infty} (h-1)^2 \rho_\infty \diff s \diff y.
	\end{align*}
	The expression between the brackets will be negative as long as the initial relative entropy is small enough, and then there is exponential convergence in relative entropy.\\

As the condition \eqref{eq:poincareestimate} depends nonlinearly on $\sigma$, it is not easy to check theoretically. Results on Poincar\'e constants \cite{M2} and bounds on $M_\infty$ provide the clarified (less sharp) condition
\[ \frac1{L^d}\norme{\Phi'}_\infty \norme{W}_{L^2(\mathbb T^d)} < \frac12. \]

The problem with this result is that, as it can be checked numerically, it applies only in cases where there is a unique stationary state (the space-homogeneous one). As only a norm of $W$ is involved, nothing can be concluded about nonlinear stability of the constant stationary state in terms of the shape of the connectivity kernel.

\subsubsection{Linear stability of the space-homogeneous stationary state}\label{sec:GC_lin}

Stability of the space-homogeneous stationary state was first investigated in \cite{CHS} by linearising the equation \eqref{eq:4PDE}. Let $\rho_\infty$ be a space-homogeneous stationary state. Denoting $\rho^\beta = \rho_\infty + u^\beta$, a linearised equation for the perturbation $(u^1,u^2,u^3,u^4)$ is,
\begin{equation}
    \tau \dfrac{\partial u^\beta}{\partial t} + \dfrac{\partial \rho_\infty}{\partial s}\dfrac{\Phi_0'}{4}\sum_{\beta'=1}^4 W^{\beta'} \ast \bar u^{\beta'} + \dfrac{\partial}{\partial s}\left[ \left(\Phi_0 - s\right) u^\beta \right] - \dfrac{\partial^2 u}{\partial s^2} = 0.
\end{equation}
Consider here an exponential basis for Fourier modes (not to be confused with the cosine and sine Fourier modes $\tilde W(k)$ above), for multi-indexes $k\in\Z^d$,
\begin{equation}
    \hat W(k) = \int_{\mathbb T^d} W(x) e^{-i 2\pi k\cdot x}\diff x.
\end{equation}
Assume that they satisfy $F(k)<\frac{\sigma}{M_\infty}$, with
\begin{equation}
    F(k) = \dfrac{\Phi_0'}{4} \hat W(k) \sum_{\beta=1}^4 e^{-2i\pi k\cdot r^\beta}.
\end{equation}
Notice that these quantities are real due to the symmetry by component of the connectivity $W$ and the shifts $r^\beta$.
If the linear perturbation $u$ is taken in $L^2(\mathbb T^d\times \R_+)$ and decomposed in plane waves as
\[ u^\beta(x,s,t) = \rho_\infty(s) \sum_{k\in \Z^d} u_k^\beta (s,t) e^{-2i\pi k\cdot x},  \]
then, we can reduce the $L^2$-relative entropy \eqref{l2relentropy} to estimate
\begin{align}\label{l2relentsinglemode}
		\int_0^{+\infty}  \rho_\infty |U|^2 \diff s \qquad \mbox{with } \qquad U = \sum_\beta \exp(-i\bk \cdot {\bf r}^\beta ) u^\beta
\end{align}
for all $k\in \Z^d$. However, computing the evolution of this linearized $L^2$-entropy for one Fourier mode \eqref{l2relentsinglemode} over the linearized dynamics leads to a dissipation of the quantity
$$
\int_0^{+\infty} \rho_\infty  |U|^2 \diff s - \frac{F(\bk)}{\sigma} \left|\int_0^\infty \rho_\infty  U \diff s\right|^2
$$
This quantity can be rewritten as 
\begin{align}\label{l2relentsinglemode2}
\int_0^{+\infty} \rho_\infty  |V|^2 \diff s + \left( 1 - M_\infty \frac{F(\bk)}{\sigma} \right) M_\infty |c(t)|^2
\end{align}
by expanding $U$ into $U = V + c(t)(s-\bar \rho_\infty )$. Here $\int_0^{+\infty} V \rho_\infty ds = 0$ and $c(t)$ is chosen such that $\int_0^{+\infty} V (s-\bar \rho_\infty ) \rho_\infty ds = 0$ leading to
\begin{equation}
c(t)=\frac1{M_\infty} \int_0^\infty \rho_\infty  U \diff s,
\end{equation}
Detailed computations in \cite[Theorem 3.6]{CHS} allowed to estimate both terms in the quantity \eqref{l2relentsinglemode2}
under the assummption \eqref{linstabcond}
leading to the following result.

\begin{theorem}[Carrillo, Holden, Solem \cite{CHS}]
    Grant Assumption \ref{as:1} and assume $d=2$. The unique spatially-homogeneous stationary state $\rho_\infty$ of \eqref{eq:4PDE} is linearly asymptotically stable in $L^2(\mathbb T^d\times \R_+)$ as long as for all multi-index $k\in \Z^d$,
\begin{equation}\label{linstabcond}
        F(k) < \dfrac{\sigma}{M_\infty}.
    \end{equation}
\end{theorem}
By taking all shifts equal to zero, $B^\beta=B$ and all components $\rho^\beta_\infty$ equal, this implies a linear stability result for \eqref{eqn:2}. It has then been proved in \cite{CRS} that the loss of linear stability perfectly matches the first point of bifurcation on the space-homogeneous stationary states curve.

\subsubsection{Nonlinear stability of the space-homogeneous stationary state}

In order to go from linear stability to nonlinear stability, a possible method is to find nonlinear quantities which would, after linearisation, yield quantities similar to the linearly decaying \eqref{l2relentsinglemode}. With this in mind, define
\begin{align*}
	E(x,t) =  \frac{1}{2}\int_0^{+\infty}  \left(\frac{\rho- \rho_\infty}{ \rho_\infty}\right)^2 \rho_\infty \diff s, \quad \mathrm{and} \quad H(x,t) =    \Phi^\delta \bar \rho^\delta,
\end{align*}
where $ \Phi^\delta = \Phi_{\bar \rho}-\Phi_0$, $ \bar \rho^\delta = \bar \rho - \bar \rho_\infty$. Linearising $\Phi^\delta$, we get
\begin{align*}
	2E-\frac{H}{\sigma} \simeq  \int_0^{+\infty} \left(\frac{\rho- \rho_\infty}{ \rho_\infty}\right)^2 \rho_\infty \diff s - \frac{\Phi'_0}{\sigma} W\ast(\bar \rho - \bar \rho_\infty) (\bar \rho - \bar \rho_\infty).
\end{align*}
Define then
	\begin{align*}
		\mathcal{Q}=\int_{\mathbb T^d} 2E-\frac{H}{\sigma} \diff x, \qquad 
		\mathcal E(t) = \int_{\mathbb T^d} \int_{0}^{+\infty} \left( \dfrac{ \rho(x,s,t) - \rho_\infty(x,s) }{ \rho_\infty(x,s) } \right)^2 \rho_\infty(x,s) \diff s \diff x.
	\end{align*}
	It is possible to show under the condition \eqref{eq:nonlinstab-condition} that there are constants $\alpha,C_1,C_2,C_3>0$ such that
	\begin{equation*}
		\frac{\diff}{\diff t} \mathcal{Q}\leqslant -\sigma C_1  \alpha \left(\alpha - C_2 \mathcal{Q}^\frac12 \right) \mathcal{Q}
	\end{equation*}
	and,
	\begin{equation}\label{eq:GC_equiv_Q_E}
		2\alpha \mathcal{E}(t) \leqslant \mathcal Q(t) \leqslant C_3\mathcal{E}(t).
	\end{equation}
    The differential inequality allows to prove nonlinear decay of $\mathcal Q$ to 0 at exponential speed, provided $\mathcal Q(0)$ is small enough. The second inequality proves that the relative entropy $\mathcal E$ and the quantity $\mathcal Q$ are equivalent. In the proof of \eqref{eq:GC_equiv_Q_E}, the condition \eqref{eq:nonlinstab-condition} appears naturally as a requirement.

\begin{theorem}[Carrillo, Roux, Solem \cite{CRS2}]\label{thm:GC_CV_2}
		Let $\Phi \in C^2(\mathbb R), \nabla \Phi\in \mathcal W^{1,\infty}(\R)$, $W \in L^2_S(\mathbb T^d)$, and $\rho_\infty$ be a spatially homogeneous stationary state. Assume for a suitably small $\alpha\in(0,1)$, that $\sigma>0$ satisfies
		\begin{align}\label{eq:nonlinstab-condition}
			\int_{\mathbb T^d} (1-\alpha)g^2(x) - \frac{\Phi^\delta_g(x)}{\sigma} g(x) \, \diff x > 0, \quad g \in L^2(\mathbb T^d),
		\end{align}
		where $$\Phi^\delta_g(x) =\Phi\big( M_\infty W\ast g (x) +W_0 \bar\rho_\infty +B\big)-\Phi_0.$$ Then, if $\mathcal E(0)$ is small enough, $\mathcal E(t)$ converges to 0 exponentially fast. 
	\end{theorem}

In order to better understand the link between condition \eqref{eq:nonlinstab-condition} and the linear stability condition, take for example the simplest case $\Phi(z)=z$ and send formally $\alpha$ to 0; then
		\[\Phi^\delta_g(x) = M_\infty W\ast g (x) + W_0 \bar\rho_\infty +B - \Phi_0 = M_\infty W\ast g (x).\]
		and condition \eqref{eq:nonlinstab-condition} becomes
		\begin{align}
			\int_{\mathbb T^d} g^2(x) >  \frac{ M_\infty}{\sigma} \int_{\mathbb T^d}\int_{\mathbb T^d} W(x-y) g(y)g(x) \, \diff y \diff x, \quad g \in L^2(\mathbb T^d),
		\end{align}
		If all Fourier modes of $W$ are negative (implying unconditional global linear stability), then for all $g\in L^2(\mathbb T^d)$,
		 \[  0 > \int_{\mathbb T^d}\int_{\mathbb T^d} W(x-y) g(y)g(x) \, \diff y \diff x \]
		 and there is also local nonlinear stability. Condition \eqref{eq:nonlinstab-condition} thus appears as a perturbation of the linear stability condition.

   In \cite{CRS2}, the asymptotic stability results for \eqref{eqn:2} are extended to the model \eqref{eq:4PDE} under additional technical assumptions.

\section{Rate models for decision making tasks}

Understanding the mechanisms underlying decision-making has been one of the fundamental questions in behavioral neuroscience for the last decades. Optimizing tasking behaviour involves the constant making of decisions between alternative choices. Neurophysiological and theoretical neurosciences have started to reveal the neural mechanisms underlying decision-making. Identifying areas of the cortex of animals involved in particular decisions has been one of the cornerstones of neuroscience. The most basic studies focused on binary decision-making tasks, i.e., the subject is asked to make a choice between two alternatives according to an experimentally defined criterion based on sensory input information. The difficulty of making a decision can be manipulated by the level of uncertainty in the discrimination of the sensory input information influenced by regulating the signal to noise ratio. 

In a wide range of such decision-making experiments, the behavioral response, i.e. performance and reaction times, can be properly described by a simple stochastic diffusion model or Fokker-Planck equation \cite{RS04}. It is plausible to think that the underlying neuronal system perform decision-making by accumulation of evidences. Indeed, electrophysiological recordings from awake behaving monkeys performing decision-making have shown that trial-averaged spiking activity of some neurons shows a gradual, nearly linear increase shortly after stimulus onset until just before the response is made \cite{RS01,RS03,SN01,RS02,HS05}. One example of these binary decision making tasks is the so-called random-dot task. The random-dot task is a two-choice visual decision-making task. In this task, the subject is confronted with a visual stimulus consisting of a field of randomly moving dots. On any one trial, a fixed fraction of the dots, determined by the experimentalist, moves coherently in one of two directions. The subject must discriminate in which direction the majority of dots are moving. Electrophysiological recordings from awake behaving monkeys performing this task have shown that trial-averaged spiking activity of neurons in the lateral intraparietal (LIP) cortex reflects the accumulation of information mentioned above \cite{SN01,RS02,HS05}. The difference in the firing rates of cells with opposing preferred directions is a likely measure of the available visual evidence for the direction of coherent motion of the moving dot stimulus, further indicating that the area LIP is involved in this decision-making task.

A computational model of the cortical circuit underlying LIP cell activity has been proposed and studied numerically \cite{W02,WW,LW}. This model was able to explain qualitatively the experimentally observed trial-averaged spiking
activity. The model consists of two populations of spiking neurons, within which interactions are mediated by excitatory synapses and between which interactions occur principally through an intermediary population of inhibitory interneurons. Sensory input may bias competition in favour of one of the populations, potentially resulting in a gradually developing decision in which neurons in the chosen population exhibit increased activity while activity in the other population is inhibited. Here, activity is measured as the global firing rate of each of the populations of spiking neurons. When the activity of one of the two populations exceeds a pre-defined threshold, then a behavioural response is generated. 

\begin{figure}[ht!]
\centering{\includegraphics[width=9cm]{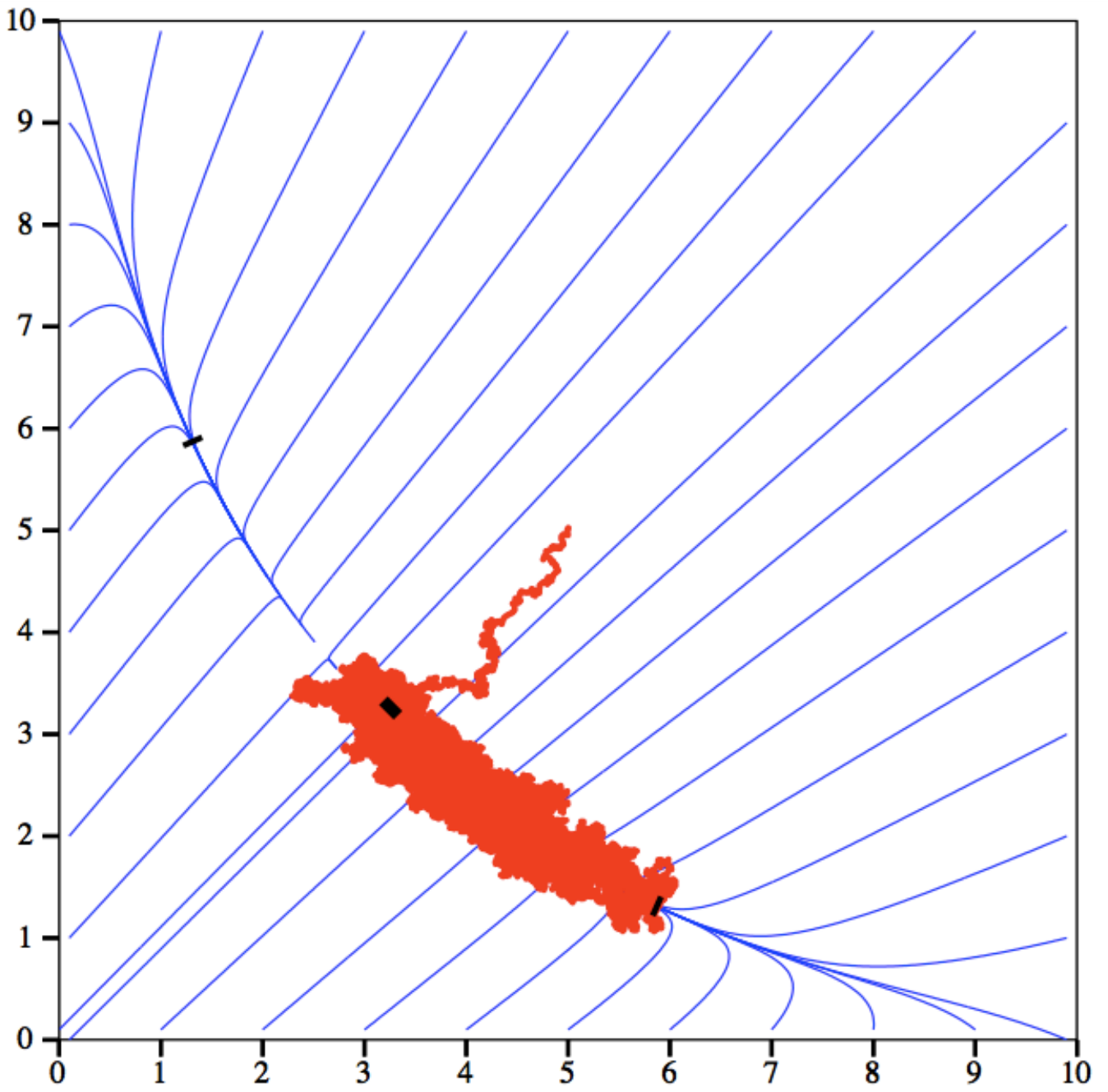}}
\caption{Dynamics in the unbiased ($\lambda_i=0$) case of \eqref{ODE} in the
deterministic system $\beta=0$ (straight lines) and the stochastic system $\beta>0$
(wiggled line). Straight lines emphasize the location of the
equilibrium points. The black
bars are the equilibrium points of the deterministic differential system. The single particle realization of the stochastic system starts its dynamics at point $(5,5)$, moves almost straight forward towards
a slow-manifold nearby the spontaneous state, and then slowly oscillates towards one of the
decision states. Carrillo, Cordier, Mancini \cite{CCM}.} \label{sde}
\end{figure}

There are two possible mechanisms by which the decision can be made which correspond to two different dynamical working points, see Figure \ref{sde}. In one mechanism, the appearance of the stimulus destabilises the so-called spontaneous state, in which all neurons show low activity. In the other mechanism, even after stimulus presentation, the network can still sustain low activity, but random fluctuations eventually drive the collective activity of the network to an activated state, in which a fraction of the neurons fires at a high rate. Let us remark that the behavioural performance of the model (i.e. the fraction of trials in which the 'correct' direction is chosen) and corresponding reaction time to make a decision match well the data from both monkeys and humans performing the moving-dot task \cite{RL}. Therefore, it may be possible to establish a direct link between the neuronal and behavioural correlates of decision making. 

In order to perform the explicit link between a neuronal circuit processing decision-making and a diffusion model able to describe the behavioural level, we consider a mean-field approximation consisting of two competing rate models. The precise model uses a Wilson-Cowan \cite{WC1} type system, describing the time evolution of two population of neurons firing rates $\nu_i$, $i=1,2$:
\begin{equation}\label{ODE}
\tau \frac{\diff \nu_i(t)}{\diff t} = - \nu_i(t) + \phi \left( \lambda_i +
\sum_{j=1,2} w_{ij} \nu_j(t) \right) + \beta \mathcal B_{i}(t), \quad i=1,2,
\end{equation}
where $\tau$ is the typical time relaxation and $\mathcal B_{i}(t)$,
$i=1,2$, represent normalized independent Brownian motions.
In \eqref{ODE} the function $\phi(x)$ has a sigmoidal shape
determining the response function of the neuron population to a
mean excitation $x$ given by $x_i(t)=\lambda_i+\sum_j w_{ij}
\nu_j$, $i=1,2$, in each population:
\begin{equation}\label{phi}
\phi(x)=\frac{\nu_c}{1+{\rm exp} (-\alpha (x/\nu_c -1))},
\end{equation}
where $\lambda_i$ are the external stimuli applied to each neuron
population and $w_{ij}$ are the connection coefficients. The
parameter $\nu_c$ represents both the maximal activity rate of the
population and the frequency input needed to drive the population
to half of its maximal activity.

Following \cite{RBW,LRLSF,DM}, we assume that neurons
within a specific population are likely to correlate their
activity, and to interact \textit{via} strong recurrent excitation with a
dimensionless weight $w_+>1$ greater than a reference baseline
value established to 1. Analogously, neurons in two different
populations are likely to have anti-correlated activity expressed
by a excitatory weight lower than the baseline, $w_- < 1$.
Furthermore, we assume that there is global feedback inhibition,
as a result of which all neurons are mutually coupled to all other
neurons in an inhibitory fashion (we will denote this inhibitory
weight by $w_I$). As a result, the synaptic connection
coefficients $w_{ij}$, representing the interaction between
population $i$ and $j$, are the elements of a $2\times 2$
symmetric matrix $W$ given by
$$
W=\left[
\begin{array}{ccc}
w_{+}-w_{I} &\phantom{t} & w_{-}-w_{I}\\
w_{-}-w_{I} &\phantom{t} & w_{+}-w_{I}
\end{array}
\right],
$$
The typical synaptic values considered are such
that $w_-<w_I<w_+$ leading to cross-inhibition and
self-excitation.

Standard stochastic calculus methods show that the probability density $p=p(t,\nu)$ of finding the neurons of both populations firing at
averaged rates $\nu=(\nu_1, \nu_2)$ at $t>0$, satisfies a Fokker-Planck equation of the form
\begin{equation}\label{FPold}
\partial_t p + \nabla \cdot \left( \left[-\nu + \Phi(\Lambda + W \cdot \nu) \right]p \right) - \frac{\beta^2}2 \Delta p =0
\end{equation}
where $\nu  \in \Omega=[0, \nu_m] \times [0,\nu_m]$,
$\Lambda=(\lambda_1,\lambda_2)$, $\Phi(x_1,x_2)=(\phi(x_1), \phi
(x_2))$, $\nabla  = (\partial_{\nu_1}, \partial_{\nu_2})$ and
$\Delta = \Delta_\nu$. We choose to complete equation
\eqref{FPold} by the following no-flux boundary conditions:
\begin{eqnarray}\label{BCold}
\left( \left[-\nu + \Phi(\Lambda + W \cdot \nu) \right] p - \frac{\beta^2}2
\nabla p \right) \cdot n =0
\end{eqnarray}
where $n$ is the outward normal to the domain $\Omega$.
Physically, these boundary conditions mean that neurons cannot
spike with arbitrarily large firing rates and thus there is a
typical maximal value of the averaged firing rate $\nu_m$ and that
the solution to \eqref{FPold} is a probability density function,
i.e.,
\begin{equation}\label{normalisation}
\int_\Omega p(t,\nu) \, \diff \nu = 1 .
\end{equation}

In order to simplify notations, we will from now on consider the
vector field $F= (F_1,F_2)$, representing the flux in the
Fokker-Planck equation:
\begin{equation}\label{flux}
F = -\nu + \Phi(\Lambda + W \cdot \nu) = \left(
\begin{array}{c}
-\nu_1 + \phi(\lambda_1 + w_{11}\nu_1 + w_{12}\nu_2) \\
-\nu_2 + \phi(\lambda_2 + w_{21}\nu_1 + w_{22}\nu_2)
\end{array}\right)
\end{equation}
then, equation \eqref{FPold} and boundary conditions \eqref{BCold} read:
\begin{equation}\label{FP}
\partial_t p + \nabla \cdot \left( F\, p  - \frac{\beta^2}2 \nabla p \right) =0
\end{equation}
\begin{equation}\label{BC}
\left( F\, p - \frac{\beta^2}2 \nabla p \right) \cdot n =0
\end{equation}
The nonlinear nature of the SDE system \eqref{ODE},
however, hinders analytical progress. For this reason, the main analysis of such noise driven
probabilistic decision-making systems remains based on numerical investigations, which are time consuming because of the need for sufficiently many trials to generate statistically meaningful data. This is an advantage of the mean-field Fokker-Planck equation due to the existence of efficient numerical solvers. Moreover, the Fokker-Planck framework \eqref{FP}-\eqref{BC} is amenable to nonlinear analysis.

\subsection{Analysis of the rate model}

The analysis of \eqref{FP}-\eqref{BC} was performed in \cite{CCM} where we refer for further details. Let us elaborate on the most important structural property of this model. Assume that $u_1$ and $u_2$ are smooth solutions of \eqref{FP}-\eqref{BC} with $u_1>0$ in $\bar\Omega$. Given any convex smooth function $H=\HH$, where $\omega=u_2/u_1$, we have
\begin{equation}\label{entropy}
\frac{\diff}{\diff t}\left[ u_1 \HH \right] =  \dfrac{\beta^2}{2}\left( \Delta \left( u_1\HH  \right) - 
u_1 H''\left( \omega \right) \left| \nabla \omega 
\right|^2 \right) - \nabla \cdot \left[F u_1 \HH \right]\,.
\end{equation}
In order to check this, let us develop the left hand side of \eqref{entropy},
using \eqref{FP} to obtain
\begin{align*}
\displaystyle \frac{\diff }{\diff t}\left[u_1 \HH \right]=& \, - \nabla \cdot \left[F u_1 \HH  \right] + u_1 F \cdot \nabla  \HH  \\
&+ \frac{u_2}{u_1}  H' \left(\omega \right) \nabla \cdot (u_1 F) -  H' \left(\omega \right) \nabla \cdot (u_2 F)   \\
&+\bb \left[  \HH  - \frac{u_2}{u_1}  H' \left(\omega \right)
\right] \Delta u_1 + \bb  H' \left(\omega \right) \Delta u_2 \, .
\end{align*}
We separate now the computation in two parts: the one concerning
first order derivatives (I) and the one concerning second order
derivatives (II). We start by computing (I). Leaving the first
term unchanged and developing the following ones, we get:
\begin{align*}
(I) = &\,-\nabla \cdot \left[F u_1 \HH  \right] + H' \left(\omega \right)F \cdot \left(\nabla u_2 - \frac{u_2}{u_1} \nabla u_1 \right) \\
&+ \frac{u_2}{u_1} H' \left(\omega \right) \nabla u_1 \cdot  F + H' \left( \omega \right) u_2\nabla \cdot F \\
&- H' \left(\omega \right) \nabla u_2 \cdot F - H' \left(
\omega \right) u_2 \nabla \cdot F= -\nabla \cdot \left[ F u_1 
\HH  \right].
\end{align*}
Concerning the second part (II), we can compute $\Delta \left( u_1 \HH  \right) $ to obtain the relation
\begin{align*}
\Delta \left( u_1 \HH  \right)- u_1 \left| \nabla \omega \right|^2 H'' \left(\omega \right) =   \left[ \HH  - \frac{u_2}{u_1} H' \left(\omega \right) \right] \Delta u_1 + H' \left(\omega \right) \Delta u_2 
\end{align*}
Multiplying by $\bb$ and collecting  (I) and (II) completes the proof of \eqref{entropy}.
Let us now define the relative entropy between two solutions as
\begin{equation}\label{Hcal}
{\mathcal H} (u_2 | u_1) = \int_\Omega u_1 \HH  \, \diff \nu\,.
\end{equation}
For any $u_1$ and $u_2$ strong solutions of \eqref{FP}-\eqref{BC},
and any $H$ smooth convex function then, the relative entropy is
decreasing in time and
\begin{equation}\label{decrease}
\frac{d}{dt} {\mathcal H}(u_2 | u_1) = -\bb {\mathcal D}(u_2 |
u_1)=-\bb \int_\Omega u_1 H'' \left(\omega\right) \left| \nabla \left( \frac{u_2}{u_1} \right) \right|^2 \, \diff \nu \leq 0 .
\end{equation}
This is simply obtained by integrating \eqref{entropy} over the domain $\Omega$ and showing that the boundary terms disappear due to \eqref{BC}.

Under an additional assumption on the flux $F$ one can prove that the stationary problem has a unique classical solution. More precisely, assuming that
\begin{equation}\label{bcf}
F\in C^1(\bar\Omega,\R^2) \,\, \mbox{ with } \,\quad F\cdot n < 0
\quad \mbox{on } \partial \Omega ,
\end{equation}
there exists an unique probability density function $u_\infty \in H^4(\Omega)$, $u_\infty(\nu)>0$ in
$\bar\Omega$ satisfying
\begin{equation} \label{meq}
\begin{cases}
 \nabla \cdot \left( F\, p  - \frac{\beta^2}2 \nabla p \right)= 0 \quad \textrm{in} \; \Omega,\\[2mm]
\displaystyle \left( Fu -\bb \nabla u\right) \cdot n =0  \quad
\textrm{on} \; \partial\Omega.
\end{cases}
\end{equation}

The consequences of the existence of this family of Liapunov functionals given in \eqref{decrease} for \eqref{FP}-\eqref{BC} have already been explored for several equations in \cite{MMP2,MMP1} and it is usually referred as the general relative entropy (GRE) inequalities. The same conclusions apply here as in those cases apply here and we just list them here refering to \cite{CCM} for more details.

\begin{corollary}[Carrillo, Cordier, Mancini \cite{CCM}]
Given $F$ satisfying \eqref{bcf} and any solution $u$ with normalized initial data $u_0$ to \eqref{FP}-\eqref{BC}, then the following
properties hold:
\begin{itemize}
\item[i)] Nonnegativity: If $u_0$ is nonnegative, then the solution $u$ of problem \eqref{FP}-\eqref{BC} is nonnegative.
\item[ii)] Contraction principle:
\begin{equation}
\int_\Omega |u(t,\nu)| \,\diff \nu \leq \int_\Omega |u_0(\nu)| \,\diff \nu.
\label{eq:contraction}
\end{equation}
\item[iii)] $L^p$ bounds, $1<p<\infty$:
\begin{equation}
\int_\Omega u_\infty(\nu)
\left|\frac{u(t,\nu)}{u_\infty(\nu)}\right|^p \,\diff \nu \leq
\int_\Omega u_\infty(\nu)
\left|\frac{u_0(\nu)}{u_\infty(\nu)}\right|^p \,\diff \nu.
\label{eq:lp}
\end{equation}
\item[iv)] Pointwise estimates:
\begin{equation}
\inf_{\nu\in\Omega} \frac{u_0(\nu)}{u_\infty(\nu)} \leq
\frac{u(t,\nu)}{u_\infty(\nu)} \leq \sup_{\nu\in\Omega}
\frac{u_0(\nu)}{u_\infty(\nu)}\, . \label{eq:bound}
\end{equation}
\end{itemize}
\label{apriori}
\end{corollary}

This corollary is a consequence of the GRE inequality \eqref{decrease} with $H(\omega) = \omega^-$ (the negative part of $\omega$), $H(s)=|s|$, $H(s)=|s|^p$, and $H(s)=(s-k)_+^2$
respectively by approximation from smooth convex functions.
Moreover, the GRE inequality gives the convergence of the solution
$u(t)$ to the stationary state $u_\infty$.

\begin{corollary}[Long time asymptotic]
Given $F$ satisfying \eqref{bcf} and any solution $u$ with
normalized initial data $u_0$ to \eqref{FP}-\eqref{BC}, then
\begin{equation}
\lim_{t\to \infty} \int_\Omega |u(t,\nu)-u_\infty(\nu)|^2 \,\diff \nu =
0 \, . \label{eq:lta}
\end{equation}
\end{corollary}

Convergence rates have been identified in terms of Poincar\'e type inequalities in \cite{CMT}.

\subsection{Reduction to a diffusion model}

One dimensional Fokker-Planck effective reductions near bifurcation points were first obtained in \cite{RL}. They proposed to locally approximate the evolution on the slow manifold at the bifurcation point by Taylor expanding the nonlinearities with respect to the slow variable based on the different time scales. In this way, they obtained an effective potential of degree 4 for the one dimensional Fokker-Planck
dynamics on the slow manifold. 

In \cite{CCM2,CCDM} a different strategy was used instead by not performing any Taylor expansion and rather approximate the slow manifold in terms of the nonlinearities defining the evolution using the slow-fast character of the dynamical system, see \cite{BG}. Once, the slow manifold is approximated then one should project the Fokker-Planck dynamics on it to obtain a one dimensional Fokker-Planck reduction possibly valid beyond the local character of the expansion in \cite{RL}. In this way, an approximated full reduced potential was obtained together with an explicit formula for the stationary state distribution on the slow manifold. This effective reduced potential covers the other two stable equilibrium points, and not only the central equilibrium point changing its stability character as in \cite{RL}. 

We illustrate this strategy in the particular example studied in \cite{CCDM} of a neuronal population model with two pools with self-excitation and cross-inhibition as in \cite{DM} related to \eqref{ODE}. The firing rates $\nu_1$ and $\nu_2$ of the neuronal networks are determined by the stochastic dynamical system
\begin{equation}\label{ODE2}
\begin{cases}
{\diff \nu_1} = \left[-\nu_1 +\phi(\lambda_1+ w_+ \nu_1 + w_I
\nu_2)\right]{\diff t} +
\beta \,\diff \mathcal{B}_t^1\\[3mm]
{\diff \nu_2} = \left[- \nu_2 + \phi(\lambda_2+ w_I \nu_1 + w_+
\nu_2)\right]{\diff t} + \beta \,\diff \mathcal{B}_t^2
\end{cases},
\end{equation}
with $\nu_1, \nu_2\geq 0$. Here the applied stimuli are, $\lambda_1 =33$,
$\lambda_2=\lambda_1-\Delta \lambda$, with the bias $\Delta \lambda\in
[0, 10^{-3}]$,  the inhibitory connectivity coefficient is $w_I=1.9$, the standard deviation of the Brownian motion is $\beta=3\cdot 10^{-3}$, and the excitatory connectivity coefficient $w_+$ is the bifurcation parameter, its range of values will be discussed later on. Moreover, the sigmoid (response)
function $\phi(z)$ is given by
\[
\phi(z)=\dfrac{\nu_c}{1+\exp(-b z+\alpha)}\,,
\]
with $\nu_c= 15$, $b= 0.25$, and $\alpha=11.1$.

\begin{figure}[ht!]
\begin{center}
\includegraphics[width=2.9in]{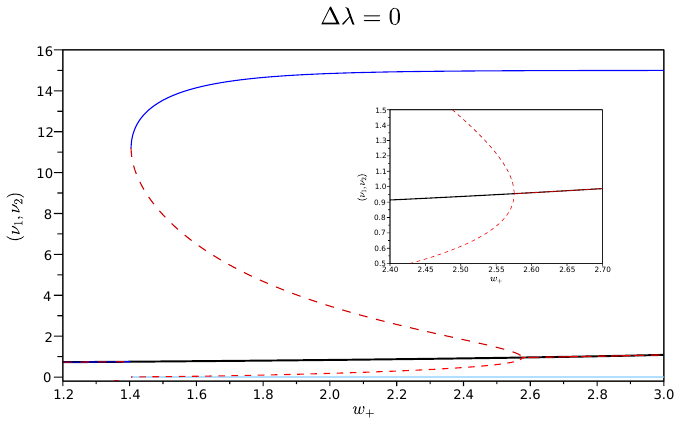}
\includegraphics[width=2.9in]{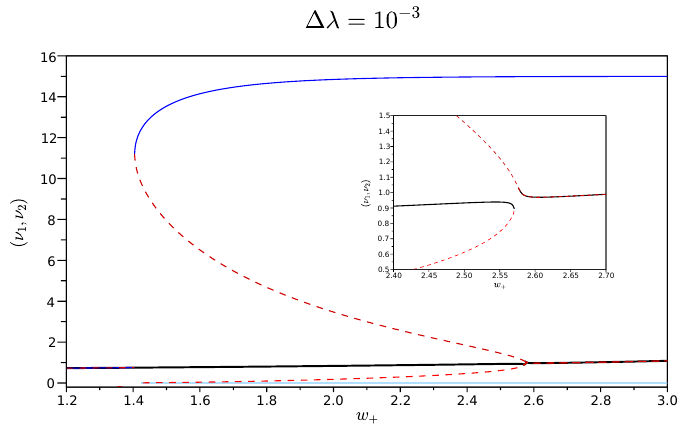}
\end{center}
\caption{Bifurcation diagrams with respect to $w_+$ for $\Delta \lambda=0$ and $10^{-3}$. The spontaneous state is the black line. There are other two stable equilibria: the blue and light-blue lines. The zoom shows the behavior at the bifurcation point $w_+ = 2.5695$. Carrillo, Cordier, Deco, Mancini \cite{CCDM}.} \label{bifurcation}
\end{figure}

Let us first consider the underlying dynamical system \eqref{ODE2} without noise. In terms of $w_+$, it presents a subcritical Hopf bifurcation whose bifurcation diagram is shown in Figure \ref{bifurcation}. For values of the excitatory coefficient around $w_+=1.4$ the system goes from a single asymptotically stable equilibria to a situation in which there are three stable (continuous lines) and two unstable equilibria (hashed lines). The second local bifurcation, where the central asymptotic equilibria disappears, happens at around $w_+=2.5695$ for $\Delta \lambda=10^{-3}$.

\begin{figure}[ht!]
\centering{\includegraphics[width=2in]{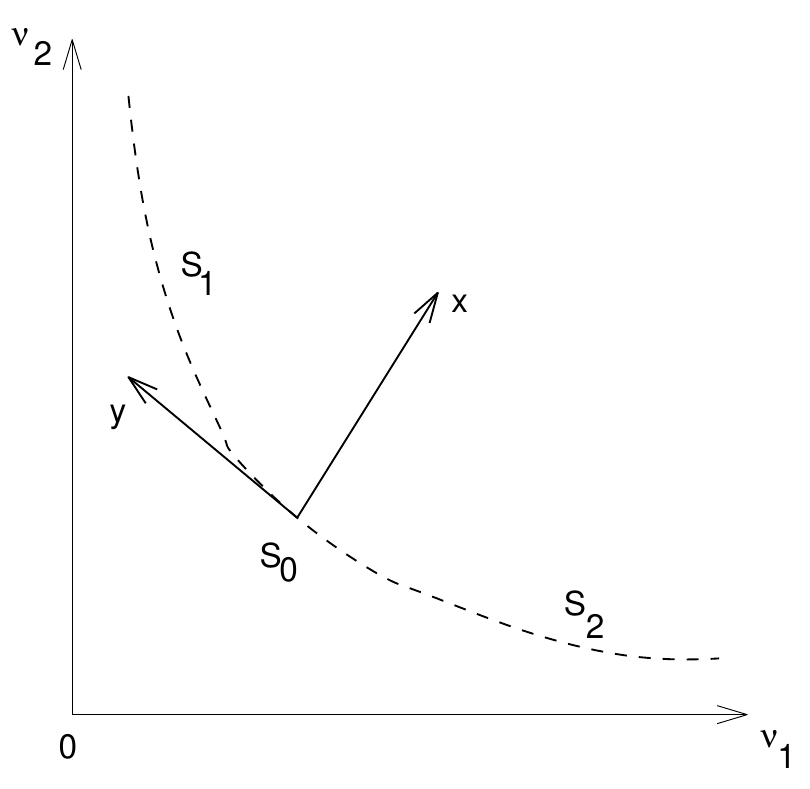}}
\caption{Change of variable from the phase space $(\nu_1, \nu_2)$
to $(x,y)$. Carrillo, Cordier, Deco, Mancini \cite{CCDM}.} \label{chang-variable}
\end{figure}

An approximation of the slow manifold joining the spontaneous $S_0$ and the decision states $S_1$ and $S_2$ is needed (Figure \ref{chang-variable}). This curve is found by introducing the linearization of the dynamical system around the spontaneous equilibria $S_0$ as a new set of variables defined by $X=P^{-1}(\nu-S_0)$, with $X=(x,y)^T$, or equivalently  $\nu = S_0 + P X$ with $\nu=(\nu_1,\nu_2)^T$. Here, $P$ is the matrix diagonalizing the jacobian of the dynamical system \eqref{ODE} at the spontaneous state $S_0$, see \cite{CCM2,CCDM} for more details. 

This change of variables, sketched in Figure \ref{chang-variable}, is a natural way to introduce slow $y$ and fast $x$ variables in the system determined by the eigenvalues of the linearization at the spontaneous state $S_0$. By rewriting the dynamical system \eqref{ODE2}
above in these new variables together with the new time variable $\tau=\varepsilon t$, we have an evolution governed by
\begin{equation*}
\begin{cases}
\varepsilon\displaystyle\frac{\diff x}{\diff \tau} = f(x,y) \\[4mm]
\displaystyle\frac{\diff y}{\diff t} = g(x,y)
\end{cases}\, ,
\end{equation*}
where $\epsilon$ is a small parameter. Assuming that the scaling ratio is zero, $\epsilon=0$, we find the implicit relation $f(x,y)=0$ that determines the approximated slow manifold $x^*(y)$. Inverting this identity to find the approximated curves restrict the set of possible values of the bifurcation parameter in order to keep positive firing rates, in the particular case chosen above, one needs $w_+\geq 1.9$ approximately.

Given the approximation of the slow manifold $x^*(y)$, we can
restrict the dynamics in \eqref{ODE2} to a single
effective Langevin equation. This equation is determined by a
potential $G(y)$ obtained from the evaluation of the dynamics over
the slow manifold approximation, leading to
$$
{\diff y} = g\left(x^*(y),y\right) {\diff t} + \beta_y \, \diff W_t^2\,,
$$
with $\beta_y$ properly obtained in terms of $\beta$ via the
change of variables, see \cite{CCM2} for details. The effective potential is determined by the relation $\partial_y G(y)=-g(x^*(y),y)$. 

\begin{figure}[ht!]
\centering{\includegraphics[width=2.9in]{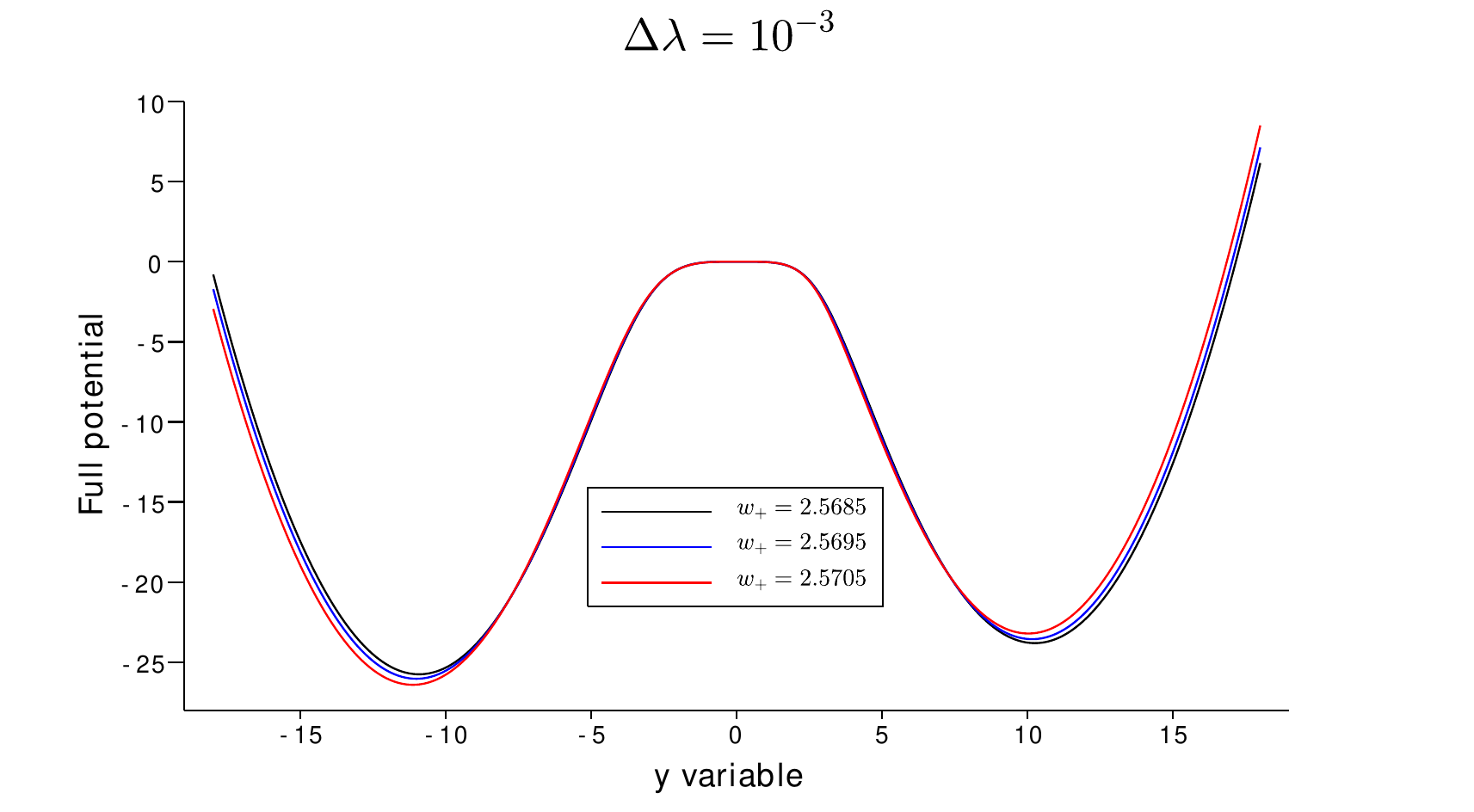}\;
\includegraphics[width=2.9in]{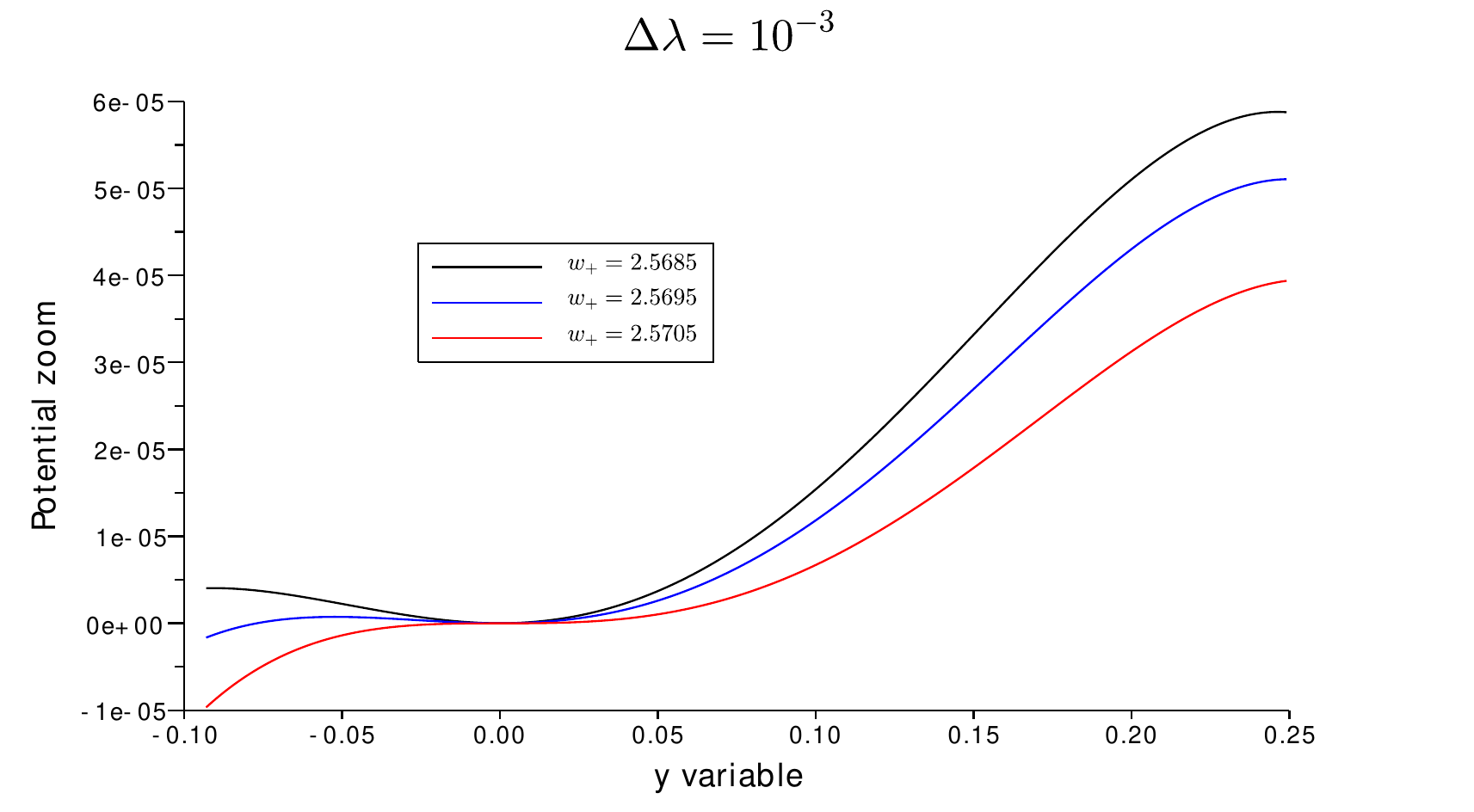}}
\caption{Potential associated to the 1D reduction. Left: $\Delta \lambda=10^{-3}$, right: zoom for $\Delta \lambda=10^{-3}$ around the spontaneous
state. Carrillo, Cordier, Deco, Mancini \cite{CCDM}.}
\label{pot1D}
\end{figure}

Returning to the example \eqref{ODE2}, the previous described reduction strategy yields an approximated potential as in Figure \ref{pot1D} for the case of the subcritical Hopf bifurcation in Figure \ref{bifurcation}. This is in contrast with the approach in \cite{RL} leading to a potential of degree 4 approximating an interval around the central equilibrium point, see \cite[Figure 2]{RL}. Our strategy is valid as long as the approximated slow manifold lies fully in the
relevant biological range of positive rates leading to a restricted range of the bifurcation parameter $w_+$.

\begin{figure}[ht!]
\centering{\includegraphics[width=2.8in]{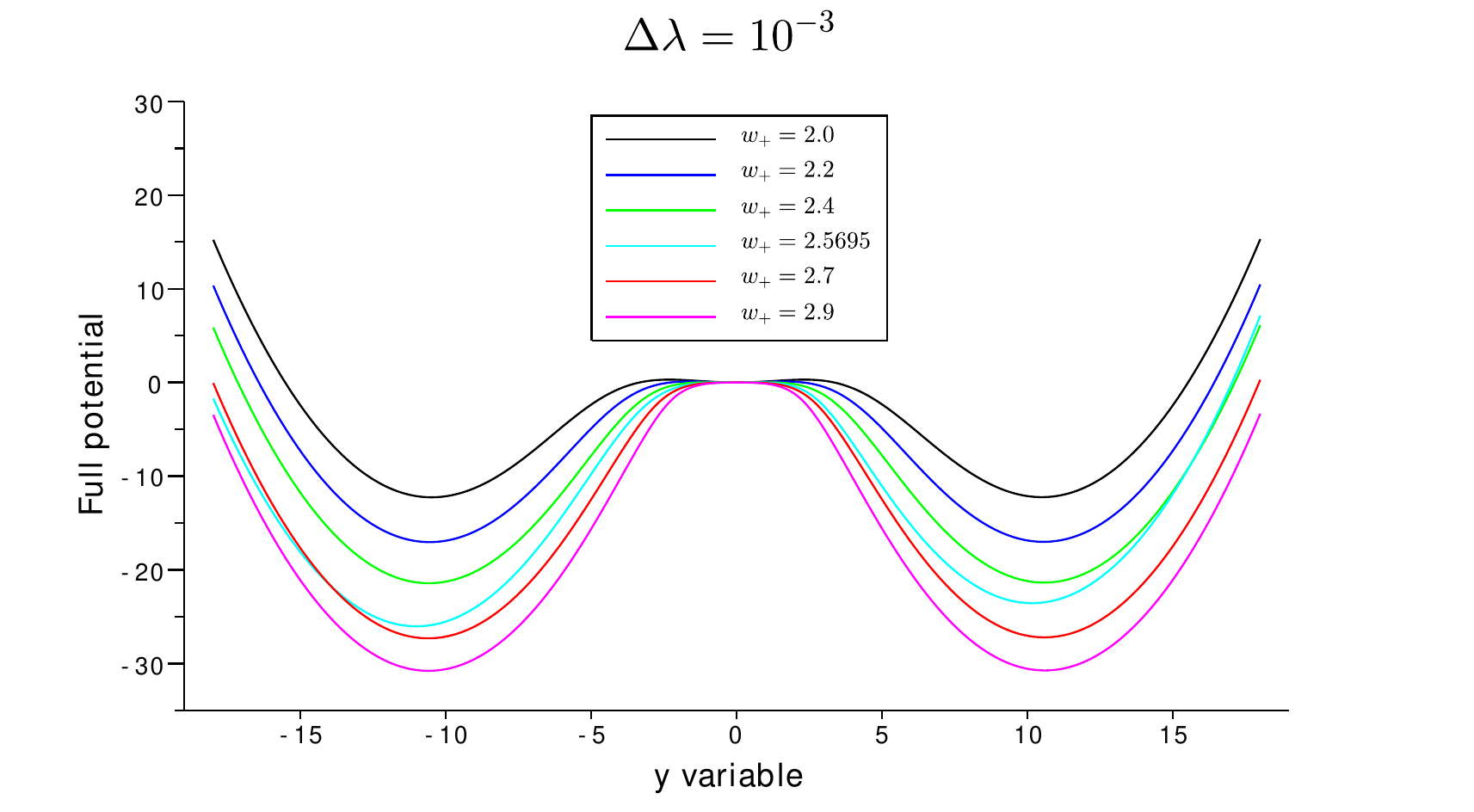}\;
\includegraphics[width=2.8in]{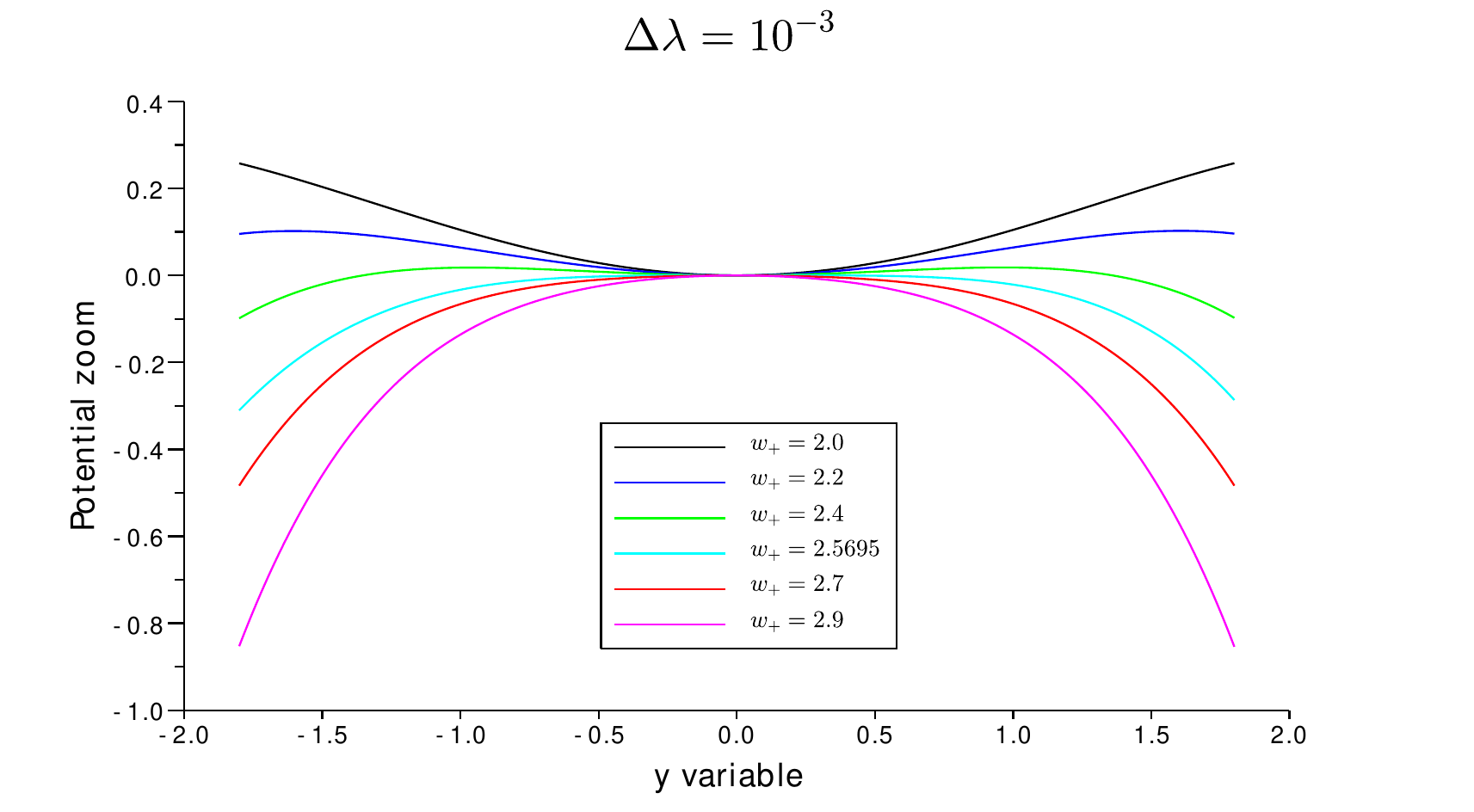}}
\caption{Full 1D potentail along the $y$ variable, for $w_+=2.0,
\ldots, 2.9$ and $\Delta \lambda=10^{-3}$ (left) and a zoom around
the spontaneous state (right). Carrillo, Cordier, Deco, Mancini \cite{CCDM}.}
\label{pot1D-all}
\end{figure}

The 1D effective computed potential $G(y)$ for the biased case ($\Delta \lambda=10^{-3}$)
with respect to the slow variable $y$ is plotted in Figure
\ref{pot1D-all} for $w_+=2.0, \ldots, 2.9$ (left), together with a
zoom at the spontaneous state (right). We note that increasing
$w_+$ beyond the second bifurcation point, the spontaneous state pass from local minimum to local maximum. 

These reduced Fokker-Planck equations carry information about reaction times and performance of the decision tasks since one can compute probabilities for these events in terms of the stationary distributions. We refer for more details to \cite{CCM2,CCDM}.

\part{Hyperbolic equations of transport type}

\section{The Time Elapsed model for neural networks}
\label{sec:TE}

Although the introduction of the model \eqref{NNLIF} was motivated by observations of self-sustained periodic activity in inhibitory networks with sparse firing, the Time Elapsed model \eqref{eq:T-E} described below, on the contrary, was elaborated from observations of periodic activity in excitatory networks with intermittent intense spiking.

Multiple neuron models have exhibited self-sustained oscillations \cite{VVCL,PPCV,PPV,KPV}, leading to think that the core mechanism relies upon a few key aspects of neural dynamics. The Time Elapsed model focuses on one key point: the probability for a neuron to spike in a time window $[t,t+dt)$ can be modelled according to the time $s$ elapsed since the last action potential, the intrinsic refractory period of the neuron and the synaptic integration of the recent history of network activity.

We will first explain the constitutive ideas of Time Elapsed neuron models and propose a set of mathematical assumptions used throughout the section. Then, we will discuss the first attempts, through direct calculations and characteristics, at constructing solutions and obtaining convergence towards a stationary state. Many of these results are obtained in simplified cases where formulae, relative entropies and decay rates can be obtained \textit{via} direct computations. Then, we present more abstract methods adapted from semigroup theory or the Doeblin-Harris machinery, which allow to prove results about existence of solutions and long-time asymptotics in more general frameworks. We then focus on the theoretical and numerical insights into the emergence and stability of periodic solutions. Finally, we discuss the rigorous derivation of this family of equations and present many variants with additional features.

\subsection{Neurons structured by time elapsed since the last discharge}

Summarising \cite{PPCV,PPV}, the article \cite{PPS} introduces and studies rigorously a hyperbolic partial differential equation of transport type modelling an infinite neural network structured by the time $s$ elapsed since the last time a neuron has emitted an action potential. If $n^0(s)$ is the initial probability density for finding in the network at time $t=0$ a neuron which did not fire since a time $s$, then this density is transported at time $1$ and neurons are continuously removed upon firing:
\begin{equation*}
    \dfrac{\partial n}{\partial t}(s,t) + \dfrac{\partial n}{\partial s}(s,t) + p(s,t)n(s,t) = 0,
\end{equation*}
where $p(s,t)$ is at time $t$ the firing rate of neurons which did not fire since a time $s$. The firing function $p(s,t) = p(s,X(t))$ can depend both on the time $s$ elapsed since the last discharge and the past history of network activity $X(t)$. This model is a nonlinear version of renewal equations \cite{GW,P,G}. The key idea for the Time Elapsed model is that when $s$ is small for a neuron, \textit{i.e.} when this neuron just fired, it is not able to fire again due to the refractory period. This refractory period has an absolute component (physical impossibility to fire) and a relative component dependent on the network activity. The Time Elapsed model aims at understanding how self-sustained oscillations can arise from the interplay between refractoriness of neurons and network activity.

The instantaneous firing rate of the network, \textit{i.e.} the flux of neurons firing at time $t$ in the network, is given by
\[ N(t) = \int_{0}^{+\infty} p(s,t)n(s,t)\diff s. \]
Then those neurons that just fired must be reset at time elapsed $s=0$ through the boundary condition $n(0,t)=N(t)$. Network activity is defined through a memory kernel represented by the measure $\mu_X$,
\[ X(t) = \int_{0}^t N(t-\tau)\diff \mu_X(\tau). \]
The transmission can be delayed (smoothly if $\mu_X(\tau)= J \alpha(\tau)\diff \tau$ or singularly with $\mu_X = J \delta_{\tau=d}$, $d > 0$, which is $X(t)=JN(t-d)$), or instantaneous ($\mu_X = J\delta_{\tau=0}$, which is $X(t) = JN(t)$). In both cases, the strength of interconnections is represented by a parameter $J\geqslant 0$.

\subsection{Mathematical setting and the choice of parameters}

Let us summarise the Time Elapsed model into one PDE system.
\begin{equation}\label{eq:T-E}\tag{TE}
	\left\{  \begin{array}{ll}
		\displaystyle \dfrac{\partial n}{\partial t}(s,t)   + \dfrac{\partial n}{\partial s}(s,t)   +  p\big(s,X(t)\big)\, n(s,t) = 0,\qquad& s,t>0,\\[0.2cm]
		\displaystyle N(t) = \int_{0}^{+\infty} p(s,X(t))\, n(s,t) \diff s, \quad X(t) = \int_{0}^{+\infty} N(t-\tau) \mu_X(\diff\tau),\ \ & t\geqslant 0,\\
		\displaystyle n(0,t) = N(t), \qquad n(s,0) = n^0(s)\geqslant 0,\qquad \int_0^{+\infty} n^0(s) \diff s = 1, & s,t\geqslant 0.
	\end{array}     \right.
\end{equation}
Assuming enough regularity, direct integration yields the conservation of mass
\begin{equation}\label{eq:T-E_mass_c}
    \forall t\geqslant 0,\qquad \int_0^{+\infty} n(s,t) \diff s = \int_0^{+\infty} n^0(s) \diff s = 1.
\end{equation}
For the time elapsed neuron model, this feature is at the core of many methods, as the use of characteristics lines in the free transport part provides handy representation formulae for $N(t)$ and $n(s,t)$.

The discharge function $(s,X)\mapsto p(s,X)$ represents the firing rate for neurons with age $s$, given that the feedback from activity in the network is $X$. In the excitatory case, it should increase with activity and with age (note that the inhibitory case, although considered later, has also been defined and studied rigorously \cite{CPST}). We expect that, because of the refractory period, neurons will have zero probability to fire between $s=0$ and some $X$-dependent age $\sigma(x)$. This firing rate should also be bounded. These properties can be summarised as
\begin{equation*}
    \dfrac{\partial p}{\partial s}\geqslant 0, \quad \dfrac{\partial p}{\partial X}\geqslant 0,\quad p\equiv 0\ \mathrm{on} \ [0,\sigma(X)], \quad p> 0 \ \mathrm{on} \ (\sigma(X),+\infty),\quad \norme{p}_\infty < +\infty.
\end{equation*}

In order to perform explicit computations and take full advantage of the simple characteristics in the transport part of \eqref{eq:T-E}, many studies make the following choice for $p(s,X)$.
\begin{hyp}[Simplified choice for $p(s,X)$]\label{as:TE_p}
    The discharge function has the shape \[ p(s,X) = \mathds 1_{s > \sigma(X)},\]
    where $\sigma:\R_+\to\R_+$ is a Lipschitz function such that $\sigma'(x)\leqslant 0$ almost everywhere and
    \[ \forall X\in\R_+,\qquad 0 < \sigma^-\leqslant \sigma(X)\leqslant \sigma^+= \sigma(0) < 1. \]
\end{hyp}

In contrast, it is possible to choose a smooth function for the discharge function $p$.
\begin{hyp}[Smooth firing function $p(s,X)$]\label{as:TE_p_smooth}
    The discharge function satisfies $p\in W^{2,\infty}(\R^2)$ and \[ \dfrac{\partial p}{\partial s}\geqslant 0,\qquad \dfrac{\partial p}{\partial X}\geqslant 0, \qquad 0 < p^- = \lim_{s\to +\infty} p(s,0) \leqslant \lim_{s,X\to +\infty} p(s,X) = p^+ < +\infty.\]
\end{hyp}
Note that Assumption \ref{as:TE_p_smooth} excludes Assumption \ref{as:TE_p} because of the higher regularity it requires. Although existence and long time behaviour can be studied under more general assumptions, many studies choose well-behaved initial conditions in order to get global-in-time bounds for the solution and simplify the computations.
\begin{hyp}[Agreeable initial condition]\label{as:TE_n_0}
    The initial condition satisfies $n^0\in L^\infty(\R_+)\cap L^1(\R_+)$ and
    \begin{equation}
       0\leqslant n^0(\cdot)\leqslant 1,\qquad  \int_0^{+\infty} n^0(s)\diff s = 1.
    \end{equation}
\end{hyp}

The impact of firing neurons on the dynamics is represented through the quantity $X(t)$, itself defined through the integration of the firing rate $N(t)$ under a measure $\mu_X$. The choice of this measure plays a critical role in the modelling goals and the behaviour of the solutions. The simplest choice is the non-delayed Time Elapsed neuron model where $\mu_X(\diff\tau)=J\delta_{\tau = 0}$ for some connectivity strength $J\geqslant 0$. This translates to
\[ X(t) = J N(t).  \]
In the same fashion, one can choose a constant delay for the feedback between firing neurons and the dynamics, with $\mu_X(\diff\tau)=J\delta_{\tau = d}$ for some $d > 0$, yielding
\[X(t) = JN(t-d).\]
In larger networks where the transmission delay and synaptic integration are not evenly distributed, a random delay with distribution $\alpha$ and average connectivity strength $J$ can be considered.
\begin{hyp}[Smoothly delayed feedback]\label{as:TE_delay}
    The measure $\mu_X$ can be written as $\mu_X = J\alpha(\tau) \diff \tau$ with $J\geqslant 0$, $\alpha\in L^q(\R_+)\cap L^1(\R_+)$ for some $q\in(1,+\infty)$ and 
    \[\int_{0}^{+\infty}\alpha(\tau)\diff \tau = 1.\]
\end{hyp}
For some results we will require stonger decay on the delay feedback.
\begin{hyp}[Exponential delay]\label{as:TE_delay_exp}
    The measure $\mu_X$ can be written as $\mu_X = J\alpha(\tau) \diff \tau$ with $J\geqslant 0$, $\alpha\in L^1(\R_+)$ and there exists $\lambda > 0$ such that
    \[\int_{0}^{+\infty}\alpha(\tau)\diff \tau = 1, \quad \int_0^{+\infty} e^{\lambda \tau} (\alpha(\tau) + |\alpha'(\tau)|)\diff \tau < +\infty .\]
\end{hyp}

Solutions are to be understood in the sense of characteristics or in the sense of distributions, as is typical for transport equations. If the distributional solution has a density, for example if it belongs in $L^1$ for all times, then it will be a weak-solution in the following sense:
\begin{multline}\label{eq:T-E_weak}
    \int_0^{+\infty} \phi(s,T) n(s,T) \diff s  = \int_0^{+\infty} \phi(s,0) n^0(s)\diff s + \int_0^T\int_0^{+\infty}\left(\dfrac{\partial\phi}{\partial t}(s,t)+ \dfrac{\partial \phi}{\partial s}(s,t)\right) n(s,t) \diff s\diff t \\\ - \int_0^T\int_0^{+\infty} \phi(s,t)p(s,X(t))n(s,t) \diff s\diff t + \int_0^T N(t)\phi(0,t)\diff t,\\ N(t) = \int_0^{+\infty} p(s,X(t)) n(s,t) \diff s, \qquad X(t) = \int_0^T N(t-\tau)  \mu_X(\diff\tau),
\end{multline}
for any test function $\phi\in C^1(\R_+\times\R_+)$ with  appropriate decay at $+\infty$.

Last, if we want to talk about solutions to \eqref{eq:T-E} in the sense of measures, we have to introduce the following notions. 

Let $(E,\mathcal A)$ be a measurable space. Denote $\mathcal M(E)$ the set of finite measures on $E$, $\mathcal M_+(E)$ the set of positive finite measures on $E$. Recall that the total variation norm $\norme{\cdot}_{\mathrm{TV}}$ on $\mathcal M(E)$ is defined by
\[  \norme{\mu}_{\mathrm{TV}} = \int_E \mu_+ + \int_E \mu_-, \]
where $\mu_+,\mu_-$ is the Hahn-Jordan decomposition of the measure $\mu$ into its positive and negative parts.
Then, $(\mathcal M(\R_+),\norme{\cdot}_{\mathrm{TV}}) $ is a Banach space and the topological dual of $C_{\to0}(\R_+)$. The weak-$*$ topology on $\mathcal M(\R_+)$ can be metrised on any tight subset of $\mathcal M(\R_+)$ with bounded total variation by the bounded Lipschitz norm
\[ \norme{\mu}_{\mathrm{BL}} = \sup_{\psi\in\mathcal W^{1,+\infty}} \int_0^{+\infty}\psi(x) \mu(\diff x). \]
Then we denote $C([0,T),\mathcal M(\R_+))$ the set of functions on $[0,T)$ with values in  $\mathcal M(\R_+)$ which are continuous for the bounded Lipschitz norm. We denote the translation semigroup generated on $(\mathcal M(\R_+),\norme{\cdot}_{\mathrm{BL}})$ by the operator $-\partial_s$, which means by abuse of notation $T_t n^0 = n(\cdot-t)$.

\begin{definition}\label{def:T-E_meas_sol}
    Consider the case $\mu_X = \delta_{\tau=0}$, \textit{i.e.} $X(t) = N(t)$. A couple $(n,N)$ is said to be a measure solution to \eqref{eq:T-E} on $[0,T)$ with initial condition $n^0\in\mathcal M_+(\R_+)$, $N^0\in\R_+$ if $n\in C([0,T),\mathcal M_+(\R_+))$, $N\in C([0,T),\R_+)$, $n(0)=n^0$, $N(0)=N^0$, and for all $t\in[0,T)$,
    \begin{equation*}\label{eq:T-E_meas_sol}
        n(t) = T_t n^0 + \int_0^t T_{t-\tau}\, \left( p(\cdot,N(\tau)) n(\cdot,\tau)\right)(s)\diff \tau + \int_0^t T_{t-\tau}\, \left( N(\tau)\delta_{s=0}\right)\diff \tau.
    \end{equation*}
    and        $N(t) = \int_0^{+\infty} p(s,N(\tau)) n(t)(\diff s)$.
\end{definition}
In the definition, we write $n(\cdot,t)$ by abuse of notation even if $n(t)\in\mathcal M_+(\R_+)$ does not always have a density. The time integrals have to be understood as Bochner integrals. It can be readily checked that the crucial mass conservation \eqref{eq:T-E_mass_c} holds for measure solutions too, which implies by positivity $\norme{n(t)}_{\mathrm{TV}}=\norme{n^0}_{\mathrm{TV}}$. Note also that
\[  \int_0^t T_{t-\tau}\, \left( N(\tau)\delta_{s=0}\right)(s)\diff \tau = \int_0^t  N(\tau)\delta_{s=t-\tau}\diff \tau = N(t-s) \mathds 1_{[0,+\infty)}(t-s), \]
and that \eqref{eq:T-E_meas_sol} can be turned into the easier equivalent form(see arguments in \cite{CY}): for all $\phi\in C^\infty_c(\R_+)$,
\begin{equation*}\label{eq:T-E_weak_measures}
    \dfrac{\diff }{\diff t}\int_0^{+\infty} \phi(s) n(t) (\diff s)  =  \int_0^{+\infty} \dfrac{\partial \phi}{\partial s}(s) n(t) (\diff s)  - \int_0^{+\infty} \phi(s)p(s,X(t))n(t) (\diff s) +  N(t)\phi(0).
\end{equation*}

In the context of measure solutions, we introduce another assumption on the firing function $p$.

\begin{hyp}[Firing function for measure solutions]\label{as:TE_p_meas}
    The firing function $p$ belongs in $W^{1,\infty}(\R_+\times\R_+)$, $ p\geqslant 0$ and there exists $L>0$ such that
    \[  \forall s,X_1,X_2\geqslant 0,\qquad |p(s,X_1)-p(s,X_2)|\leqslant L|X_1-X_2|.  \]
\end{hyp}
\begin{hyp}[Preparation for decay in measures]\label{as:TE_p_meas_2}
    The firing function $p\in W^{1,+\infty}(\R_+\times\R_+)$ satisfies, for some constants $\sigma, p^-,p^+>0$,
    \[ \forall s,X\geqslant 0,\qquad  \dfrac{\partial p}{\partial s}(s,X)\geqslant 0,\qquad p^- \mathds 1_{[\sigma,+\infty)} \leqslant p(s,X) \leqslant p^+.  \]
\end{hyp}

Before proceeding with the description of the numerous analytical results on \eqref{eq:T-E}, let us emphasise that two equivalent ways of making a weak-nonlinearity assumption are used in the literature: one can choose a Lipschitz condition on $p$ like $\norme{\sigma'}_\infty$ small enough in Assumption \ref{as:TE_p} or $L$ small enough in Assumption \ref{as:TE_p_meas}, and then $X$ is not scaled as in $X(t) = N(t)$ or $X(t) = \int_0^t N(t-\tau)\alpha(\tau)\diff \tau$, $\int_0^{+\infty} \alpha(\tau)\diff \tau = 1$; or one could choose $J$ small enough with $X(t) = J N(t)$ or $J$ in Assumption \ref{as:TE_delay} or Assumption \ref{as:TE_delay_exp}. This switch between a Lipschitz condition on the discharge function $p$ and a smallness condition on the connectivity parameter $J$ allows to simplify computations depending on the method. We try, for each reference described, to stay compatible with the global framework proposed above for the Time Elapsed model, which may require some mental gymnastic on the part of the reader.

\subsection{Global existence and well-posedness}

\subsubsection{Uniformly bounded solutions}
\label{sec:TE_exist}

In \cite{PPS}, global-in-time solutions are constructed \textit{via} a fixed point argument. Let us choose the smoothly delayed assumption \ref{as:TE_delay} and fix some time $T > 0$. Given some external input $Y(t)$ in place of the nonlinear $X(t)$, there is a unique solution $n$ with firing rate $N$ to the linear problem

\begin{equation}\label{eq:T-E_linear}
	\left\{  \begin{array}{ll}
		\displaystyle \dfrac{\partial n}{\partial t}(s,t)   + \dfrac{\partial n}{\partial s}(s,t)   +  p\big(s,Y(t)\big)\, n(s,t) = 0,\qquad& s,t>0,\\[0.2cm]
		\displaystyle N(t) = \int_{0}^{+\infty} p(s,Y(t))\, n(s,t) \diff s, \ & t\geqslant 0,\\
		\displaystyle n(0,t) = N(t), \qquad n(s,0) = n^0(s)\geqslant 0,\qquad \int_0^{+\infty} n^0(s) \diff s = 1, & s,t\geqslant 0.
	\end{array}     \right.
\end{equation}
Given such a solution, we can define the functional
\begin{equation*}
	\left\{  \begin{array}{lccl}
		F:&L^\infty([0,T])&\to& L^\infty([0,T])\\
        &Y&\mapsto&\displaystyle J\int_0^{+\infty} \alpha(\tau) N(t-\tau) \diff \tau.
	\end{array}     \right.
\end{equation*}
Consider two external inputs $Y_1,Y_2\in L^\infty([0,T])$, for $T$ small enough but independent from $n^0$. Because the differential operator is speed one transport, it is possible to exploit the characteristics:
for all $s\geqslant t$, the initial condition is propagated along characteristics while being dampened by the loss of firing neurons,
\begin{equation}\label{eq:T-E_char1}
    n(s,t) = n^0(s-t) e^{-\int_0^t p(\tau+s-t,Y(\tau))\diff \tau}, \qquad s\geqslant t.
\end{equation}
On the contrary, when $s\leqslant t$, the neuron have already fired and reset and come from the boundary condition at $s=0$, which can be computed as
\begin{equation}\label{eq:T-E_char2}
    n(s,t) = N(t-s) e^{-\int_0^s p(\tau,Y(\tau+t-s))\diff \tau}, \qquad s\leqslant t.
\end{equation}
Combining \eqref{eq:T-E_char1} and \eqref{eq:T-E_char2} gives the representation formula
\begin{multline}
    N(t) = \int_0^t p(s,Y(t))N(t-s)e^{-\int_0^s p(\tau,Y(\tau+t-s))\diff \tau}\diff s\\ + \int_t^{+\infty} p(s,Y(t))n^0(s-t) e^{-\int_0^t p(\tau+s-t,Y(\tau))\diff \tau}\diff s.
\end{multline}
Applying this formula for both $Y_1$ and $Y_2$, it is then possible to decompose the difference $F(Y_1)-F(Y_2)$ into many terms, each of those one can control by the norm $\norme{Y_1-Y_2}_\infty$. It is possible to choose $T$ small enough independently of $n^0$ such that
\begin{equation*}
    \norme{F(Y_1)-F(Y_2)}_\infty \leqslant \dfrac12 \norme{Y_1-Y_2}_\infty.
\end{equation*}
Hence, by applying the Banach fixed point theorem, there is a unique fixed-point for $F$, which yields a unique local-in-time solution of \eqref{eq:T-E}. Since $T$ is independent from $n^0$, the procedure can be iterated to extend the local solution into a global-in-time solutions. Uniform bounds can be derived along the way.

\begin{theorem}[Pakdaman, Perthame, Salort \cite{PPS}]\label{thm:T-E_exist}
    Grant Assumptions \ref{as:TE_p} and \ref{as:TE_n_0} and \ref{as:TE_delay} (simple $p$, good $n^0$, smooth delay). There exists a unique solution $n\in C(\R_+,L^1(\R_+))$ to \eqref{eq:T-E} and for all $s,t\geqslant 0$,
    \begin{equation*}
        0\leqslant n(s,t)\leqslant 1,\quad 0\leqslant N(t)\leqslant 1,\qquad J(1-\sigma^+)\leqslant X(t)\leqslant J.
    \end{equation*}
\end{theorem}

\begin{remark}
    The authors of \cite{PPS} claim that, if $p$ is taken to be Lipschitz continuous, then the same proof works when the smooth delay Assumption \ref{as:TE_delay} is replaced by a fixed delay assumption $\mu = J \delta_{\tau = d }$, $d\geqslant 0$, $J\geqslant 0$. The case $d=0$ yields the non-delayed Time Elapsed model where $X(t) = J N(t)$.
\end{remark}

\subsubsection{Solutions in the sense of measures}
\label{sec:TE_meas_sol}

In \cite{CY}, measures solution in the sense of definition \ref{def:T-E_meas_sol} are constructed for well-behaved firing functions $p$ and initial conditions with small enough total variation norm, in the special case $X(t) = N(t)$.

Following \cite{CCC,PPS}, a fixed point argument is applied in an appropriate space. Fix $C>0$ and consider
\[  \mathcal X = \{ n\in C([0,T], \mathcal M_+(\R_+))\ | \ n(0) = n^0, \ \forall t\in[0,T],\ \norme{n(t)}_{\mathrm{TV}}\leqslant C \},\]
endowed with the norm $\norme{n}_{\mathcal X} = \sup_{t\in[0,T]} \norme{n(t)}_{\mathrm{TV}}$. Beware that the continuity in $C([0,T], \mathcal M_+(\R_+))$ is not with respect to the total variation norm $\norme{\cdot}_{\mathrm{TV}}$ but according to the bounded Lipschitz norm $\norme{\cdot}_{\mathrm{BL}}$. The fixed point is sought for the operator $\Psi:\mathcal X\to\mathcal X$ defined as
\[  \Psi[n](t) =   T_t n^0 + \int_0^t T_{t-\tau}\, \left( p(\cdot,N(\tau)) n(\cdot,\tau)\right)\diff \tau + \int_0^t T_{t-\tau}\, \left( N(\tau)\delta_{s=0}\right)(s)\diff \tau,\]
where
\[  N(t) = \int_0^{+\infty} p(s,N(t))n(t)(\diff s). \]
An important step is to check that $N$ is uniquely defined, allowing the operator $\Psi$ to be well-defined as well.
\begin{lemma}[Ca\~nizo, Yoldaş \cite{CY}]
    Let $n\in \mathcal M_+(\R_+)$. Grant Assumption \ref{as:TE_p_meas} with $L$ satisfying
    \begin{equation} \label{eq:T-E_L_cond_1} L\leqslant \frac{1}{4\norme{n}_{\mathrm{TV}}}.  \end{equation}
    Then there exists a unique solution $N\geqslant 0$ to the fixed-point problem
    $ N = \int_0^{+\infty} p(s,N) n(\diff s). $
    Given $n_1,n_2\in \mathcal M_+(\R_+)$ satisfying \eqref{eq:T-E_L_cond_1}, the associated $N_1,N_2$ satisfy
    \[ |N_1-N_2| \leqslant \dfrac{\norme{p}_\infty}{1-L\norme{n_1}_{\mathrm{TV}}} \norme{n_1-n_2}_{\mathrm{TV}}. \]
\end{lemma}

Using this lemma, the mapping $\Psi$ is well-defined and can be proved to be a contraction on $\mathcal X$ for $T$ small enough, yielding an existence result in the space of finite measures.

\begin{theorem}[Ca\~nizo, Yoldaş \cite{CY}]
    Let $n^0\in \mathcal P(\R_+)$. Grant Assumption \ref{as:TE_p_meas} with 
    \begin{equation} \label{eq:T-E_L_cond_2}L\leqslant \frac{1}{4\norme{n^0}_{\mathrm{TV}}}. \end{equation}
    Consider the non-delayed case $\mu_X=\delta_{\tau=0}$, \textit{i.e.} $X(t) = N(t)$. There exists a unique global-in-time measure solution $(n,N)$ to \eqref{eq:T-E} in the sense of Definition \ref{def:T-E_meas_sol}. If $n^{0,1},n^{0,2}\in \mathcal P(\R_+)$ satisfying \eqref{eq:T-E_L_cond_2}, then the associated solutions $n_1$ and $n_2$ satisfy
    \begin{equation*}
        \forall t\geqslant 0, \qquad \norme{n_1(t)-n_2(t)}_{\mathrm{TV}} \leqslant e^{4\norme{p}_\infty t} \norme{n^{0,1}-n^{0,2}}_{\mathrm{TV}}.
    \end{equation*}
\end{theorem}
We enforce $\int n^0 = 1$ in this theorem for the sake of consistency with our formulation \eqref{eq:T-E} but this assumption is unnecessary in \cite{CY}. 

The Lipschitz condition on $p$ could be translated onto a connectivity parameter $J$ taking $X(t) = J N(t)$. Note that this smallness condition on $L$ is important, as will be seen in Section \ref{sec:TE_periodic} where solutions with discontinuous firing rate are constructed using the non-uniqueness of $N(t)$ for large $L$ (translated to large $\norme{\sigma'}_\infty$ in the context of Assumption \ref{as:TE_p}). See also the Lipschitz assumption on $\sigma$ in Theorem \ref{thm:T-E_CV1}.

\subsubsection{Improved existence proof including the strongly inhibitory case}

In the article \cite{STV}, the existence results in \cite{PPS} and \cite{CY} are revisited and improved. They proceed by considering a cutoff and solutions to an approximate problem, then going back to the full \eqref{eq:T-E}.

In the non-delayed case, they obtain:
\begin{theorem}[Sepulveda, Torres, Villada \cite{STV}]
    Assume $\mu_X = \delta_{\tau=0}$, \textit{i.e.} $X(t)=N(t)$, $p\in W^{1,+\infty}(\R_+\times\R_+)$ and let $\gamma = \sup_{s,X\geqslant 0} \partial_X p (s,X)$. For all initial condition $n^0\in L^1(\R_+)\cap \mathcal P(\R_+)$ such that $\gamma \norme{n^0} < 1$, there exists a unique solution $n\in C(\R_+,L^1(\R_+))$, $N\in C(\R_+)$ to \eqref{eq:T-E}. 
\end{theorem}
It is claimed in \cite{STV} that the regularity can be relaxed to include step functions such as Assumption \ref{as:TE_p} or the case $p(s,X) = \phi(X)\mathds 1_{s > \sigma}$ used in \cite{CPST}. A key advantage of this result is that, because there is no absolute value or norm in the definition of $\gamma$, the condition $\gamma = \sup_{s,X\geqslant 0} \partial_X p (s,X)$ is always satisfied when $\gamma\leqslant 0$, which included the strongly inhibitory framework that is studied for example in \cite{CPST}. Similarly, for the model \eqref{NNLIF}, solutions are always global-in-time in the inhibitory case as described in Section \ref{sec:NNLIF}.

In the distributed delayed case, they obtain a more general existence result:
\begin{theorem}[Sepulveda, Torres, Villada \cite{STV}]
    Assume $\mu_X = \alpha(\tau)\diff \tau$, $\alpha\in L^1(\R_+)$, and $p\in W^{1,+\infty}(\R_+\times\R_+)$. For all initial condition $n^0\in L^1(\R_+)\cap \mathcal P(\R_+)$, there exists a unique solution $n\in C(\R_+,L^1(\R_+))$, $N,X\in C(\R_+)$ to \eqref{eq:T-E}. 
\end{theorem}
Note that the delayed assumption is slightly more general than Assumption \ref{as:TE_delay} used in Theorem \ref{thm:T-E_exist} as the delayed kernel is not assumed to be in some $L^q$, $q > 1$.

\subsection{Study of the stationary equation}

Under either Assumption \ref{as:TE_delay} or for fixed delay $\mu_X = J \delta_{\tau=d}$, $J,d\geqslant 0$, a stationary firing rate implies that
$  X_\infty = J N_\infty. $
Then stationary states are solutions of
\begin{equation*}
    \dfrac{\diff n_\infty}{\diff s}(s) = - p(s,X_\infty) n_\infty(s),\qquad X_\infty = J \int_0^{+\infty} p(s,X_\infty) n_\infty(s)\diff s,\qquad \int_0^{+\infty} n_\infty(s)\diff s = 1,
\end{equation*}
which has solutions of the form
\begin{equation*}
    n_\infty(s) = \dfrac{e^{- \int_0^s p(\eta,X_\infty)\diff \eta }}{\int_0^{+\infty} e^{- \int_0^s p(\eta,X_\infty)\diff \eta }\diff s}.
\end{equation*}
Under Assumption \ref{as:TE_p} for $p$, it simplifies down to the full characterisation
\begin{equation}\label{eq:T-E_st_st}
    n_\infty(s) = \dfrac{e^{ - (s-\sigma(X_\infty))_+ }}{1+\sigma(X_\infty)},\qquad X_\infty = J ( 1+\sigma(X_\infty) ).
\end{equation}
The number of stationary states thus depends on the number of solutions $X_\infty$ in \eqref{eq:T-E_st_st}; there can be several of them for some choices of $\sigma$. For small enough $J$, a Banach fixed point argument yields uniqueness. For larger $J$, the criterion is more involved.

\begin{theorem}[Pakdaman, Perthame, Salort \cite{PPS}]
    Grant Assumption \ref{as:TE_p} and assume either \ref{as:TE_delay} or $\mu_X = J \delta_{\tau = d}$, $J,d\geqslant 0$. 
    \begin{itemize}
        \item If $J$ is small enough, there is a unique stationary state of \eqref{eq:T-E}.
        \item There exists a unique stationary state of \eqref{eq:T-E} if $J$ is large enough and
    \begin{equation*}
        \limsup_{J\to+\infty} \sup_{X>\frac{J}{1+\sigma^+}}  \dfrac{J|\sigma'(X)|}{(1+\sigma(x))^2} < 1.   
    \end{equation*}
    \end{itemize} 
\end{theorem}

If we choose a smoother function $p$, going for Assumption \ref{as:TE_p_smooth}, we can also obtain existence of stationary states and uniqueness for small connectivities.
\begin{theorem}[Mischler, Weng \cite{MW}]
    Grant Assumption \ref{as:TE_p} and assume either \ref{as:TE_delay_exp} or $\mu_X = J \delta_{\tau = 0}$, $J\geqslant 0$. Then, there exists at least one stationary state $(n_\infty,N_\infty, X_\infty)$ of \eqref{eq:T-E} such that for some constant $C > 0$,
    \[  \forall s\geqslant 0,\qquad 0 \leqslant n_\infty(s) \leqslant C e^{-\frac{p^-}{2}s}, \quad \left|\dfrac{\partial n_\infty(s)}{\partial s}\right| \leqslant C e^{-\frac{p^-}{2}s}. \]
    There exists $J_0>0$ small enough such that for all $J\in[0,J_0)$, there is a unique stationary state to \eqref{eq:T-E}.
\end{theorem}

This result has been extended to the strongly nonlinear case in \cite{MQW}. In the context of measures solutions, stationary states can be found in the sense of Definition \ref{def:T-E_meas_sol}.

\begin{theorem}[Ca\~nizo, Yoldaş \cite{CY}]
    Grant Assumptions \ref{as:TE_p_meas} and \ref{as:TE_p_meas_2} and assume $\mu_X = \delta_{\tau=0}$, \textit{i.e.} $X(t) = N(t)$. If
    \begin{equation*}
        L< \dfrac{1}{\left( \frac{\sigma^2}{2} + \frac{\sigma}{p^-} + \frac{1}{(p^-)^2} \right) (p^+)^2},
    \end{equation*}
    then there exists a unique stationary state of \eqref{eq:T-E} in the space of probability measures $\mathcal P(\R_+)$.
\end{theorem}
Further insights into existence and uniqueness of measure stationary states can be obtained \textit{via} the Doeblin-Harris method, as described later in Subsection \ref{sec:TE_Doeblin}

\subsection{Asymptotic relaxation to a stationary state for the step function case}

We first discuss here the long-time asymptotic results under Assumption \ref{as:TE_p}, when the firing function $p(s,X) = \mathds 1_{s > \sigma(X)}$ is a step function, or when we choose a firing function $p(s,X) = \phi(X)\mathds 1_{s > \sigma}$ which separates into a network influence $\phi$ and a constant refractory period $\mathds 1_{s > \sigma}$.

\subsubsection{Convergence in relative entropy}

\paragraph{Linear case:} let us first showcase the method on the simpler linear case, which corresponds, under Assumption \ref{as:TE_p}, to $\sigma(X)\equiv \sigma\in(0,1)$ constant, and we write for clarity $\bar p(s) = p(s,X) = \mathds 1_{s>\sigma}$. Then, there is a unique stationary state to \eqref{eq:T-E} which can be written as
\begin{equation}
    n_\infty(s) = \dfrac{e^{-\int_0^s \bar p(\tau) \diff \tau}}{\int_0^{+\infty} e^{-\int_0^s \bar p(\tau) \diff \tau}}\diff s.
\end{equation}
Note that this function $n_\infty$ is constant on $[0,\sigma]$ and then is exponentially decaying.

Taking inspiration from the general methods in \cite{MMP1,MMP2}, we want to find some function $\Psi$ and some $\mu>0$ such that
\begin{equation}\label{eq:T-E_entropy}
    \dfrac{\diff }{\diff t} \int_0^{+\infty} \Psi(s) |n(s,t)-n_\infty(s)|\diff s \leqslant -\mu \int_0^{+\infty} \Psi(s) |n(s,t)-n_\infty(s)|\diff s.
\end{equation}
First, we can notice that by linearity, $\tilde n = n-n_\infty$ still satisfies \eqref{eq:T-E}, and then it is also the case for $|\tilde n|$. By conservation of mass, $\int_0^{+\infty} \tilde n(s,t)\diff s = 0$ for all $t\geqslant 0$, which allows the computation, for any $\nu\in(0,1)$,
\[ |\tilde n(0,t)| = |\tilde N(t)| = \left| \int_0^{+\infty} \bar p(s) \tilde n(s)\diff s \right| = \left| \int_0^{+\infty} (\bar p(s)-\nu) \tilde n(s)\diff s \right|. \]
Using  integration by part (with a density argument on $n^0$ with compact support if necessary),
\begin{align*}
    \dfrac{\diff }{\diff t} \int_0^{+\infty} \Psi(s) |n(s,t)-n_\infty(s)|\diff s = & \int_0^{+\infty}\left( \dfrac{\diff \Psi}{\diff s}(s) - \bar p(s) \Psi(s) \right)|n(s,t)-n_\infty(s)|\diff s\\
    & - \Psi(0)\left| \int_0^{+\infty} (\bar p(s)-\nu) \tilde n(s)\diff s \right|
\end{align*}
Hence, it is possible to prove \eqref{eq:T-E_entropy} by choosing $\Psi$ as a solution of
\begin{equation}\label{eq:T-E_Psi}
    \dfrac{\diff \Psi}{\diff s}(s) - (\bar p(s) -\mu) \Psi(s) = - \Psi(0) |\bar p(s) - \nu|.
\end{equation}
It is possible to prove that there is a choice of $\mu,\nu$ and $\Psi(0)=1$ such that there exists a non-negative solution $\Psi$ to this problem.

\begin{theorem}[Pakdaman, Perthame, Salort \cite{PPS}]
    Grant Assumptions \ref{as:TE_p} and \ref{as:TE_n_0} with $\sigma(X)\equiv \sigma\in(0,1)$ constant. There exists $\mu,B>0$ and a non-negative function $\Psi:\R_+\to\R_+$ such that $\Psi(0)=1$, $\Psi(\cdot)\geqslant B$, $\Psi$ grows as $e^{\mu s}$ as $s\to +\infty$ and if $\int_0^{+\infty} \Psi(s)n^0(s)\diff s +\infty$, then for all $t\geqslant 0$,
    \begin{equation}
        \int_0^{+\infty} \Psi(s) |n(s,t)-n_\infty(s)|\diff s  \leqslant e^{-\mu t}  \int_0^{+\infty} \Psi(s) |n^0(s)-n_\infty(s)|\diff s.
    \end{equation}
\end{theorem}

\paragraph{Weakly-nonlinear and strongly nonlinear case:}

For $J\neq 0$ small enough (and $\sigma$ not constant), it is still possible to construct the function $\Psi$ as a solution of \eqref{eq:T-E_Psi}, which can be interpreted as spectral gap property.

Denoting $f(t) = \int_0^{+\infty} \Psi(s) |n(s,t)-n_\infty(s)|\diff s $, it is possible to write a differential inequality of the form
\[  \dfrac{\diff f}{\diff t}(t) \leqslant -\mu f(t) + C \norme{\sigma'}_{L^\infty([J(1-\sigma^+),J])} |X(t) - X_\infty|.  \]
The difficulty is to control the difference $|X-X_\infty|$ which depends on past values of $N(t)$. However, according to Theorem \ref{thm:T-E_exist}, we have the control
\[ J(1-\sigma^+)\leqslant X(t)\leqslant J, \qquad J(1-\sigma^+)\leqslant X_\infty \leqslant J. \]
Hence $ \frac{\diff f}{\diff t}\leqslant -\mu f + C\norme{\sigma'}_{L^\infty([J(1-\sigma^+),J])} J$. Under either a smallness assumption on $J$ or taking $J$ large and assuming that $J\norme{\sigma'}_{L^\infty([J(1-\sigma^+),J])} $ is small enough when $J$ is large enough, we can control the size of $f$ and feed this control into an iteration procedure for further refined estimation of $|X-X_\infty|$.

\begin{theorem}[Pakdaman, Perthame, Salort \cite{PPS}]\label{thm:T-E_CV_PPS1}
    Grand Assumptions \ref{as:TE_p}, \ref{as:TE_n_0} and \ref{as:TE_delay}. If either $J$ is small enough or if $J$ is large enough and the quantity
    \[  \limsup_{J\to+\infty}\sup_{X\in[J(1-\sigma^+),J]} J|\sigma'(X)| \]
    is finite and small enough, there exists $\mu,B>0$ and a non-negative function $\Psi:\R_+\to\R_+$ such that $\Psi(0)=1$, $\Psi(\cdot)\geqslant B$, $\Psi$ grows as $e^{\mu s}$ as $s\to +\infty$. For any initial condition such that $\int_{0}^{+\infty}\Psi(s)|n^0(s)-n_\infty(s)|\diff s < +\infty$, the solution of \eqref{eq:T-E} satisfies
    \begin{equation}
        \lim_{t\to +\infty} \int_{0}^{+\infty}\Psi(s)|n(s,t)-n_\infty(s)|\diff s = 0,
    \end{equation}
    where $n_\infty$ is the unique stationary state.
\end{theorem}

Without more information on $\alpha$ in Assumption \ref{as:TE_delay}, it is not possible to conclude about the rate of convergence to the stationary states. However, in the simpler non-delayed case $\mu = J \delta_{\tau = 0}$, \textit{i.e.} $X(t)=JN(t)$, \cite{PPS} proves under similar hypotheses that if $J$ is small enough, for some positive $\beta$,
\begin{equation}
         \int_{0}^{+\infty}\Psi(s)|n(s,t)-n_\infty(s)|\diff s \leqslant e^{-\beta t} \int_{0}^{+\infty}\Psi(s)|n^0(s)-n_\infty(s)|\diff s.
    \end{equation}

\subsubsection{Convergence of the firing rate}\label{sec:TE_CV_1}

Because of the simple characteristic flow of the transport operator $\partial_t + \partial_s$, it is also possible to work directly at the level of the firing rate $N(t)$. If we know that $N(t)$ converges to a stationary firing rate $N_\infty$, then we can conclude that the density $n$ will also become stationary as a solution to a linear equation with external force $X_\infty$. Following the approach of \cite{PPS2} we choose the simplified non-delayed case $\mu = J \delta_{\tau = 0}$ with $J=1$, wich is equivalent to $X(t) = N(t)$. Note that the general case $J\neq 1$ can be retrieved by a suitable modification of $\sigma$. The connectivity parameter is thus transmitted onto the derivative of $\sigma$, which we will assume bounded:
\[ 0\leqslant -\sigma'(X)\leqslant \norme{\sigma'}_\infty.  \]

\paragraph{Linear case:} Consider again the linear case $\sigma(X)\equiv \sigma\in(0,1)$ constant, under Assumption \ref{as:TE_p}. Note first that then the unique stationary state has the firing rate
\[ N_\infty = \dfrac{1}{1+\sigma}. \]
A fundamental idea is that for $t\geqslant s$, the initial condition is gone through the free speed one transport at point $(s,t)$ and if $s\leqslant \sigma$ there is no firing yet. Hence
\[ n(s,t) = N(t-s),\qquad t\geqslant s, \ s\leqslant \sigma. \]
Then the conservation of mass can be rewritten as \[1 = \int_0^{+\infty} n(s,t)\diff s = \int_0^{\sigma}n(s,t)\diff s + \int_\sigma^{+\infty}\mathds 1_{s\geqslant \sigma}n(s,t)\diff s,\] which yields
\begin{equation}\label{eq:T-E_Nmethod}
    N(t) + \int_{t-\sigma}^t N(s)\diff s = 1 , \qquad t\geqslant \sigma.
\end{equation}
Because the stationary firing rate satisfies this identity and because, by Theorem \ref{thm:T-E_exist}, $0\leqslant N,N_\infty\leqslant 1$, a contraction can be obtained by iteration:
\begin{equation}
    N(t)-N_\infty = \int_{t-\sigma}^t \big(N_\infty- N(s)\big)\diff s \implies \forall t\geqslant k\sigma,\quad |N(t)-N_\infty| \leqslant \sigma^k.
\end{equation}
Looking further at \eqref{eq:T-E_Nmethod}, we can also see that $N(t)$ oscillates around $N_\infty$ while exponentially converging towards it.

\paragraph{Weakly-nonlinear case:} The same procedure can be applied when $\sigma$ is not constant but varies slowly $\norme{\sigma'}_\infty < 1$. The idea is to obtain and use again a formula like \eqref{eq:T-E_Nmethod}. First, note that $N(t)$ is defined by
\begin{equation}\label{eq:T-E_fixed}
		 N(t) = \int_{\sigma(N(t))}^{+\infty} n(s,t)\diff s
\end{equation}
In the case $\norme{\sigma'}_\infty < 1$, there is a unique solution to this problem. As we will see later, it is not the case in other parameter regimes and non-uniqueness of $N$ can induce oscillations.
From this fixed point formula we have \begin{equation}\label{conservation_nonlin}   \forall t\in\R_+, \qquad N(t)+ \int_{0}^{\sigma(N(t))} n(s,t) \diff s = 1.   \end{equation}
Using \eqref{conservation_nonlin} and the equation for $n$, it is possible to prove that $N$ is a Lipschitz function. It then has a almost everywhere derivative that satisfies
\begin{equation}\label{hardcore}
		 		N'(t) \Big[   1  + \sigma'\big(  N(t)   \big)  n\big(  \sigma(N(t)) ,  t  \big)  \Big] = - N(t) + n\big(  \sigma(N(t)) ,  t  \big) \qquad \mathrm{a.e.\ in\ }  [\sigma^+,+\infty).
\end{equation}
One can then obtain $
	N'(t)\sigma'\big(  N(t) \big) \leqslant \norme{\sigma'}_\infty$,
which leads to the nonlinear version of \eqref{eq:T-E_Nmethod} after some tedious developments involving characteristics:
\begin{equation}\label{eq:T-E_Nmethod2}
	N(t) + \int_{t-\sigma(  N(t) )}^t N(s) \diff s = 1.
\end{equation}
Application of the Banach fixed point theorem to \eqref{eq:T-E_Nmethod2} proves that there is a unique stationary firing rate $N_\infty$. The argument above for the linear case can then be adapted to get the contraction
\[  |N(t) - N_\infty| \leqslant \left(\dfrac{\sigma^+}{1-\norme{\sigma'}_\infty N_\infty}\right)^k. \]

\begin{theorem}[Pakdaman, Perthame, Salort \cite{PPS2}]\label{thm:T-E_CV1}
    Grant Assumptions \ref{as:TE_p} and \ref{as:TE_n_0}. Assume that $\mu = \delta_0$, \textit{i.e.} $X(t) = N(t)$, and $\norme{\sigma'}_\infty < 1$ and 
    \[ \sigma(x) < \sigma^+ < 1-\norme{\sigma'}_\infty N_\infty. \]
    Then there exists a unique stationary firing rate $N_\infty$ and $\nu > 0$ such that the solution of \eqref{eq:T-E} satisfies $\frac{\diff \sigma(N(t))}{\diff t} < 1$ and 
    \[  \forall t\geqslant 2\sigma, \qquad |N(t)-N_\infty| \leqslant e^{-\nu t}. \]
\end{theorem}

\subsubsection{Reduction to a delay differential equation}
\label{sec:TE_delay}

Consider the non-delayed case $X(t)=N(t)$, and take now a firing function with the form $p(s,X) = \varphi(X)\mathds 1_{s > \sigma}$, where $\sigma>0$ is a constant and $\varphi$ a positive bounded function $0 < p^-\leqslant \varphi\leqslant p^+$ (then, by definition of $N(t)$, $0\leqslant N(t)\leqslant p^+$). The constant $\sigma>0$ is a fixed refractory period. When $\varphi$ is non-decreasing, the network is excitatory as more activity implies more discharges after the refractory time is elapsed. On the contrary, if $\varphi$ is non-increasing, the network influence is inhibitory: more activity in the network diminishes the quantity of neurons undergoing discharges.

Following \cite{CPST}, let us define the auxiliary function $\Psi:\R_+\to\R_+$ and its derivative by
\[  \psi(X) = \dfrac{X}{\varphi(X)},\qquad \psi'(X) = \dfrac{\varphi(X)-X\varphi'(X)}{\varphi^2(X)}.  \]
As we will see, complex dynamics can occur only when $\psi'$ changes sign on $[0,p^+]$.
\begin{itemize}
    \item In the inhibitory case $\varphi'\leqslant 0$, we always have $\psi'(X)>0$ for $X\in[0,p^+]$. 
    \item In the excitatory case $\varphi'\geqslant 0$, we can distinguish between a weakly-excitatory case if $\psi'(X)>0$ and a strongly-excitatory case when it changes sign on $[0,p^+]$. A sufficient condition for being in the weakly-excitatory case is for example $0 \leqslant \varphi'(X)\leqslant \frac{p^-}{p^+}$
\end{itemize}

The main idea to exploit the properties of $\psi$ in the study of \eqref{eq:T-E} is to reduce the full system to a delayed differential equation on the firing rate $N(t)$ that involves $\psi'$. This can be done using mass conservation and the methods of characteristics.

\begin{lemma}[C\'aceres, Perthame, Salort, Torres \cite{CPST}]\label{lm:TE_delay}
    Assume $\mu_X = \delta_{\tau=0}$, \textit{i.e.} $X(t)=N(t)$, $p(s,X) = \varphi(X)\mathds 1_{s > \sigma}$, where $\varphi\in C^1(\R_+)$, $\sigma>0$ and $0 < p^-\leqslant \varphi\leqslant p^+$. Then, given a solution to \eqref{eq:T-E} with initial condition $n^0$,
    \begin{itemize}
        \item for $0 < t < \sigma$,
        \begin{equation}
            \int_0^t N(\tau)\diff \tau + \int_0^{\sigma-t} n^0(s)\diff s + \psi\big(N(t)\big) = 1,
        \end{equation}
        and if $N$ is differentiable on $(0,\sigma)$,
        \begin{equation}
            \psi'\big(N(t)\big) N'(t) = -N(t) + n^0(\sigma-t);
        \end{equation}
        \item for $\sigma < t$,
        \begin{equation}\label{eq:T-E_integral}
            \int_{t-\sigma}^t N(\tau)\diff \tau + \psi\big(N(t)\big) = 1,
        \end{equation}
        and if $N$ is differentiable on $(\sigma,+\infty)$,
        \begin{equation}\label{eq:T-E_delay}
            \psi'\big(N(t)\big) N'(t) = - N(t) + N(t-\sigma).
        \end{equation}
    \end{itemize}
\end{lemma}
Then, \cite{CPST} proves that from a solution to the integral equation \eqref{eq:T-E_integral} and with compatibility conditions, it is possible to construct a solution to the full system \eqref{eq:T-E}.

\begin{theorem}[C\'aceres, Perthame, Salort, Torres \cite{CPST}]\label{thm:T-E_extension}
    Assume $\mu_X = \delta_{\tau=0}$, \textit{i.e.} $X(t)=N(t)$, $p(s,X) = \varphi(X)\mathds 1_{s > \sigma}$, where $\phi\in C^1(\R_+)$, $\sigma>0$ and $0 < p^-\leqslant \phi\leqslant p^+$. Let $N\in L^\infty(\R_+)$ such that $\psi\circ N\in C([\sigma,+\infty))\cap C^1((0,\sigma))$. Then if
    \begin{itemize}
        \item for all $t\in (0,\sigma)$, $\tilde n(t):= N(\sigma-t) + \frac{\diff \psi(N)}{\diff t}(\sigma-t) \geqslant 0$
        \item the function $N$ is solution to \eqref{eq:T-E_integral} on $[\sigma,+\infty)$.
    \end{itemize}
    for any initial condition $n^0$ of mass 1 such that for all $s\in(0,\sigma)$, $n^0(s)=\tilde n(s)$, there exists a unique solution to \eqref{eq:T-E} with firing rate $N$.
\end{theorem}

This theorem is proved by inputing $N$ in the linear system \eqref{eq:T-E_linear} and checking that the solution is also a solution to the nonlinear problem \eqref{eq:T-E}.

From the delay differential equation formulation, we can conclude that in both the inhibitory and the weakly-excitatory case, the firing rate $N(t)$ oscillates an infinite number of times. Indeed, if we reason by contradiction and assume that the firing rate is, for example, increasing over an interval of length $\sigma$, then by \eqref{eq:T-E_delay},
\[\psi'\big(N(t)\big) N'(t) = - N(t) + N(t-\sigma) < 0.\]
Using that $\psi'\big(N(t)\big)>0$, we get $N'(t) < 0$, which is a contradiction with $N$ increasing. This way, \cite{CPST} proves that $\psi'>0$ on $[0,p^+]$ implies the existence of an infinite quantity of local maxima and minima. 

Then, it can be proved that the amplitude of oscillations is actually decreasing: define $I_k = [k\sigma,(k+1)\sigma]$, then
\[  \max_{t\in I_{k-1}} N(t) \geqslant \max_{t\in I_k} N(t),\qquad \min_{t\in I_{k-1}} N(t) \leqslant \min_{t\in I_k} N(t).  \]
Last, defining the sequence $N_n(t) = N(t+n)$, by compactness there exists a subsequence converging to a limit $N_\infty(t)$ which defines a $C^1$ bounded function of time. Using $\psi'>0$ and Lemma \ref{lm:TE_delay}, $N_\infty$ is proved to be a constant function and in fact the firing rate of the only possible stationary state. Then, it it proved in \cite{CPST} that $N(t)\to N_\infty$ and $n(\cdot,t)\to n_\infty$.

\begin{theorem}[C\'aceres, Perthame, Salort, Torres \cite{CPST}]\label{thm:T-E_CV2}
    Assume $\mu_X = \delta_{\tau=0}$, \textit{i.e.} $X(t)=N(t)$, $p(s,X) = \varphi(X)\mathds 1_{s > \sigma}$, where $\varphi\in C^1(\R_+)$, $\sigma>0$ and $0 < p^-\leqslant \phi\leqslant p^+$. Assume that for all $X\in[0,p^+]$, $\psi'(X) > 0$. Then, there exists a unique stationary state $(n_\infty,N_\infty)$ to \eqref{eq:T-E} and for any solution to \eqref{eq:T-E}, $N(t)$ converges to $N_\infty$ when $t\to+\infty$.
\end{theorem}

As we have seen, the method relies on the sign of $\psi'$. We will see later in \ref{sec:TE_periodic} that when the system is not inhibitory nor weakly-excitatory, the ocillations can still happen but are not necessarily dampened towards the stationary state, as explicit periodic solutions can be constructed.

However, even in the strongly excitatory case, converge to a (now not necessarily unique) stationary state can happen: \cite{CPST} proves that under technical conditions on a function $N\in C([0,\sigma])$, a solution $(n,N)$ that prolongs $N$ on $[0,+\infty)$ can be constructed using Theorem \ref{thm:T-E_extension} and the extended $N$ monotonically converges towards a stationary state $N_\infty$. Numerical simulations in \cite{CPST} show that it happens when the solution remains in an area where the sign of $\psi'$ is constant. For example, they choose
\[  \varphi(X) = \frac{1}{1+e^{-9X+3.5}},\qquad \sigma =\frac12.  \]
In this case, $\psi$ changes sign and there are three stationary states $N_\infty^1,N_\infty^2,N_\infty^3$. Even more surprising, all of them can attract solutions from the same initial density. Choose
 $n^0(s) = \frac12 e^{-(s-1)_+}$.
Then, both the theoretical and the numerical solution are not uniquely defined, because the equation
\[  N(0) = \varphi\big(N(0)\big) \int_\sigma^{+\infty} n^0(s)\diff s \]
has three different solutions $N_1(0)$, $N_2(0)$ and $N_3(0)$. Each of them leads to a different numerical solution $n_1,n_2,n_3$ converging to $N_\infty^1,N_\infty^2,N_\infty^3$ respectively. Two of them converge while oscillating and the third does so monotonically. For the same $\phi$ and $\sigma$, some other initial conditions exhibit jumps on the firing rate. Those jumps preserve the continuity of $\psi(N(t))$. 

Contrary to the case of Assumption \ref{as:TE_p} described before, the speed of convergence is unknown in Theorem \ref{thm:T-E_CV2}; \cite{CPST} conjectures that the speed is exponential.

\subsection{Relaxation to a stationary state for a smooth discharge function}

In the first works \cite{PPS,PPS2} on the model \eqref{eq:T-E}, asymptotic convergence is studied in the case where the firing function $(s,X)\mapsto p(s,X)$ has the simple form $p(s,X) = \mathds 1_{s \geqslant \sigma(X)}$. This choice allows for explicit computations, predictable characteristics paths and an explicit construction for relative entropy decay. In later works \cite{MW, CY}, other methods are introduced which allow to prove weakly-nonlinear stability of a stationary state for other families of firing functions $p$.

\subsubsection{Abstract semigroup approach for the weakly-nonlinear problem}
\label{sec:TE_MW}

In \cite{MW}, the strategy of \textit{perturbation of semigroup} initiated in \cite{MM} and previously used in a neuroscience context in \cite{MQT} was applied to the Time Elapsed neuron model. It allows to prove a general weakly-nonlinear stability result.
\begin{theorem}[Mischler, Weng \cite{MW}]\label{thm:T-E_MW}
    Grant Assumption \ref{as:TE_p_smooth}. Assume either $\mu_X = J\delta_{\tau=0}$, $J\geqslant 0$, \textit{i.e.} $X(t) = JN(t)$ or Assumption \ref{as:TE_delay_exp}. Then, there exists $J_0 > 0$, $\nu>0$, $C\geqslant 1$, $\eta>0$ such that for all $J\in[0,J_0)$, for all initial condition $n^0\in L^1(\R_+)\cap\mathcal P(\R_+)$ and satisfying $\norme{n^0-n_\infty}_{L^1}<\frac{\eta}{J}$,
    \[ \forall t\geqslant 0, \qquad \norme{n(\cdot,t)-n_\infty(\cdot)}_{L^1} \leqslant C e^{-\nu t}.  \]
\end{theorem}

The method proceeds as follows.

\paragraph{Case without delay.} Consider the case $\mu_X = J \delta_0$, which yields $X(t) = JN(t)$.

Let us write the stationary state perturbation
\[  n(s,t) = n_\infty(s) + m(s,t), \qquad N(t) = N_\infty + M(t). \]
Then the linearised system for the pair $(m,M)$ is
\begin{equation}\label{eq:T-E_lin_nodelay}
    \left\{\begin{array}{l}
           \dfrac{\partial m}{\partial t} + \dfrac{\partial m}{\partial s} + p(s,JN_\infty)m + J \dfrac{\partial p}{\partial X}(s,JN_\infty)n_\infty M = 0,\\
           \displaystyle M(t) = \int_0^{+\infty} \left( p(s,JN_\infty)m(s,t) + J\dfrac{\partial p}{\partial X}(s,JN_\infty)n_\infty(s)M(t) \right) \diff s,
           \\
           m(0,t) = M(t), \quad m(s,0) = m^0(s).
    \end{array}\right.
\end{equation}
Define
\[ \kappa = J\int_0^{+\infty} \dfrac{\partial p}{\partial X}(s,JN_\infty)n_\infty(s) \diff s. \]
For $J$ small enough, $\kappa < 1$, and solving for $M(t)$ in the second equation of \eqref{eq:T-E_lin_nodelay} yields
\[  M(t) = \dfrac{1}{1-\kappa} \int_0^{+\infty} p(s,JN_\infty) m(s,t) \diff s  = : \mathcal M_J\big[m(\cdot,t)\big].  \]
Defining this linear operator $\mathcal M_J$ allows to rewrite the linearised system with a linear operator acting on $m$:
\[ \dfrac{\partial m}{\partial t} = L_J\, m := -\dfrac{\partial m}{\partial s} - p(s,JN_\infty)m - J \dfrac{\partial p}{\partial X}(s,JN_\infty)n_\infty \mathcal M_J[m].  \]
A semigroup $S_{L_J}(t)$ with generator $L_J$ acting on $ L^1(\R_+)$ is then classically associated to \eqref{eq:T-E_lin_nodelay} and it has the domain
\[  D(L_J) = \{ m\in W^{1,1}(\R_+) \ | \ m(0) = \mathcal M_J(m)  \}. \]
The linearised system can also be rewritten as
\begin{equation}
    \dfrac{\partial m}{\partial t} = \Lambda_J\, m := -\dfrac{\partial m}{\partial s} - p(s,JN_\infty)m - J \dfrac{\partial p}{\partial X}(s,JN_\infty)n_\infty \mathcal M_J[m] + \delta_{s=0} \mathcal M_J[m],
\end{equation}
allowing to define a semigroup $S_{\Lambda_J}(t)$ with generator $\Lambda_J$ acting on the space of bounded Radon measures
\[ \mathcal X = M^1(\R_+) = \{ m\in C_{\to0}(\R)'\ |\ \supp(m)\subset \R_+ \} \]
endowed with the weak-$*$ topology $\sigma(M^1,C_{\to0})$. It can be checked by duality that ${S_{\Lambda_J}}_{|L^1(\R_+)} = L_J$.

A core idea of \cite{MW} is to split the operator $\Lambda_J$ as
\begin{equation}
    \Lambda_J = \mathcal A_J + \mathcal B_J, \qquad \mathcal A_J\, m := \nu_J \mathcal M_J[m],\qquad \mathcal B_J\, m := -\dfrac{\partial m}{\partial s} - p(s,JN_\infty) m,
\end{equation}
where $\gamma_J = \delta_{s=0} + J \partial_X p(s,J N_\infty)n_\infty$.
The authors of \cite{MW} then propose to adapt the strategy in \cite{MS2,MS3} relying on Weyl's theorem and the Spectral Mapping theorem. However, an obstacle is that $\Lambda_J$ does not generate a strongly continuous semigroup on $\mathcal X$, but only a weakly-$*$ continuous semigroup. Among diverse ways to circonvent the issue, \cite{MW} chooses to rely on the splitting and the Banach structure of $L^1(\R_+)\subset \mathcal X$.

\begin{definition}[Hypodissipativity \cite{GMM}]\label{def:T-E_hyp}
    A closed operator $L$ on a Banach space $(E,\norme{\cdot})$ with dense domain $D(L)$ is $\beta$-hypodissipative if there exists a norm $|||\cdot|||$ on $E$ equivalent to $\norme{\cdot}$ such that
    \[  \forall x\in D(L),\ \exists \phi\in F_{|||\cdot|||}(x), \quad \mathcal Re\pscal{\phi}{(L-\beta)\,x}\leqslant 0 \]
    where
    \[ F_{|||\cdot|||}(x) = \{ \phi\in E' \ | \ \pscal{\phi}{x} = |||x|||^2 = |||\phi|||_{E'} \}. \]
\end{definition}
An operator $L:D(L)\to E$ being hypodissipative is equivalent to the growth estimate
\[  \norme{S_L(t)}_{\mathscr B(E)} \leqslant M e^{\beta t},    \]
on the semigroup $S_L$ associated to $L$, where $\mathscr B(E)$ is the space of linear and bounded operators on $E$.

\begin{lemma}[Mischler, Weng \cite{MW}]
    Grant Assumption \ref{as:TE_p_smooth}. Denote $\beta_0=-\frac{p^-}{2}$. For all $J\geqslant 0$,
    \begin{itemize}
        \item $\mathcal A_J\in \mathscr B(W^{-1,1}(\R_+),\mathcal X)$.
        \item $\mathcal B_J$ is $\beta_0$-hypodissipative in both $L^1(\R_+)$ and $\mathcal X$.
        \item The family of operators $S_{\mathcal B_J} \ast \mathcal A_J S_{\mathcal B_J}$ satisfies
        \[  \forall \beta\in (\beta_0,+\infty),\qquad  \norme{(S_{\mathcal B_J} \ast \mathcal A_J S_{\mathcal B_J})(t)}_{\mathscr B(\mathcal X,\mathcal Y)} \leqslant C_\beta\, e^{\beta t}, \]
        where $C_\beta > 0$, $\mathcal Y = BV(\R_+)\cap L^1_1(\R_+)$, $BV$ is the space of functions with bounded variation and $L^1_1$ the Lebesgue space with weight $s\mapsto <s>$.
    \end{itemize}
\end{lemma}

Applying this lemma, adapting Weyl's theorem and the Spectral Mapping theorem \cite{MS2,MS3} and taking advantage of the splitting introduced above, the article \cite{MW} characterises the spectrum $\Sigma(\Lambda_J)$ of $\Lambda_J$ on $\mathcal X$.

\begin{theorem}[Mischler, Weng \cite{MW}]\label{thm:T-E_spectral_noD}
    Grant Assumption \ref{as:TE_p_smooth}. Denote $\beta_0=-\frac{p^-}{2}$. For all $J\geqslant 0$, the operator $\Lambda_J$ is the generator of a weekly-$*$ continuous semigroup $S_{\Lambda_J}$ acting on $\mathcal X$ endowed with the weak-$*$ topology. There exist a finite rank projector $\Pi_J$ which commutes with $S_{\Lambda_J}$, $j\in\N$ and complex numbers
    \[  \xi_1,\dots,\xi_j \in \Delta_{\beta_0} : = \{ z\in\mathbb C\ | \ \mathcal Re z > \beta_0 \}, \]
    such that $\Sigma({\Lambda_J}_{|\Pi_J\mathcal X})\cap \Delta_{\beta_0} = \{ \xi_1,\dots,\xi_j\}$ and for any $\beta\in(\beta_0,+\infty)$, there exists $C_\beta$ such that the remainder semigroup satisfies
    \[\forall t\geqslant 0, \quad \norme{S_{\Lambda_J} (I_\mathcal X-\Pi_J)}_{\mathscr B(\mathcal X)} \leqslant C_\beta e^{\beta t}.\]
\end{theorem}

In the linear case $J=0$, the semigroup is positive and \cite{MW} proves the following version of the Krein-Rutman theorem: there exists $\beta < 0$, $C > 0$ such that $\Sigma(\Lambda_0)\cap \Delta_{\beta} = \{0\}$ and for any $m^0\in \mathcal X$ of average 0,
\[ \forall t\geqslant 0,\qquad \norme{S_{\Lambda_0}(t) m^0}_\mathcal X \leqslant Ce^{\beta t}\norme{m^0}_\mathcal X. \]
Then, it is possible to use the stability theory for semigroups \cite{K,MM,MT,T} to apply a perturbation argument and prove stability for the weakly connected linearised system.

\begin{theorem}[Mischler, Weng \cite{MW}]
    Grant Assumption \ref{as:TE_p_smooth}. There exists $J_0$, $\beta < 0$, $C > 0$ such that for all $J\in[0,J_0)$, $\Sigma(\Lambda_J)\cap \Delta_{\beta} = \{0\}$ and for any $m^0\in \mathcal X$ of average 0,
\[ \forall t\geqslant 0,\qquad \norme{S_{\Lambda_0}(t) m^0}_\mathcal X \leqslant Ce^{\beta_0 t}\norme{m^0}_\mathcal X. \] 
\end{theorem}
For the weakly-nonlinear case, consider the functional
\[  \Phi[m,N] = \int_{0}^{+\infty} p(s,JN) - N.    \]
Using the $W_1$ Wasserstein distance on the set $\mathcal P(\R_+)$ of probability measures, \cite{MW} proves (similar to the key uniqueness result in \cite{PPS2} under Assumption \ref{as:TE_p}) that for $J$ small enough, there is a function $\phi_J: \mathcal P(\R_+)\to \R$, Lipschitz continuous for the weak topology of probability measures, such that $N = \phi_J(m)$ is the unique solution to $\Phi(m,N) = 0$, and thus the unique stationary firing rate associated to the perturbation $m$. Then, using this result, the decay for the linearised system and a Grönwall estimate, it is possible to prove that the perturbation $m$ decays exponentialy fast in $L^1(\R_+)$ for the nonlinear problem, hence proving Theorem \ref{thm:T-E_MW} in the non-delayed case.

\paragraph{Case with delay.} Consider the case of a delayed interaction with an exponentially decaying regular time kernel as in Assumption \ref{as:TE_delay_exp}. The strategy described above for the non-delayed case can be adapted using a modified functional framework. First, to write tbe linearised system as a time-autonomous problem, introduce a local time $\tau$ and, given the perturbation $M(t)$ of the firing rate $N(t) = N_\infty + M(t)$, transport the value in time through
\[  \dfrac{\partial v}{\partial t} + \dfrac{\partial v}{\partial \tau} = 0,\qquad v(0,t) = M(t), \qquad v(\tau,0) = 0.  \]
The solution is  $v(\tau,t) = M(t-\tau)\mathds 1_{0\leqslant \tau \leqslant t}$.
The variation $x$ of $X=X_\infty + Y(t)$ is characterised by
\[  Y(t) = \mathcal D[v]:= J\int_{0}^{+\infty} v(\tau,t) \alpha(\tau)\diff \tau \]
and we have, given the perturbation $n(s,t) = n_\infty(s) + m(s,t)$, the linearised equation
\begin{equation}
    M(t) = \int_0^{+\infty} p(s,X_\infty) m(s,t) \diff s + Y(t) J \int_0^{+\infty} \dfrac{\partial p}{\partial X}(s,X_\infty) n_\infty(s) \diff s,
\end{equation}
which can be rewriten in operator form
\begin{equation}
    M(t) = \mathcal O_J[m,v] : = \mathcal N_J[m] + \kappa \mathcal D[v], \qquad \kappa = \int_0^{+\infty} \dfrac{\partial p}{\partial X}(s,X_\infty) n_\infty(s) \diff s.
\end{equation}
This leads to rewritting the linearised system in the delayed case in the form
\begin{equation}
    \dfrac{\partial }{\partial t} \left(\begin{matrix} m\\v
    \end{matrix}\right) = \mathcal L_J \left(\begin{matrix} m\\v \end{matrix}\right),
\qquad
    \mathcal L_J \left(\begin{matrix} m\\v \end{matrix}\right) = \left(\begin{matrix} 
    -\partial_s m - p(s,X_\infty)m - \partial_X p(s,X_\infty) n_\infty \mathcal D[v]\\
    -\partial_\tau v
    \end{matrix}\right),
\end{equation}
where the domain of $\mathcal L_J$ is
\[  D(\mathcal L_J) = \{ (m,v)\in W^{1,1}(\R_+)\times W^{1,1}(\R_+,\omega \} \ | \ m(0) = v(0) = \mathcal O_J[m,v] \},\quad \omega(s) = e^{-\lambda s}, \]
with $\lambda > 0$ from Assumption \ref{as:TE_delay_exp}. The semigroup $S_{\mathcal L_J}$ acts on $L^1(\R_+)\times L^1(\R_+,\omega)$. As in the non-delayed case, let us introduce the larger space $\mathcal X = M^1(\R_+)\times M^1(\R_+,\omega)$, and the extended operator $\Lambda_J$ defined as
\begin{equation}
    \Lambda_J \left(\begin{matrix} m\\v \end{matrix}\right) = \left(\begin{matrix} -\partial_s m - p(s,X_\infty) m - \partial_X p(s,X_\infty)n_\infty \mathcal D[v] + \delta_{s = 0} \mathcal O_J[m,v]  \\ - \partial_\tau v + \delta_{\tau = 0} \mathcal O_J[m,v] \end{matrix}\right).
\end{equation}
It is again possible to decompose it in the form $\Lambda_J = \mathcal A_J + \mathcal B_J$ and to obtain a semigroup satisfying ${S_{\Lambda_J}}_{|L^1(\R_+)\times L^1(\R_+,\omega)} = S_{\mathcal L_J}$. From there, it is possible to adapt the arguments in the non-delayed case and to prove the other part of Theorem \ref{thm:T-E_MW}. Let us mention the main intermediary results about the operator $\Lambda_J$ and its spectrum, as they are useful in themselves.
\begin{lemma}[Mischler, Weng \cite{MW}]
     Grant Assumptions \ref{as:TE_p_smooth} and \ref{as:TE_delay_exp}. Let $\beta_0 = \max(-\frac{p^-}{2},-\lambda) < 0$.
     \begin{itemize}
         \item $\mathcal A_J\in \mathscr B(W^{-1,1}(\R_+)\times W^{-1,1}(\R_+,\omega),\mathcal X)$.
        \item $\mathcal B_J$ is $\beta_0$-hypodissipative in both $L^1(\R_+)\times L^1(\R_+,\omega)$ and $\mathcal X$.
        \item The family of operators $S_{\mathcal B_J} \ast \mathcal A_J S_{\mathcal B_J}$ satisfies
        \[  \forall \beta\in (\beta_0,+\infty),\qquad  \norme{(S_{\mathcal B_J} \ast \mathcal A_J S_{\mathcal B_J})(t)}_{\mathscr B(\mathcal X,\mathcal Y)} \leqslant C_\beta\, e^{\beta t}, \]
        where $C_\beta > 0$, $\mathcal Y = BV(\R_+)\cap L^1_1(\R_+)\times BV(\R_+,\omega)\cap L^1_1(\R_+,\omega)$.
     \end{itemize}
\end{lemma}
Using this lemma, it is possible to obtain the same spectral result as in the non-delayed case.
\begin{theorem}[Mischler, Weng \cite{MW}]
     Grant Assumptions \ref{as:TE_p_smooth} and \ref{as:TE_delay_exp}. The conclusion of Theorem \ref{thm:T-E_spectral_noD} holds with $\beta_0 = \max(-\frac{p^-}{2},-\lambda) < 0$.
\end{theorem}

\subsubsection{Back to non-smooth firing functions and the strong connectivity case}

In \cite{MQW}, the method of \cite{MW,W} and ideas from \cite{PPS,PS2} are combined with refined spectral theory \cite{MS2}, leading to convergence results for more general firing functions $p$ and allowing to consider the strong connectivity regime. It namely extends the strong connectivity convergence result of \cite{PPS} described above in Theorem \ref{thm:T-E_CV_PPS1} and extends the weakly-nonlinear convergence results of Theorem \ref{thm:T-E_MW} to less regular firing functions like those satisfying Assumption \Ref{as:TE_p}.

Denote $L^1_q(\R_+)$ the space of measurable functions $f$ such that $\norme{f}_{L^1_q(\R_+)} := \int_0^{+\infty}(1+s^q)|f|\diff s < +\infty$ and $L^1_w$ the space $L^1(\R_+)$ endowed with the weak topology $\sigma(L^1,L^\infty)$.

Let us define a generalised weak nonlinearity condition.
\begin{hyp}\label{as:TE_MQW1}
    There exists $L>0$ small enough such that for all $x_0 > 0$ there exists $J_0 > 0$ such that
        \begin{equation*}
            \forall J\in(0,J_0),\ \forall X_1,X_2\in\left[0,\frac{x_0}{J}\right),\qquad \int_0^{+\infty} |p(s,X_1)-p(s,X_2)| \diff s \leqslant \frac{L}{J}|X_1-X_2|.
        \end{equation*}
\end{hyp}
Similarly we have the strong nonlinearity condition
\begin{hyp}\label{as:TE_MQW2}
    There exists $L>0$ small enough such that for all $x_\infty > 0$ there exists $J_\infty > 0$ such that
        \[  \forall J\in(J_\infty,+\infty),\ \forall X_1,X_2\in\left(\frac{x_\infty}{J},+\infty\right),\qquad \int_0^{+\infty} |p(s,X_1)-p(s,X_2)| \diff s \leqslant \frac{L}{J}|X_1-X_2|. \]
\end{hyp}
Note that these assumptions cover many of the firing functions in Assumption \ref{as:TE_p} used in \cite{PPS,PPS2}.

\begin{theorem}[Mischler, Qui\~ninao, Weng \cite{MQW}]\label{thm:T-E_MQW1}
    Assume $p\geqslant 0$, $\partial_s p\geqslant 0$, $\partial_X p\geqslant 0$, \[0 < \lim_{s\to+\infty} p(s,0)  \leqslant \lim_{s,X\to+\infty} p(s,X) < +\infty,\] and for all $s\geqslant 0$, $\int_0^s p(w,\cdot)\diff w \in C(\R_+)$. Then for all $q > 0$, for all $n^0\in L^1_q(\R_+)\cap L^\infty(\R_+)\cap \mathcal P(\R_+)$, if one of the following three conditions hold:
    \begin{enumerate}
        \item $\mu_X = J \alpha$ with $\alpha \in L^1(\R_+)$ satisfying $\int_0^{+\infty} \alpha(\tau)e^{\delta \tau}\diff \tau < +\infty$;
        \item $\mu_X = J \delta_{s=0}$ with $J>0$ small enough and Assumption \ref{as:TE_MQW1} holds;
        \item $\mu_X = J \delta_{s=0}$ with $J>0$ large enough, Assumption \ref{as:TE_MQW2} holds and $\int_0^{+\infty}p(s,0)n^0(s)\diff s > 0$;
    \end{enumerate}
     there exists a global-in-time weak solution $n\in C(\R_+,L^1_w(\R_+))\cap L^\infty(\R_+,L^1_q(\R_+))\cap L^\infty(\R_+^2)$, $N,X\in C(\R_+)$ to \eqref{eq:T-E}. When case 2. or 3. is assumed, this solution is unique.
\end{theorem}
The reason the solution is not proved to be unique in case 1. is that a Schauder fixed point theorem is used. Cases 2. and 3. yield uniqueness as they rely upon the Banach fixed point theorem.

The proof of existence of stationary states $(n_\infty,N_\infty,X_\infty)$ with $n_\infty\in W^{1,+\infty}(\R_+)$ in \cite{MW} is also extended to more general firing functions $p$ in \cite{MQW}. When $J$ is small enough or large enough, they prove that there exists a unique stationary state.

\begin{theorem}[Mischler, Qui\~ninao, Weng \cite{MQW}]\label{thm:T-E_MQW2}
    Grant the assumptions of Theorem \ref{thm:T-E_MQW1}. There exists $J_0 > 0$ small enough and $J_\infty > 0$ large enough, $\nu >0$ and $C\geqslant 1$ such that if $J\in (0,J_0)\cup(J_\infty,+\infty)$, the solution to \eqref{eq:T-E} satisfies
    \begin{equation}
        \forall t\geqslant 0, \qquad \norme{n(\cdot,t)-n_\infty}_{L^1} \leqslant Ce^{-\nu t}.
    \end{equation}
\end{theorem}

\subsubsection{Convergence using the Doeblin-Harris method}
\label{sec:TE_Doeblin}

In \cite{CY}, a general convergence result is achieved in the space of probability measures using the Doeblin-Harris machinery. After the seminal \cite{H3}, these methods have been revisited between 2000 and 2010, namely in \cite{G,S,HM}, and they have since been applied to a variety of linear and nonlinear partial differential equations; see \cite{Y} for a review. We start by reproducing here the introduction provided by \cite{CY}.

\begin{definition}[Stochastic operators and semigroups]
    Let $(E,\mathcal A)$ be a measurable space. We call a stochastic operator any linear operator $S:\mathcal M(E)\to\mathcal M(E)$ such that
    \begin{equation}
        \forall \mu \in \mathcal M(E),\qquad \int_E S\mu = \int_E \mu, \qquad  \mu\geqslant 0\implies S\mu \geqslant 0.
    \end{equation}
    A stochastic semigroup is a semigroup $(S_t)_{t\geqslant 0}$ of operators $S_t:\mathcal M(E)\to\mathcal M(E)$ such that for all $t\geqslant 0$, $S_t$ is a stochastic operator. An equilibrium or stationary state of a stochastic semigroup is a finite measure $\mu_\infty\in \mathcal M(E)$ such that for all $t\geqslant 0$, $S_t\mu_\infty = \mu_\infty$.
\end{definition}

Equivalently, stochastic operators on $E$ preserve the space of probability measures $\mathcal P(E)$.

\begin{definition}[Doeblin condition \cite{G}]\label{def:Doeblin}
    Let $(E,\mathcal A)$ be a measurable space. We say that a stochastic operator $S$ satisfies the Doeblin condition if there exists $\alpha\in(0,1)$ and a probability measure $\nu\in \mathcal P(E)$ such that,
    \begin{equation}
        \forall \mu \in \mathcal P(E),\qquad S\mu \geqslant \alpha\nu.
    \end{equation}
\end{definition}

\begin{theorem}[Doeblin's theorem \cite{G}]
    Let $(E,\mathcal A)$ be a measurable space and $(S_t)_{t\geqslant 0}$ a stochastic semigroup on $\mathcal M(E)$. If there exists $t_0>0$ such that $S_{t_0}$ satisfies the Doeblin condition with constant $\alpha\in(0,1)$, then the semigroup has a unique equilibrium $n_\infty \in \mathcal P(E)$, and
    \begin{equation}
        \forall n^0\in\mathcal P(E),\qquad \norme{S_t(n^0-n_\infty)}_{\mathrm{TV}} \leqslant \dfrac{1}{1-\alpha} e^{-\lambda t}\norme{n^0-n_\infty}_{\mathrm{TV}}
    \end{equation}
    where $\lambda = -\frac{1}{t_0}\log(1-\alpha)>0$.
\end{theorem}
In order to apply Doeblin's theorem, \cite{CY} focuses first on the linear problem $p(s,X) = \bar p(s)$. Then, as seen in Subsections \ref{sec:TE_exist} and \ref{sec:TE_CV_1} for a simpler $\bar p$, the method of characteristics provides explicit estimates on the solution and allow to construct solutions without any smallness assumption. This extends to the framework of measure solutions: in the linear case, there exists a global-in-time measure solution to \eqref{eq:T-E} for any probability $n^0\in\mathcal P(\R_+)$ and any $L$ in Assumption \ref{as:TE_p_meas}. Using the linear characteristics and performing lower bounds on them, the Doeblin condition can be proved in the form
\[  n(2\sigma) \geqslant  p^- \sigma e^{-2p^+ \sigma}\mathds 1_{(0,\sigma)} \]

\begin{lemma}[Ca\~nizo, Yoldaş \cite{CY}]
    Grant Assumptions \ref{as:TE_p_meas} and \ref{as:TE_p_meas_2}. Assume for all $s,X\geqslant 0$, $p(s,X) = \bar p(s)$. Then there exists a stochastic semigroup $(S_t)_{t\geqslant 0}$, such that for all $n^0\in \mathcal P(\R_+)$, $t\mapsto S_t n^0$ is the unique measure solution to \eqref{eq:T-E}. Moreover, $(S_t)_{t\geqslant 0}$ satisfies the Doeblin condition with $t_0 = 2\sigma$ and $\alpha = p^- \sigma e^{-2p^+ \sigma}$, $\nu = \frac{1}{\sigma} \mathds 1_{(0,\sigma)}$.
\end{lemma}
Then, Doeblin's theorem allows to straightforwardly conclude that the linear Time Elapsed equation has a spectral gap in the space of probability measures.

For the weakly-nonlinear case, the idea is to rewrite and decompose the evolution problem as
\begin{equation*}
    \dfrac{\partial n}{\partial t} = \mathcal L_N n = \mathcal L_{N_\infty}n + h
\end{equation*}
with the operator
\begin{equation*}
    \mathcal L_N n = -\dfrac{\partial n}{\partial s} - p(\cdot,N) n + \delta_{s=0} \int_0^{+\infty} p(u,N) n(\diff u).
\end{equation*}
and the perturbation
\begin{equation*}
    h = [p(s,N_\infty) - p(s,N)] n + \delta_{s=0} \int_0^{+\infty} [p(s,N) - p(s,N_\infty)] n(\diff u)
\end{equation*}
It can be proved that under a smallness condition on $L$ (condition \eqref{eq:T-E_L_cond_3} below),
\begin{equation}\label{eq:T-E_stab}
    \norme{h(t)}_{\mathrm{TV}} \leqslant \frac{2p^+L}{1-L}\norme{n(t) - n_\infty}_{\mathrm{TV}}, \qquad \int_0^{+\infty} h(t)(\diff s) = 0
\end{equation}
The second step is to use the Duhamel formula
\begin{equation*}
    n(t) = S_t n^0 + \int_0^t S_{t-\tau} h(\cdot,\tau) \diff \tau,
\end{equation*}
which yields after subtracting the stationary state and taking the total variation norm
\begin{equation*}
    \norme{n(t)-n_\infty}_{\mathrm{TV}} \leqslant \norme{S_t(n^0-n_\infty)}_{\mathrm{TV}} + \norme{ \int_0^t S_{t-\tau} h(\cdot,\tau) \diff \tau}_{\mathrm{TV}}.
\end{equation*}
Then, using the linear decay result for the first term in the right-hand side, \eqref{eq:T-E_stab} and Grönwall's inequality for the second part of the right-hand side, one can obtain the following result.

\begin{theorem}[Ca\~nizo, Yoldaş \cite{CY}]
    Let $n^0\in\mathcal P(\R_+)$. Grant Assumptions \ref{as:TE_p_meas} and \ref{as:TE_p_meas_2} and $\mu_X = \delta_{\tau=0}$ with
    \begin{equation} 
    \label{eq:T-E_L_cond_3}L\leqslant \frac{1}{4\norme{n^0}_{\mathrm{TV}}}. 
    \end{equation}
    Then, there exists a unique measure stationary solution $(n_\infty,N_\infty)$ in the sense of Definition \ref{def:T-E_meas_sol} and there exists constants $C\geqslant 1,\lambda>0$ independent of $n^0$ such that the solution to \eqref{eq:T-E} satisfies
    \begin{equation}
        \forall t\geqslant 0,\qquad \norme{n(t)-n_\infty}_{\mathrm{TV}}\leqslant Ce^{-\lambda t}\norme{n^0-n_\infty}_{\mathrm{TV}}.
    \end{equation}
\end{theorem}

The constants $C$ and $\lambda$ in this result are constructive; in \cite{CY}, explicit formulae are provided.

\subsection{Self-sustained oscillations}
\label{sec:TE_periodic}

\subsubsection{Explicit construction in the step function case}

Theorem \ref{thm:T-E_CV1} can be considered sharp in some sense, as there are analytic examples of periodic solutions when its assumptions are not satisfied.

Consider again the non-delayed case $X(t) = N(t)$ and the step firing function $p(s,X) = \mathds 1_{s > \sigma(X)}$ as in Assumption \ref{as:TE_p}. Given the firing rate $N(t)$, the system \eqref{eq:T-E} becomes linear. Hence, a way to construct a periodic solution is to find a suitable $T$-periodic function $t\mapsto N(t)$ such that the solution to the linear hyperbolic boundary-value problem
\begin{equation}
    \dfrac{\partial n}{\partial t} + \dfrac{\partial n}{\partial s} + p(s,N(t)) = 0,\qquad n(0,t) = N(t), \qquad n(s,0) = n(s,T), 
\end{equation}
is a solution to the full system \eqref{eq:T-E}. Because there is no access to $n$, the mass conservation has to be translated as a property on $N$ like in \eqref{eq:T-E_Nmethod2}. More precisely, one can prove:

\begin{theorem}[Pakdaman, Perthame, Salort \cite{PPS2}]\label{thm:T-E_periodic}
    Grant Assumption \ref{as:TE_p} and assume $\mu_X = \delta_0$, \textit{i.e.} $X(t) = N(t)$. Let $N\in L^\infty(\R)$ a periodic function such that, in the sense of distributions,
    \begin{equation}\label{eq:T-E_condition_periodic}
        \dfrac{\diff }{\diff t}(\sigma\circ N) \leqslant 1,\qquad  N(t)+ \int_{0}^{\sigma(N(t))} N(t-s) \diff s = 1.
    \end{equation}
    Then there exists a periodic solution of \eqref{eq:T-E} with firing rate $N$.
\end{theorem}

\begin{figure}[ht!]
	\begin{center}
		\includegraphics[width=300pt]{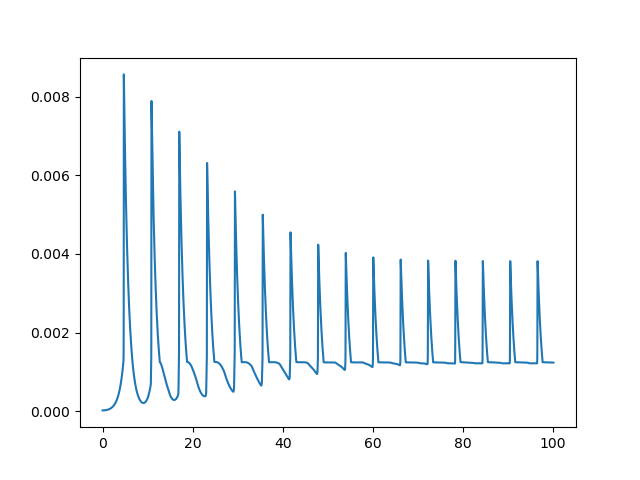}
	\end{center}
	\caption{Numerical simulation of the firing rate $N(t)$ of a solution to \eqref{eq:T-E} with $\mu_X = J \delta_{\tau=0}$, $J=20$, $\sigma$ as in \eqref{eq:T-E_sigma2} with $\alpha=3$, $p(s,X) = \mathds1_{s\geqslant \sigma(X)}$, $n^0(s) = e^{-s}$, original simulation by us, numerical scheme from \cite{PPS2}.\label{fig:T-E_periodic}}	
\end{figure}

Proving existence of periodic solutions to \eqref{eq:T-E} is then reduced to finding periodic functions satisfying \eqref{eq:T-E_condition_periodic}. Such functions can have discontinuity points. Then, $\frac{\diff }{\diff t}(\sigma\circ N)$ in the sense of distributions means that $\sigma(N)$ has a negative jump, which is a positive discontinuity jump for $N$ given that $\sigma$ is non-increasing by Assumption \ref{as:TE_p}. The reason for these discontinuities is the possibility for non-uniqueness of $N$ in the fixed point formulation \eqref{eq:T-E_fixed} when $\norme{\sigma'}_\infty$ exceeds one. Those positive jumps are a feature of the system and a core mechanism in non-delayed periodic solutions to \eqref{eq:T-E}, as can be seen on Figure \ref{fig:T-E_periodic} and the explicit analytical examples below.

The article \cite{PPS2} proposes two kinds of explicit construction using Theorem \ref{thm:T-E_periodic}.

First, choose a number $n$ of discontinuities and an increasing sequence $0=\alpha_0 < \alpha_1<\cdots<\alpha_{n+1}$. Then, define the quantities
\begin{align*}
    N^+_n = \dfrac{1}{1+\sum_{j=1}^{n-1}(e^{\alpha_{j+1}-\alpha_j}-1 ) + \alpha_{n+1}-\alpha_n  },\qquad N^+_i = e^{\alpha_{i+1}-\alpha_i} N_n^+,\ \ 0\leqslant i \leqslant n-1, 
\end{align*}
and the periodic function
\begin{equation}\label{eq:T-E_N_periodic_1}
    N(t) = \left\{
    \begin{array}{ll}
         N_i^+ e^{\alpha_i-t},\qquad&  t\in (\alpha_i,\alpha_{i+1}),\ \ 0 \leqslant i \leqslant n-1 \\[0.1cm]
         N_n^+, &  t\in (\alpha_n,\alpha_{n+1}),\\[0.1cm]
         N(t-k\alpha_{n+1}),& t\in [k\alpha_{n+1},(k+1)\alpha_{n+1}),\ \ k\in \N,\ k>0.
    \end{array}
    \right.
\end{equation}
If we define 
\[  \underline N = N_n^+, \qquad \overline N = \sup_{0\leqslant i\leqslant n} N_i^+,\qquad \sigma_0 = \alpha_{n+1}-\alpha_1 + \log\left(\dfrac{N_0^+}{\overline N}\right), \]
and the Lipschitz function
\begin{equation}\label{eq:T-E_sigma1}
    \sigma(X) = \left\{
    \begin{array}{ll}
         \sigma_0 - \log(\underline N) + \log(\overline N),\qquad& X\in \left[\,0,\underline N\,\right], \\[0.1cm]
        \sigma_0 - \log(x) + \log(\overline N),&  X\in \left[\,\underline N,\overline N\,\right], \\[0.1cm]
         \sigma_0, &  X\in \left[\,\overline N,+\infty\,\right],
    \end{array}
    \right.
\end{equation}
there exists a periodic solution to \eqref{eq:T-E} with firing rate $N$ as defined in \eqref{eq:T-E_N_periodic_1}

The authors of \cite{PPS2} construct a second kind of periodic solution. Choose a number $\alpha>0$ and set
\[  \underline N = \dfrac{1}{2e^\alpha - 1 } ,\qquad \overline N =  \dfrac{e^\alpha}{2e^\alpha - 1 }.  \]
Similar to \eqref{eq:T-E_sigma1} we use the Lipschitz function
\begin{equation}\label{eq:T-E_sigma2}
    \sigma(X) = \left\{
    \begin{array}{ll}
         2\alpha,\qquad& X\in \left[\,0,\underline N\,\right], \\[0.1cm]
        \alpha - \log(x) + \log(\overline N),&  X\in \left[\,\underline N,\overline N\,\right], \\[0.1cm]
         \alpha, &  X\in \left[\,\overline N,+\infty\,\right].
    \end{array}
    \right.
\end{equation}
Then, choose any $q\in(0,\alpha)$. There is a solution $\delta = \delta(q)\in [q+\alpha,2\alpha]$ to the problem
\[   e^{\alpha - \delta} \left(-2\alpha +1 + e^\alpha + \delta\right) + e^{-\alpha} \left( \delta-q-\alpha+e^q-1\right) = 2  \]
and there exists a solution of \eqref{eq:T-E} with firing rate the $2\alpha$-periodic function defined on $[0,2\alpha]$ by
\begin{equation}\label{eq:T-E_N_periodic_2}
    N(t) = \left\{
    \begin{array}{ll}
         \overline N e^{-t},\qquad&  t\in [0,\alpha], \\[0.1cm]
         e^q \underline N e^{-t+\alpha}, &  t\in (\alpha,q+\alpha],\\[0.1cm]
         \underline N,& t\in (q+\alpha,\delta(q)],\\[0.1cm]
         \overline N e^{-t+\alpha} \left(t-2\alpha+e^\alpha\right),& t\in (\delta(q),2\alpha].
    \end{array}
    \right.
\end{equation}
Note the particular case $q=\alpha$, $\delta(q) = 2\alpha$.

These numerous analytical examples and further numerical insights \cite{PPS,PPS2} show that for fixed parameters there can be a great variety of periodic solutions. In the delayed cases, the firing rate may be periodic without exhibiting discontinuities. Characterising the stability profile of all those co-existing periodic solutions is an open question.

\subsubsection{Construction from a delayed differential equation}

Let us now come back to the framework of Subsection \ref{sec:TE_delay}, which is the non-delayed case $\mu_X = \delta_{\tau=0}$ ($X(t)=N(t)$), and the firing function $p(s,X) = \varphi(X)\mathds 1_{s > \sigma}$, $\sigma>0$, $0 < p^-\leqslant \phi\leqslant p^+$. Recall the auxiliary functions
\[  \psi(X) = \dfrac{X}{\varphi(X)},\qquad \psi'(X) = \dfrac{\varphi(X)-X\varphi'(X)}{\varphi^2(X)}.  \]
We have seen already that when $\psi'>0$ on $[0,p^+]$ (inhinitory case or weakly-excitatory case) there is only convergence to the stationary state.

However, in the strongly excitatory case, which is when $\psi'$ changes sign, explicit examples of periodic solutions are described in \cite{CPST}. The main idea is to use the reconstruction theorem \ref{thm:T-E_extension} to obtain a periodic solution to \eqref{eq:T-E} from a well-chosen function $N$ on $[0,\sigma]$.

The first kind of solution \cite{CPST} finds is piecewise constant $\sigma$-periodic solutions.

\begin{theorem}[C\'aceres, Perthame, Salort, Torres \cite{CPST}]\label{thm:T-E_delay_per1}
    Assume $\mu_X = \delta_{\tau=0}$, \textit{i.e.} $X(t)=N(t)$, $p(s,X) = \phi(X)\mathds 1_{s > \sigma}$, where $\phi\in C^1(\R_+)$, $\sigma>0$ and $0 < p^-\leqslant \phi\leqslant p^+$. Assume $\psi'$ changes sign on $[0,p^+]$. Let $N_1,N_2\in(0,p^+]$ such that $N_1\neq N_2$ and $\psi(N_1)=\psi(N_2)$. If there exists $\beta\in(0,\sigma)$ such that
    \begin{equation}\label{eq:T-E_integral_constant}  \beta N_1 + (\sigma-\beta) N_2 + \psi(N_1) = 1,  \end{equation}
    then there exists a $\sigma$-periodic solution to \eqref{eq:T-E} with firing rate defined on $[0,\sigma)$ by
    \begin{equation}
        N(t) = \left\{
        \begin{matrix}
        N_1,&\qquad t\in[0,\beta),\\
        N_2,&\qquad t\in[\beta,\sigma).
        \end{matrix}
        \right.
    \end{equation}
\end{theorem}

Note that the condition \eqref{eq:T-E_integral_constant} is equivalent to satisfying \eqref{eq:T-E_integral} for this piece-wise-constant stationary state. When such a $\beta$ exists, there is necessarily at least one stationary state \[N_\infty\in[\min(N_1,N_2),\max(N_1,N_2)].\] A similar construction can be made with more than one jump in $[0,\sigma)$, provided there are values $N_1,\dots,N_k$ such that $\psi(N_1) =\dots = \psi(N_k)$. In the special case where there exists an interval $[a,b]\subset[0,p^+]$ where $\phi$ is linear, \textit{i.e.} $\phi(X) = CX$ on $[a,b]$, and if
\[  a\sigma < 1-\frac1C < b\sigma, \]
then the same method allows to construct a smoother $\sigma$-periodic solution with firing rate
\[  N(t) = \dfrac1\sigma \left(1-\frac1C\right) + \beta u(t). \]
where $u(t)$ is any $\sigma$-periodic function such that $\int_0^\sigma u(t)\diff t = 0$ and the constant $\beta>0$ is chosen small enough depending upon $u$ to keep $N(t)$  in $[a,b]$.\\

A second, more involved construction is proposed in \cite{CPST}, which allows to obtain piece-wise smooth $2\sigma$-periodic solutions to \eqref{eq:T-E}. Similar to what we have described for the delayed NNLIF model in Section \ref{sec:NNLIF_delay}, it is possible to use arguments from \cite{HT} in the study of the delay equation 
\begin{equation}\label{eq:T-E_delay2}
    \dfrac{\diff \psi\circ N }{\diff t}(t) = - N(t) + N(t-\sigma)
\end{equation}
that we find in Lemma \ref{lm:TE_delay}. Assume that $\psi$ is strictly convex around a local minimum $\underline N$. Then we can choose $N^-<\underline N < N^+$ such that $\psi(N^-)=\psi(N^+)$, $\psi$ decreasing on $[N^-,\underline N]$ and increasing on $[\underline N, N^+]$.
Then the idea of \cite{CPST} is to seek a periodic solution of \eqref{eq:T-E_delay2} satisfying
\begin{itemize}
    \item $N\in(N^-,\underline N)$ and $N'<0$ on $(0,\sigma)$ ; $N\in(\underline N, N^+)$ and $N'<0$ on $(\sigma,2\sigma)$;
    \item $N(0) = N^+$, $N((2\sigma)^-)= N^-$;
    \item $\psi(N(\sigma^-)) = \psi(N(\sigma^+))$.
\end{itemize}
Then, the solution on $[\sigma,2\sigma)$ is given by the solution on $[0,\sigma)$ as an initial condition in the DDE \eqref{eq:T-E_delay2}. Such a procedure can be reduced to finding a fixed-point of a well-chosen functional $T:\mathcal C\to \mathcal  C$ defined over the Banach space
\[  \mathcal  C = \{ N \in C[0,\sigma] | N(0) = N^+,\ N\ \mathrm{non-increasing},\ \psi\circ N \geqslant \psi(N^+) - (N^+-N^-)\sigma \}. \]
Given $N\in \mathcal C$, the corresponding extension on $[\sigma,2\sigma)$ with the DDE \eqref{eq:T-E_delay2} is chosen as a solution $M$ to the backward ODE problem
\[  \dfrac{\diff \psi\circ M }{\diff t}(t) = - M(t) + N(t-\sigma),\qquad M(2\sigma) = N^-.  \]
Then, this defines a new portion $L$ of the solution on $[2\sigma, 3\sigma]$ as a solution to the forward ODE
\[  \dfrac{\diff \psi\circ L }{\diff t}(t) = - L(t) + M(t-\sigma),\qquad L(2\sigma) = N^+.  \]
The natural choice to get a periodic solution of \eqref{eq:T-E_delay2} as a fixed point of $T$ is then to define
\[  T[N](t) = L(t+2\sigma).  \]
When $\sigma$ is small enough, it can be proved that $T$ is a contraction on $\mathcal C$. Provided additional technical assumptions, a solution to the full system \eqref{eq:T-E} can be reconstructed from this periodic solution of \eqref{eq:T-E_delay2} \textit{via} Theorem \ref{thm:T-E_extension}.

\begin{theorem}[C\'aceres, Perthame, Salort, Torres \cite{CPST}]\label{thm:T-E_delay_per2}
    Assume $\mu_X = \delta_{\tau=0}$, \textit{i.e.} $X(t)=N(t)$, $p(s,X) = \phi(X)\mathds 1_{s > \sigma}$, where $\phi\in C^1(\R_+)$, $\sigma>0$ and $0 < p^-\leqslant \phi\leqslant p^+$. Assume $\psi'$ changes sign on $[0,p^+]$ and there exist $\underline N,\varepsilon>0$ such that $\psi$ is strictly convex on $(\underline N-\varepsilon,\underline N+\varepsilon)$ and $\underline N$ is a local minimum of $\psi$. Then,
    \begin{itemize}
        \item If $\sigma$ is small enough, there exists a $2\sigma$-periodic solution $N$ of \eqref{eq:T-E_delay2} such that $\psi\circ N\in W^{1,\infty}$, $N$ is decreasing on $(0,\sigma)$ and on $(\sigma,2\sigma)$ with a discontinuity at $t=\sigma$.
        \item If in addition $\psi(\underline N) < 1$ and $\psi(\underline N \pm \varepsilon) > 1$, then there exists a solution of \eqref{eq:T-E} with firing rate $N$.
    \end{itemize}
\end{theorem}

\begin{figure}[ht!]
	\begin{center}
		\includegraphics[width=330pt]{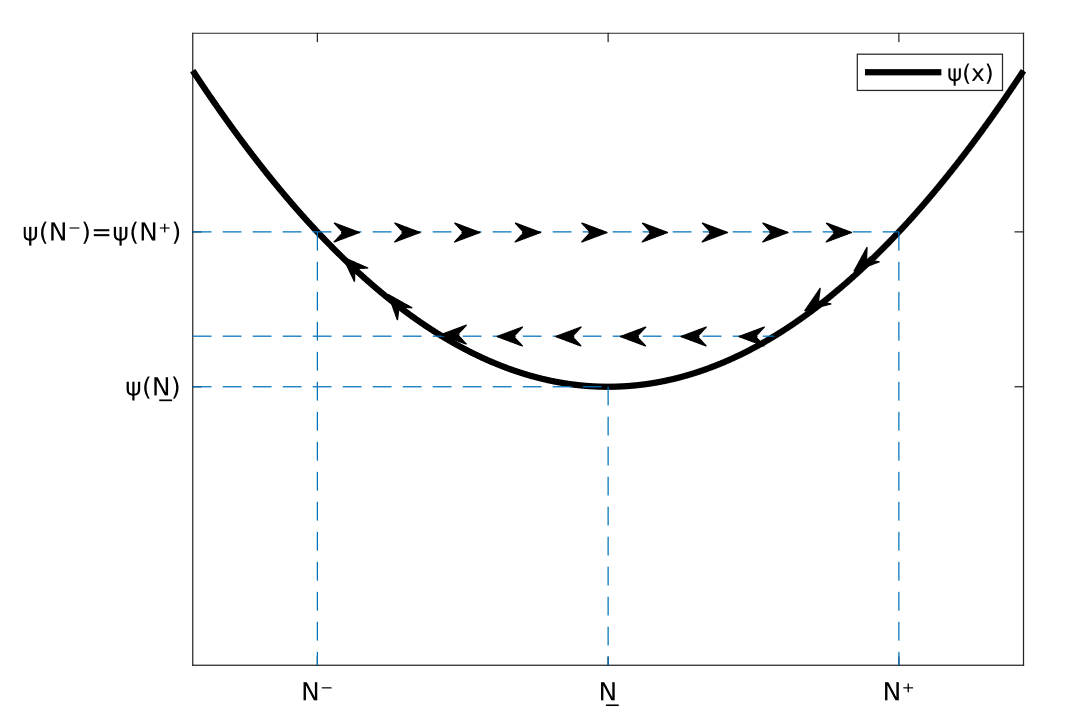}
	\end{center}
	\caption{Construction of a periodic firing rate yielding a solution to \eqref{eq:T-E} from a local minimum of $\psi$ as in Theorem \ref{thm:T-E_delay_per2}. The arrows represent the path of N(t) over a cycle of length $2\sigma$. The horizontal lines are the jump discontinuities at $t = \sigma$ and $t=2\sigma$; note that they preserve the continuity of $\psi(N)$. C\'aceres, Perthame, Salort, Torres \cite{CPST}.\label{fig:T-E_periodic2}}	
\end{figure}

Note that given the fixed point operator is defined through a succession of ODEs, depending on the regularity of $\psi$, more regularity can be proved on the solution $N$ on $(0,\sigma)$ and $(\sigma, 2\sigma)$, meaning that this solution can be in fact smooth by part around the discontinuity point $t=\sigma$. The same procedure allows to construct solutions from a local maximum of $\psi$, but this time the solution will be increasing by part.

Numerical simulations in \cite{CPST} confirm that different types of periodic solutions can be attractive. The solution to \eqref{eq:T-E} converges to a periodic solution with constant by part firing rate for the choice
\[ \varphi(X) = \max(\min(1.6 X,1),0.25),\qquad \sigma = 1,\qquad n^0(s)=e^{-s}. \]
With the choice
\[ \varphi(X) = \frac{10X^2}{1+X^2},\qquad \sigma = 1,\qquad n^0(s)=e^{-(s-1)}\mathds 1_{s > 1}. \]
the solution converges towards a periodic solution with a jump discontinuity and two non-constant paths in each period.

Note that only periodic solutions of period a multiple of $\sigma$ have been constructed so far. Given that multiple periodic solutions are possible for a same choice of $\psi$ and $\sigma$, \cite{CPST} conjectures that the stable solutions are the one with few jump discontinuities. The speed of convergence towards the stable periodic solutions is unknown.

\subsubsection{Further numerical evidence with theoretical guarantees on the scheme}

In the article \cite{STV}, numerical schemes are carefully constructed for the model \eqref{eq:T-E} in both the delayed and non-delayed cases. They propose first order upwind schemes based on a finite volume framework. These schemes are proved to be stable, consistent and convergent for small enough compactly supported initial conditions in the space $\mathrm{BV}(\R_+)$ and $W^{1,+\infty}$ firing functions $p$.

The improved and theoretically validated numerical simulations in \cite{STV} confirm the previous observations of stable periodic solutions \cite{PPS2,CPST} and exhibit some other cases. For example, periodic solutions are obtained in the delayed strongly inhibitory case in the framework of \cite{CPST} with singular delay $d=\frac12$, with
\[ p(s,X(t)) = \varphi(X)\mathds 1_{s > \sigma},\quad X(t)=N\left(t-\frac12\right),\quad \phi(N) = e^{-9N},\quad \sigma=\frac12,\quad n^0(s) = \frac12 e^{-(s-1)^+}.  \]
These oscillations are also observed in delayed settings.
Moreover, new delayed excitatory oscillations are shown in \cite{STV} with detailed information about the numerical density $n(s,t)$.

\subsection{Rigorous derivation from microscopic models}
\label{sec:TE_microscopic}

In \cite{CCDR}, a rigorous connection is made between point processes used by neuroscientists to fit spike data (Poisson, Wold, Hawkes) and the macroscopic Time Elapsed model for neural networks.

Consider a locally finite point process $(T_z)_{z\in\Z}$ on $(\R,\mathcal B(\R))$, \textit{i.e.} a random set of points such that there is almost surely a finite number of points in any bounded interval. For a measurable set $A\in\mathcal B(\R)$, $N_A\in\N\cup\{+\infty\}$ denotes the number of points in $A$. Then, $t\mapsto N_t := N_{(0,t]}$, defined on $\R_+$, is the counting process associated to $N$. The natural filtration of $N$ is the family $(\mathcal F_t)_{t\in\R}$ of sigma-algrebras $\mathcal F_t = \sigma(N\cap(-\infty,t])$ and the predictable filtration is $(\mathcal F_{t^-}^N)_{t\in\R}$ defined by $\mathcal F_{t^-}^N = \sigma(N\cap(-\infty,t))$. The former contains all the information given by the point process at times $t$ and the latter only the information strictly before $t$, thus does not contain information about whether or not $t$ is a point. In most cases, a stochastic intensity $\lambda(t,\mathcal F_{t^-}^N)$ can be associated to the point process $(N_t)_{t\geqslant 0}$. Then, $\lambda(t,\mathcal F_{t^-}^N)\diff t$ represents the probability to have a new point in $[t,t+\diff t)$ given the past points.

Those point processes here correspond to spike times. Then, $T_{N_t}$ is the last discharge for the neuron that occurred at or before time $t$. The time since the last discharge is then $S_t = t-T_{N_t}\geqslant 0$, and for all $t\geqslant 0$, $S_t$ is $\mathcal F_t^N$ measurable. However, as the conditional intensity cannot use knowledge exactly at $t$, and thus cannot use $S_t$, we also need the variable $S_{t^-} = t-T_{N_{t^-}}$ which is $\mathcal F_{t^-}^N$ measurable; $T_{N_{t^-}}$ is the time of the last spike emitted strictly before $t$. Given this framework, the firing function $p(s,X(t))$ and the stochastic intensity $\lambda(t,\mathcal F_{t^-}^N)$ play the same role: they give the firing rate of neurons with time $s$ elapsed since the last discharge given the full information about strictly past activity.
The article \cite{CCDR} then proposes three different point processes to model the time since the last discharge.
\begin{itemize}
    \item Poisson processes: if the Poisson process is homogeneous, which is $\lambda(t,\mathcal F_{t^-}^N)=\lambda$ independent of $t$ and of the past activity, it corresponds to $p(s,X(t))= \lambda$ as
    \[  \mathbb P(S_{t^-} > s) = \mathbb P(N_{[t-s,t)}=0) = e^{-\lambda s},  \]
    which corresponds to the stationary distribution $n_\infty(s) = \lambda e^{-\lambda s}$ in \eqref{eq:T-E} with $p(s,X(t))= \lambda$. The Poisson process can then be chosen inhomogeneous by letting $\lambda(t)$ depend on $t$. If we want to take into account the refractory period, the intensity has to depend on age $\lambda(t,\mathcal F_{t^-}^N)=f(S_{t^-})$. It is then equivalent to the setting $p(s,X(t)) = \bar p(s)$, which is the linear case of \eqref{eq:T-E}. This renewal process enjoys i.i.d. inter-spike intervals (ISI).
    \item Wold processes: consider $A_t^1$ the time between the last discharge and the penultimate discharge. Then $A_t^1= T_{N_{t^-}}-T_{N_{t^-}-1}$ is a Wold process characterised by $\lambda(t,\mathcal F_{t^-}^N)=f(S_{t^-},A_t^1)$ \cite{L2}. This would correspond to $p(s,X(t)) = p(s,a)$ where a is the time elapsed since the penultimate discharge. This extension of \eqref{eq:T-E} has later been investigated at the PDE level in \cite{PST} (see also Subsection \ref{sec:TE_double} below). One could define in the same way the intervals $A_t^k$ between the last $k^{th}$ and $(k+1)^{th}$ discharges, which was also to be investigated at the PDE level in \cite{DPQSZ} (see Subsection \ref{sec:TE_multiple}
    \item Hawkes processes: let us consider now a univariate Hawkes process
    \[ \lambda(t,\mathcal F_{t^-}^N) = \mu + \int_{-\infty}^{t^-} h(t-\tau)N(\diff \tau), \]
    with $\mu$ the spontaneous rate and $h$ an interaction kernel supported in $\R_+$ such that $\int_0^{+\infty} h(\tau)\diff \tau < 1$. Under mild assumptions the second part can be rewritten as $\int_{0}^{t^-} h(t-\tau)N(\diff \tau)$ and represents the influence of the network activity $X(t)$ in $p(s,X(t))$. Although we expect that $p(s,X(t)) = \mu + X(t)$ in this framework, \cite{CCDR} proves that it does not hold and provides more specific rigorous information about the link between these Hawkes processes and the model \eqref{eq:T-E}.
\end{itemize}
Rigorous links between these point processes and the model \eqref{eq:T-E} are proved in \cite{CCDR} after technical considerations beyond the scope of this review. Their main tool is Ogata's thinning algorithm \cite{LS3,O}.

In \cite{C}, the approach of \cite{CCDR} is generalised and a rigorous link between generalised Hawkes processes and the Time Elapsed model is proved.

Another microscopic approach has been proposed for the Time Elapsed model: in \cite{Q} a jump process interacting through a set of random global activity variables is proved to enjoy propagation of chaos and to converge toward \eqref{eq:T-E}. The convergence rate is quantified by use of a coupling method.

\subsection{Variants of the Time Elapsed model}

\subsubsection{Inhomogeneous distribution of the refractory period}

In many cases, the refractory period can be divided into an absolute refractory period (physiological impossibility to fire) and a relative refractory period (the firing threshold is increased), ranging on average from 0.5 to 1 ms for the former and around 10ms for the later. Different types of nerve cells can have different absolute and relative refractory periods \cite{GC}.

A possible extension of the Time Elapsed neuron model is then to add an other structuring variable $r\in\R$, and to choose a firing function $p(r,s,X)$ where for each $r\in\R$ the refractory period is different. The article \cite{KPS} proposed the model
\begin{equation}\label{eq:T-E_var_inhomogeneous}
	\left\{  \begin{array}{l}
		\displaystyle \dfrac{\partial n}{\partial t}(r,s,t)   + \dfrac{\partial n}{\partial s}(r,s,t)   +  p\big(r,s,X(t)\big)\, n(r,s,t) = 0,\\[0.3cm]
		\displaystyle N(r,t) = \int_{0}^{+\infty} p(r,s,X(t))\, n(r,s,t) \diff s, \quad X(t) = J \int_{-\infty}^{+\infty} N(r,t)\diff r,\ \ \\[0.3cm]
		\displaystyle n(r,0,t) = N(r,t), \qquad n(r,s,0) = n^0(r,s)\geqslant 0,\qquad  \int_{-\infty}^{+\infty}\int_0^{+\infty} n^0(r,s) \diff s = 1.
	\end{array}     \right.
\end{equation}
For the firing function, they make the choice
\begin{equation}\label{eq:T-E_var_inh_1}
    p(r,s,X) = \left\{
    \begin{matrix}
        a(r)\qquad& s \leqslant \sigma(r,X),  \\
        b(r)\qquad& s > \sigma(r,X),
    \end{matrix}
    \right.\quad\qquad a(r)<b(r),
\end{equation}
with 
\begin{equation}\label{eq:T-E_var_inh_2}
    0<a_0=\inf_{r\in\R} a(r),\quad b_0 = \sup_{r\in\R} b(r) < +\infty, \quad \sigma\in L^\infty(\R,C^1_b(\R_+)),\quad \partial_X \sigma \leqslant 0.
\end{equation}
The quantity $\sigma(r,X)$ represents the refractory period for neurons of type $r$ under global network activity $X$. The firing rate during the refractory period is taken to be a small value $a(r)$ and not 0, in part to account for stochastic variability in firing and in part for purely mathematical reasons.

From the initial condition $n^0$, the distribution in $r$ is
$\rho(r) = \int_0^{+\infty} n^0(r,s)\diff s$; it is preserved over time by the equation. Let us then assume 
\begin{equation}\label{eq:T-E_var_inh_4}
    0 \leqslant n^0\leqslant b(r) \rho(r).
\end{equation}
The article \cite{KPS} claims that standard arguments and the previous work \cite{PPS} allow to construct a distributional solution $n\in C(\R_+,L^1(\R))$ such that
    $0\leqslant n(r,s,t) \leqslant b(r) \rho(r)$.
Stationary states are then of the form
\begin{equation}\label{eq:T-E_var_inh_6}
    n_\infty = \dfrac{\rho(r)}{ \dfrac{1}{a(r)} - \left( \dfrac{1}{a(r)}-\dfrac{1}{b(r)} \right)e^{-a(r) \sigma(r,X_\infty))}  }.
\end{equation}
There exists at least one stationary network activity $X_\infty$ defining a stationary state $(n_\infty,N_\infty,X_\infty)$, and \cite{KPS} proves that if the average connectivity strength $J$ is small enough, the stationary network activity, and thus the stationary state, is unique.

Then, similar to the standard model \eqref{eq:T-E}, it is possible to prove weakly-nonlinear stability of the (then unique) stationary state.

\begin{theorem}[Kang, Perthame, Salort \cite{KPS}]
    Assume the conditions \eqref{eq:T-E_var_inh_1}, \eqref{eq:T-E_var_inh_2} and \eqref{eq:T-E_var_inh_4}. If there exists $\beta>0$ such that
    \begin{equation}\label{eq:T-E_var_inh_7} \forall r\in\R,\qquad  a(r) - \dfrac{2J b(r) \norme{ \int_{-\infty}^{+\infty} \rho(r)b(r)(b(r)-a(r))\partial_X\sigma(r,\cdot)\diff r}_\infty}{1-J\norme{ \int_{-\infty}^{+\infty}b(r)(b(r)-a(r))\partial_X\sigma(r,\cdot)\diff r}_\infty} \geqslant \beta > 0, \end{equation}
    then for all $t\geqslant 0$,
    \[ \int_{-\infty}^{+\infty}\int_0^{+\infty} |n(r,s,t) - n_\infty(r,s)|\diff s\diff r \leqslant e^{-\beta}\int_{-\infty}^{+\infty}\int_0^{+\infty} |n^0(r,s) - n_\infty(r,s)|\diff s\diff r. \]
\end{theorem}
Note that \eqref{eq:T-E_var_inh_7} will be satisfied for $J$ small enough given \eqref{eq:T-E_var_inh_2}.

\subsubsection{Adaptation and fatigue of neurons modelled as a fragmentation kernel}

An important aspect of neuronal dynamics which is overlooked in the models \eqref{eq:T-E} and \eqref{eq:T-E_var_inhomogeneous} is that the refractory period of neurons can depend on their past activity. Following previous literature on fragmentation-type equations \cite{DG3,PR,M3,CLDDMP,P}, the article \cite{PPS3} proposes to incorporate a fragmentation kernel into the Time Elapsed model.
\begin{equation}\label{eq:T-E_var_AF_kernel}
	\left\{  \begin{array}{l}
		\displaystyle \dfrac{\partial n}{\partial t}(s,t)   + \dfrac{\partial n}{\partial s}(s,t)   +  p\big(s,N(t)\big)\, n(s,t) = \int_0^{+\infty} K(s,a) p(a, N(t))n(a,t) \diff a,\\[0.3cm]
		\displaystyle N(t) = \int_{0}^{+\infty} p(s,N(t))\, n(s,t) \diff s, \\[0.3cm]
		\displaystyle n(0,t) = 0, \qquad n(s,0) = n^0(s)\geqslant 0,\qquad  \int_0^{+\infty} n^0(s) \diff s = 1.
	\end{array}     \right.
\end{equation}
The kernel $K(s,a)$, written as a function by abuse of notation, can be any positive measure with total mass 1 on $\R_+\times\R_+$. Its values represent the probability density to reset to the state $s$ a neuron which fired while being in a state $a$. The variable $s$ then becomes a more general state variable not necessarily related to actual time elapsed since the last discharge. Note that the $a$-independent choice $K = \delta_{s = 0}$ corresponds\footnote{In such case, the condition $n(0,t)=0$ has to be overlooked or replaced by $n(0,t)=N(t)$.} to system \eqref{eq:T-E} with $X(t) = N(t)$. An other simple choice would be $K = \delta_{ s = \psi(a) }$,
which means that upon firing a neuron is deterministically reset to a state determined as a function $\psi(a)$ of its previous state $a$. Taking $(s,a)\mapsto K(s,a)$ as a smooth function instead adds randomness to which state neurons are reset to.

For the function $p$, \cite{PPS3} chooses a regularisation of Assumption \ref{as:TE_p}.

\begin{hyp}[Structural assumptions on the firing function]\label{as:TE_var1}
    Assume there exists $p_M>0$ such that for all $N$,
    \begin{equation*}
        0\leqslant p(s,N)\leqslant p_M,\quad \norme{\dfrac{\partial p}{\partial s}}_\infty + \norme{\dfrac{\partial^2 p}{\partial s^2}}_\infty + \norme{\dfrac{\partial^2 p}{\partial s\partial N}}_\infty + \int_0^{+\infty} \left| \dfrac{\partial p}{\partial s}(s,N) \right|\diff s  < +\infty,
    \end{equation*}
    and a bounded measurable function $\sigma:(0,+\infty)\to(0,\sigma^+]$ such that $p(s,N) \equiv p_M$ for $s\geqslant \sigma(N)$ and $p(s,N)<p_M$ for $s<\sigma(N)$.
\end{hyp}

Regarding the fragmentation kernel $K$, the natural properties to enforce are to not reset neurons at a larger time elapsed than just before they fire, and to preserve the mass of the equation (redistribute the firing neurons without creation nor annihilation). Also, following \cite{LP}, an important assumption is that neurons firing a short time after the previous discharge suffer from fatigue and are reset at lower value $s$, but neurons which didn't fire since a long time are less fatigued and are reset at larger values of $s$, thus enjoying a shorter refractory period. A last assumption is that the change of state is significant when neurons fire. All these assumptions are summarised in mathematical form in Assumption \ref{as:TE_var2} hereafter.
\begin{hyp}[Structural assumptions on the fragmentation kernel]\label{as:TE_var2}
    Assume that the measure $K$ satisfies, by abuse of notation
    \begin{equation}
        K(s,a) \geqslant 0,\qquad K(s,a) = 0,\ s > a, \qquad \int_0^a K(s,a)\diff s = 1.
    \end{equation}
    and, there exists $\theta \in [0,1)$ such that
    \begin{equation}
        -\dfrac{\partial}{\partial a} \int_0^s K(w,a) \diff w \geqslant 0,\qquad \dfrac{\partial}{\partial a} \int_0^s s K(s,a) \diff s \leqslant \theta < 1.
    \end{equation}
\end{hyp}

As we have seen for the standard model \eqref{eq:T-E}, it is not to be taken for granted that the firing rate $N(t)$ is uniquely defined and multiple solutions for $N$ can be the driving force for the arising of periodic solutions. Smallness assumptions can however enforce uniqueness both for the time-dependent firing rate $N(t)$ and the stationary firing rate $N_\infty$. Since \eqref{eq:T-E_var_AF_kernel} does not contain a connectivity parameter, a weak-nonlinearity assumption has to be expressed as a condition on the variations of $p$. If there exists $\eta$ small enough such that
\begin{equation}\label{eq:TE_var_as1}
     \norme{\dfrac{\partial p}{\partial N}}_\infty < \eta, 
\end{equation}
then the firing rate $N(t)$ is uniquely defined and \cite{PPS3} proves that there exists a unique stationary state $(n_\infty,N_\infty)$. We can then define the functions and quantities $p_\infty(s) = p(s,N_\infty)$, $\sigma_\infty=\sigma(N_\infty)$, and
\begin{equation}
    \mathcal B_\infty = e^{\int_0^{\sigma_\infty} p_\infty(s) \diff s } \left[ \sigma_\infty \int_0^{\sigma_\infty} \left| p_\infty'(s) \right|\diff s + \theta \int_0^{\sigma_\infty} p_\infty(s)\diff s \right] > 0.
\end{equation}
Consider then the following assumption on the stationary state
\begin{equation}\label{eq:TE_var_as2}
    \mathcal B_\infty < 1 , \quad -\dfrac{e^{\int_0^{\sigma_\infty} p_\infty(s) \diff s }}{1-\mathcal B_\infty} \int_0^{\sigma_\infty}\dfrac{\partial}{\partial a}\int_0^s K(w,a) \diff w\diff a  - \int_{\sigma_\infty}^{+\infty}\dfrac{\partial}{\partial a} \int_0^s K(w,a) \diff w\diff a < p_M,
\end{equation}
which has to be understood as a smallness assumption on $\sigma$, and thus a smallness assumption on $p$.
Following the approach of \cite{LP}, \cite{PPS3} introduces the quantities
\begin{equation*}
    m(s,t) = n(s,t) - n_\infty(s),\quad M(s,t) = \int_0^s m(w,t)\diff w, \quad J(s,t) = \dfrac{\partial M}{\partial t}(s,t).
\end{equation*}
and applies a similar method to prove convergence in the linear case.
\begin{theorem}[Pakdaman, Perthame, Salort \cite{PPS3}]
    Grant Assumptions \ref{as:TE_var1} and \ref{as:TE_var2}. Assume \eqref{eq:TE_var_as2}, $p(s,N)=p(s)$ and $\sigma(N)=\sigma$ are independent of $N$ and there exists $\lambda > 0$
    \begin{equation*}
         \sup_{u\leqslant\sigma} p(u) \left( \int_0^u K(s,u) e^{\lambda(u-s)}\diff s - 1 \right) < \lambda
    \end{equation*}
    There exists $b>0$ and a function $P(s)>b$ such that if
    \begin{equation*}
        \int_0^{+\infty} P(s)|M(s,0)| \diff s + \int_0^{+\infty}P(s)|J(s,0)| \diff s < +\infty,
    \end{equation*}
    there exists $C > 0$ and $\nu > 0$ such that
    \begin{equation*}
        \int_0^{+\infty} P(s)|n(s,t)-n_\infty(s)|\diff s \leqslant C e^{-\nu t} \int_0^{+\infty}P(s)\Big(|M(s,0)|+|J(s,0)|+\mathds 1_{s\leqslant \sigma}|m(s,0)|\Big)\diff s.
    \end{equation*}
\end{theorem}
This linear result then allows to prove a weakly-nonlinear result.
\begin{theorem}[Pakdaman, Perthame, Salort \cite{PPS3}]
    Grant Assumptions \ref{as:TE_var1} and \ref{as:TE_var2}. Assume \eqref{eq:TE_var_as1} with $\eta$ small enough and \eqref{eq:TE_var_as2}.
    There exists $b>0$ and a function $P(s)>b$ such that if
    \begin{equation}
        \int_0^{+\infty} P(s)|M(s,0)| \diff s + \int_0^{+\infty}P(s)|J(s,0)| \diff s < +\infty,
    \end{equation}
    there exists $C > 0$ and $\nu > 0$ such that
    \begin{multline*}
        \int_0^{+\infty} |n(s,t)-n_\infty(s)|\diff s + \int_0^{+\infty} P(s)|M(s,t)|\diff s + |N(t) - N_\infty| + |N'(t)|\\ \leqslant C e^{-\nu t}\left( 1 + \int_0^{+\infty}P(s)|M(s,0)|\diff s \right).
    \end{multline*}
\end{theorem}
The method of \cite{CY} presented in Subsections \ref{sec:TE_meas_sol} and \ref{sec:TE_Doeblin} for the standard model \eqref{eq:T-E} also allows to prove existence and weakly-nonlinear stability of the stationary state in the space of probability measures. Besides Assumptions \ref{as:TE_p_meas} and \ref{as:TE_p_meas_2} already introduced, we need, similar to Assumption \ref{as:TE_var2}, to ensure that the mass is preserved over time, that the kernel is a positive measure and that neurons are reset in a physical way. We now require a stronger condition, which is that for any state $u$ before firing, a least some neurons are reset in a fixed neighbourhood of $s=0$. We sum up these requirements as 
\begin{hyp}[Fragmentation kernel for measure solutions] \label{as:TE_var3}
For all $u\geqslant 0$, $K(\cdot,u)$ is a probability measure supported on $[0,u]$. There exists $\varepsilon>0$, $\sigma>0$ and $0 < \delta < \sigma$ such that
\[ \forall u\geqslant \sigma,\qquad  K(\cdot,u)\geqslant \varepsilon \mathds 1_{[0,\delta]}. \]
\end{hyp}

Using this assumption, \cite{CY} extends the existence of measure solution to \eqref{eq:T-E_var_AF_kernel} with a similar fixed-point method.

\begin{theorem}[Ca\~nizo, Yoldaş \cite{CY}]
    Let $n^0\in \mathcal P(\R_+)$. Grant Assumptions \ref{as:TE_p_meas} and \ref{as:TE_var3} with 
    \begin{equation} \label{eq:T-E_L_cond_22}L\leqslant \frac{1}{4\norme{n^0}_{\mathrm{TV}}}. \end{equation}
    There exists a unique global-in-time measure solution $(n,N)$ to \eqref{eq:T-E_var_AF_kernel}. If $n^{0,1},n^{0,2}\in \mathcal M_+(\R_+)$ satisfying \eqref{eq:T-E_L_cond_22}, then the associated solutions $n_1$ and $n_2$ satisfy
    \begin{equation}
        \forall t\geqslant 0, \qquad \norme{n_1(t)-n_2(t)}_{\mathrm{TV}} \leqslant e^{4\norme{p}_\infty t} \norme{n^{0,1}-n^{0,2}}_{\mathrm{TV}}.
    \end{equation}
\end{theorem}

Then, the linear and weakly-nonlinear proofs involving Doeblin's theorem can also be extended with similar arguments, leading to
\begin{theorem}[Ca\~nizo, Yoldaş \cite{CY}]
    Grant Assumptions \ref{as:TE_p_meas}, \ref{as:TE_p_meas_2} and \ref{as:TE_var3} with same constant $\sigma$ and $L$ small enough. Then, there exists a unique stationary state $(n_\infty,N_\infty)$ to \eqref{eq:T-E_var_AF_kernel}, there exist constants $C\geqslant 1$, $\lambda > 0$ such that for all probability density $n^0$, if $L\leqslant\frac{1}{4\norme{n^0}_{\mathrm{TV}}}$, then
    \begin{equation}
        \forall t\geqslant 0,\qquad \norme{n(t)-n_\infty}_{\mathrm{TV}} \leqslant Ce^{-\lambda t} \norme{n^0-n_\infty}_{\mathrm{TV}}.
    \end{equation}
\end{theorem}

In \cite{PPS3}, the strongly-nonlinear case is studied numerically with the choice $K(s,u) = \delta_{s= \frac u2}$, where a neuron is reset to half the time elapsed before its discharge. They observe that this fragmentation kernel has a smoothing effects on the periodic solutions constructed and observed for the standard model (which is $K(s,u) = \delta_{s=0}$): there is a unique continuous periodic firing rate attracting all solutions. Periodic solutions with bursting pattern can be observed for some parameters. Note that the choice $K(s,u) = \delta_{s= \frac u2}$ is incompatible with the kernel assumption \ref{as:TE_var3}. To the best of our knowledge, the numerical exploration of periodic solutions in the framework of \cite{CY} has not been performed yet.

\subsubsection{Adaptation and fatigue modelled as an additional variable}

In \cite{FS}, another point of view is proposed for the modelling of adaptation and fatigue of neurons. Instead of resetting the neurons through a fragmentation kernel, they keep the physical meaning of the time elapsed since the last discharge $s$ and add a separate leaky memory variable $m\geqslant 0$. It allows to take into account with more flexibility phenomena like short-term plasticity and spike-frequency adaptation \cite{KTS,TIMGF,TEFDB}.

\begin{figure}[ht!]
	\begin{center}
		\includegraphics[width=300pt]{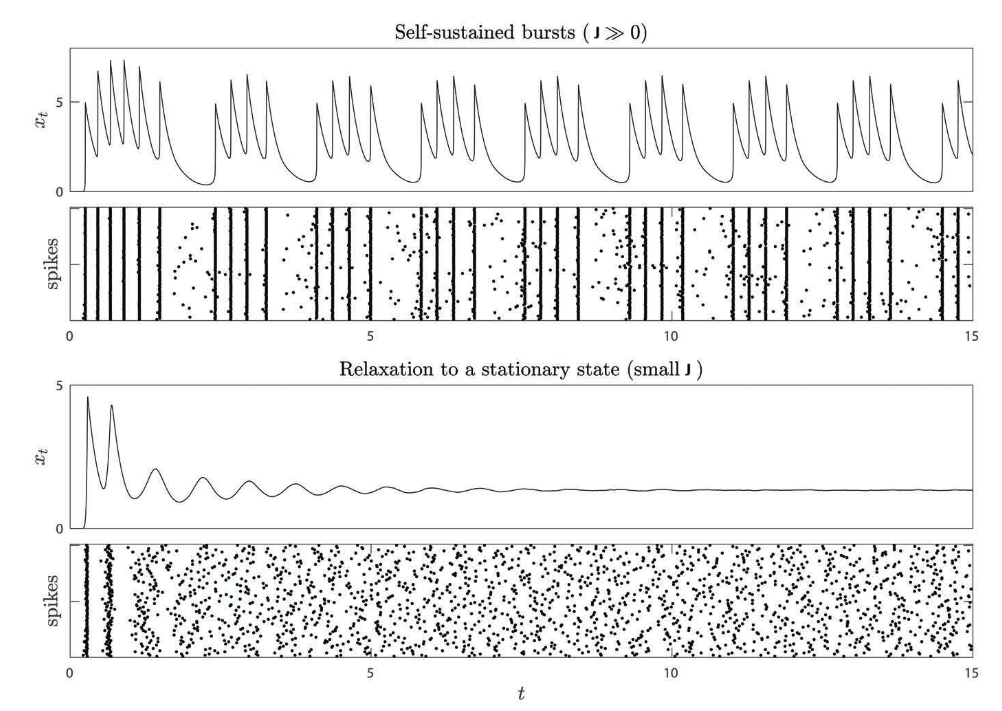}
	\end{center}
	\caption{ Self sustained bursts for large $J$ and relaxation to stationary state for small $J$, seen as the network activity $X(t)$ in model \eqref{eq:T-E_var_FS} or the stochastic spikes in the underlying particle system \cite{S4}. Fonte  S\'anchez, Schmutz \cite{FS}.\label{fig:T-E_periodic3}}	
\end{figure}

In \cite{S4}, the following model is obtained rigorously through a mean-field limit of age and memory dependent Hawkes processes.
\begin{equation}\label{eq:T-E_var_FS}
	\left\{  \begin{array}{l}
		\displaystyle \dfrac{\partial n}{\partial t}   + \dfrac{\partial n}{\partial s} + \dfrac{\partial }{\partial m}\left(-\lambda m\, n \right)  +  p\big(s,m,X(t)\big)\, n = 0,\qquad\\[0.2cm]
		\displaystyle n(0,m,t)  = \mathds 1_{m\geqslant \gamma(0)} \left| (\gamma^{-1})'(m) \right| \int_{0}^{+\infty} p(s,\gamma^{-1}(m),X(t))\, n(s,\gamma^{-1}(m),t) \diff s , \\
        \displaystyle X(t) = J\int_0^t\int_{0}^{+\infty}\int_{0}^{+\infty} h(t-\tau,s,m)p(s,m,X(\tau)) n(s,m,\tau) \diff m\diff s\diff \tau,\\
		\displaystyle n(s,m,0) = n^0(s,m)\geqslant 0,\qquad \int_0^{+\infty}\int_0^{+\infty} n^0(s,m) \diff m\diff s = 1.
	\end{array}     \right.
\end{equation}
As in the previous Time Elapsed models we presented, $p$ is the firing function, $X(t)$ the network activity at time $t$ and $J$ the connectivity strength. The network activity is computing according to a memory kernel $h$ and the history of all the neurons spiking between time 0 and $t$; this variant of the Time Elapsed model is then delayed. 
If $h(\tau,s,m)=h(\tau)$, the network activity simplifies down to the familiar
\[  X(t) = \int_0^t h(t-\tau) \mathcal N(\tau)\diff \tau, \qquad \mathcal N(t) = \int_{0}^{+\infty}\int_{0}^{+\infty} p(s,m,X(\tau)) n(s,m,\tau) \diff m\diff s. \]
However, more general biologically relevant forms are possible for $h$, for example:
\begin{itemize}
    \item Spike frequency adaptation: Take $\gamma(m)=m+\Gamma$ for some constant $\Gamma > 0$, $h(t,s,m) = \tilde h(t)$ and $p(s,m,t)=\tilde p(\eta(s)-m+X)$ with $\tilde p$ non-decreasing and $\eta$ a bounded refractory Kernel such that $\eta(+\infty)=0$. This can represent the fact that neurons spiking at high frequency will undergo spike frequency adaptation and be unable to spike for some time \cite{BH2}. With these parameters, the model exhibits bursting at the population level as can be seen on Figure \ref{fig:T-E_periodic3}.
    \item Short term synaptic depression: when neurons fire a lot in recent past, it can induce short term depression at the level of the synapse \cite{TPM}, it can be modelled by restricting $m\in(0,1)$ and choosing $h(t,s,m)=\tilde h(t)(1-m)$, $\gamma(m)=1-v+vm$ for some fixed $v\in(0,1)$. Neurons with a leaky memory $m$ closed to 1 will have very little impact on the network activity $X(t)$.
\end{itemize}
The resetting of neurons depends on an increasing $C^1$-diffeomorphism $\gamma:\R_+\to(0,+\infty)$ called the jump mapping. A neuron firing with state $(s,m)$ is reset to $(0,\gamma(m))$, and then the boundary condition at $(0,m)$ is written in term of $\gamma^{-1}$.

It can be formally checked and rigorously proved that system \eqref{eq:T-E_var_FS} preserves the mass. Define the weight $\omega(s,m) = m + 1$, consider the weighted space
 $L^1_+(\omega)$
and the set of hypotheses
\begin{hyp}[For well-posedness]\label{as:TE_var_4}
    There exists $L$ such that for all $(s_1,m_1,X_1),(s_2,m_2,X_2)\in \R+\times(0,+\infty)\times\R_+$,
    \[ |p(s_1,m_1,X_1)-p(s_2,m_2,X_2)| \leqslant L( |s_1-s_2|+|m_1+m_2|+|X_1-X_2|). \]
    There exists a positive bounded function $\Gamma$ such that for all $m > 0$, $\gamma(m) = m + \Gamma(m)$.
\end{hyp}
Then, similar to the strategy in Subsection \ref{sec:TE_exist} two applications of the Banach fixed point theorem to an appropriate functional yields existence of solutions in the weighted space $L_+^1(\omega)$.
\begin{theorem}[Fonte S\'anchez, Schmutz \cite{FS}]
    Grant Assumption \ref{as:TE_var_4}. There exists constants $\alpha$ and $b\in\R$ such that for any probability density $n^0\in L^1_+(\omega)$, there exists a unique weak-solution $n\in C(\R_+,L_+^1(\omega))$ to \eqref{eq:T-E_var_FS} and
    \[\forall t\geqslant 0, \quad \norme{n(\cdot,t)}_{L^1_+(\omega)}\leqslant \norme{n^0}_{L^1_+(\omega)} e^{-\alpha t} + \frac{b}{\alpha}(1-e^{-\alpha t}).\]
\end{theorem}
Note that well-posedness is also proved in $L^1$ not weighted in \cite{FS}, but the weighted version is required for the Doeblin-Harris method.

\begin{hyp}[On firing and resetting]\label{as:TE_var_5}
    There exists $\sigma>0$, $p_0>0$ such that if $s>\sigma$, $p(s,\cdot,\cdot)\geqslant p_0$. There exists $C_\gamma\in(0,1]$ such that $C_\gamma\leqslant \gamma'\leqslant 1$. The function $\bar h:(s,m)\mapsto \int_0^{+\infty} h(t,s,m)\diff t$ is bounded.
\end{hyp}

\begin{theorem}[Fonte S\'anchez, Schmutz \cite{FS}]
    Grant Assumptions \ref{as:TE_var_4} and \ref{as:TE_var_5}. There exist at least one stationary state $(n_\infty,X_\infty)$ to \eqref{eq:T-E_var_FS}. If $|J|$ is small enough, the stationary state is unique.
\end{theorem}

Then, similar to \cite{CY} but adding the weighted framework and using Harris' theorem instead of Doeblin theorem, \cite{FS} proves that when $J$ is small enough, there is global convergence towards the unique stationary state.

\begin{theorem}[Fonte S\'anchez, Schmutz, \cite{FS}]
    Grant Assumptions \ref{as:TE_var_4} and \ref{as:TE_var_5}. Assume there exists $C_1,c_2>0$ such that $h(t,\cdot,\cdot)\leqslant C_1e^{-c_2 t}$. For any $A>0$, if $|J|$ is small enough with respect to $A$, then there exists $C\geqslant 1$, $\nu > 0$ such that for all initial condition satisfying $\norme{n^0}_{L^1_+(\omega)}\leqslant A$, the solution to \eqref{eq:T-E_var_FS} satisfies
    \[ \norme{n(\cdot,t)-n_\infty}_{L^1_+(\omega)} +|X(t)-X_\infty|\leqslant C e^{-\nu t}\left( 1 + \norme{n^0-n_\infty}_{L^1_+(\omega)} \right) . \]
\end{theorem}
A similar result can be proved when instead of $\gamma(m)=m+\Gamma(m)$ we assume that $\gamma$ is bounded.

\subsubsection{Time elapsed since the last two discharges}
\label{sec:TE_double}

Another direction for extending the standard Time Elapsed model is to incorporate another variable $a$ representing, for a neuron with time $s$ elapsed since the last discharge, the time elapsed since the penultimate discharge. Such a model is proposed and explored in \cite{PST}. Multiple renewal equations had previously been explored in \cite{FP}. Define
\[ \Omega = \{ (s,a)\in\R^2 \ | \ 0\leqslant s\leqslant a  \}. \]
The multiple time elapsed model reads
\begin{equation}\label{eq:T-E_double}
	\left\{  \begin{array}{l}
		\displaystyle \dfrac{\partial n}{\partial t}(s,a,t)   + \dfrac{\partial n}{\partial s}(s,a,t) + \dfrac{\partial n}{\partial a}(s,a,t)  +  p\big(s,a,X(t)\big)\, n(s,a,t) = 0,\qquad\\[0.2cm]
		\displaystyle N(a,t) = \int_{a}^{+\infty} p(s,a,X(t))\, n(a,u,t) \diff u, \quad X(t) = \int_{0}^{+\infty} N(a,t) \diff a,\\
		\displaystyle n(0,a,t) = N(a,t), \qquad n(s,a,0) = n^0(s,a)\geqslant 0,\qquad \int_0^{+\infty}\int_s^{+\infty} n^0(s,a) \diff a\diff s = 1.
	\end{array}     \right.
\end{equation}

As for the standard Time Elapsed model, we have conservation of mass for weak solutions
\begin{equation}
    \int_0^{+\infty}\int_s^{+\infty} n^0(s,a) \diff a\diff s = \int_0^{+\infty}\int_s^{+\infty} n(s,a,t) \diff a\diff s =  1.
\end{equation}
Note that if $p(s,a,X) = p(s,X)$ is independent of $a$, then the quantity $m(s,t) = \int_s^{+\infty} n(s,a,t)\diff a$ solves \eqref{eq:T-E}. The model \eqref{eq:T-E_double} can hence exhibit all the behaviours of the standard Time Elapsed model, like periodic solutions and the existence of multiple stationary states.

By studying the linear problem and then applying a fixed point argument, \cite{PST} proves the well-posedness of this model in the weakly-nonlinear regime.
\begin{theorem}[Perthame, Salort, Torres, \cite{PST}]
    Assume that $n^0\in L^1(\Omega)$ is a probability density and that $p\in W^{1,\infty}(\Omega\times\R)$. If $\norme{\partial_X p}_\infty < 1 $, there exists a unique solution of \eqref{eq:T-E_double} satisfying
    \[ n \in C(\R_+,L^1(\Omega))\cap L^\infty(\Omega\times\R_+),\quad N\in C(\R_+,L^1(\R_+))\cap L^\infty(\R_+^2),\quad X\in C(\R_+)\cap L^\infty(\R_+).\]
\end{theorem}
Then, for the linear system $p(s,a,X) = p(s,a)$, they prove the following entropy dissipation: for all convex function $H\geqslant 0$ with $H(0)=0$, for all solution $n$ to the linear system,
\begin{equation*}
    \dfrac{\diff}{\diff t}\int_0^{+\infty}\int_s^{+\infty} H\left(\frac{n(s,a,t)}{n_\infty(s,a)} \right) n_\infty(s,a) \diff a \diff s = - D_H[n](t) \leqslant 0,
\end{equation*}
where $(n_\infty,N_\infty,X_\infty)$ is the unique linear stationary state and
\begin{equation*}
    D_H[n](t) = \int_0^{+\infty}\int_s^{+\infty} p(s,a)H\left(\frac{n(s,a,t)}{n_\infty(s,a)} \right) n_\infty(s,a) \diff a \diff s - \int_s^{+\infty} H\left(\frac{N(a,t)}{N_\infty(a)} \right)N_\infty(a)\diff a.
\end{equation*}
However, a requirement for obtaining this entropy dissipation is that the stationary state satisfies $n_\infty>0,N_\infty>0$ which is not true in general (it can be proved to be false for example with $p(s,a,X) = \mathds 1_{s > 1}$). Moreover, it is not easy to prove an appropriate Poincaré-like inequality to bound $D_H[n]$. For these reasons, \cite{PST} turns to the Doeblin theory used in \cite{CY} and described in Subsection \ref{sec:TE_Doeblin}. In the linear case $p(s,a,X) = \bar p(s,a)$, they prove that there exists a unique stationary state $(n_\infty,N_\infty,X_\infty)$ of \eqref{eq:T-E_double} and the solution $n$ satisfies
    \begin{equation}
        \norme{n(t)-n_\infty}_{L^1} \leqslant \dfrac{1}{1-\beta} e^{-\lambda t}\norme{n^0-n_\infty}_{L^1} ,
    \end{equation}
    with $\beta = \frac12 p_0^2 \sigma^2 e^{-3p_\infty\sigma}$ and $\lambda = -\frac{1}{3\sigma}\log(1-\beta)>0$.
In the weak-connectivity regime $\partial_X p$ small, there is also a unique stationary state for the nonlinear problem. The linear result can be extended to the weakly-nonlinear case.

\begin{theorem}[Perthame, Salort, Torres, \cite{PST}]
    Assume that $p$ is smooth, $\norme{\partial_X p}_\infty$ is small enough and that there exists $p_0>0,p_\infty>0,\sigma>0$ such that $p_0\mathds 1_{a>s\geqslant\sigma}\leqslant p\leqslant p_\infty$. Then, there exists a unique stationary state $(n_\infty,N_\infty,X_\infty)$ of \eqref{eq:T-E_double} and constants $C\geqslant 1$, $\lambda > 0$ such that if $n^0\in L^1(\Omega)$ is a probability density, the solution $n$ to \eqref{eq:T-E_double} satisfies
    \begin{equation}
        \forall t\geqslant 0,\qquad \norme{n(t)-n_\infty}_{L^1} \leqslant C e^{-\lambda t}\norme{n^0-n_\infty}_{L^1}.
    \end{equation}
\end{theorem}

Numerical exploration in \cite{PST} shows that \eqref{eq:T-E_double} also displays solutions with jump discontinuities in the strongly-nonlinear regime. It was observed that convergence to a periodic solution in the case $p(s,a,X)=p(s,X)$ is replaced by steady asymptotics when the dependency in $a$ is added. Complex dynamics elicited by the $a$ variable has yet to be observed and studied. One can hope to find a reduction to a 2-delays equation similar to the findings of \cite{CPST} for \eqref{eq:T-E}.

\subsubsection{Convergence towards an infinite renewal model}
\label{sec:TE_multiple}

The idea of taking into account both the time since the last and the penultimate discharges is pushed further in \cite{DPQSZ}: they consider for each neuron $N$ variables $s_1,\dots,s_N$ representing the time elapsed since the last $N$ discharges. Then, the firing function can depend upon more information about the past history of the neurons. With the aim of sending $N$ towards $+\infty$, a particular structure is proposed, which gives less and less weight to later discharges
\begin{equation}\label{eq:T-E_var_infinite_1}
    p_N(s_1,\dots,s_N) = \sum_{i=1}^N \phi_i(s_1,\dots,s_N),\qquad \phi_i > 0,\qquad \sum_{i=1}^{+\infty} \norme{\phi_i}_\infty \leqslant a_+ <  +\infty.
\end{equation}
Note that, for the sake of mathematical tractability, $p_N$ does not depend on the network activity $X(t)$. The model for $N$ time elapsed variables can be written as
\begin{equation}\label{eq:T-E_var_N}
    \left\{
    \begin{array}{l}
         \displaystyle \dfrac{\partial n_N}{\partial t} + \sum_{i=1}^{N} \dfrac{\partial n_N}{\partial s_i} + p_N( s_1, \dots, s_N) n_N= 0.  \\
          \displaystyle n_N(0,s_2,\dots,s_N,t) = \int_0^{+\infty} p_N(s_2, \dots, s_N,u) n_N(s_2,\dots,s_N,u,t) \diff u,\\
          \displaystyle n_N(s_1,\dots,s_N,0) = n^0(s_1,\dots,s_N)\geqslant 0,\qquad \int_{\mathcal C_N} n_N(s_1,\dots,s_N) \diff s_1\dots\diff s_N = 1. 
    \end{array}
    \right.
\end{equation}
where we define the domains
\begin{equation}
    \mathcal C_N = \{ 0\leqslant s_1\leqslant \dots,\leqslant s_N  \},\qquad \mathcal C_\infty = \{0\leqslant s_1\leqslant\dots\leqslant s_N\leqslant\dots\}.
\end{equation}
As we have mentioned in Subsection \ref{sec:TE_microscopic}, there is a connection between \eqref{eq:T-E_var_N} and Wold processes. However, there is currently, to our knowledge, no framework for Wold processes with an infinite number of variables.

In order to improve readability, especially when $N\to+\infty$, let us use $[s]_N$ to denote $(s_1,\dots,s_N)\in \mathcal C_N$ and $[s]_\infty=(s_1,\dots,s_N,\dots)\in \mathcal C_\infty$. Consider also the shift operator $\tau$ defined by
\[  \tau(s_1,s_2,\dots,s_N) = (0,s_1,s_2,\dots,s_{N-1}),\qquad \tau(s_1,s_2,\dots,s_N,\dots) = (0,s_1,s_2,\dots,s_N,\dots). \]
A weak-solution to \eqref{eq:T-E_var_N} can be defined with the shift operator: $n\in C(\R_+,L^1(C_N))$ if for all $\psi\in C_b^1([0,T]\times C_N)$,
\begin{multline}
    \int_0^T \int_{\mathcal C_N} n_N([s]_N,t) \left[ -\dfrac{\partial \psi}{\partial t} - \nabla\cdot \psi + \big(\psi([s]_N,t) - \psi(\tau([s]_N,t) \big)p_N([s]_N) \right] \diff [s]_N \diff t\\
    = - \int_{\mathcal C_N}  \psi(\tau([s]_N,T)n_N([s]_N,T) \diff [s]_N + \int_{\mathcal C_N} \psi(\tau([s]_N,0)n_N([s]_N,0)\diff [s]_N
\end{multline}
As for the model with $N=2$ in \cite{PST} and other Time Elapsed models \cite{CY,FS}, the Doeblin-Harris method allows to prove convergence towards a unique stationary state.
\begin{theorem}[Dou, Perthame, Qi, Salort, Zhou \cite{DPQSZ}]
    Assume $0 < a_- \leqslant p_N \leqslant a_+$, with constants $a_- >0,a_+ > 0$. Then, for all appropriate initial condition, for all $[s]_N\in \mathcal C_N$, $n_N([s]_N,t)\geqslant a_-^N e^{-s_N a_+}$. There exist a unique stationary state $n_{N,\infty}$ and constants $c_N >0,\lambda_N >0$ such that
    \[ \forall t\geqslant 0, \qquad \norme{n_N(t)-n_{N,\infty}}_{L^1(\mathcal C_N)} \leqslant c_n e^{-\lambda_N t}\norme{n_N^0-n_{N,\infty}}_{L^1(\mathcal C_N)}. \]
    Moreover, the following tightness estimate holds for the stationary state:
    \begin{equation}
        \int_{\mathcal C_N} \sigma_N([s]_N) n_{N,\infty} \diff [s]_N \leqslant \frac{2}{a_-},\qquad \sigma_N([s]_N) := \sum_{i=1}^{N} \dfrac{s_i}{2^i}.
    \end{equation}
\end{theorem}
The constants $c_N,\lambda_N$ are explicitly constructed in \cite{DPQSZ}. Unfortunately, they converge towards 0 when $N$ goes to $+\infty$. Thus, it is not possible to use such decay estimates for the limit$N\to+\infty$. Another difficulty is to write a rigorous formulation for the limit of \eqref{eq:T-E_var_N}. Formally, it would be an infinite renewal equation with variables $[s]_\infty\in C_\infty$ in the form
\begin{equation}
    \displaystyle \dfrac{\partial n_\infty}{\partial t} + \sum_{i=1}^{+\infty} \dfrac{\partial n_\infty}{\partial s_i} + p_\infty( [s]_\infty) n_\infty= 0.
\end{equation}
But because there is no "last" variable, it is not easy to define a solution, or simply to prescribe the boundary condition at $s_1=0$. Therefore, \cite{DPQSZ} proposes to study the $K$-marginales of the solution to this limiting infinite renewal equation, similar to BBGKY hierarchies in classical kinetic theory \cite{CIP}.

First, let us write the marginals for the solution $n_N$ to the problem with a finite number $N$ of variables $s_1,\dots,s_N$. For $K\leqslant N$, denote $n_N^{(K)}$ the first $K$ marginals of $n_N$:
\begin{equation}
    n_N^{(K)}([s]_K,t) = \int_{s_K\leqslant s_{K+1}\leqslant\dots\leqslant s_N} n_N([s]_N,t) \diff s_{K+1} \dots\diff s_{N}.
\end{equation}
Note that $n_N^{(K)} = \int_{s_K}^{+\infty} n_N^{(K+1)} \diff s_{K+1}$.
\begin{definition}[Consistent sequence] A sequence of measures $(\mu_N)_{N\geqslant 1}$ is said to be consistent if for all $K < N$, the marginals satisfy $\mu_K = \mu_N^{(K)}$. A measure $\mu_\infty$ on $C_\infty$ is consistent if for all $K\geqslant 1$, $\mu_\infty^{(K)} = \int_{s_K}^{+\infty} \mu_\infty^{(K+1)} \diff s_{K+1}$.
\end{definition}
Given the coupling term
\begin{equation}
   E_N^{(K)}([s]_K,t) = \sum_{i=K+1}^N \int_{s_{K+1}}^{+\infty}\dots \int_{s_N=0}^{+\infty} \phi_i([s]_i) n_N^{(i)}([s]_i,t)\diff [s]_{K+1,i}, 
\end{equation}
the marginals $n_N^{(K)}$ satisfy
\begin{align}
    &\dfrac{\partial n_N^{(K)}}{\partial t} + \sum_{i=1}^{K} \dfrac{\partial n_N^{(K)}}{\partial s_i} + p_K( [s]_K) n_N^{(K)} + E_N^{(K)} = 0,\\\nonumber& n_N^{(K)}(0,[s]_{2,K},t) = \int_0^{+\infty} \Big[p_K([s]_{2,K},u) n_N^{(K)}([s]_{2,K},u,t) + E_N^{(K)}([s]_{2,K},u,t)\Big]\diff u
\end{align}
It leads to the infinitely coupled hierarchy system
\begin{equation}\label{eq:T-E_var_infinite}
    \left\{
    \begin{array}{l}
         \displaystyle \dfrac{\partial n_\infty^{(K)}}{\partial t} + \sum_{i=1}^{K} \dfrac{\partial n_\infty^{(K)}}{\partial s_i} + p_K( [s]_K ) n_\infty^{(K)} + E_\infty^{(K)} = 0.  \\
        \displaystyle n_\infty^{(K)}(0,s_2,\dots,s_K,t) = \int_0^{+\infty} \Big[p_K([s]_{2,K},u) n_\infty^{(K)}([s]_{2,K},u,t)+E_\infty^{(K)}([s]_{2,K},u,t)\Big] \diff u,\\
        n_\infty^{(K)}([s]_K,0) = n_\infty^{0,(K)}([s]_K).
    \end{array}
    \right.
\end{equation}
The coupling system is obtained by passing to the limit in $E_N^{(K)}$, which yields
\begin{equation}
    E_\infty^{(K)}([s]_K,t) = \sum_{i=K+1}^{+\infty} \int_{s_{K+1}=0}^{+\infty}\dots \int_{s_{i}=0}^{+\infty} \phi_i([s]_i) n_\infty^{(i)}([s]_i,t)\diff [s]_{K+1,i}.
\end{equation}
\begin{theorem}[Dou, Perthame, Qi, Salort, Zhou \cite{DPQSZ}]Assume \eqref{eq:T-E_var_infinite_1}. Consider a sequence of initial conditions $n_N^0$ of mass 1 on $\mathcal C_N$ such that there is a consistent family $(n_\infty^{0,(K)})_K$ such that
    \[ \lim_{N\to+\infty}\sup_K \norme{n_N^{0,(K)} - n_\infty^{0,(K)}} = 0. \]
    \begin{itemize}
        \item For all $T > 0$ and $K\in\N$, $(n_N^{(K)})_{N\in\N}$ is a Cauchy sequence in $C([0,T],L^1(\mathcal C_K))$ and thus has a consistent limit $n_\infty^{(K)}\in C([0,T],L^1(\mathcal C_K))$.  
        \item $E_N^{(K)}$ converges towards $E_\infty^{(K)}$ in $C([0,T],L^1(\mathcal C_K))$ as $N\to+\infty$.
        \item $(n_\infty^{(K)})_{K\in\N}$ is the unique consistent weak solution of the hierarchy \eqref{eq:T-E_var_infinite}.
    \end{itemize}
\end{theorem}
Such a solution to the hierarchy \eqref{eq:T-E_var_infinite} is obtained as a local-in-time limit on $[0,T]$ of the models with finite number of variables. Under additional smallness assumptions, \cite{DPQSZ} also proves a uniform-in-time, using the Doeblin method. It allows to approximate the stationary solution of the infinite renewal equation by the stationary state $n_{N,\infty}$ for large $N$:
\[\lim_{N\to+\infty} \norme{n_{N,\infty} - n_{\infty,\infty}^{(N)}}_{L^1(\mathcal C_N)} = 0. \]
Under smallness assumption on the functions $\phi_i$ in \eqref{eq:T-E_var_infinite_1}, they obtain uniqueness of the stationary state $n_{\infty,\infty}$ of the hierarchy.

Last, \cite{DPQSZ} proposes a further interpretation of the limiting infinite renewal equation by introducing solutions as measures on $\mathcal C_\infty$ which are weak limits of the solutions for finite $N$. The time-dependent infinite-dimensional measure $n_\infty(t) \in \mathcal P(\mathcal C_\infty)$ can be built using the Kolmogorov extension theorem. This notion of solution is proved to be equivalent to the hierarchy \eqref{eq:T-E_var_infinite}. It allows in particular to use the Monge-Kantorovich distance to write and prove a long-term convergence result for the infinite system.

\subsubsection{Spatial extension and evolution of the connectivity}

The standard Time Elapsed model does not take into account spatial variability of the connectivity $J$ nor evolution of the connectivity in time. Following previous spatial extensions in Fitzhugh-Nagumo models \cite{C4,C5} (see also Section \ref{sec:FN}) and foreshadowing spatial extension in Fokker-Planck models (see Section \ref{sec:GC}), \cite{ST2} proposes a spatial extension of the Time Elapsed model where neurons at locations $x\in\Omega$ and $y\in\Omega$ interact through the time-dependent connectivity $w(x,y,t)$.
\begin{equation}\label{eq:T-E_var_spatial}
	\left\{  \begin{array}{l}
		\displaystyle \dfrac{\partial n}{\partial t}(x,s,t)   + \dfrac{\partial n}{\partial s}(x,s,t)   +  p\big(s,X(x,t)\big)\, n(x,s,t) = 0,\\[0.2cm]
		\displaystyle N(x,t) = \int_{0}^{+\infty} p(s,X(t))\, n(x,s,t) \diff s,\\
        \displaystyle X(t) = \int_\Omega w(x,y,t) N(y,t)\diff y + I(x,t),\\
        \dfrac{\partial w}{\partial t}(x,y,t) = -w(x,y,t) + \gamma G[N(x,t),N(y,t)] , \qquad w(x,y,0)=w^0(x,y).\\
		\displaystyle n(0,t) = N(t), \qquad n(s,0) = n^0(s)\geqslant 0,\qquad \int_0^{+\infty} n^0(s) \diff s = 1.
	\end{array}     \right.
\end{equation}
Neurons at location $x$ are stimulated by a network activity $X(x,t)$ which is determined as an integral on the spatial domain $\Omega$ of the firing rate $N(y,t)$ with weight given by the connectivity function $w(x,y,t)$, plus an external input $I(x,t)$ which is a given function. The connectivity kernel $w$ evolves according to its own equation, depending on the firing rates at locations $x$ and $y$. The smooth function $G:\R^2\to\R$ is the learning rule for the connectivity, and $\gamma>0$ is the learning rate. We can take for example a representation of Hebbian learning like
$G[N_1,N_2] = N_1 N_2$. The following assumption is made for simplicity
\begin{equation}\label{eq:T-E_spatial_G}
    \norme{G}_\infty + \norme{\nabla G}_\infty \leqslant 1.
\end{equation}
There is no loss of generality as $\gamma$ can be scaled to enforce it.
With a fixed point argument, \cite{ST2} proves that the spatially extended problem is well-posed in the weak-interaction regime for small enough initial condition.
\begin{theorem}[Salort, Torres, \cite{ST2}]
    Assume \eqref{eq:T-E_spatial_G}, $p\in W^{1,+\infty}(\R_+\times\R)$ and that there exists $p_0,p_\infty,\bar\sigma$ such that $p_0\mathds 1_{s > \bar\sigma}\leqslant p \leqslant p_\infty$. For all initial condition $n^0\in C_b(\Omega,L^1(\R_+)),w^0\in C_b(\Omega\times\Omega)$ with $n^0$ of total mass 1, if 
    \[  |\Omega| \norme{\int_0^{+\infty} n^0(\cdot,s)\diff s}_\infty\norme{\dfrac{\partial p}{\partial X}}_\infty \max\left(\norme{w^0}_\infty,\gamma\right) < 1 ,\]
    then the system \eqref{eq:T-E_var_spatial} has a unique solution such that $n\in C_b(\Omega \times \R_+ ,L^1(\R_+))$ is non-negative, $N\in C_b(\Omega\times\R_+)$ and $w\in C_b(\Omega\times\Omega\times\R_+)$.
\end{theorem}
The result can also be proved for a larger class of firing functions $p$. A fixed point argument also allows \cite{ST2} to prove that, given an initial spatial distribution of neurons $g(x) = \int_0^{+\infty} n^0(\cdot,s)\diff s$, there exists a unique continuous stationary state in the weak-connectivity and small learning rate regime (it can be related to the homogeneous brain assumption for the grid cells model in Section \ref{sec:GC}). Note that since the connectivity kernel is evolving in time, a kernel converging towards a small function in $L^\infty$ norm means that the nonlinearity through $X(x,t)$ will be small too, hence the uniqueness of the $n_\infty$ stationary density. 

Following \cite{MMP1,MMP2,P,PPS,ST2} proves the convergence towards the unique stationary state in the weak-connectivity regime via a relative entropy method.
\begin{theorem}[Salort, Torres, \cite{ST2}]
    Assume \eqref{eq:T-E_spatial_G}, $p\in W^{1,+\infty}(\R_+\times\R)$ and that there exists $p_0,p_\infty,\bar\sigma$ such that $p_0\mathds 1_{s > \bar\sigma}\leqslant p \leqslant p_\infty$. For all initial condition $n^0\in C_b(\Omega,L^1(\R_+)),w^0\in C_b(\Omega\times\Omega)$ with $n^0$ of total mass 1, if $\gamma,\norme{\partial_X p}_\infty$ are small enough with respect to $n^0,w^0$, then there exists $C\geqslant 1$, $\nu > 0$ such that
    \[ \forall t\geqslant 0,\quad \norme{n(\cdot,\cdot,t)-n_\infty}_{L^1} + \norme{w(\cdot,\cdot,t)-w_\infty}_{L^1} \leqslant C e^{-\nu t} \left(\norme{n^0-n_\infty}_{L^1} + \norme{w^0-w_\infty}_{L^1}\right),\]
    and $X\to X_\infty$, $N\to N_\infty$ both exponentially fast.
\end{theorem}
Then, \cite{ST2} builds on the Doeblin theory in \cite{CY} and improves the relative entropy result, proving under similar conditions that
\[ \norme{n(\cdot,\cdot,t)-n_\infty}_{L^\infty_x L^1_s} + \norme{w(\cdot,\cdot,t)-w_\infty}_{\infty} \leqslant C e^{-\nu t} \left(\norme{n^0-n_\infty}_{L^\infty_x L^1_s} + \norme{w^0-w_\infty}_{\infty}\right).\]
The article \cite{ST2} also provides theoretical insigths into the effect of large inputs $I(x,t)$ and slow learning dynamics compared to the evolution of network activity. For the later, they consider a slow timescale represented by a parameter $\varepsilon>0$ and they change the first equation of the system with
\[  \varepsilon\partial_t n + \partial_s n + p(s,X(x,t)) n = 0. \]
Then, \cite{ST2} rigorously proves that when $\max(\norme{w^0},\gamma)$ is small enough, the solution $n_\varepsilon$ converges in $L^1$ towards a solution to the limiting problem
\[  \partial_s n + p(s,X(x,t)) n = 0, \]
and that this limiting problem is well-posed in the weakly-nonlinear regime. They also prove asymptotic convergence towards a unique stationary state for the weakly-nonlinear limiting problem.

The theoretical analysis of the strongly-nonlinear regime is yet unexplored for this model, but \cite{ST2} provides numerical simulations showing that a non-constant external input $I(x)$ can lead to asymptotic convergence towards spatial patterns. When $I$ is constant, convergence towards spatially uniform stationary states are observed even for strong nonlinearities. The numerical experiments are performed with a Hebbian learning rule as explicited above.

\subsubsection{Master equation and link to the NNLIF model}

A series of papers \cite{DHT,DHT2,DHT3} proposed a connection between the NNLIF model described throughout Section \ref{sec:NNLIF} and the Time Elapsed model.

First, in the linear case, consider the auxiliary equation
\begin{equation}\label{eq:link}
    \dfrac{\partial q}{\partial s} + \dfrac{\partial }{\partial v}[-v q ] - a\dfrac{\partial q}{\partial v^2} = 0, \qquad q(V_F,s) = q(-\infty,s)=0, \qquad q(v,0) = \delta_{v=V_R},
\end{equation}
posed for voltage $v\in(-\infty,V_F]$ and time since the last discharge $s\in[0,+\infty)$. The solution $q$ represents, when a neuron fires and is reset at voltage $v=V_R$, how its probability distribution in voltage evolves until the next spike. In the limit $s\to +\infty$, the neuron has more and more probability to have fired a second time, as represented by the decay of $q(v,s)$ towards 0 in age. The quantity
\[   S(s) = \int_{-\infty}^{V_F} q(v,s)\diff v  \]
is the survival probability at age $s$ since the last discharge. It has been used in \cite{IDL} for parameter estimation of the linear NNLIF parameters.

\begin{theorem}[Dumont, Henri, Tarniceriu, \cite{DHT,DHT2}]
    Let $\rho$ a solution to \eqref{NNLIF} with $b=0$ from an initial condition $\rho^0$, $n$ a solution to \eqref{eq:T-E} with $\mu_X = 0$ from an initial condition $n^0$ and $q$ a solution of \eqref{eq:link}. If
    \begin{equation*}
        \forall v\in(-\infty,V_F],\qquad \rho^0(v) = \int_{0}^{+\infty} \dfrac{q(s,v)}{\int_{-\infty}^{V_F}q(v,a)\diff a} n^0(s) \diff s,
    \end{equation*}
    then
    \begin{equation*}
        \forall (v,s,t)\in(-\infty,V_F]\times \R_+\times\R_+,\qquad \rho(v,t) = \int_{0}^{+\infty} \dfrac{q(s,v)}{\int_{-\infty}^{V_F}q(v,a)\diff a} n(s,t) \diff s.
    \end{equation*}
    Moreover, for all $t$ such that $0 < s < t$, 
    \begin{equation*}
        n(s,t) = \dfrac{\partial}{\partial s}\int_{-\infty}^{V_F} q(v,s) \rho(v,t-s)\diff v.
    \end{equation*}
    and for $s\geqslant t$,
    \begin{equation*}
        \int_t^{+\infty} n(s,t) \diff s = \int_{-\infty}^{V_F} q(v,s)\rho^0(v)\diff v.
    \end{equation*}
\end{theorem}
This correspondence can be interpreted as a Bayes rule: the kernel in the integral transform represents the probability density for a neuron at electric potential $v$ given that it did not spike at time $s$ since the last discharge.

In \cite{DHT3}, a unifying approach is proposed in the nonlinear case $b\neq 0$, $\mu_X\neq 0$. The idea is to consider the stochastic system
\[ \diff V_t = ( - v + bN(t))\diff t + \sigma \diff B_t, \qquad \diff S_t = \diff t,  \]
with $N(t)$ the firing rate of the network. When a neuron has voltage $V_{t^-}=V_F$, it is reset at $V_{t^+}=V_R, S_{t^+} = 0$.
The two-dimensional kinetic Fokker-Planck equation for the density, that we also call the master equation is
\begin{equation}
    \left\{
    \begin{array}{l}
         \dfrac{\partial M}{\partial t} +\dfrac{\partial M}{\partial s} + \dfrac{\partial }{\partial v}\left[  (-v+bN(t))M \right] - a \dfrac{\partial M}{\partial v^2} = 0,  \\
         \displaystyle N(t) = - a \int_0^{+\infty} \dfrac{\partial M}{\partial v}(V_F,s,t)\diff s,\\
         M(V_F,s,t ) = M(-\infty,s,t) = M(v,+\infty,t)=0, \quad M(v,0,t) = N(t)\delta_{v=V_R},\\
         M(v,s,0) = M^0(v,s)\geqslant 0, \quad \int_0^{+\infty}\int_{-\infty}^{V_F} M^0(v,s)\diff v\diff s = 1.
    \end{array}
    \right.
\end{equation}
where $a=\frac{\sigma^2}{2}$.
Integrating the density $M$ with respect to age $s$ yields a solution $\rho$ to \eqref{NNLIF}:
\[   \rho(v,t) = \int_0^{+\infty} M(v,s,t)\diff s.  \]
A solution $n$ to \eqref{eq:T-E} with $p(s,t) = \frac{-a\partial_v q(V_F,s,t)}{\int_{-\infty}^{V_F}q(v,s,t)\diff v}$ can also be recovered by computing
\[ M(v,s,t) = \dfrac{q(v,s,t)}{\int_{-\infty}^{V_F}q(v,s,t)} n(s,t), \]
where $q$ solves the first passage time problem
\[  \dfrac{\partial q}{\partial t} +\dfrac{\partial q}{\partial s} + (-v+bN(t))\dfrac{\partial q}{\partial v} - a \dfrac{\partial q}{\partial v^2} = 0, \quad q(v,0,t) = \delta_{v=V_R}, \quad q(V_F,s,t) = 0. \]
Both these fact are rigorously proved in the linear case $b = \mu_X = 0$ in \cite{DHT3}. Formal integration show they hold in the nonlinear case $b\neq 0$, but then, the strong coupling through $N(t)$ between the equations for $q$ and $M$ has to be taken into account. Also, the reduction to \eqref{eq:T-E} is not standard as the firing function does not depends straightforwardly on the network activity.

\section{Jump based Leaky Integrate \& Fire model}

An issue with many of the models we have presented so far is that they rely on a mean-field assumption of an infinite number of all-to-all connected neurons. For example, in the NNLIF model in Section \ref{sec:NNLIF}, the K neurons in the particle system make jumps of size $\frac{b}{K}$ upon receiving a spike. In some contexts, it can be useful to keep jumps at a non-vanishing size during the mean-field limit. The jump-based Leaky Integrate \& Fire model is an example of this idea.

We will discuss how it has been obtained by physicists, then propose a mathematical setting, discuss well-posedness theory in both function and measure spaces, and talk about the occurence of finite time blow-up in strongly excitatory networks. Last, we will explain how this model can be connected back to the NNLIF model through heuristic arguments in the vanishing jumps limit.

\subsection{Macroscopic description of deterministic IF neurons}

Like for the derivation of the model \eqref{NNLIF} in Section \ref{sec:NNLIF}, let us start from a (now renormalised) Lapicque Integrate and Fire model
\begin{equation}
    \dfrac{\diff V}{\diff t}(t) = -V(t) + h\sum_{j\in\N} \delta_{t=t_{j}},\qquad V(t) \geqslant V_F \implies V(t^+) = V_R,
\end{equation}
where a representative neuron of voltage $V(t)$ receives spikes from other neurons at times $t_j$. The quantity $h$ represents an average jump in voltage upon receiving a spike. 

Instead of assuming Poisson discharges and replacing the incoming spikes by an Ornstein-Uhlenbeck process of same mean and variance like in the NNLIF derivation, another popular choice \cite{OKS,KOS,GK} to obtain a macroscopic description consistent with behaviours in particle simulations \cite{OBH} is to formally consider a large number of $N$ replicas of the population of neurons and to take an average over all the replicas. If for the $k^{th}$ replica, there are $n_k(v,t)\diff v$ neurons in a voltage element $[v,v+\diff v)$, then the macroscopic probability density $\rho$ is constructed as
\[ n(v,t) = \lim_{N\to+\infty}\frac{1}{N}\sum_{N\to +\infty} n_k(v,t),\qquad \rho(v,t) = \dfrac{n(v,t)}{\int_{-\infty}^{V_F} n(w,t)\diff w}.  \]
The initial conditions for the replicas are randomly drawn but there is no source of noise in the time evolution of each of them, leading to a hyperbolic equation with no intrinsic noise source. Note an important distinction between the model \eqref{NNLIF} and this approach: because we do not take the mean-field limit for a number $N$ of all-to-all connected neurons going to $+\infty$, we do not need to introduce a mean-field scaling $\frac{h}{N}$ in the process. The macroscopic model obtained with the replica ansatz then has finite-sized, non-vanishing jumps when excitatory spikes are received at the macroscopic level, leading to a jump-based structure in the PDE.

After constructing each quantity involved through averages over the replicas, \cite{OKS} obtains the evolution equation
\[\dfrac{\partial \rho}{\partial t}(v,t) + \dfrac{\partial}{\partial v}\left[-v\rho(v,t)\right] + \sigma(t) \left[\rho(v,t)-\rho(v-h,t)\right] = r(t)\delta_{V_R},\]
where the function $\sigma(t)$ is the deterministic flux of excitatory spikes reaching the population and $r(t)$ is the quantity of neurons spiking and being reset at time $t$. Then, $\sigma(t)$ is itself a sum of the external stimulation $\sigma_0\geqslant 0$ from outside the network and the internal excitatory stimulation $J_E r(t)$ coming from neurons in the network spiking, with $J_E\geqslant 0$ the excitatory connectivity strength.

The quantity $\rho(v,t)-\rho(v-h,t)$ represents the neurons in a window $[v-h,v]$ which are excited by the spikes. It can be better understood in divergence form
\[   \sigma(t) \left[\rho(v,t)-\rho(v-h,t)\right] = \dfrac{\partial}{\partial v} \left(\sigma(t)\int_{v-h}^{v} \rho(w,t)\diff w\right). \]
In particular, the flux $r(t)$ of neurons spiking at time $t$ corresponds to the part of the network making the $h$-sized jump from $V_F-h$ to $V_F$:
\[  r(t) = \sigma(t)\int_{V_F-h}^{V_F} \rho(w,t)\diff w. \]

Recurrent inhibition can be taken into account by adding the term
$ I(t) [\rho(v+h)-\rho(v)]$, to the right-hand side of the equation, where $I(t) = J_I r(t)$ is the inhibitory current and $J_I\geqslant 0$ is the inhibitory connectivity parameter. A comparison between the numerical simulation of the PDE \eqref{eq:jLIF} and a Monte-Carlo simulation of a particle system can be found in \cite{DH}.

Instead of the instantaneous transmission $\sigma(t) = \sigma_0 + J_E r(t)$, $I(t) = J_I r(t)$ a delayed transmission can be considered like for the model \eqref{eq:T-E}:
\[ \sigma(t) = \sigma_0 + J_E \int_0^t r(t-\tau) \alpha(\tau)\diff \tau, \qquad I(t) = J_I \int_0^t r(t-\tau) \alpha(\tau)\diff \tau, \qquad \int_{0}^{+\infty} \alpha(\tau)\diff \tau = 1. \]
Since there is no noise (no diffusion) an the deterministic field $-v$ does not send the mass in negative voltage space, the deterministic case $J_E\geqslant 0$, $J_I=0$ can be equivalently posed on $[0,V_F]$ by choosing an initial condition supported in $[0,V_F]$.

\subsection{Mathematical setting}
Following several references on the problem 
\cite{DH,DH2,DG}, we make the simplifying choice $V_F=1$ and sum up the jump-based Leaky Integrate and Fire model as the closed system
\begin{equation}\label{eq:jLIF}\tag{jLIF}
    \left\{
    \begin{array}{l}
        \dfrac{\partial \rho}{\partial t} + \dfrac{\partial}{\partial v}\left[-v\rho\right] + \sigma(t) \left[\rho(v,t)-\rho(v-h,t)\right] = I(t)\left[\rho(v+h,t)-\rho(v,t)\right] + r(t)\delta_{v=V_R} \\
        \displaystyle r(t) = \sigma(t)\int_{1-h}^{1} \rho(v,t) dv,\\ \displaystyle \sigma(t) = \sigma_0(t) +J_E\int_0^t r(t-\tau) \alpha(\tau)\diff\mu_\sigma,\qquad I(t) = J_I\int_0^t r(t-\tau) \alpha(\tau)\diff\mu_\sigma\\
        \displaystyle \rho(1,t) = 0,\qquad \rho(v,0) = \rho^0(v)\geqslant 0,\qquad \int_{-\infty}^{1} \rho^0(v,t)\diff v = 1.
    \end{array}
    \right.
\end{equation}
The abstract measure\footnote{Beware that contrary to our convention for \eqref{eq:T-E} we do not include the connectivity strength $J$ into the transmission measure; in both cases we chose the most suitable notation to encompass the largest proportion of references on the topic.} $\mu_\sigma$ allows to represent both instantaneous transmission with $\mu_\sigma = \delta_{\tau=0}$ and smoothly delayed transmission with $\mu_\sigma = \alpha(\tau)\diff \tau$. Note that when there is instantaneous transmission, $\sigma(t)$ solves the closed equation
\[
    \sigma(t) = \dfrac{\sigma_0}{1- J_E\int_{1-h}^{1}\rho(v,t)dv}.
\]
In most cases, solutions are to be understood in the weak sense: for all suitable $C^1$ test function $v\mapsto\phi(v)$, we require in the sense of distributions
\begin{multline}\label{eq:jLIF_weak}
    \dfrac{\diff }{\diff t}\int_{-\infty}^1 \phi(v)\rho(v,t)\diff v + \int_{-\infty}^1 v\dfrac{\diff \phi}{\diff v}(v)\rho(v,t)\diff v + \sigma(t) \int_{-\infty}^1 [\phi(v)-\phi(v+h)]\rho(v,t)\diff v\\
    +\sigma(t) \int_{1-h}^1 \rho(v,t)\phi(v+h)\diff v = I(t)\int_{-\infty}^1 [\phi(v-h)-\phi(v)]\rho(v,t)\diff v + \phi(V_R) r(t).
\end{multline}
Note that this formulation is obtained by extending $\rho$ and $\phi$ by 0 past $V_F=1$ in the computations. Taking $\phi(v)=1$ in \eqref{eq:jLIF_weak} yields mass conservation:
\[  \int_{-\infty}^1 \rho(v,t)\diff v = \int_{-\infty}^1\rho^0(v)\diff v = 1.  \]

When working with measure solutions, we use the same formalism as for \eqref{eq:T-E} in Section \ref{sec:TE}. We denote $\mathcal M(E)$ the set of finite measures on a measurable space $E$, $\mathcal M_+(E)$ the set of positive finite measures, $\mathcal P(E)$ the set of probability measures and $\norme{\cdot}_{\mathrm{TV}}$ the total variation norm.

\begin{definition}[measure solutions]\label{def:jLIF_meas} Let $T > 0$, we say that a family $(\rho_t)_{t\geqslant 0}$ of $\mathcal P([0,1])$ is a solution to the excitatory ($J_I=0$, $J_E\geqslant 0$) system \eqref{eq:jLIF} on $[0,T)$ with initial datum $\rho_0=\rho^0$ if
\begin{itemize}
    \item $t\mapsto \sigma(t):=\frac{\sigma_0(t)}{1-J_E \rho_t([1-h,1])}$ is positive and locally integrable;
    \item $t\mapsto \rho_t$ is weak-* continuous on $[0,T)$;
    \item for all $\phi\in C^1([0,1])$, for all $t\in[0,T)$,
    \begin{multline*}
        \int_0^1 \phi(v)\diff \rho_t = \int_0^1 \phi(v)\diff \rho_0\\
        + \int_0^t\int_0^1 \left(-v\dfrac{\diff \phi}{\diff v}(v) + \sigma(s)\Big[\phi(v+h)\mathds 1_{[0,1-h]}(v)-\phi(v)+\phi(V_R)1_{[1-h,1]}(v)\Big]\right)\diff \rho_s \diff s.
    \end{multline*}
\end{itemize}
\end{definition}

Define the operators $\mathcal A_\sigma$ and $\mathcal B$ acting on $C^1([0,1])$ by
\begin{equation*}
    \mathcal A_\sigma \phi (v)= -v\dfrac{\diff \phi}{\diff v}(v) + \sigma \mathcal B\phi (v),\qquad \mathcal B \phi (v) = \phi(v+h)\mathds 1_{[0,1-h]}(v)-\phi(v)+\phi(V_R)1_{[1-h,1]}(v).
\end{equation*}
Then, denoting $\rho_t \phi = \int_0^1\phi(v)\diff \rho_t$, the last part of the definition of a measure solution reduces to
\[  \rho_t \phi = \rho_0\phi + \int_0^t  \rho_s \mathcal A_{\sigma(s)} \phi \diff s.  \]

\subsection{Local and global well-posedness of the jump equation}

\subsubsection{Case with distributed delay}

Consider the excitatory case $J_E \geqslant 0$ and $J_I=0$. Let us first define a capped problem using a function $D(x)=\min(M,\max(0,x))$ for some constant $M\geqslant 0$:
\begin{equation}\label{eq:jLIF_capped}
    \dfrac{\partial \rho}{\partial t}(v,t) + \dfrac{\partial}{\partial v}\left[-v\rho(v,t)\right] + D(\sigma(t)) \left[\rho(v,t)-\rho(v-h,t)\right] = D(\sigma(t))\int_{1-h}^{1} \rho(w,t) dw\ \delta_{V_R}.
\end{equation}
Then, a solution to \eqref{eq:jLIF_capped} is a solution to \eqref{eq:jLIF} as long as $0 \leqslant \sigma(t)\leqslant M$.

Like for \eqref{eq:T-E}, the method to prove well-posedness is to first consider the linear problem with $J = 0$. The solution can be constructed \textit{via} a fixed-point argument in the Banach space $\mathcal X=C([0,T],L^1(0,1))$ endowed with it's natural Banach norm on the mapping $\Lambda:q\to\rho$ where $\rho$ is the solution to
\begin{equation*}
    \dfrac{\partial \rho}{\partial t}(v,t) + \dfrac{\partial}{\partial v}\left[-v\rho(v,t)\right] + D(\sigma(t)) \left[\rho(v,t)-q(v-h,t)\right] = D(\sigma(t))\int_{1-h}^{1} q(w,t) dw\ \delta_{V_R},
\end{equation*}
with same initial and boundary condition as \eqref{eq:jLIF}. For small enough $T$, $\Lambda$ is a contraction mapping and the argument can be iterated to obtain a global-in-time solution to the linear problem.

Then, the solution to the (capped) nonlinear problem can be obtain with a second fixed point argument in $C([0,T])$ applied to the mapping
\begin{equation}
    F:\sigma\mapsto \sigma_0 + J \int_0^t G[\sigma](t-\tau)\alpha(\tau)\diff \tau, \qquad G[\sigma](t)=D(\sigma) \int_{1-h}^1 \rho(w,t)\diff w,
\end{equation}
with $\rho$ solution to the linear problem with external input $\sigma$.

Last, by deriving local in time bounds of the form
\[ r(t) \leqslant \norme{\sigma_0}_\infty e^{J\norme{\alpha}_\infty t}, \]
the cap $D$ can be removed on any finite time interval $[0,T]$ by letting $M$ be large enough.

\begin{theorem}[Dumont, Henry \cite{DH}]
    Assume $\mu_\sigma = \alpha(\tau)\diff \tau$, with $\alpha\in L^\infty(\R_+)\cap L^1(\R_+)$ and $\int_0^{+\infty}\alpha(\tau)\diff \tau = 1$, $\sigma_0\in C(\R_+)\geqslant 0$, $J_I=0$ and $J_E\geqslant 0$ (excitatory case). For all $\rho^0\in L^1_+(0,1)$ of unit mass, there exists a unique solution $\rho\in C(\R_+, L^1_+(0,1)$ to \eqref{eq:jLIF}, and it satisfies
    \begin{equation*}
        r(t) \leqslant \norme{\sigma_0}_\infty e^{J\norme{\alpha}_\infty t}.
    \end{equation*}
    If $J_E < 1$, it satisfies the improved bound
    \begin{equation*}
        r(t)\leqslant \dfrac{\norme{\sigma_0}_\infty}{1-J_E}.
    \end{equation*}
\end{theorem}

In the inhibitory case, the same strategy applies, but now the equation has to be considered on $(-\infty,1]$ instead of just $[0,1]$.

\begin{theorem}[Dumont, Henry \cite{DH}]
    Assume $\mu_\sigma = \alpha(\tau)\diff \tau$, with $\alpha\in L^\infty(\R_+)\cap L^1(\R_+)$ and $\int_0^{+\infty}\alpha(\tau)\diff \tau = 1$, $\sigma_0\in C(\R_+)\geqslant 0$, $J_I\geqslant 0$ and $J_E = 0$ (inhibitory case). For all $\rho^0\in L^1_+(-\infty,1)$ of unit mass, there exists a unique solution $\rho\in C(\R_+, L^1_+(-\infty,1)$ to \eqref{eq:jLIF}.
\end{theorem}

\subsubsection{Instantaneous transmission}

As was already noticed within a less rigorous formalism in \cite{OKS}, when $\sigma(t) = \sigma_0+J\sigma(t)$, finite time blow-up can happen. It is easier to see in the reformulation
\[ \sigma(t) = \dfrac{\sigma_0}{1- J\int_{1-h}^{1}\rho(w,t)dw},\]
as when $J\int_{1-h}^{1}\rho(w,t)dw$ approaches 1, $\sigma(t)$ ceases to exist in the classical sense.

However, since $\rho$ is a probability density, as long as $J<1$, we should have $J\int_{1-h}^{1}\rho(w,t)dw<1$ and the firing rate $\sigma(t)$ should be well-defined.

\begin{theorem}[Dumont, Henry \cite{DH}]
    Assume $\mu_\sigma = \delta_{\tau=0}$, $\sigma_0\in C(\R_+)\geqslant 0$, $J_I=0$ and $J_E\in [0,1)$ (weakly-excitatory case). For all $\rho^0\in L^1_+(0,1)$ of unit mass, there exists a unique solution $\rho\in C(\R_+, L^1_+(0,1)$ to \eqref{eq:jLIF}, and it satisfies
    \begin{equation*}
        r(t)\leqslant \dfrac{\norme{\sigma_0}_\infty}{1-J_E}.
    \end{equation*}
\end{theorem}

The proof follows the same pattern as in the delayed case, with some shortcuts. However, for a general excitatory connectivity $J_I=0$, $J_E\geqslant 0$, as the solution may blow-up (see later), it is only possible to obtain a local-in-time solution. It can be sought in the space of admissible functions
\[\mathcal X = \left\{ \rho\in L^1_+(0,1)\ | \  \int_0^1\rho(v)\diff v = 1,\quad J_E \int_{1-h}^{1} \rho(v)\diff v < 1 \right\}.\]
Then, it is possible to apply another fixed point procedure to a version of \eqref{eq:jLIF} with a mollified impact of $\sigma$ and obtain the following result.
\begin{theorem}[Dumont, Henry \cite{DH2}]
    Assume $\sigma_0\in C_b(\R_+,\R_+)$, $J_I=0$, $J_E\geqslant 0$ (excitatory case). For all $n^0\in  L^1_+(0,1)$ of unit mass, there exists a unique maximal solution $\rho\in C([0,T^*),L^1_X(0,1))$ to \eqref{eq:jLIF}.
\end{theorem}

The inhibitory case also proceeds similarly, but there is no issue with finite time blow-up even when $J_I$ is large, allowing to achieve unconditional global-in-time existence.

\begin{theorem}[Dumont, Henry \cite{DH}]
    Assume $\mu_\sigma = \delta_{\tau=0}$, $\sigma_0\in C(\R_+)\geqslant 0$, and $J_I\geqslant 0$, $J_E=0$ (inhibitory case). For all $\rho^0\in L^1_+(-\infty,1)$ of unit mass, there exists a unique solution $\rho\in C(\R_+, L^1_+(-\infty,1))$ to \eqref{eq:jLIF}.
\end{theorem}

\subsubsection{Solutions in the sense of measures}

Consider the case $\mu_\sigma$ (instantaneous transmission) and $J_I=0,J_E\in[0,1)$ (weakly-excitatory). Following the duality method applied in \cite{G} on a toy example, \cite{DH} constructs global-in-time solutions to \eqref{eq:jLIF}. See also the construction of measure solutions to the system \eqref{eq:T-E} in Subsection \ref{sec:TE_meas_sol} which lies in the same toolbox and was performed a year earlier in \cite{CY}. If $\psi(s,t,v)$ satisfies
\[  \frac{\partial \psi}{\partial t} = v\frac{\partial \psi}{\partial v} + \dfrac{\sigma_0}{1-J_E \mu_0 \psi(0,s,\cdot)} \left[ \psi(s,t,v)-\psi(s,t,v+h)\mathds1_{[0,1-h)}(v) -\psi(s,t,V_R)\mathds1_{[1-h,1]}(v)\right], \]
with terminal condition $\psi(t,t,v) = \mathds 1_{[1-h,1]}(v)$, then $\rho_t ([1-h,1]) = \rho_0 \psi(0,t,\dot) $. Solving the dual equation allows to obtain a Duhamel formula and to desingularise the problem, proceeding similarly to the construction of the linear semigroup below.

However, because of the discontinuity introduced by terms like $\mathds1_{[1-h,1]}$, and because of the need to solve for $\psi$ in the space of continuous functions in order to get to $\mathcal M([0,1])$ by duality, the existence proof involves intermediate mollified $\psi_n$ equations and a lot of technicalities that we omit in our account. 

\begin{theorem}[Dumont, Gabriel \cite{DG}]
    Assume $\mu_\sigma = \delta_{\tau=0}$ (instantaneous transmission), $J_I=0$, $J_E < 1$ (weakly-excitatory case). Then for all $\rho_0\in \mathcal P([0,1])$, there exists a unique global-in-time measure solution to \eqref{eq:jLIF} in the sense of Definition \ref{def:jLIF_meas}.
\end{theorem}

\subsection{Finite time blow-up in the strongly-excitatory regime}

In \cite{DH2}, it is proved that finite time blow-up does occur for the excitatory \eqref{eq:jLIF}. Unlike in the model \eqref{NNLIF}, blow-up does not arise in weakly-connected networks. However, the universal blow-up property for high enough connectivity still holds.

\begin{theorem}[Dumont, Henry \cite{DH2}]
    If $\mu_\sigma = \delta_{\tau=0}$ (instantaneous transmission) and
    \[ J_I=0, \qquad J_E \geqslant \dfrac{1-V_R}{h} + 1 ,\qquad h \sigma_0 > 1,  \]
    then for all initial condition $\rho^0\in L^1_+(0,1)$ with unit mass, the solution $\rho$ to \eqref{eq:jLIF} cannot be global-in-time. The maximal time of existence $T^*$ satisfies
    \begin{equation}
        T^*\leqslant \dfrac{1}{h\sigma_0-1},\qquad \limsup_{t\to T^*} r(t) = +\infty.
    \end{equation}
\end{theorem}

\begin{remark}
    Under those conditions, there is in particular no stationary state to \eqref{eq:jLIF}.
\end{remark} 

The proof proceeds similar to the one of Theorem \ref{ExploPert}, by considering the weak formulation with test function $\phi(v)=e^{\mu v}$, for a constant $\mu > 0$ large enough (not to be confused with the measure $\mu_\sigma$). The exponential moment $M_\mu(t) = \int_0^1 \rho(v,t)e^{\mu v}\diff v$ then satisfies, after some computations and using the hypotheses,
\[  \dfrac{\diff M_\mu}{\diff t} (t)\geqslant K_\mu M_\mu(t) +  \mu r(t)\left( Jh M_\mu(t) - \frac{e^{1+h}-e^{V_R}}{\mu}  \right). \]
for some constant $K_\mu>0$. We have
\[  \dfrac{\diff M_\mu}{\diff t} (t)\geqslant \mu r(t)\left( Jh M_\mu(t) - \frac{e^{1+h}-e^{V_R}}{\mu}  \right). \]
So if $M_\mu(0)>\frac{e^{1+h}-e^{V_R}}{\mu}$, then $M_\mu$ is non-decreasing and
\[  \dfrac{\diff M_\mu}{\diff t} (t)\geqslant K_\mu M_\mu(t), \]
which leads to a contradiction with the natural bound $M_\mu(t)\leqslant e^{\mu} \int_0^1\rho(v,t)\diff v = e^{\mu}.$
The problem is reduced to prove that for any initial condition there exists a time $t^*$ such that $M_\mu(0)>\frac{e^{1+h}-e^{V_R}}{\mu}$, which happens due to the strength of the external drift $\sigma_0$ in the assumptions.

\subsection{Long time behaviour}

In \cite{DG}, the long time asymptotic of solutions in the linear and weakly-connected excitatory regimes are investigated for \eqref{eq:jLIF}. Like in \cite{CY}, the Doeblin method allows to do so in the space of measures and to obtain quantitative global asymptotic stability of a unique invariant measure. Since we already described some parts of this procedure, we will not enter all technical details; see Subsection \ref{sec:TE_Doeblin}.

\subsubsection{Study of the linear asymptotic}

Starting with the spectral gap in the linear case, \cite{DG} builds on the method of \cite{G} based on the dual equation
\[\dfrac{\partial f}{\partial t} + v\dfrac{\partial f}{\partial v} + \sigma_0 f = \sigma_0[f(v+h,t)\mathds 1_{[0,1-h)}(v)+f(V_R,t)\mathds 1_{[1-h,1]}(v),\]
which they prove to be well-posed in the space of continuous functions. Then, it is possible to define a positive conservative semigroup $(M_t)_{t\geqslant 0}$ acting on $C([0,1])$ such that $M_t f^0 = f(\cdot,t)$. It allows in turn to define by duality a semigroup $(M_t)_{t\geqslant 0}$ acting on $\mathcal M([0,1]) = C([0,1])'$, which is also positive and conservative (because of the duality construction, this semigroup acts left). It preserves the set of probabilities $\mathcal P([0,1])$ and is a contraction in total variation norm.

After some technicalities involving regularisation, \cite{DH} proves that this semigroup provides solutions to the linear problem. They also prove that it satisfies Doeblin's condition, which leads, applying Doeblin theorem, to the following result.

\begin{theorem}[Dumont, Gabriel \cite{DG}]
     Assume $\mu_\sigma = \delta_{\tau=0}$ (instantaneous transmission), $\sigma_0(t)\equiv\sigma_0>0$ is constant and $J_I=J_E=0$ (linear case). Then the equation \eqref{eq:jLIF} generates a positive contraction left-acting weak-* continuous semigroup $(M_t)_{t\geqslant 0}$ on $\mathcal M([0,1])$. The set of probability measures $\mathcal P([0,1])$ is invariant under $(M_t)_{t\geqslant 0}$. It admits a unique invariant probability measure $\rho_\infty$ and for all $\rho^0\in\mathcal P([0,1])$,
     \begin{equation*}
         \norme{\rho^0 M_t  - \rho_\infty}_{\mathrm{TV}} \leqslant e^{-\nu (t-t_0)}\norme{\rho^0 - \rho_\infty}_{\mathrm{TV}},
     \end{equation*}
     where 
     \[ t_0 = \log\left(\frac{4}{h}\right)>0,\qquad  \nu = -\dfrac{\log\left(1-\frac{\sigma_0}{2}\left(\frac{h}{4}\right)^{\sigma_0}\right)}{\log\left(\frac{4}{h}\right)} >0.   \]
\end{theorem}
\subsubsection{Nonlinear stability and existence of stationary states}

Existence of a nonlinear stationary state in the excitatory case with instantaneous transmission is equivalent to finding a constant $\sigma_\infty$ such that
\[  \sigma_\infty = \dfrac{\sigma_0}{1-J_E \rho^{\sigma_\infty}([1-h,h]), } \]
where $\rho^{\sigma_\infty}$ is the invariant measure to the linear problem with external input $\sigma_\infty$. This done in \cite{DG} by defining the mappings $F,G$ by
\[  F(\sigma) = \rho^\sigma([1-h,h]),\qquad G(\sigma) = \frac1{J_E}\left( 1-\frac{\sigma_0}{\sigma} \right), \]
and by looking for a solution to $F(\sigma_\infty)=G(\sigma_\infty)$.

First, one can prove that for fixed constant external $\sigma$, the invariant measure is absolutely continuous with respect to Lebesgue measure and enjoys some regularity;
\[  v\rho_\infty^\sigma \in W^{1,\infty}([0,V_R)\cup(V_R,1]),\qquad \forall v\in(0,1),\ \ 0 < \rho_\infty^\sigma (v) \leqslant \min\left( \frac{\sigma v^{\sigma-1}}{h^\sigma} , \frac{\sigma}{v}\right). \]
Note that then the nonlinear stationary state also enjoys such regularity. From this regularity and the decomposition of the measure formulation with $\mathcal A_\sigma$ and $\mathcal B$, \cite{DG} proves asymptotic properties on $F$:
\[ \lim_{\sigma\to 0} F(\sigma) = 0,\qquad \lim_{\sigma\to+\infty} F(\sigma) = \dfrac{1}{1+ \left\lfloor\frac{1-V_R}{h}\right\rfloor}. \]

Combining these results, using the Duhamel formula and the contraction of the linear semigroup $(M_t)_{t\geqslant 0}$, and working some more yields
\begin{theorem}[Dumont, Gabriel \cite{DG}]\label{thm:LIF_CV}
     Assume $\mu_\sigma = \delta_{\tau=0}$ (instantaneous transmission), $\sigma_0(t)\equiv\sigma_0>0$ is constant and $J_I=0$, $J_E\geqslant 0$ (excitatory case).
    \begin{itemize}
        \item There exist at least two stationary states to \eqref{eq:jLIF} under the condition
        \[  J_E >1+ \left\lfloor\frac{1-V_R}{h}\right\rfloor,\qquad \sigma_0 <\frac{1-h}{4J_E} . \]
        \item There exists at least one stationary states to \eqref{eq:jLIF} under the condition \[  J_E < 1+ \left\lfloor\frac{1-V_R}{h}\right\rfloor.\] 
        \item If
        \[ J_E < \left(5-2\sqrt6\right)\left(\frac{h}{4}\right)^{\sigma_0+1},  \]
        there exists a unique stationary state $\rho_\infty$, and there are constants $t_0\geqslant 0,\nu>0$ such that for all initial condition $\rho_0 \in \mathcal P([0,1])$ the solution to \eqref{eq:jLIF} satisfies,
        \begin{equation*}
            \norme{\rho_t-\rho_\infty}_{\mathrm{TV}} \leqslant e^{-\nu(t-t_0)}\norme{\rho_0-\rho_\infty}_{\mathrm{TV}}. 
        \end{equation*}
    \end{itemize}
\end{theorem}

\subsection{Heuristic link to the NNLIF model}

As noted in the physics literature \cite{OKS} and explained in detail in \cite{DH}, the singular structure of the jump mechanism in \eqref{eq:jLIF} can be approximated by a Taylor expansion when the jumps of size $h$ are small enough, leading to a diffusive model which happens to be the model \eqref{NNLIF} from Section \ref{sec:NNLIF}. Consider the excitatory case $J_I= 0, J_E\geqslant 0$. Denote $b_0 + b N(t) = h\sigma(t)$ for some $b_0 > 0, b\geqslant 0$. If $h$ is small enough, we can write the second order expansion
\begin{align*}  \sigma(t) [\rho(v-h,t)-\rho(v,t)] \simeq\ &  \sigma(t) h\dfrac{\partial \rho}{\partial v} - \sigma(t)\frac{h^2}{2}\dfrac{\partial^2 \rho}{\partial v^2}\\
\simeq\ & (b_0+b N(t))\dfrac{\partial \rho}{\partial v} - (b_0+b N(t)) \frac{h}{2}\dfrac{\partial^2 \rho}{\partial v^2},
\end{align*}
and, using that $\rho(1,t)=0$
\begin{align*}
    r(t) = \sigma(t) \int_{1-h}^1 \rho(v,t)\diff v \simeq \ & \sigma(t) h \rho(1,t) - \sigma(t)\frac{h^2}{2} \dfrac{\partial \rho}{\partial v}(1,t)\\
    \simeq\ & - (b_0+b N(t))\frac{h}{2} \dfrac{\partial \rho}{\partial v}(1,t).
\end{align*}
Note then that for instantaneous transmission $\mu_\sigma = \delta_{\tau=0}$,
\[ b_0+b N(t) = h\sigma(t) = h(\sigma_0 + J_E r(t)) =h\sigma_0+ h J_E r(t), \]
which allows to identify $b_0 = \sigma_0 h$, $N(t)=r(t)$, $b = hJ_E$.
Plugging these approximations into \eqref{eq:jLIF}, denoting 
\[a\big(N(t)\big) = \frac{h}{2}\left( b_0 + b N(t)\right),\] and rescaling the voltage $v$ so as to make a term $b_0$ disappear in the $-v$, eventually yields
\begin{equation*}
    \dfrac{\partial \rho}{\partial t} + \dfrac{\partial }{\partial v}\big[ (-v + bN(t) )\rho \big] - a\dfrac{\partial^2 \rho}{\partial v^2} = N(t) \delta_{v=V_R}.
\end{equation*}
with
\[  N(t) = - a\big(N(t)\big) \dfrac{\partial \rho}{\partial t}(V_F,t). \]
It is exactly the model \eqref{NNLIF} in the case of intrinsic noise. Note that the right-hand side in the condition $J_E < \left(5-2\sqrt6\right)\left(\frac{h}{4}\right)^{\sigma_0+1}$ for unconditional asymptotic stability in Theorem \ref{thm:LIF_CV} goes to 0 as $h\to0$, consistent with the fact that the limiting NNLIF model exhibits blow-up for any positive connectivity.

\section{Some other hyperbolic models}

We provide here a short description of some other hyperbolic PDEs used in mathematical neuroscience.

\subsection{Leaky Integrate and Fire with random firing}
\label{sec:LIFrand}

Another hyperbolic model that was rigorously obtained as a mean-field limit \cite{FL,DGLP} is the Leaky Integrate and Fire with random firing. At the particle level, the neurons of potential $V_t^1,\dots,V_t^N$ evolve according to the system of SDEs
\begin{multline}\label{LIFRand_sto_particles}
    V_t^{i} = V_0^i + \int_0^t (-V_s^i+j_0)\diff s + \frac{J}{N} \sum_{j\neq i} \int_0^t \int_0^{+\infty} \mathds 1_{z\leqslant \phi_\varepsilon(V_{s^-}^j)} N^j(\diff s,\diff z)\\
    - \int_0^t \int_0^{+\infty} V_{s^{-}}^i\mathds 1_{z\leqslant \phi_\varepsilon(V_{s^-}^i)} N^i(\diff s,\diff z),
\end{multline}
where the $N^i(\diff s,\diff z)$ are independent Poisson measures with intensity $\diff s\diff z$, $v\mapsto\phi_\varepsilon(v)$ is a firing function, $J>0$ a connectivity parameter, $j_0>0$ an external excitatory input. In the limit $N\to+\infty$, the system converges towards
\begin{equation}\label{LIFRand_sto}
    V_t = V_0^i + \int_0^t (-V_s+j_0)\diff s + J  \int_0^t \mathbb E\left[  \phi_\varepsilon(X_s) \right]\diff s
    - \int_0^t \int_0^{+\infty}  V_{s^-}\mathds 1_{z\leqslant \phi_\varepsilon(V_{s^-})} N(\diff s,\diff z),
\end{equation}
whose law is a solution of
\begin{equation}\label{LIFRand}
	\left\{\begin{array}{l}
		\displaystyle \dfrac{\partial \rho}{\partial t} + \dfrac{\partial}{\partial v}\Big[(-v + j_0 + J r(t))\rho\Big]  + \phi_\varepsilon(v)\rho  = 0  ,\\ 
		\displaystyle r(t)= \int_{0}^{+\infty} \phi_\varepsilon(v)\rho(v,t)\diff v,\\
		\displaystyle \rho(0,t) = \dfrac{ r(t)}{j_0+J r(t)},\qquad 
		\displaystyle \int_{0}^{+\infty} \rho^0(v)\diff v = 1, \quad  \rho(v,0)=\rho^0(v) \geqslant 0,
	\end{array}\right.
\end{equation}
The model \eqref{LIFRand} shares similarities both with the LIF with jump \eqref{eq:jLIF} and the NNLIF model with random firing \eqref{NNLIFRand}.

\begin{figure}[ht!]
	\begin{center}
		\includegraphics[width=430pt]{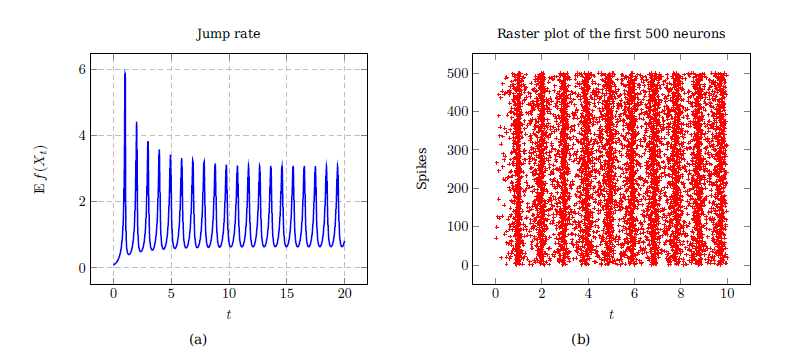}
	\end{center}
	\caption{Periodic solution in the LIF model with random firing. Left is the firing rate $N(t)$ in the PDE \eqref{LIFRand}, right represents spikes from the mean-field particle system \eqref{LIFRand_sto}. Cormier, Tanré, Veltz \cite{CTV2}.  \label{fig:CVT2}}	
\end{figure}

In \cite{CTV}, the model is studied at the level of the stochastic mean-field equation \eqref{LIFRand_sto}. They prove that for small enough $J$ there is a unique invariant measure and that the solutions converge towards it when $t\to+\infty$. In \cite{CTV2}, it is rigorously proved that the model \eqref{LIFRand} undergoes a Hopf birfurcation in terms of the parameter $\varepsilon$ for well-chosen $\phi_\varepsilon$, which proves the existence of time-periodic solutions (see Figure \ref{fig:CVT2}). For both studies \cite{CTV,CTV2}, a crucial tool is a Volterra equation related to the system. Further results about long-time behaviour and bistability were recently proved in \cite{C8}.

Note that \eqref{LIFRand} was also studied from the PDE point of view in the linear case in \cite{PS2} as the noiseless limit of \eqref{eq:VC_RF_lim_1} below in Subsection \ref{sec:VC_RF}. Conditions are given for linear instability of the stationary state in the nonlinear case, which can be related to the rigorously obtained periodic solutions in \cite{CTV2}.

\subsection{Continuity equation for pulse coupled oscillators}

A simple model for interacting neurons are the so-called pulse coupled oscillators first introduced in the context of modelling cardiac cells \cite{P4}. It consists in a system of ordinary differential equations (ODE) coupled by singular interaction. Consider $V^1(t),\dots,V^n(t)$ the potential at time $t$ of $n$ neurons. Each of them evolves according to
\[  \dfrac{\diff V^i}{\diff t}(t) = f\big(V^i(t)\big),  \]
and when $V^i(t^-) \geqslant V_F$, $V^i(t^+) = V_R = 0$ and the other neurons receive a kick $V^j(t^+) = V^j(t^-) + \frac{J}{n}$, $j\neq i$, $J\geqslant 0$. This is similar to the construction of \eqref{NNLIF} in Section \ref{sec:NNLIF}, but in a fully deterministic framework. This deterministic particle system has been investigated in depth, with the milestone \cite{MS4} establishing simple conditions for synchronisation or relaxation to the stationary state. A common strategy is to rescale the system from voltage variable $v\in[0,V_F]$ to phase variable $\phi\in[0,\Phi_F]$ with a constant vector field and a modified jump function:
\[  \dfrac{\diff \phi^i}{\diff t}(t) = 1,\qquad \phi^i(t^-) \geqslant V_F \implies \phi^i(t^+)=0, \ \phi^j(t^+) = \phi^j(t^-) + \frac{1}{n} K(\phi^j), \ j\neq i.  \]
A formal mean-field limit \cite{AV,K2} yields the continuity equation
\begin{equation}\label{eq:PC}
	\left\{\begin{array}{l}
		\displaystyle \dfrac{\partial \rho}{\partial t} + \dfrac{\partial}{\partial \phi}\Big[(1+K(\phi) N(t))\rho\Big]  = 0  ,\\[0.3cm] 
        N(t)= (1+K(\Phi_F) N(t))\rho(\Phi_F,t),\\[0.3cm]
		(1+K(\Phi_F) N(t))\rho(\Phi_F,t) = (1+K(0) N(t))\rho(0,t),\\[0.3cm]
		\displaystyle \int_{0}^{\Phi_F} \rho^0(\phi)\diff \phi = 1, \quad  \rho(\phi,0)=\rho^0(\phi) \geqslant 0.
	\end{array}\right.
\end{equation}
Notice that as in \eqref{eq:jLIF}, we have for the firing rate $N$ the closed expression
\begin{equation*}
    N(t)=\frac{\rho(t,\Phi_F)}{1-K(\Phi_F)\rho(t,\Phi_F)}.\quad
\end{equation*}
which indicates the model can blow-up in finite time.

This nonlinear mean-field problem has scarcely been investigated from the mathematical perspective. Criteria for the existence of stationary states and partial results about the synchronisation-relaxation dichotomy have been provided first in \cite{MS}. In particular, they show, similar to the findings for the particle system \cite{MS4}, that synchronisation occurs for increasing $K$ and relaxation to a stationary state for decreasing $K$ and adequate initial conditions. However, their results require many technical assumptions. Another approach has been proposed in \cite{CDRZ}: consider the repartition function $F(\phi,t) = \int_0^\phi \rho(\psi,t)\diff \psi$ and its left-continuous pseudo-inverse
\begin{equation}
    Q(\eta)=\inf \{\phi\geq0: F(\phi)\geq \eta\},\qquad \eta\in[0,1].
\end{equation}
The function $Q$ is also known as the quantile function of $\rho$. Consider the time dilation $\diff \tau = N(t)\diff t$; then, $Q(\eta,\tau)$ solves
\begin{equation}\label{eq:PC_Q}
\left\{\begin{array}{llll}
     &\partial_{\tau}Q+\partial_{\eta}Q=\dfrac{1}{N(\tau)}+{K}(Q), &&\tau>0,\ \eta\in(0,1),\\[0.3cm]
     &Q(\tau,0)=0, &&\tau>0,\\[0.2cm]
     &\dfrac{1}{N(\tau)}=\partial_{\eta}Q(\Phi_F,\tau)-{K}(\Phi_F), &&\tau>0,\\[0.3cm]
     &Q(\eta,0)=Q^0(\eta), &&\eta\in[0,1].
     \end{array}\right.
\end{equation} 
Note that this kind of time dilation has been used in previous \cite{DZ,TW,ST} and subsequent \cite{DPSZ2} works. The simpler drift and the decoupled nonlinearities $1/N$ and $K(Q)$ in \eqref{eq:PC_Q} have allowed \cite{CDRZ} to propose generalised versions of the results in \cite{MS} and to pave the road for continuation of solutions to \eqref{eq:PC} after blow-up.

\part{Kinetic Fokker-Planck equations}

\section{The kinetic Voltage-Conductance model}

We discuss here some kinetic voltage-conductance descriptions of neural networks based on the Moris-Lecar model for individual neurons. The degenerate diffusion, nonlinearity, and peculiar boundary conditions in these models make them difficult to study analytically, and many questions are still widely open.

We start by describing the historical construction of kinetic voltage-conductance models at the network level. Then we propose a mathematical framework, and we discuss how the already difficult linear case can be tackled with methods from kinetic theory or the use of probabilistic tools like Doeblin-Harris techniques. We then turn to the nonlinear case and discuss both \textit{a priori} estimates relying on direct integration and moment estimates, and a recent, more abstract approach which allows to prove the existence of solutions for weakly-connected networks. Finally, we address the issues around numerical simulations of this model and present some asymptotic limits and variants.

\subsection{Fluctuations at the level of the conductance}

The models we have discussed so far are overlooking an important aspect of the dynamics of individual neurons: the evolution of the membrane potential $v$ of individual cell depends on underlying variables with faster timescales. In large and mostly inhibitory networks observed over long periods of time, the voltage or activity descriptions provide good approximations, which is why the delayed inhibitory NNLIF model in Section \ref{sec:NNLIF} or the activity-based grid cell model in Section \ref{sec:GC} are good at capturing macroscopic behaviours like self-sustained oscillations and the hexagonal firing pattern of grid cells.

However, as we have seen, those models are less precise for the modelling of excitatory ensembles of neurons. For example, in the primary visual cortex V1, so-called \textit{simple} neurons can respond in a linear fashion to frequency in visual stimulation, when others, called \textit{complex}, respond nonlinearly, see \cite{HW}. It is believed that this \textit{simple} or \textit{complex} nature is not an individual characteristic of those cells but rather an emerging property from the network they are embedded in \cite{CNA}. The behaviour of the \textit{complex} type of V1 neuron has been obtained in a model with a strong excitation; The simple type has been found to be better represented by inhibitory models; see \cite{TSMS}.

More generally, accurately describing the fluctuation-driven dynamics of highly nonlinear networks of neurons with slow-postsynaptic responses requires to work at smaller scales for space and time. Around 2005, in the pioneering works \cite{CTSM,CTRM,RCM-PNAS05}, based upon previous literature about deterministic and stochastic representation of coarsed grained neural networks (see \textit{e.g.} \cite{FB,KMS,NT1,NT2,OKS,OKKS,BH}), Cai, McLaughlin, Rangan, Shelley and Tao developped an accurate kinetic description for excitatory and excitatory-inhibitory ensembles of neurons.

Let us focus on the fully excitatory case. They start with a coarse-grained patch of $N$ all-to-all connected excitatory neurons described by their voltage $V$ and conductance $G$: assume each neuron labelled by $j\in\{1,\dots,N\}$ evolves according to the Integrate \& Fire equations
\begin{align*}
	\tau \dfrac{\diff V^j}{\diff t}\, &=\, -(V^j-V_L) - G^j ( V^j - V_E ), \\
	\sigma_E \dfrac{\diff G^j}{\diff t}\, &=\, - G^j + S_{\mathrm{ext}} \sum_{l=1}^{+\infty}\delta(t-T_j^l) + \frac{S_E}{N}\sum_{k=1, k\neq j}^{N}\sum_{l=1}^{+\infty} \delta(t-t_k^l),
\end{align*} 
where $V_E$ is the excitatory reversal potential, $\tau$ is the leakage time constant, $\sigma_E$ the time constant of the synaptic receptors, $S_{\mathrm{ext}}$ the strength of connections from outside the network, $S_E$ the strength of synapses inside the network, $T_j^l$ is the time of the $l^{th}$ exterior spike coming on neuron $j$ and $t_k^l$ is the time of the $l^{th}$ spike coming from the neuron $k$ onto the neuron $j$. We can compare with the simpler Lapicque model \eqref{voltage1} that was used in \cite{BH} to derive the model \eqref{NNLIF}: here the incoming current is replaced by a more complex dynamics involving the conductance of the cell membrane.

When neurons reach the voltage firing threshold $V_F$, which is $V(t^-) = V_F$, they are reset to the reset potential $V_R$, which is $(V(t^+),G(t^+))=(V_R,G(t^-))$. From each individual neuron $k$, the output train statistics are not Poisson in general. However, the inputs \textit{received} from the network by each neuron $j$ is a spike train summed over many output spike trains from the other neurons. If we assume that these incoming spikes are asymptotically independent from each other, the total input from the network can be assumed to be Poisson \cite{CA} or \cite[pp. 549--606]{L2}. Applying this simplification and considering the descrete probability density
\[  \rho^N(v,g,t) = \mathbb E \left[ \dfrac{1}{N} \sum_{j=1}^{N} \delta\big(v =V^j(t)\big)\delta\big(g =G^j(t)\big) \right], \] 
where $\mathbb E $ is the expectation with respect to all possible sets of $N$ independent realisations of the exterior Poisson spikes $(T_j^l)_{j\in\{1,\dots,N\},l\in\N}$, \cite{CTSM,CTRM} derive by a moment closure approach a $(2+1)-$dimensional kinetic Fokker-Planck equation for the density $\rho(v,g,t)$ in the limit $N\to+\infty$. Up to rescaling and assuming that
$ 0 = V_R = V_L <  V_F <  V_E$,  
it writes
\begin{equation}\label{eq:Kinetic}
	\dfrac{\partial \rho}{\partial t} + \dfrac{\partial }{\partial v} \left[ J_v(v,g) \rho(v,g,t) \right] + \dfrac{\partial }{\partial g}\left[ J_g(g, \mathcal N(t),t) \rho(v,g,t) \right] - \dfrac{a(\mathcal N(t),t)}{\sigma_E} \dfrac{\partial^2 \rho }{\partial g^2} = 0,
\end{equation}
where the fluxes are defined by
\begin{equation}\label{eq:Kineticflux}
	J_v(V_F,t) = -g_L v + g(V_E - v), \qquad J_g(g, \mathcal N,t) = \dfrac{1}{\sigma_E}\left( S_{\mathrm{ext}} \nu(t) + S_E \mathcal N - g \right),
\end{equation}
the firing rate $\mathcal N(t)$ and the intensity of the synaptic noise $a(\mathcal N(t),t)$ are characterised by
\begin{equation}\label{eq:KineticN}
	\mathcal N(t) = \int_0^{+\infty} N(g,t) \diff g, \quad N(g,t) = J_v(V_F,g) \rho(V_F,g,t),
\end{equation}
\begin{equation}\label{eq:KineticA}
	a(\mathcal N(t),t) =\dfrac{1}{2\sigma_E} \left(S_{\mathrm{ext}}^2\nu(t) + \dfrac{S_E}{N_E} \mathcal N(t) \right) .
\end{equation}
and the system is supplemented with initial and boundary constraints
\begin{equation}\label{eq:KineticBC}
	J_g (0,\mathcal N(t),t) \rho(v,0,t) - \dfrac{a(\mathcal N(t),t)}{\sigma_E} \dfrac{\partial \rho}{\partial g}(v,0,t) = 0\qquad \mathrm{and} \qquad \mathrm{BC\ \eqref{eq:BCv}}
\end{equation}
and 
\begin{equation}\label{eq:KineticIC}
	\rho(v,g,0)=\rho^0(v,g) \geqslant 0, \qquad \int_{0}^{V_F}\int_0^{+\infty} \rho^0(v,g)\diff g\diff v = 1.
\end{equation}
Let us comment on the boundary conditions \eqref{eq:KineticBC}. The first part is necessary for the conductance to stay positive: we impose a dynamic no-flux condition on $g=0$. We could also impose this condition at $g=+\infty$, but in order to avoid technical complications, we will simply assume fast decay at $+\infty$ for $\rho$ and $\partial_v \rho$. The second boundary condition \eqref{eq:BCv}, that we have not described yet, is more tricky. There is no diffusion in $v$, so the fluxes entering or exiting through the boundary at $v=0$ or $v=V_F$ depend only on $ J_v \rho$. The model prescribes that nothing comes in the domain from $\{V_F\}\times (0,+\infty)$ and anything that goes out through this boundary line must come in through the line  $\{ 0 \}\times (0,+\infty)$. Looking at the drift $J_v$, we see that it vanishes on a curve\footnote{This will be an issue for mathematical analysis later.} 
\[g = \dfrac{g_L v}{V_E-v}\]
and is always non-negative when 
\[g\geqslant g_F := \dfrac{g_L V_F}{V_E-V_F}.
\]
Below this critical conductance ($g\leqslant g_F$), the fluxes are coming inside the domain from both sides, so they are set to 0 (if nothing comes out at $v=V_F$, nothing should come in at $v=0$). Hence, we impose the following flux conditions in the $v$ direction:
\begin{equation}\label{eq:BCv}
	\left\{\begin{array}{ll}
		\rho(v,g,t) = \rho(V_F,g,t) = 0,\qquad& g\in (0,g_F], \\
		J_v(0,g) \rho(0,g,t) = J_v(V_F,g)\rho(V_F,g,t)=N(g,t),\qquad \qquad& g\in(g_F,+\infty).
	\end{array}\right.
\end{equation}

Given the profusion of parameters, let us also recapitulate about their meaning.
\begin{itemize}
	\item The exterior environment (other networks, visual stimuli, etc) acts through the given function $\nu(t)>  0$, that we will later take as constant, amplified by the synaptic strength $S_{ext} > 0$.
	\item The internal parameters are the average synaptic strength $S_E > 0$, the normalisation factor $N_E$, the leak conductance $g_L$ and the decay rate $\sigma_E$ for the conductance.
    \item The neurons are driven towards the leak potential $V_L=0$ with conductance $g_L$ and towards the excitatory reversal potential $V_E > V_F$ with conductance $g$; when a neuron reaches the potential $V_F$ on its way to $V_E$, it fires (see the relaxed model \eqref{eq:VC_var_ML} below for a smoother approach). 
\end{itemize}
Although we will not discuss it further, note that \cite[Section 5]{CTRM} proposes a model for excitatory-inhibitory networks. It consists in introducing two separate conductance variables $g_E$ and $g_I$, two excitatory and inhibitory reversal potentials $V_I < V_E$ and a $v$-flux defined as
\[  J_v(v,g_E,g_I) = g_L (V_L-v) + g_E(V_E-v) + g_I(V_I-v). \]

\subsection{Notation for mathematical analysis}\label{subsec:not}

Consider, performing renormalisations and taking constant external inputs, the fluxes and noise term
\[  J_v(v,g) = -g_Lv+g(V_E-v),\quad J_g(g,\mathcal N) = \dfrac{ s_0 + S_E \mathcal N(t)-g}{\sigma_E},\quad a(\mathcal N) = \dfrac{ a_0 +  a_1 S_E^2\, \mathcal N}{\sigma_E} ,  \]
with $g_L,V_E,s_0,a_0, \sigma_E>0$ positive parameters and $a_1,S_E>0$ non-negative connectivity coefficients. Recall the crucial $g$ threshold
\[g_F := \dfrac{g_L V_F}{V_E-V_F},\]
where the boundary condition switches. We can sum up the Voltage-Conductance kinetic model as
\begin{equation}\label{eq:VC}\tag{V-C}
    \left\{
    \begin{array}{l}
         \dfrac{\partial \rho}{\partial t} + \dfrac{\partial }{\partial v} \Big[ \big(-g_Lv+g(V_E-v)\big) \rho \Big] + \dfrac{\partial }{\partial g}\left[ \dfrac{ s_0 + S_E \mathcal N(t)-g}{\sigma_E} \rho \right] - \dfrac{a(\mathcal N(t))}{\sigma_E} \dfrac{\partial^2 \rho }{\partial g^2} = 0, \\[0.3cm]
         \displaystyle \mathcal N(t) = \int_0^{+\infty} N(g,t) \diff g, \quad N(g,t) = J_v(V_F,g) \rho(V_F,g,t),\\[0.3cm]
        \rho(0,g,t) = \rho(V_F,g,t) = 0,\quad g\leqslant g_F, \quad
		J_v(0,g) \rho(0,g,t) = J_v(V_F,g)\rho(V_F,g,t), \\[0.3cm]
        J_g (0,\mathcal N(t)) \rho(v,0,t) - \dfrac{a(\mathcal N(t))}{\sigma_E} \dfrac{\partial \rho}{\partial g}(v,0,t) = 0,\\[0.3cm]
        \displaystyle\rho(v,g,0)=\rho^0(v,g) \geqslant 0, \qquad \int_{0}^{V_F}\int_0^{+\infty} \rho^0(v,g)\diff g\diff v = 1.
    \end{array}
    \right.
\end{equation}

In most of the results about \eqref{eq:VC}, the solutions have to be understood in the distributional sense, or in the following weak sense first introduced in \cite{DKXZ}: consider the test function space
\[  \Phi = \left\{ \phi \in C_c^\infty([0,V_L]\times\R_+\times\R_+)\ \ | \ \ \phi(0,g,t) = \phi(V_F,g,t),\ \dfrac{\partial g}{\partial g}(v,0,t) = 0 \right\}.  \]
then a weak-solution $\rho\in L^\infty(\R_+,L^1_+([0,V_L]\times\R_+))$ with firing rate $\mathcal N \in L^1_{\mathrm{loc}}(R_+)$ satisfies
\begin{multline*}
    \int_{0}^{V_F}\int_{-\infty}^{+\infty} \phi(v,g,T) \rho(v,g,T) \diff v\diff g = \int_{0}^{V_F}\int_{-\infty}^{+\infty} \phi(v,g,0) \rho(v,g,0) \diff v\diff g\\
    + \int_0^T\int_{0}^{V_F}\int_{-\infty}^{+\infty}\left[\dfrac{\partial \phi}{\partial t} + J_v\dfrac{\partial \phi}{\partial v} + J_g(g,\mathcal N(t))\dfrac{\partial \phi} {\partial g} + \frac{a(\mathcal N(t))}{\sigma_E}\dfrac{\partial \phi}{\partial g^2}\right] \rho(v,g,t) \diff v\diff g \diff t.
\end{multline*}
Note the that the test functions can have any values $\phi(v,0,t)$, which is crucial to ensure the boundary condition $J_g (0,\mathcal N(t)) \rho(v,0,t) - \tfrac{a(\mathcal N(t))}{\sigma_E} \partial_g \rho(v,0,t) = 0$.

Most of the computations in articles about \eqref{eq:VC}, for example in \cite{PS, PS2, KPS,DPSZ}, are done in the sense of distributions but would work in function spaces with this notion of weak solutions. It relieves a large part of the notational burden to just integrate and differentiate everything \textit{a priori} and then prove \textit{a posteriori} higher $L^q$ or Sobolev regularity from the propagation of finiteness in norms.

\subsection{The linear Voltage-Conductance model}

Consider the simplified case of a linear equation $\eqref{eq:VC}$ by taking $S_E=0$; this case has been studied first in \cite{PS}, later improved in \cite{DPSZ} and recently generalised to weak interaction in \cite{FM}. Then, there is no interaction between neurons. Yet, because of the degenerate structure of the differential operator, the vanishing $v$-flux and the complex boundary conditions, the mathematical analysis of solutions requires some care.

Indeed, consider a stationary state $(\rho_\infty,\mathcal N_\infty)$ of the linear problem: contrary to the stationary states of \eqref{NNLIF}, \eqref{eqn:2} or \eqref{eq:T-E}, there is no semi-implicit formulation or parabolic regularity to ensure they are more than $L^1$ functions. They have to be constructed in a more abstract way, akin to the stationary solutions of \eqref{eq:jLIF}. As we will explain, it is possible to obtain more regularity for linear stationary states by carefully performed iteration methods \cite{PS,DPSZ}.

\begin{theorem}[Dou, Perthame, Salort, Zhou \cite{DPSZ}]\label{thm:VC_DPSZ1}
	Assume $S_E=0$. There exists a unique stationary state $\rho_\infty \in L^\infty([0,V_F]\times\R_+)$ to \eqref{eq:VC}.
\end{theorem}
Although we include it in the theorem for clarity, existence was supposed a priori in \cite{PS,DPSZ} and proved in the linear case first in \cite{DKXZ} and then in \cite{FM} which extends it to the weakly-nonlinear case. 

\subsubsection{\textit{A priori} properties of linear stationary states}

Denote
\[ a = \frac{a_0}{\sigma_E},\quad M(g) = e^{-\frac{(g-s_0)^2}{2a}},\quad Z = \int_{0}^{+\infty} M(g)\diff g. \]
The stationary equation for \eqref{eq:VC} can be rewritten
\[  \dfrac{\partial }{\partial v}[J_v(v,g)\, \rho_\infty(v,g)] - \dfrac{a}{\sigma_E} \dfrac{\partial}{\partial g}\left[ M(g) \dfrac{\partial}{\partial g}\left(\dfrac{\rho_\infty(v,g)}{M(g)} \right)  \right] = 0. \]

By direct integration and simple computations,
\begin{equation}
    \int_0^{V_F}\rho_\infty(v,g)\diff v = \dfrac{M(g)}{Z}.
\end{equation}

Then, by multiplying the equation by $\log \rho_\infty$ and integrating, it is possible to prove: 
\begin{lemma}[Perthame, Salort \cite{PS}]\label{lm:VC_PS}
	Assume $S_E=0$. There exist positive constants $K_0, K_1$ depending on $\mathcal N_\infty$ such that
	\begin{equation*}
		\int_{0}^{V_F}\int_{0}^{+\infty} \dfrac{1}{\rho(v,g)}\left(\dfrac{\partial \rho}{\partial g}(v,g)\right)^2\diff g\diff v \leqslant K_0,\qquad \int_{0}^{V_F}\int_{0}^{+\infty} e^{\frac{g^2}{8a}}\left|\dfrac{\partial \rho}{\partial g}(v,g)\right|\diff g\diff v \leqslant K_1. 
	\end{equation*}
\end{lemma}

From there, the strategy of \cite{PS} was to prove that $J_v \rho_\infty$ was lying in the Besov space
\[ B_\infty^{\frac12,1}((0,V_R)\times\R_+) = \left\{f\in L^1 \ | \sup_{h\leqslant 1} h^{-\frac12} \norme{f(\cdot+h)-f(\cdot)}_{L^1} \right\}.  \]
Then, by Besov injection, they could get $J_v \rho_\infty \in L^q$ for all $1\leqslant q\leqslant \frac43$. Unfortunately, because $J_v$ vanishes on a curve in the domain, they could not get better than $\rho_\infty\in L^q$ for $1\leqslant q < \frac87 $.

A significant improvement was proposed in \cite{DPSZ}. They bypass the Besov spaces by writting
\[\rho_\infty(v,g) = - \int_g^{+\infty}\dfrac{\partial \rho_\infty}{\partial g}\diff g \leqslant \int_0^{+\infty}\left|\dfrac{\partial \rho_\infty}{\partial g}(v,g')\right|\diff g',  \]
and performing the direct bound
\[  \int_0^{V_F} \sup_{g\in\R_+} \rho_\infty(v,g)\diff v \leqslant \int_0^{V_F} \int_0^{+\infty}\left|\dfrac{\partial \rho_\infty}{\partial g}(v,g')\right|\diff g' \leqslant K_1.\]
After more computations and using Lemma \ref{lm:VC_PS}, they build a set of useful estimates.
\begin{lemma}[Dou, Perthame, Salort, Zhou \cite{DPSZ}]
	Assume $S_E=0$. There is a constant $K^*$ such that
	\begin{equation}
		\norme{\rho_\infty}_{L^1([0,V_F];L^{\infty}(\R_+))} \leqslant K_1, \quad \norme{ J_v^2\, \rho_\infty }_{L^\infty([0,V_F]; L^1(\R_+)) }\leqslant K^*, \quad  \norme{J_v\, \rho_\infty}_{L^2} \leqslant K_1 K^*  
	\end{equation}
  and for all $q\in (0, \tfrac{4}{3})$,
  \begin{equation}\label{eq:VC_improve}
  	 \int_{0}^{V_F}\int_{0}^{+\infty} (g+g_L)^2 \rho(v,g)^q \diff g \diff v < +\infty.
  \end{equation}
\end{lemma}
The direct $\sup_g$ bound is possible due to the smoother parabolic $g$ structure of the equation, embedded in the constant $K_1$ of Lemma \ref{lm:VC_PS}. It is not possible to do the same in $v$, which is why \cite{DPSZ} has to derive the weaker $\sup_v$ bound
\[ \norme{ J_v^2\, \rho }_{L^\infty([0,V_F]; L^1(\R_+)) } = \sup_{v\in[0,V_F]}\int_0^{+\infty} J_v(v,g)^2\rho_\infty(v,g)\diff g \leqslant K^*,  \]
before getting to the $L^2$ bound on $J_v\, \rho_\infty$. Note that this constitutes a slight improvement from $J_v \rho_\infty \in L^{\frac43-}$ in \cite{PS} to $J_v \rho_\infty \in L^{2}$ in \cite{DPSZ}. After carefully overcoming the curve of zeros of $J_v$, it allows to go from $\rho_\infty \in L^{\frac87-}$ to $(g+g_L)\rho_\infty\in L^{\frac43-}$. 

This does not sound impressive, but the difference allows to make more progress. Using the Alikakos iteration method \cite{A3} starting from $L^{\frac43-}$ and using specific weights, \cite{DPSZ} manages to get to $\rho_\infty\in L^\infty$. Uniqueness of linear stationary states had already been obtained at hand in \cite{PS}.

\subsubsection{Asymptotic stability in the linear case \textit{via} an entropy method}

Let $H\in C^2$ a convex function. The entropy dissipation for solutions of \eqref{eq:VC} in the linear case $S_E=0$ writes
\begin{equation*}
    \dfrac{\diff}{\diff t} \int_0^{V_F}\int_0^{+\infty} H\left( \dfrac{\rho}{\rho_\infty} \right) \rho_\infty\diff v\diff g \leqslant - \int_0^{V_F}\int_0^{+\infty} H''\left(\dfrac{\rho}{\rho_\infty}\right) \left[\dfrac{\partial }{\partial g}\left(\dfrac{\rho}{\rho_\infty}\right)\right]^2\rho_\infty\diff v\diff g \leqslant 0.
\end{equation*}
Provided $\frac{\rho^0}{\rho_\infty}\in L^\infty$, following the method of \cite[Chapter 8]{P2} and using the regularity in Theorem \ref{thm:VC_DPSZ1}, one can deduce from the entropy dissipation a bound of the form
\begin{equation}\label{eq:VC_bound}  C^0_- \leqslant \rho(v,g,t)\leqslant C^0_+,   \end{equation}
with $C^0_-,C^0_+>0$ depending on $\frac{\rho^0}{\rho_\infty}$. The dissipation not being coercive, it is not possible to get exponential decay in relative entropy. Instead, \cite{DPSZ} uses a compactness argument by defining
\begin{equation*}
    p_k(v,g,t) = \rho(v,g,t+k),\quad h_k(v,g,t) = \frac{\rho(v,g,t+k)}{\rho_\infty(v,g,t+k)}.
\end{equation*}
By \eqref{eq:VC_bound}, $h_k$ is bounded in $L^\infty$. Thus, there is a subsequence $h_{k(n)}$ converging weakly towards some $h_\infty$ when $n\to+\infty$.

With the choice $H(x) = (x-1)^2$, the entropy dissipation allows first to prove that $h_\infty(v,g,t)$ does not depend on $g$, and then that $\partial_t h_\infty + \partial_v h\infty = 0$, implying by mass conservation that $h_\infty \equiv 1$. In turn, this allows to prove time compactness for $h_k$ and $p_k$, and to conclude using Aubin-Lions lemma and the structure of the equation that $p_k$ converges towards $\rho_\infty$ almost everywhere and in $L^q_{\mathrm{loc}}$, $1\leqslant q < +\infty$.

The last step is then to use these converging sequences to prove the decay of the relative entropy.

\begin{theorem}[Dou, Perthame, Salort, Zhou \cite{DPSZ}]
	Assume we are in the linear case $S_E=0$. Assume there exists $C>0$ such that $\rho^0\leqslant C \rho_{\infty}$. Then, for all time $t\geqslant 0$, $\rho(\cdot,\cdot,t)\leqslant C \rho_{\infty}$ and
	\[  \lim_{t\to +\infty} \int_{0}^{V_F} \int_{0}^{+\infty} \left( \dfrac{\rho(v,g,t)-\rho_\infty(v,g)}{\rho_\infty(v,g)}  \right)^2 \rho_\infty(v,g)  \diff g\diff v   = 0. \]
\end{theorem}

Note however that in the linear case, the entropy method does allow to prove exponential convergence to a stationary state if we first integrate the voltage variable. Denote \[M_\infty(g) = \frac{M(g)}{Z}, \qquad \varphi(g,t) = \int_0^{V_F} \rho(v,g,t)\diff v. \]
Then, \cite{PS} computes,
\begin{equation*}
    \dfrac{\diff}{\diff t} \int_0^{+\infty} \left(\frac{\varphi(g,t)-M_\infty(g)}{M_\infty(g)}\right)^2M_\infty(g)\diff g\leqslant - \frac{a}{\sigma_E} \int_0^{+\infty} \left[\dfrac{\partial}{\partial g}\left(\frac{\varphi(g,t)}{M_\infty(g)}\right)\right]^2 M_\infty(g)\diff g,
\end{equation*}
which by a Poincaré-like inequality and standard computations yields for some $\nu > 0$,
\begin{equation*}
    \int_0^{+\infty} \left(\frac{\varphi(g,t)-M(g)}{M_\infty(g)}\right)^2M_\infty(g)\diff g\leqslant e^{-\nu t}\int_0^{+\infty} \left(\frac{\varphi(g,0)-M(g)}{M_\infty(g)}\right)^2M_\infty(g)\diff g.
\end{equation*}

\subsubsection{Probabilistic reformulation and exponential ergodicity}
\label{sec:VC_Harris}

There are two issues in the approach from \cite{PS,DPSZ} described so far. First, the knowledge of the properties of the statonary state is \textit{a priori} and this is not enough to prove existence. Second, the convergence to a prospective stationary state is obtained by compactness, and there is thus no rate of convergence. To circumvent the difficulty of the degenerate diffusion in \eqref{eq:VC}, \cite{DKXZ} proposes to use a probabilistic reformulation (see also the similar approach in \cite{HMP} described below in Subsection \ref{sec:VC_ML}) in order to take advantage of the ergodicity of the microscopic trajectories.

Consider again $S_E=0$ and fix for the sake of simplicity $a_0 = \sigma_E = s_0 = V_F = 1$. Denote $(V_t,G_t)$ the solution of the system of stochastic differential equations
\begin{equation}\label{eq:VC_lin_sto}
    \left\{
    \begin{array}{l}
         \diff V_t = J_v(V_t,G_t)\diff t,\\
         \diff G_t = -(G_t-1)\diff t + \sqrt{2} \diff B_t + \diff L_t,\\
         V_{t^-} = 1 \implies V_t = 0,
    \end{array}
    \right.
\end{equation}
whith $(B_t)_{t\geqslant 0}$ a standard Brownian motion and $L_t$ the local time preventing the conductance $G_t$ from becoming negative (see \cite[Chapter 6]{RY} or the discussion on a similar issue in Section \ref{sec:GC}).

Since we are in the linear case, there is no feedback on $G_t$ from the behaviour of $V_t$. Hence, as we can see in \eqref{eq:VC_lin_sto}, the equation for $G_t$ is autonomous and constitutes a standard reflected Ornstein-Uhlenbeck SDE, whose solution exist and enjoy many good properties. Thus, following \cite{DIRT2}, \cite{DKXZ} considers the sequence of jump times of $V_t$,
\[ \tau_k = \inf \{ t>\tau_{k-1} \ | \ \lim_{s\to t^-} V_s = 1  \}, \]
and proves that $\mathbb P(\lim_{k\to+\infty} \tau_k) = 1$, which implies that there exists a trajectory-wise almost surely unique solution to \eqref{eq:VC_lin_sto}.

Then, \cite{DKXZ} proceeds to prove that the law of a solution to \eqref{eq:VC_lin_sto} is a weak solution to the linear \eqref{eq:VC}, a key difficulty being the boundary behaviour, they prove separately that
\[  \mathbb P(V_t=0)=\mathbb P(V_{t^-}=1) = \mathbb P(G_t=0) = 0. \]
Given the hypoellipticity of the linear operator defining the PDE and the stochastic process, there is enough interior regularity to ensure that the law of $(V_t,G_t)$ is absolutely continuous with respect to Lebesgue measure, which allows to build a link to weak solutions of the PDE using Itô's formula for semi-martingales with jumps \cite[Theorem 2.33]{P3}.

Last, they prove that the $v$-boundary condition is satisfied for weak solutions in the sense that for test functions $\psi\in C^\infty_c(\R_+^2)$,
\[  \lim_{v\to 0^+} \int_{-\infty}^{+\infty}\int_{-\infty}^{+\infty} J_v(v,g)\rho(v,g,t) \psi(g,t) \diff g\diff t = \lim_{v\to 1^-} \int_{-\infty}^{+\infty}\int_{-\infty}^{+\infty} J_v(v,g)\rho(v,g,t) \psi(g,t) \diff g\diff t .  \]

To prove the exponential convergence of the density associated to $(V_t,G_t)$, \cite{DKXZ} uses the Doeblin-Harris method. See our presentation of Doeblin's theorem in the context of the Time Elapsed neuron model in Subsection \ref{sec:TE_Doeblin}. For \eqref{eq:VC_lin_sto}, Doeblin theorem is not enough. Hence, like for example \cite{FS} for the system \eqref{eq:T-E_var_FS}, they use Harris theorem \cite{HM}, which requires Doeblin's minoration condition (Definition \ref{def:Doeblin}) on a compat subset only, but as a trade-off requires on top a Lyapunov control on the stochastic semigroup $(S_t)_{t\geqslant 0}$:
\[  S_t \omega \leqslant \alpha_1 \omega + \alpha_2,  \]
for some weight function $\omega(v,g)\geqslant 0$ and constants $\alpha_1\in (0,1)$ and $\alpha_2 \geqslant 0$. The idea is that the Lyapunov function sends back the mass towards the compact subset where the Doeblin condition is applied, which is exploiting the ergodicity of the process and the good properties at a location it visits by ergodicity. Given both Doeblin's condition on a compact subset and the Lyapunov control with weight $\omega$, Harris theorem says that the stochastic semigroup admits a unique invariant measure and that for some $\beta\geqslant 0$, for all initial probability densities, $S_t \rho^0$ converges exponentially fast to the invariant measure in the weighted total variation norm
\[ \norme{\mu}_{\mathrm{TV}(\beta\omega)} = \int_{[0,1)\times \R_+} (1+\beta\omega)\diff |\mu| . \]

The domain being bounded in the $v$ direction, \cite{DKXZ} constructs the Lyapunov weight $\omega(v,g) = \omega(g)$ as a function of $g$ only and proves that $\omega(g) = (g-1)^2$ satisfies the Lyapunov condition. For Doeblin's condition, in some area $g\leqslant R$, they exploit the interaction between the deterministic dynamical system, sending the particles towards the curve where $J_v(v,g)=0$, and the Brownian noise making them move away. A careful construction yields a measure $\nu$ satisfying the minoration condition.

\begin{theorem}[Dou, Kong, Xu, Zhou \cite{DKXZ}]
    Let $\omega(g) = (g-1)^2$. Assume the linear case $S_E=0$. There exists a unique invariant measure $\rho_\infty$ and constants $C\geqslant 1,\lambda>0,\beta>0$ such that for all $\rho^0\in \mathcal P([0,1)\times\R_+)$, there exists a unique solution to \eqref{eq:VC_lin_sto} associated for all $t>0$ to a density $\rho_t$ which is a weak solution to \eqref{eq:VC} and satisfies
    \begin{equation*}
        \forall t\geqslant 0,\qquad \norme{\rho_t-\rho_\infty}_{\mathrm{TV}(\beta\omega)} \leqslant C e^{-\lambda t} \norme{\rho_0-\rho_\infty}_{\mathrm{TV}(\beta\omega)} .
    \end{equation*}
    
\end{theorem}

\subsection{\textit{A priori} estimates on solutions in the nonlinear case}

\subsubsection{Existence of nonlinear stationary state}

Given the existence of a linear stationary state, it is possible to investigate the existence of nonlinear stationary state by a fixed point on the stationary firing rate $\mathcal N_\infty$. Consider the functional $\Psi:\R_+\to\R_+$ defined by
\begin{equation*}
    \mathcal N\mapsto \Psi(\mathcal N) = V_E \int_0^{+\infty} \rho_\infty^\mathcal N(0,g)\diff g,
\end{equation*}
where $\rho_\infty^\mathcal N$ is the unique stationary state from the linear system with $\mathcal N$ as a constant external input. It is continuous and $\Psi(0)>0$. Hence, to prove the existence of a fixed point (and hence of a nonlinear stationary firing rate yielding a stationary state) it is enough to find $\mathcal N_1$ such that $\Psi(\mathcal N_1) < \mathcal N_1$. To do so, \cite{PS} proves the bounds
\[  -V_F + (V_E-V_F) \left( g_\infty + \sqrt{\frac{a_\infty}{2\pi}}e^{-\frac{g_\infty^2}{2a_\infty}} \right) \leqslant \Psi(\mathcal N)\leqslant V_E\left(  g_\infty+ \sqrt{\frac{2a_\infty}{\pi}}e^{-\frac{g_\infty^2}{2a_\infty}} \right),  \]
with
 $g_\infty = s_0 + S_E \mathcal N$, $a_\infty = a_0 + S_E^2 \mathcal N$ .
Analysing the explicit $\mathcal N$ dependency in the bounds allows to prove
\begin{theorem}[Perthame, Salort \cite{PS}]
    If
    $ \frac{V_E}{V_F} S_E < 1$, 
    there exists at least one stationary state to \eqref{eq:VC}.
    However, if
    \[ \dfrac{V_E-V_F}{V_F} S_E>1,\qquad \mathrm{and} \qquad (V_E-V_F)s_0 > V_F^2, \]
    then there is no stationary state to \eqref{eq:VC}.
\end{theorem}
The $L^\infty$ regularity from Theorem \ref{thm:VC_DPSZ1} translates to these nonlinear stationary states, as they are linear stationary states with a modified $\tilde s_0 = s_0 + S_E \mathcal N_\infty$.

\subsubsection{Propagation of moments and global bounds on the firing rate}

For any $k\in\N$, define the moments
\[ \Psi(t) = \int_{0}^{V_F}\int_0^{+\infty} v\rho(v,g,t)\diff g\diff v ,  \]
\[ H_k(t) = \int_{0}^{V_F}\int_{0}^{+\infty} g^k \rho(v,g,t) \diff g\diff v\quad \mathrm{and}\quad f(t) = \int_{0}^{V_F} \rho(v,0,t)\diff v.\]
Then, providing these quantities exist, multiplying \eqref{eq:Kinetic} by the appropriate weight and integrating, yields the moment equations
\begin{equation}
	\begin{array}{rcl}
		\sigma_E\dfrac{\diff H_1}{\diff t}(t)&=&  - H_1(t) +  s_0 + S_E \mathcal N(t)  + a(t)f(t),\\[6pt]
		\sigma_E\dfrac{\diff H_2}{\diff t}(t)&=& - 2 H_2(t) + 2\left[ s_0 + S_E \mathcal N(t) \right]H_1(t) + 2a(t), \\[6pt]
		\sigma_E\dfrac{\diff H_k}{\diff t}(t)&=& - k H_k(t) + k \left[ s_0 + S_E \mathcal N(t) \right]h_{k-1}(t)  + k(k-1)a(t)h_{k-2}(t), \quad k\geqslant 3\\[6pt]
		\dfrac{\diff \Psi}{\diff t}(t)&=& \displaystyle - g_L \Psi(t) - \int_{0}^{V_F}\int_{0}^{+\infty} gv \rho(v,g,t)\diff g\diff v + V_E H_1(t) - V_F \mathcal N(t).
	\end{array}
\end{equation}

Note that since $\rho$ is a probability density, $H^1(t)^2 \leqslant H_2(t)$. Define
\[  \omega = \dfrac{1}{\sigma_E}\left(\dfrac{V_E}{ V_F}\left[ S_E + a_1 S_E^2  \right]  - 1 \right).\]
From this system of differential equations for the moments and using the signs, \cite{PS} proves that if the initial moments are finite, it is propagated in time for any solution of the nonlinear problem \eqref{eq:VC}.

\begin{theorem}[Perthame, Salort \cite{PS}]
	Let $\rho$ be a solution to \eqref{eq:VC}. Assume that $H_k(0) < +\infty$ for all $k\in \N$. Then there exists $C > 0$ such that
	\begin{equation}
		H_1(t) \leqslant C\max(1,e^{\omega t}),\qquad \int_0^t \mathcal N(s)\diff s \leqslant C\max(1,e^{\omega t}),
	\end{equation}
	and for all $T>0$ there exists $\tilde C(T)>0$ such that for all $t\in[0,T]$,
	\begin{equation}
		H_k(t) \leqslant \tilde C(T),\qquad \int_0^T\int_0^{+\infty} (1+g)^{k-1} N(g,t)\diff g\diff t \leqslant \tilde C(T).
	\end{equation}
\end{theorem}

This allows to prove that the firing rate is unbounded in the high connectivity regime. Using the moment equations, we can write
\[  \dfrac{\diff \Psi}{\diff t}(t) \geqslant (V_E-V_F)H^1(t)-g_L\Psi(t)-V_F\mathcal N,\qquad \dfrac{\diff H^1}{\diff t}(t) \geqslant \dfrac{1}{\sigma_E}\big( s_0 + S_E\mathcal N(t) - H^1(t)\big).  \]
Multiplying the first inequality by $\frac{S_E}{V_F\sigma_E}$ and merging the constants into $\bar C, \zeta>0$, it leads to
\[ \dfrac{\diff}{\diff t} \left( \dfrac{S_E}{V_F\sigma_E}\Psi + H^1\right) \geqslant  \zeta\left( \dfrac{S_E}{V_F\sigma_E}\Psi + H^1\right) - \bar C. \]
For large enough $H^1(0)$ it implies a uniform lower bound on $\tfrac{S_E}{V_F\sigma_E}\Psi + H^1$, which in turn yields exponential growth of $\mathcal N$ by re-injecting the lower bound in the differential inequality for $\Psi$.

\begin{theorem}[Perthame, Salort \cite{PS}]
	Assume that
	\begin{equation}
		\zeta := \dfrac{1}{\sigma}\left(\dfrac{S_E(V_E-V_F)}{V_F}-1\right) > 0.
	\end{equation}
	Then, if $H_1(0)$ is large enough, there exists $C>0$ such that for any solution to \eqref{eq:VC},
	\begin{equation}
		\int_0^t \mathcal N(s)\diff s \geqslant C e^{\zeta t}.
	\end{equation}
\end{theorem}


Given that several over neuroscience models like \eqref{NNLIF}, \eqref{eq:jLIF} or \eqref{eq:PC} exhibit finite time blow-up at the level of the firing rate, that reciprocally some of them are proved to have global-in-time solutions as long as the firing rate stays in $L^\infty_{\mathrm{loc}}$ \cite{CGGS,CRSS,CDRZ}, a natural question is whether or not the firing rate $\mathcal N(t)$ of \eqref{eq:VC} is locally uniformly bounded. The article \cite{PS} provides the following partial answer.

\begin{theorem}[Perthame, Salort \cite{PS}]
	Let $q\geqslant 2$, and $\varepsilon > 0$. If
    \[  \int_0^{V_F} \int_0^{+\infty} (1+g)^{2q+\varepsilon - 1}\rho^0(v,g)^q \diff v\diff g < +\infty ,\qquad \int_0^{V_F} \int_0^{+\infty} (1+g)^2\rho^0(v,g)\diff v\diff g <+\infty,  \]
    then for all $T\in\R_+$, any solution to \eqref{eq:VC} satisfies
	\begin{equation}
		\int_0^T \mathcal N(t)^q\diff t < +\infty.
	\end{equation}
\end{theorem}

This is a $L^q_{\mathrm{loc}}$ bound for $q\geqslant 2$ and not a $L^\infty_{\mathrm{loc}}$ bound, but it should be enough to prove that blow-up cannot happen, provided a suitable well-posedness framework. Note that \cite{RS} has proved that an $L^q$ bound, $q > 2$, on the firing rate $N(t)$ is enough\footnote{The proof in \cite{RS} is done for $q=3$ for the sake of clarity, but the argument works for any $q>2$, and does not work for $q=2$ or lower.} to ensure global existence of solutions to \eqref{NNLIF}.

\subsection{Existence theory in the weak interaction regime}

As we have seen, \cite{DKXZ} has provided a framework for existence of weak solutions to \eqref{eq:VC} in the linear case by relying upon a probabilistic setting for the corresponding system of SDEs. Until very recently, there was no result about existence of time-dependent solutions for the nonlinear case. In \cite{FM}, a partial answer is proposed: weak solutions are constructed in the case of weak interactions (small enough $S_E$).

\subsubsection{Ultracontractivity estimates}

One way to handle the structure of \eqref{eq:VC} is by using the hypocoercivity structure of the kinetic Fokker-Planck equation \cite{V,BCMT}. However, the article \cite{FM} goes another way: following \cite{CM,FGM}, they prove ultracontractivity estimates and rely upon the Doeblin-Harris method (see Subsections \ref{sec:TE_Doeblin} and \ref{sec:VC_Harris}), in order to construct a framework for well-posedness.

Consider, for external source terms $\alpha,\beta\in L^\infty(0,+\infty)$, the linear kinetic Voltage-Conductance equation
\begin{equation}\label{eq:VC_lin_FM}
    \dfrac{\partial \rho}{\partial t} = \mathscr L_{\alpha,\beta}\, \rho,\qquad \mathscr L_{\alpha,\beta}\, \rho = \dfrac{\partial }{\partial v}[J_v(v,g) \rho] + \dfrac{\partial }{\partial g}[ (\beta(t) - g)\rho ]  + \alpha(t)\dfrac{\partial^2 \rho}{\partial g^2},
\end{equation}
with boundary conditions as in \eqref{eq:VC}. Taking either $\omega(v,g) = (1+g^2)^\frac{k}{2}$, $k\in\N$, $k > 1$ or $\omega(v,g)=e^{\alpha g}$, $\alpha > 0$ as a weigth function, \cite{FM} formally integrates the equation and proves that there is another equivalent weigth $\tilde \omega$ and a constant $ \kappa\geqslant 0$ depending only on $\omega,\alpha,\beta$ satisfying the \textit{a priori} estimate:
\begin{equation}\label{eq:VC_MF_first}
       \norme{\rho(\cdot,\cdot,t)}_{L^p(\tilde \omega)} \leqslant e^{\kappa t} \norme{\rho^0}_{L^p(\tilde \omega)}. 
\end{equation}
Then, using the penalisation of the $v$-boundary $\gamma(v) = \min(v,V_F-v)$ and building step by step from the proof of \eqref{eq:VC_MF_first}, they prove that for the exponential weigth $\omega(v,g)=e^{\alpha g}$, $\alpha > 0$, solutions to the linear equation \eqref{eq:VC_lin_FM} satisfy
\begin{equation}\label{eq:VC_MF_second}
    \int_0^T \phi^2\int_{0}^{V_F} \int_{0}^{+\infty} \left[\dfrac{1}{\delta^{\frac1{27}}}+\sqrt{1+g^2}\right]\rho\, \omega\, \diff v\diff g\diff t \leqslant C \int_{0}^T (\phi\phi'_+ +\phi^2 )\int_{0}^{V_F} \int_{0}^{+\infty} \rho^2\omega^2 \diff t\diff v\diff g,
\end{equation}
with $C\geqslant 1$ a constant and any test function $\phi\in C^1_c((0,T))$.

From there, \cite{FM} considers $\bar \rho = \rho\, \phi(t) \chi(v) \omega(g)$, $\phi\in C^1_c((0,T))$, $\chi\in C^1_c((0,V_F))$ and performs a regularisation argument \cite[Appendix]{V} in order to get a bound of $\int_0^T \norme{\bar \rho}_{L^{\frac52}}^2\diff t$ in function of $\int_0^T\norme{\rho}_{L^2(\omega)}^2\diff t$ and the regularisation functions. Combining with \eqref{eq:VC_MF_second}, they conclude that there exists $r>2$ such that for any $\phi\in C^1((0,T))\cap C_{\to 0} ((0,T))$, $T\in(0,1)$,
\begin{equation}\label{eq:VC_MF_third}
    \int_0^T \phi^2 \norme{\rho(\cdot,\cdot,t)}_{L^r(\omega)} \diff t \leqslant \tilde C \norme{(\phi+|\phi'|)\rho\omega}_{L^2([0,V_F]\times\R_+\times\R_+)}.
\end{equation}
Using Nash's argument \cite{N}, it leads to the bound
\[ \forall t\in(0,1),\qquad \norme{\rho(\cdot,\cdot,t)}_{L^2(\omega)}\leqslant \frac{1}{t^{\nu_1}} \norme{\rho^0}_{L^1(\omega)},  \]
for some $\nu_1>0$. Last, they build from this result using the Nash-Moser iteration procedure and estimate \eqref{eq:VC_MF_first} in order to get a full ultracontractivity estimate for the linear system \eqref{eq:VC_lin_FM}.

\begin{theorem}[Fonte Sanchez, Mischler \cite{FM}]\label{thm:VC_ultracontractive}
    Assume $a_1=1$. Let $\omega(v,g)=e^{\alpha g}$, $\alpha > 0$ and $\rho^0\in L^2(\omega)$. Then there exists $C\geqslant 1,\nu > 0$ and $\kappa\in[0, C(1+\norme{\alpha}_\infty+\norme{\beta}_\infty)]$ such that any solution $\rho$ to \eqref{eq:VC_lin_FM} satisfies the ultracontractivity estimate 
    \[  \norme{\rho(\cdot,\cdot,t)}_{L^\infty(\omega)} \leqslant C \dfrac{e^{\kappa(t-s)}}{(t-s)^\nu} \norme{\rho(\cdot,\cdot,s)}_{L^1(\omega)}. \]
\end{theorem}

\subsubsection{Well-posedness in a weighted $L^\infty$ framework}

Adapting classical methods from kinetic theory, \cite{FM} proves that the linear system \eqref{eq:VC_lin_FM} with external inputs $\alpha,\beta$ admits for all $T>0$ a unique weak-solution $\rho\in C([0,T],L^2(\omega))\cap L^2((0,T)\times(0,V_F); H^1(0,+\infty))$ and that this solution admits a trace $\gamma(\rho)\in L^2(\Gamma_{v,t}\,,\,\omega^2 (1+g^2)^{-\frac12}J_v^2)$ on the $(v,t)$-boundary $\Gamma_{v,t} = (0,T)\times (\{0\}\cup\{V_F\})\times \R_+$. Using a kinetic Harnack inequality in the spirit of \cite{GIMV}, those solutions to the linear problem are also proved to enjoy a higher $C^\eta$ regularity in the strict interior of the time-voltage-conductance domain, for some $\eta > 0$.

All this work allows to prove that for small enough $S_E$ with respect to $\rho^0$, the mapping that takes an external firing rate $\mathcal N_1$ defining $\alpha,\beta$ and returns the firing rate $\mathcal N_2$ defined from the trace on the $(v,t)$-boundary of $\rho$ solution to \eqref{eq:VC_lin_FM} enjoys enough regularity to apply to Schauder-Tychonov fixed-point theorem. This provides a solution to the nonlinear problem.

\begin{theorem}[Fonte Sanchez, Mischler \cite{FM}]
    Assume $a_1=1$. Let $\omega(v,g) = (1+g^2)^\frac{k}{2}$, $k\in\N$, $k > 1$ or $\omega(v,g)=e^{\alpha g}$, $\alpha > 0$. There exists $\eta(\omega) > 0$ such that for all initial condition $\rho^0\in L^\infty(\omega)$, if $(S_E+S_E^2)\norme{\rho^0}_{L^\infty(\omega)} < \eta$, there exists at least one solution
    \[  \rho \in L^\infty(\R_+,L^\infty(\omega))\cap C(\R_+,L^2(\omega)) \]
    to \eqref{eq:VC} which conserves mass and positivity.
\end{theorem}

Last, using a reformulation of the Doeblin-Harris method \cite{HM} in a Banach lattice, and taking advantage of the ultracontractivity estimate in Theorem \ref{thm:VC_ultracontractive}, they prove the best result known so far on the existence and regularity of nonlinear stationary states.

\begin{theorem}[Fonte  Sanchez, Mischler \cite{FM}]
    Assume $a_1=1$. Let $\omega_0(v,g) = (1+g^2)^\frac{k}{2}$, $k\in\N$, $k > 1$. There exists $\eta(\omega_0) > 0$ such that if $S_E < \eta(\omega_0)\norme{\omega_0^{-1}}_{L^1} $, for either $\omega=\omega_0$ or $\omega(v,g)=e^{\alpha g}$, $\alpha > 0$, there exists a unique stationary state
    \[  \rho_\infty\in L^\infty(\omega)\cap C((0,V_F)\times (0,+\infty)).  \]
    It is linearly exponentially asymptotically stable in $L^p(\omega)$ for all $p\in[1,+\infty]$.
\end{theorem}
Like in many other results obtained from Doeblin-Harris approaches, the constants are all constructive in the linear exponential contraction estimates.

\subsection{Numerical simulations of the Voltage-Conductance equation}

Numerical methods have also been introduced in order to explore the different features of the voltage-conductance kinetic equations \eqref{eq:VC}. In \cite{CCT} a numerical method based on Chang-Cooper schemes as in \cite{BCD}, see subsection \ref{subsec:numnnlif}, was developed and validated against Monte Carlo simulations. These numerical methods were used to showcase the appearance of bi-stability in this model showing the appearance of two steady states with low and high activity respectively with the appearance of a phase transition leading to complicated dynamics. They also observe the existence of periodic solutions or damped oscillations in the firing rate depending on the values of $S_E$. It is an open problem to show the existence of periodic solutions in this model even if numerically observed by multiple authors using different numerical approaches. We refer to \cite{CCT} for more details on the numerical experiments and particular values in the parameter regime leading to periodic solutions. Finally, the authors look at the validity limits of moment closures of the kinetic-conductance model introduced in \cite{RCT}. This moment-closure system reads as
\begin{eqnarray}
\partial_t \rho_v(t,v)& = &
\partial_v\left\{
\left[ \left(\frac{v-V_R}{\tau} \right) +
\mu_1(t,v)\left(\frac{v-V_E}{\tau} \right) \right] \rho_v(t,v)
\right\} \label{ke-v1}
\\
\partial_t\mu_1(t,v) &  = &
-\frac{1}{\sigma_E}\left[ \mu_1(t,v)-s_0 - S_E \mathcal N(t) \right] +
\frac{a(\mathcal N(t))}{\sigma_E\rho_v(t,v)}
\partial_v \left[
\left(\frac{v-V_E}{\tau} \right) \rho_v(t,v) \right] \nonumber
\\
& + & \left[ \left(\frac{v-V_R}{\tau} \right)
 +
\mu_1(t,v) \left(\frac{ v-V_E}{\tau} \right) \right]
\partial_v \mu_1(t,v).
\label{ke-v2}
\end{eqnarray}
where $\rho_v(t,v)$ and $\rho_v(t,v) \mu_1(t,v)$ are the first two moments in $g$ of the kinetic distribution $\rho(t,v,g)$. This closure assumption was numerically verified to be valid in steady states away from the boundaries in the voltage variable.

\subsection{The random firing voltage-conductance model and one variable reductions}

Because of the difficulty of the complex $v$-boundary condition and the singular firing at $V_F$, a strategy to simplify the study of the model \eqref{eq:VC} is to consider a softer firing mechanism \cite{PS2,KPS}. Instead of letting the neurons fire as they cross the threshold potential $V_F$, they can be allowed to go into the space $[V_F,V_E)$, with again $V_E$ the excitatory reversal potential where the flux $J_v$ is always inward in the domain. The neurons with voltage $V_F \leqslant v < V_E$ will then fire randomly with high probability $\phi_F(v)$, which translates into a removal term $\phi_F(v)\rho(v,g,t)$ in the deterministic framework. The firing rate is then the quantity of neurons over all the $(v,g)$ domain which are removed from the equation:
\[  \mathcal N(t) = \int_0^{+\infty}N(g,t)\diff g = \int_0^{+\infty} \int_0^{V_E} \Phi_F(v)\rho(v,g,t) \diff v\diff g. \]
Note that this approach has also been implemented into the model \ref{NNLIF} \cite{CP,DPSZ2} (see also Subsection \ref{sec:NNLIF_RF}).
The model writes
\begin{equation}\label{eq:VC_RF}
    \left\{
    \begin{array}{l}
         \dfrac{\partial \rho}{\partial t} + \dfrac{\partial }{\partial v} \Big[ \big(-g_Lv+g(V_E-v)\big) \rho \Big] + \dfrac{\partial }{\partial g}\left[ \dfrac{ s_0 + S_E \mathcal N(t)-g}{\sigma_E} \rho \right] - \dfrac{a}{\sigma_E} \dfrac{\partial^2 \rho }{\partial v^2} +\phi_F \rho= 0, \\[0.3cm]
         \displaystyle \mathcal N(t) = \int_0^{+\infty} N(g,t) \diff g, \quad N(g,t) = \int_0^{V_E}\phi_F(v) \rho(v,g,t)\diff v,\\[0.3cm]
        J_v(0,g) \rho(0,g,t) = N(g,t), \quad \rho(V_E,g,t) = 0, \quad
        J_g (0,\mathcal N(t)) \rho(v,0,t) = \dfrac{a}{\sigma_E} \dfrac{\partial \rho}{\partial g}(v,0,t) ,\\[0.3cm]
    \displaystyle\rho(v,g,0)=\rho^0(v,g) \geqslant 0, \qquad \int_{0}^{V_E}\int_0^{+\infty} \rho^0(v,g)\diff g\diff v = 1,
    \end{array}
    \right.
\end{equation}
where we simplify $a(\mathcal N) \equiv a>0$ and still denote
\[ J_v(v,g) = -g_Lv+g(V_E-v),\qquad J_g(g,\mathcal N) = \dfrac{ s_0 + S_E \mathcal N(t)-g}{\sigma_E}.\]
Note that we do not have anymore the complex switch at $g=g_F$ on the boundary $v=V_F$: it is replaced by a homogeneous Dirichlet boundary condition at $v=V_E$. The firing function is chosen as $\phi_F\in W^{1,+\infty}(0,V_E)$ and such that $\phi_F\geqslant 0$, $\phi_F'\geqslant 0$.

\subsubsection{Fast conductance limit}
\label{sec:VC_RF}

In \cite{PS2}, the Voltage-Conductance model with random firing \eqref{eq:VC_RF} is used to derive rigorously a voltage-only model through the limit of fast conductance, which corresponds in our setting to the limit of vanishing conductance time constant $\sigma_E\to 0$.

Let us denote $\rho^\varepsilon$ a solution obtained for $\sigma_E = \varepsilon>0$. The main idea is that when $\varepsilon$ goes to $0$, in the formal limit $\rho^\varepsilon$ converges to a solution of
\[  \dfrac{\partial }{\partial g}\Big[(s_0 + S_E\mathcal N(t) - g)\rho] - a \dfrac{\partial \rho}{\partial g^2}=0,\qquad  (s_0 + S_E\mathcal N(t))\rho(v,0,t) - a\dfrac{\partial \rho}{\partial g}(v,0,t)=0. \]
The solution can then be written in terms of a $v$-dependent part and a normalised Maxwellian
\begin{equation}
    \rho(v,g,t) = n(v,t) M(g,\mathcal N(t)),\qquad M(v,\mathcal N(t))= \frac{e^{-\frac{(s_0+S_E\mathcal N(t) - g)^2}{2a}} }{\int_0^{+\infty} e^{-\frac{(s_0+S_E\mathcal N(t) - g)^2}{2a}} \diff g}.
\end{equation}
The quantity $n(v,t)$ should then be the limit of
\[  n^\varepsilon(v,t) = \int_0^{+\infty} \rho^\varepsilon(v,g,t)\diff g.  \]
Integrating \eqref{eq:VC_RF}, we find that this quantity $n^\varepsilon$ is solution to
\begin{equation}\label{eq:VC_RF_eps}
    \left\{
    \begin{array}{l}
         \displaystyle\dfrac{\partial n^\varepsilon}{\partial t} + \dfrac{\partial }{\partial v} \left[ \left(-g_Lv+(V_E-v)\int_0^{+\infty}g \frac{\rho^\varepsilon(v,g,t)}{n^\varepsilon(v,t)}\diff g\right)\, n^\varepsilon(v,t) \right]  + \phi_F\, n^\varepsilon(v,t)= 0, \\[0.3cm]
         \displaystyle \mathcal N(t) = \int_0^{V_E}\phi_F(v) \, n^\varepsilon(v,t)\diff v,\\[0.3cm]
        \displaystyle  V_E\int_0^{+\infty}g \frac{\rho^\varepsilon(v,g,t)}{n^\varepsilon(v,t)}\diff g\, n^\varepsilon(0,t) =\mathcal  N(t), \quad n^\varepsilon(V_E,t) = 0.
        \end{array}
    \right.
\end{equation}
The idea of \cite{PS2} for performing a rigourous slow-fast limit in \eqref{eq:VC_RF_eps} is to rewrite the $v$-drift as
\begin{equation}
    -g_Lv+(V_E-v)\int_0^{+\infty}g \frac{\rho^\varepsilon(v,g,t)}{n^\varepsilon(v,t)}\diff g = G_\varepsilon(v,t) (V_\varepsilon(v,t) - v),
\end{equation}
with
\begin{equation}
    G_\varepsilon = \int_0^{+\infty} (g_L+g)\frac{\rho^\varepsilon(v,g,t)}{n^\varepsilon(v,t)}\diff g \diff g,\qquad V_\varepsilon = \dfrac{1}{G_\varepsilon(v,t)}\int_0^{+\infty} g V_E \frac{\rho^\varepsilon(v,g,t)}{n^\varepsilon(v,t)}\diff g.
\end{equation}
The rewritting $G_\varepsilon (V_\varepsilon - v)$ can then be seen as a conductance acting on the difference between the current potential and an equilibrium potential $g_{\mathrm{eq}}(t)(v_{\mathrm{eq}}(t)-v)$, and the idea is to take advantage of the inequalities $0 < V_\varepsilon \leqslant V_E$ and $G_\varepsilon \geqslant g_L$.

The quotient $\rho^\varepsilon/n$ should then converge towards $M(v,\mathcal N(t))$. In the limit $\varepsilon\to0$, \cite{PS2} proves that
\[ \mathcal N_\varepsilon \to N(t),\qquad  G_\varepsilon \to G := g_L + \int_0^{+\infty} g M(g,N(t))\diff g,\qquad V_\varepsilon \to \frac{1}{G} V_E \int_0^{+\infty} g M(g,N(t))\diff g.   \]
The limit equation will then be
\begin{equation}\label{eq:VC_RF_lim_1}
    \left\{
    \begin{array}{l}
         \displaystyle\dfrac{\partial n}{\partial t} + \dfrac{\partial }{\partial v} \left[ \left(-g_Lv+(V_E-v)\int_0^{+\infty}g M(v,N(t))\diff g\right)\, n(v,t) \right]  +\phi_F\, n(v,t)= 0, \\[0.3cm]
         \displaystyle N(t) = \int_0^{V_E}\phi_F(v) \, n(v,t)\diff v,\\[0.3cm]
        \displaystyle  V_E\int_0^{+\infty}g M(v,N(t))\diff g\, n(v,t) =  N(t), \qquad n(V_E,t) = 0,\\[0.3cm]
        \displaystyle n(v,0) = n^0(v) := \int_0^{+\infty} \rho^0(v,g)\diff g.
        \end{array}
    \right.
\end{equation}
Note that it is possible to compute explicitly
\[ \int_0^{+\infty} g M(g,N(t))\diff g = \dfrac{\sqrt{a} e^{- \frac{N(t)^2}{2a}}}{\sqrt{\frac{\pi}{2}}\left(1+\mathrm{erf}\left(\frac{s_0+N(t)}{\sqrt{2a}}\right)\right)} .  \]
\begin{theorem}[Perthame, Salort \cite{PS2}]
    Assume for all $k\in\N$, $\int_0^{+\infty}\int_0^{V_E}g^k\rho^0\diff g\diff v < +\infty$, $\phi_F\in W_+^{1,+\infty}(0,V_E)$. Then, there exists $C(k),C>0$ such that for all $\varepsilon>0$, for all $t\geqslant 0$,
    \[ \int_0^{+\infty}\int_0^{V_E} g^k\rho^\varepsilon(v,g,t)\diff g\diff v \leqslant C(k),\qquad \norme{\mathcal N}_\infty + \norme{\mathcal N'}_\infty \leqslant C.  \]
    In the weak topology of bounded measures,
    $  \rho^\varepsilon \rightharpoonup n(v,t) M(g,\mathcal N(t))$,
    where $n$ is solution to \eqref{eq:VC_RF_lim_1}
\end{theorem}
For the solutions to the limit equation \eqref{eq:VC_RF_lim_1}, \cite{PS2} also proves the $L^\infty_{\mathrm{loc}}$ bound
\[  \sup_{t\in [0,T]} \norme{n(\cdot,t)}_\infty \leqslant C(T) \norme{n^0}_\infty.  \]

Note that in the formal low-noise limit $a\to 0$, the Maxwellian converges to a Dirac mass and its average to an affine function of the firing rate:
\[  M(g,N(t))\to \delta_{s_0+S_E N(t)}, \qquad  \int_0^{+\infty} g M(g,N(t))\diff g\to s_0 + S_E N(t).  \]
The main equation of \eqref{eq:VC_RF_lim_1} then becomes
\[ \dfrac{\partial n}{\partial t} + \dfrac{\partial }{\partial v} \left[ \left(-g_Lv+(s_0+S_E N(t))(V_E-v)\right)\, n(v,t) \right]  +\phi_F(v)\, n(v,t) = 0, \]
which is almost identical to \eqref{LIFRand} studied in \cite{CTV,CTV2,PS2} and obtained as a mean-field limit in \cite{FL,DGLP}.

\subsubsection{Fast voltage limit}

If we now consider the case where it is not the conductance which evolves way faster but the voltage, then the solution should concentrate in $v$ around the curve of solutions to $J_v = 0$, which means that $v$ will converge to
\[  V^*(g) = \dfrac{g V_E}{g+g_L}. \]

Consider for simplicity $\sigma_E=0$. The fast firing regime consists in replacing the first equation in \eqref{eq:VC_RF} by
\begin{equation}
     \dfrac{\partial \rho}{\partial t} + \dfrac{1}{\varepsilon}\dfrac{\partial }{\partial v} \Big[ \big(-g_Lv+g(V_E-v)\big) \rho \Big] + \dfrac{\partial }{\partial g}\left[ (s_0 + S_E \mathcal N(t)-g) \rho \right] - a_0 \dfrac{\partial^2 \rho }{\partial v^2} +\phi_F \rho= 0,
\end{equation}
with $\varepsilon > 0$ the scaling parameter. The article \cite{KPS} thus proposes that in the infinitely fast voltage regime, the equation for
\[ n(g,t) = \int_0^{V_E} \rho(v,g,t)\diff v  \]
is the parabolic equation
\begin{equation}\label{eq:VC_RF_lim_2}
    \left\{
    \begin{array}{l}
         \dfrac{\partial n}{\partial t} + \dfrac{\partial }{\partial g}\left[ ( s_0 + S_E \mathcal N(t)-g) n \right] - a_0 \dfrac{\partial^2 n }{\partial g^2} = 0, \\[0.3cm]
         \displaystyle \mathcal N(t) = \int_0^{+\infty} \phi_F(V^*(g)) n(g,t)\diff g,\\[0.3cm]
        J_g (0,\mathcal N(t)) n(0,t) - a_0 \dfrac{\partial n}{\partial g}(0,t) = 0,\quad
    \displaystyle n(g,0)=n^0(g) :=\int_{0}^{V_E}\rho^0(v,g)\diff v.
    \end{array}
    \right.
\end{equation}
For technical reasons, the proof in \cite{KPS} is performed on a version of the equation with flux $J_g$ bounded from above uniformly. Using careful estimates on \eqref{eq:VC_RF}, first for the hyperbolic case $a=0$ and then for the case $a > 0$, they prove rigorously that the limit equation in the fast firing regime is an equation of the form \eqref{eq:VC_RF_lim_2}.

\subsection{Other variants of the voltage-conductance model}

Besides the random firing voltage-conductance model we have described, some other simplifications have been proposed in order to find which features are key for \eqref{eq:VC} to exhibit complex behaviours. Although we will not describe them in detail, note the 2-conductances model \cite{DZ2} and the kinetic Fokker-Planck equation with jumps \cite{SS2} which both share many similarities with \eqref{eq:VC}.

\subsubsection{Full understanding of a toy model and the importance of the voltage drift}

In order to understand self-organized oscillations in \eqref{eq:VC}, one can make two simplifications to isolate the difficulty of the firing rate being computed as a moment evaluated at a pointwise value of the voltage. One difficulty stems from the velocity field in $v$, i.e., $J_v(v,g)=-g_Lv+g(V_E-v)$ as in \eqref{eq:VC}. The velocity field can not be written in a separable form $J_v(v,g)=f(v)h(g)$, in contrast to the velocity field for the position variable in classical kinetic models. This non-separable velocity field prohibits the use of many analysis tools. To understand the consequences of a non-separable velocity field in $v$ and $g$, we consider a simpler velocity field in $v$ neglecting the leaky mechanism, $J_v(v,g)=g(V_E-v)$, which is further reduced to  $J_v(v,g)=g$ for the sake of simplicity, see \cite{CDZ} for further generalizations. Moreover, it leads to the complicated boundary condition \eqref{eq:BCv}, which makes it difficult to analyse even in the steady state case \cite{DPSZ}. Second, we extend the domain of $g$ from $\mathbb{R}^+$ to $\mathbb{R}$. 
Then we get the equation for the probability density function $p(t,v,g)$, denoted by $p$ in order to make the distinction to \eqref{eq:VC} apparent, as
\begin{equation}\label{eq:nonlinear-toy}
\partial_t p+g\partial_v p+\partial_g\bigl((g_{\text{in}}(t)-g)p-a(t)\partial_gp\bigr)=0,\quad v\in(0,V_F),\ g\in\mathbb{R},\ t>0,
\end{equation}
and the boundary condition in voltage $v$ is simplified to 
\begin{equation}\label{nonlinear-toy-bc}
p(t,0,g)=p(t,V_F,g),\quad \forall t>0,g\in\mathbb{R}.
\end{equation}
For the conductance $g$, we no longer need a boundary condition at $g=0$ thanks to the extension of the domain.
For the initial data, we assume it is a probability density function on $(0,V_F)\times(-\infty,+\infty)$, denoted as
\begin{equation}\label{nonlinear-toy-IC}
    p(0,v,g)=p^0(v,g),\quad v\in(0,V_F),g\in\mathbb{R}.
\end{equation}
The expressions of $g_{\text{in}}(t)$ and $a(t)$ depend on the firing rate $N(t)$ in the same way as in the original model given by
	\begin{equation}\label{eq:para-general}
	g_{\text{in}}(t)=g_0+g_1N(t),\ a(t)=a_0+a_1N(t),\quad g_0,g_1,a_0,a_1>0.
	\end{equation}
We define the firing rate $N(t)$ in a similar manner as the integration of flux at voltage $v=V_F$, over $g>0$, i.e.,
\begin{equation}\label{eq:nonliner-toy-fire}
	N(t):=\int_{0}^{\infty}gp(t,V_F,g)dg.
\end{equation}
	Though simplified, this model \eqref{eq:nonlinear-toy}-\eqref{eq:nonliner-toy-fire} inherits the nonlinear mechanism of the original model \eqref{eq:VC} -- dependence of $g_{\text{in}}(t),a(t)$ on  $N(t)$ as the firing rate $N(t)$ still only depends on the flux at one voltage value $v=V_F$.

The separability of the voltage flux makes possible an analysis via Fourier modes for classical solutions since the boundary conditions in $v$ now simplify to periodic boundary conditions. The main outcome is showing that the long time asymptotics are dictated by the zeroth-order model and then by homogeneous in $v$ solutions to \eqref{eq:nonlinear-toy}-\eqref{eq:nonliner-toy-fire}. Let us consider the Fourier expansion in $v$ direction
	\begin{equation}\label{fourier-v}
\begin{aligned}
 p(t,v,g)&=\frac{1}{V_F}\sum_{k=-\infty}^{+\infty}p_k(t,g)e^{ikv\frac{2\pi}{V_F}},\quad \text{where}\\
p_k(t,g):&=\int_{0}^{V_F}p(t,v,g)e^{-ikv\frac{2\pi}{V_F}}dv,\quad k\in\mathbb{Z}.
\end{aligned}
	\end{equation}

\begin{theorem}[Carrillo, Dou, Zhou \cite{CDZ}]\label{thm:decay-v}
		Suppose $p(t,v,g)$ is a {classical} solution of \eqref{eq:nonlinear-toy} and $p_0(t,g)$ is the zeroth mode in $v$ defined in \eqref{fourier-v}. Then $p(t,v,g)$ converges to the $v$-homogeneous mode $\frac{1}{V_F}p_0(t,g)$ exponentially. Precisely, we have
		\begin{equation}
		\begin{aligned}
		\left\|p(t,v,g)-\frac{1}{V_F}p_0(t,g)\right\|_{L^1((0,V_F)\times\mathbb{R})}&\leq 2\frac{e^{-\left(\tfrac{2\pi}{V_F}\right)^2d(t)}}{1-e^{-\left(\tfrac{2\pi}{V_F}\right)^2d(t)}},
		\end{aligned}
		\end{equation}
		where $d(t)$ is given by 
		\begin{equation}\label{def-dt}
		    d(t):=a_0\left(t-2\frac{e^t-1}{e^t+1}\right)> 0,\quad t>0.	
		\end{equation}
	\end{theorem}

 The second simplification extending the conductance to the whole line allows for a deep analysis of the $v$-homogeneous solutions to \eqref{eq:nonlinear-toy}-\eqref{eq:nonliner-toy-fire}. In fact, explicit computations show the existence of gaussian solutions with means and variances in $g$ depending only on time. The $v$-homogeneous solutions to \eqref{eq:nonlinear-toy}-\eqref{eq:nonliner-toy-fire} satisfy
\begin{equation}\label{eq:reduce-p}
\partial_t p_0=\partial_g\bigl((-g_{\text{in}}(t)+g)p_0+a(t)\partial_gp_0\bigr),\quad g\in\mathbb{R},\ t>0,
\end{equation} 
where $g_{\text{in}}(t)$ and $a(t)$ still given by \eqref{eq:para-general}. The $v$-homogeneous equation \eqref{eq:reduce-p} admits the following Gaussian type solution
 	\begin{equation}\label{eq:Gaussian anstaz}
 	p(t,g)=\frac{1}{\sqrt{2\pi c(t)}}\exp\left(-\frac{(g-b(t))^2}{2c(t)}\right),
 	\end{equation}where the mean $b(t)$ and the variance $c(t)$ satisfy the following autonomous ODE
 	\begin{equation}\label{eq:system-1}
 	\begin{aligned}
 	\frac{db(t)}{dt}&=g_0+{g_1}N(b(t),c(t))-b(t),\\
 	\frac{dc(t)}{dt}&=2a_0+2{a_1}N(b(t),c(t))-2c(t).
 	\end{aligned}
 	\end{equation}
 	Here $N(b,c)$ is a function of $(b,c)$ which denotes the firing rate of such a Gaussian type solution: 
 	\begin{equation}\label{def-N-bc}
 	N(b,c):=\frac{1}{V_F}\int_{0}^{+\infty}g\frac{1}{\sqrt{2\pi c}}\exp\left(-\frac{(g-b)^2}{2c}\right)dg\geq 0,\quad b\in\mathbb{R},\,c>0.
 	\end{equation}
 
 This combined with Theorem \ref{thm:decay-v} leads to the following final long time asymptotics result for solutions to \eqref{eq:nonlinear-toy}-\eqref{eq:nonliner-toy-fire}.

 \begin{theorem}[Carrillo, Dou, Zhou \cite{CDZ}]\label{thm:full-convergence}
For initial data $p^0(v,g)$ with $\int_{\mathbb{R}}\int_0^{V_F}|g|p^0(v,g)dvdg<\infty$, the classical solution of the simplified voltage-conductance  model \eqref{eq:nonlinear-toy}-\eqref{eq:nonliner-toy-fire} satisfy:
\begin{enumerate}
\item When $g_1/V_F\geq 1$, the firing rate $N(t)$ goes to infinity as time evolves, i.e.,
\begin{equation}
    N(t)\rightarrow+\infty,\quad \text{as }t\rightarrow+\infty.
\end{equation}
    \item 	When $0<g_1/V_F<1$, as time evolves the density function $p(t,v,g)$ converges to the unique steady state
	\begin{equation}\label{steady-vgPDE}
	p^*(v,g):=\frac{1}{V_F}\frac{1}{\sqrt{2\pi c^*}}\exp\left(-\frac{(g-b^*)^2}{2c^*}\right),
	\end{equation} in $L^1((0,V_F)\times\mathbb{R})$. Here $(b^*,c^*)$ is the unique steady state of the ODE system \eqref{eq:system-1}.
\item When $g_1/V_F\rightarrow 1^{-}$, we have $b^*\rightarrow+\infty$, {which means in $g$ direction the center of the steady state \eqref{steady-vgPDE} moves towards infinity.}
\end{enumerate}
\end{theorem}

By introducing a parameter $\varepsilon>0$, modelling the timescale ratio between the conductance $g$ and the voltage $v$, we can consider the fast conductance limit $\varepsilon\rightarrow0^+$. One can derive a limit model governing the $v$ marginal density, whose long time behavior can also be characterized clearly as in the previous theorem, see \cite{CDZ} for further details. Depending on $g_1$ and the $L^{\infty}$ norm of initial data, the solution either blows up in finite time or globally exists in the form of periodic solutions. The authors only find  periodic solutions in the fast conductance limit showing that the other neglected difficulties might be relevant for the appearance of periodic solutions in the full model \eqref{eq:VC}. 

\subsubsection{Wassertein contraction for a Moris-Lecar model}
\label{sec:VC_ML}

In \cite{HMP}, a kinetic Moris-Lecar \cite{ML} model is proposed for the modelling of a population of neurons described by voltage and conductance. Instead of letting the neurons fire at a threshold $V_F$ and resetting them, the voltage domain is extended from $[V_L,V_F]$ to $[V_L,V_E]$, where $V_L,V_E$ are again the leak and the excitatory reversal potential, satisfying here $0<V_L< V_E$. The model reads
\begin{equation}\label{eq:VC_var_ML}
    \left\{
    \begin{array}{l}
         \dfrac{\partial \rho}{\partial t} + \dfrac{\partial \rho}{\partial v}\Big[(g_L(V_L-v)+g(V_E-v))\rho\Big] + \gamma\dfrac{\partial \rho}{\partial g}\Big[(G(v)-g)\rho\Big] - a \dfrac{\partial^2 \rho }{\partial g^2} = 0,  \\
          \rho(V_L,g,t)=\rho(V_E,g,t) = 0,\\
          (G(v)-g)\rho(v,0,t) - a\dfrac{\partial \rho}{\partial g}(v,0,t) = 0,\\
          \displaystyle \rho(v,g,0) = \rho^0(v,g)\geqslant 0, \qquad \int_{V_L}^{V_E} \int_{0}^{+\infty} \rho^0(v,g)\diff v\diff g = 1,
    \end{array}
    \right.
\end{equation}
whith $g_L>0$ the leak conductance, $\gamma>0$ a timescale parameter for the conductance, $a>0$ the strength of noise on the conductance and $G\in C^2([V_L,V_E])>0$ a bounded positive function. Contrarily to \eqref{eq:VC}, this PDE system is linear, and, most important, it does not have the issue of a complex nonlocal boundary condition, which was a technical issue even in linear \eqref{eq:VC}.

The method \cite{HMP} proposes for the study of \eqref{eq:VC_var_ML} is to focus on the underlying stochastic process
\begin{equation}\label{eq:VC_var_ML_sto}
    \left\{
    \begin{array}{l}
         \diff V_t = \big[g_L(V_L-V_t) + g_E(V_E-V_t)\big]\diff t,\\
         \diff G_t = \gamma (G(V_t)-G_t)\diff t + \sqrt{2a} \diff B_t + \diff L_t,
    \end{array}
    \right.
\end{equation}
with $(B_t)_{t\geqslant 0}$ a standard Brownian motion and $(L_t)_{t\geqslant 0}$ the local time preventing the conductance $G_t$ from becoming negative (see \cite[Chapter 6]{RY} or the discussion on a similar issue in Section \ref{sec:GC}).

Denote $(P_t)_{t\geqslant 0}$ the semigroup associated to \eqref{eq:VC_var_ML_sto}, which acts left on measures and right on test functions as seen in other duality arguments in this document. The article \cite{HMP} uses the following coupling method: choose two $C^2$ Lipschitz non-negative functions $\alpha,\beta$ such that $\alpha^2+\beta^2=1$, and consider the coupled system
\begin{equation}
    \left\{
    \begin{array}{l}
         \diff v_t = \big[g_L(V_L-v_t) + g_E(V_E-v_t)\big]\diff t,\\
         \diff g_t = \gamma (G(v_t)-g_t)\diff t + \sqrt{2a} \big[\alpha_t \diff B_t + \beta_t\diff B_t'\big] + \diff L_t,\\
         \diff v_t' = \big[g_L(V_L-v_t') + g_E(V_E-v_t')\big]\diff t,\\
         \diff g_t' = \gamma (G(v_t')-g_t')\diff t + \sqrt{2a} \big[\alpha_t \diff B_t - \beta_t\diff B_t'\big] + \diff L_t',
    \end{array}
    \right.
\end{equation}
where $(B_t)_{t \geqslant 0},(B'_t)_{t \geqslant 0}$ are standard Brownian motions,
\begin{equation*}
    \alpha_t =\alpha(|g_t-g_t'|^2), \qquad  \beta_t = \beta(|g_t-g_t'|^2),
\end{equation*}
and $L,L'$ are the local times at 0 of $g_t,g_t'$. Since by Levy's characterisation the processes
\[  W_t = \int_0^t \big[\alpha_t \diff B_t + \beta_t\diff B_t'\big] ,\qquad W_t' = \int_0^t \big[\alpha_t \diff B_t - \beta_t\diff B_t'\big]   \]
are Brownian motions, $(v_t,g_t)$ and $(v_t',g_t')$ are also solutions to \eqref{eq:VC_var_ML_sto}. The regularity involved allows to write Itô formulae for the coupling system.

When $\alpha = 1, \beta= 0$ both equations are subject to the same noise source and the coupling is called synchronous. Using this coupling, \cite{HMP} proves that if  $(V_E-V_L)\norme{G'}_\infty < g_L$,
then $(v_t,g_t)$ and $(v_t',g_t')$ subject to the same Brownian noise are subject to a deterministic contraction: there exist $C\geqslant 1,\lambda > 0$ such that
\[  |(v_t,g_t)-(v_t',g_t')|\leqslant C e^{-\lambda t}|(v_0,g_0)-(v_0',g_0')| \]
This allows to prove both that there exists a unique measure stationary state $\rho_\infty\in\mathcal P([V_L,V_E]\times\R_+)$ to \eqref{eq:VC_var_ML} when $a> 0$ and that when $a=0$ the deterministic system converges exponentially fast towards a unique equilibirum $(V_\infty,G_\infty)$.

By constructing a mirror coupling, \cite{HMP} proves that there is also a stochastic contraction for \eqref{eq:VC_var_ML_sto}. Notwithstanding the regularity assumptions not being satisfied by this choice, the idea of a mirror coupling is to choose $\beta(s) = \mathds 1_{s > 0}$ and the corresponding $\alpha$ in order to have $\alpha^2 + \beta^2=1$. A regularised mirror coupling is crafted in \cite{HMP} and used to prove a very general Wasserstein contraction for \eqref{eq:VC_var_ML_sto} and hence for \eqref{eq:VC_var_ML}. More precisely, denote for any $p\geqslant 1$, $\mathcal P_p([V_L,V_E]\times\R_+)$ the set of probability measures on $[V_L,V_E]\times \R_+$ with finite $p^{th}$ moment and $\mathcal W_p$ the $L^p$-Wasserstein distance on this space. Then,
\begin{theorem}[Herda, Monmarché, Perthame \cite{HMP}]
    There are constants $\lambda>0,C \geqslant 1$ such that for any initial distributions $\nu,\mu\in \mathcal P_1([V_L,V_E]\times\R_+)$,
    \begin{equation*}
        \forall t\geqslant 0, \qquad \mathcal W_1(\nu P_t,\mu P_t) \leqslant C e^{-\lambda t} \mathcal W_1(\nu,\mu).
    \end{equation*}
    There exists a unique invariant measure $\rho_\infty $ and for all $\nu\in P([V_L,V_E]\times\R_+)$, $\nu P_t$ converges weakly to $\rho_\infty$. For all $p > q \geqslant 1$, if $\rho_\infty, \nu \in \mathcal P_p([V_L,V_E]\times\R_+)$, then $\mathcal W_q(\nu P_t,\rho_\infty)$ converges to 0 as $t\to+\infty$.  
\end{theorem}

The method can be extended to a system of finitely many neurons of potentials and conductances $V_t^1,G_t^1,\dots,V_t^N,G_t^N$ interacting through their voltage \textit{via} the quantities $G_i(V_t^1,\dots,V_t^N)$, with $G_i\in C^1([V_L,V_E])>0$. The system of $N$ neurons is then
\begin{equation}\label{eq:VC_var_ML_sto2}
    \left\{
    \begin{array}{l}
         \diff V_t^i = \big[g_L(V_L-V_t^i) + g_E(V_E-V_t^i)\big]\diff t,\\
         \diff G_t^i = \gamma (G_i(V_t^1,\dots,V_t^N)-G_t^i)\diff t + \sqrt{2a} \diff B_t^i + \diff L_t^i.
    \end{array}
    \right.
\end{equation}

Denote $\mathcal W_{1,l_1}$ the Wasserstein-1 distance associated to
\[   \norme{(v,g)+(w,h)}_1 = \sum_{i=1}^N \left( |v^i-w^i| + |g^i-h^i|\right).  \]

\begin{theorem}[Herda, Monmarché, Perthame \cite{HMP}]\label{thm:VC_var_ML_N}
    Assume there exists $K>0$ such that \[\max_{1\leqslant i,j \leqslant N}\norme{\partial_{v_j} G_i}_\infty\leqslant K.\]
    There are constants $\lambda,\eta>0,C \geqslant 1$ such that for any initial distributions $\nu,\mu\in \mathcal P_1([V_L,V_E]\times\R_+)$,
    \begin{equation*}
        \forall t\geqslant 0, \qquad \mathcal W_{1,l_1}(\nu P_t,\mu P_t) \leqslant C e^{-(\lambda-\eta) t} \mathcal W_{1,l_1}(\nu,\mu).
    \end{equation*}
    If $\lambda > \eta$, there exists a unique invariant measure $\rho_\infty^N $ whose moments are all finite and for all $\nu\in P([V_L,V_E]^N\times\R_+^N)$, $\mathcal W_p(\nu P_t,\rho_\infty^N)$ converges to 0 as $t\to+\infty$.  
\end{theorem}

The constants $C,\lambda,\eta$ being constructed explicitly in \cite{HMP}, it is in fact possible to determine when $\lambda > \eta$ and to play on the functions and parameters to ensure it is the case.

If the control $K > 0$ over the norms $\norme{\partial_{v_j} G_i}$ is uniform when $N\to+\infty$, the work done on the linear model and the system of $N$ neurons allow to pass to the mean-field limit and prove propagation of chaos. The article \cite{HMP} considers the choice
\begin{equation*}
    G_i(V^1,\dots,V^N) = H_0(V^i) + \dfrac{1}{N-1} \sum_{j\neq i} H_1(V^i,V^j).
\end{equation*}
with $H_0\in C^1([V_L,V_E],(0,+\infty))$ and $H_1 \in C^1([V_L,V_E]^2,(0,+\infty))$. Then $|\partial_{v_i}G_i|\leqslant \norme{H_0'}_\infty +\norme{\nabla H_1}_\infty$ and for $j\neq i$, $|\partial_{v_j} G_i|\leqslant \frac{1}{N-1} \norme{\nabla H_1}_\infty$.

It is then possible to take the constants as independent of $N$ in Theorem \ref{thm:VC_var_ML_N}, allowing the negative part $\lambda$ of the decay rate to be independent from $N$ as well. Using this strategy, \cite{HMP} proves rigorous propagation of chaos: if we consider the initial conditions distributed according to $(\rho^0)^{\otimes N}$ and denoting for $k\leqslant N$ by $\rho^{(k,N)}_t$ the law of $(V_t^1,G_t^1,\dots,V_t^k,G_t^k)$, then
\[  \mathcal W_2\left(\rho_t^{(k,N)},\rho(t)^{\otimes k}\right)\leqslant \sqrt{\dfrac{k}{N}} C_0 e^{c t},  \]
for some constants $C_0,c > 0$ and where $\rho$ is the solution to the mean-field system
\begin{equation}\label{eq:VC_var_ML_meanfield}
    \left\{
    \begin{array}{l}
         \dfrac{\partial \rho}{\partial t} + \dfrac{\partial \rho}{\partial v}\Big[J_v(v,g)\rho\Big] + \gamma\dfrac{\partial \rho}{\partial g}\Big[(H_0(v)+ H_1\star\rho(v)-g)\rho\Big] - a \dfrac{\partial^2 \rho }{\partial g^2} = 0,  \\[0.3cm]
          \rho(V_L,g,t)=\rho(V_E,g,t) = 0,\\
          (H_0(v)+ H_1\star\rho(v)-g)\rho(v,0,t) - a\dfrac{\partial \rho}{\partial g}(v,0,t) = 0,\\
          \displaystyle \rho(v,g,0) = \rho^0(v,g)\geqslant 0, \qquad \int_{V_L}^{V_E} \int_{0}^{+\infty} \rho^0(v,g)\diff v\diff g = 1,
    \end{array}
    \right.
\end{equation}
where $J_v(v,g)=g_L(V_L-v)+g(V_E-v)$ and
\[  H_1\star \rho(v) = \int_{V_L}^{V_E}\int_{0}^{+\infty} H_1(v,w) \rho(w,g)  \diff g \diff w.  \]

As proved in \cite{HMP}, the coupling method still applies to this now nonlinear PDE problem: there exist $C,\theta, \lambda>0$ such that if $\rho_1^0,\rho_2^0\in\mathcal P([V_L,V_E]\times\R_+)$ yield two solutions $\rho_1,\rho_2$ of \eqref{eq:VC_var_ML_meanfield}, then
\begin{equation*}
    \mathcal W_1\left(\rho_1(t),\rho_2(t)\right) \leqslant C e^{-(\lambda-\theta) t} \mathcal W_1\left(\rho_1^0,\rho_2^0\right). 
\end{equation*}
Again, $\lambda$ and $\theta$ are constructed explicitly and can be compared, providing when $\lambda > \theta$ the existence and uniqueness of a stationary state and a weakly-nonlinear stability result for \eqref{eq:VC_var_ML_meanfield}.

\section{Kinetic FitzHugh-Nagumo models}
\label{sec:FN}

One of the most comprehensive descriptions of individual neurons and the ionic exchanges leading to action potentials is the celebrated Hodgkin-Huxley model \cite{HH}. Due to its complexity and high computational cost, the FitzHugh-Nagumo model has been proposed \cite{F,NAY} as a simplification which preserves many important features while remaining mathematically simple and numerically tractable at the network scale.

Consider a neuron of voltage $v(t)$ and adaptation variable $w(t)$, whose voltage evolution responds to an external current $I_{\mathrm{ext}}(t)$, an intrinsic voltage dynamics $N(v)$ and is dampened by the adaptation variable $w$. The FitzHugh Nagumo system writes
\begin{equation}\label{eq:FN_basic}
    \left\{\begin{array}{rcl}
    \dfrac{\diff v}{\diff t} &=\ &N(v) - w + I_{\mathrm{ext}}(t),\\[0.3cm]
    \dfrac{\diff w}{\diff t} &=\ &bv-aw+c,
\end{array}\right.
\end{equation}
whith $a,b $ positive parameters and $c\in\R$ that we will often set to 0. A widespread choice \cite{NAY,FS2,BFFT, T2} for the nonlinearity $N$  is a cubic function like $N(v) = v - v^3$.

We first explain how this model for one neuron can be used to derive a mean-field equation for a whole network. Then we discuss existence of solutions, asymptotic convergence to stationary states and the arising of self-sustained oscillations. Last, we cover some results about spatial extensions of the model and the reduction towards a macroscopic limit.

\subsection{From individual FitzHugh-Nagumo neurons to a mesoscopic system}

Consider a network of $K$ interconnected neurons each described by the system \eqref{eq:FN_basic} and subject to a standard Brownian noise $(B_t^i)_{t\geqslant 0}$:
\begin{equation}\label{eq:FN_network_sto}
    \left\{\begin{array}{rcl}
     \diff V_t^i &=\ &\displaystyle \Big(N(V_t^i) - W_t^i + \sum_{j=1}^{N} J_{ij}(V_t^i - V_t^j) +  I_{\mathrm{ext}}(t) \Big)\diff t + \sigma \diff B_t^i ,\\[0.3cm]
     \diff W_t^i &=\ &\displaystyle ( b V_t^i - a W_t^i )\diff t,
\end{array}\right.
\end{equation}
where we set here the cubic nonlinearity $N(v) = v(v-\lambda)(1-v)$, $\lambda > 0$ and assume for clarity $\sigma^2 = 2$, $I(t) \equiv I_0 \in\R$. The coefficients $J_{ij}$ represent the connections between neurons $i$ and $j$.

Propagation of chaos and description of the mean-field limit for system \eqref{eq:FN_network_sto} have been rigorously obtained in \cite{BFFT} and slightly improved in \cite[Appendix A]{MQT}; see also \cite{BFT,CL2,LS4,C4} for other approaches on mean field limits for systems of FitzHugh-Nagumo neurons. In the limit $K\to+\infty$, the density $\rho(v,w)$ of the limiting process is solution to the nonlinear kinetic equation
\begin{equation}\label{eq:FN} \tag{FhN}
    \left\{\begin{array}{l}
     \displaystyle \dfrac{\partial \rho}{\partial t} + \dfrac{\partial }{\partial w}\Big( (bv-aw)\rho\Big) + \dfrac{\partial }{\partial v}\Big( B_\varepsilon\big(v,w,\mathcal J(t)\big) \rho \Big) - \dfrac{\partial^2 \rho}{\partial v^2} = 0,\\[0.3cm]
     \displaystyle \mathcal J(t) = \int_{-\infty}^{+\infty}\int_{-\infty}^{+\infty} v \rho(v,w,t)\diff v\diff w,\\
     B_\varepsilon (v,w, \mathcal J) = v(v-\lambda)(1-v) + w - \varepsilon  (v-\mathcal J) + I_0,\\
     \displaystyle \rho(v,w,0) = \rho^0(v,w)\geqslant 0,\qquad \int_{-\infty}^{+\infty} \int_{-\infty}^{+\infty}\rho^0(v,w)\diff v\diff w = 1.
\end{array}\right.
\end{equation}
The parameter $\varepsilon$ represents the average connectivity and depends upon scaling assumptions on the $J_{ij}$ in the mean-field limit. When it is small enough, we consider the model to be weakly-nonlinear.

Consider the weights $\omega(v,w) = 1+\frac{v^2}{2}+\frac{w^2}{2}$ and the spaces $L^1(\omega)$ and $L^1(\omega^\frac12)$ of functions which are integrable when multiplied by $\omega$ or $\omega^\frac12$. We take from \cite{MQT} the following notion of weak-solution.

\begin{definition}
    We say that $\rho\in C(\R_+,L^1(\omega^\frac12))$ is a global-in-time weak solution to \eqref{eq:FN} if
    \begin{itemize}
        \item for almost every $t\geqslant 0$, $\rho(\cdot,\cdot,t)\geqslant 0$ and
        \[  \int_{-\infty}^{+\infty}\int_{-\infty}^{+\infty} \rho(v,w,t) \diff v\diff w = \int_{-\infty}^{+\infty}\int_{-\infty}^{+\infty}\rho^0(v,w)\diff v\diff w = 1 ; \]
        \item for all $\phi\in C^1(\R_+,C_c^\infty(\R^2))$ and all $t\geqslant 0$,
        \begin{multline*}
            \int_{-\infty}^{+\infty}\int_{-\infty}^{+\infty}\phi(v,w,t) \rho(v,w,t)\diff v\diff w = \int_{-\infty}^{+\infty}\int_{-\infty}^{+\infty}\phi(v,w,t) \rho^0(v,w)\diff v\diff w\\
            +\int_0^t \int_{-\infty}^{+\infty}\int_{-\infty}^{+\infty}\left[\dfrac{\partial \phi}{\partial t}+(bv-aw)\dfrac{\partial \phi}{\partial w}+B_\varepsilon\big(v,w,\mathcal J(s)\big)\dfrac{\partial \phi}{\partial v}+\dfrac{\partial^2 \phi}{\partial v^2}\right] \rho(v,w,s)\diff v\diff w\diff s.
        \end{multline*}
    \end{itemize}
\end{definition}

\subsection{Existence of weak solutions}

Let us introduce the space
\[ L^1\log L^1 = \left\{ \rho\in L^1(\R^2)\ | \ \rho\geqslant 0,\ \mathscr H(\rho) < +\infty  \right\},\qquad \mathscr H(\rho) = \int_{\R^2}\rho(v,w)\log \rho(v,w)\diff v\diff w, \]
and for some $\kappa > 0$, the weight function $m(v,w) = e^{\kappa(\omega-1)} = e^{\kappa(\frac{v^2}{2}+\frac{w^2}{2})}$. Define also the space with additional $v$-derivatives
\[  H^2_v = \left\{ f\in H^1(\R^2)\ | \ \dfrac{\partial^2 f}{\partial v^2}\in L^2(\R^2) \right\}, \qquad \norme{f}_{H^2_v} = \norme{f}_{H^1} + \norme{\dfrac{\partial^2 f}{\partial v^2} }_{L^2}. \]
and its natural weighted version $H^{2}_v(m)$ where all norms are taken with the weigth $m$.

For fixed external input $\mathcal J$, \eqref{eq:FN} is the evolution equation $\partial_t \rho = \mathcal Q_{\mathcal J,\varepsilon}\,  \rho$ where the operator $\mathcal Q_{\mathcal J,\varepsilon}$ is defined by
\[  \mathcal Q_{\mathcal J,\varepsilon}\, \rho = \dfrac{\partial }{\partial w}[(aw-bv)\rho] + \dfrac{\partial }{\partial v}[ B_\varepsilon(v,w,\mathcal J)\rho ] + \dfrac{\partial^2 \rho}{\partial v^2 }.  \]

First, \cite{MQT} proves by direct calculations an estimate of the form
\[  \dfrac{\diff }{\diff t}\int_{\R^2} \rho(v,w,t)\omega(v,w) \diff v\diff w \leqslant K_1 - K_2  \int_{\R^2} \rho(v,w,t)\omega(v,w) \diff v\diff w  \]
which by application of Grönwall lemma and use of the control $|\mathcal J(t)|\leqslant 2\norme{\rho(t)}_{L^1(\omega)}$ ensures that any solution of \eqref{eq:FN} satisfies, for appropriate initial conditions, the \textit{a priori} bounds
\[  \norme{\rho(\cdot,\cdot,t)}_{L^1(\omega)} \leqslant \max\left(C,\norme{\rho^0}_{L^1(\omega)}\right),\qquad \sup_{t\geqslant 0} |\mathcal J(t)| \leqslant C,   \]
for some constant $C$ depending on $\rho^0$ and the parameters. Then, for fixed $\mathcal J$, they derive a bound the operator $\mathcal Q_{\mathcal J,\varepsilon}$ and deduce a similar bound for $\rho$ in the space $L^1(m)$. This first step allows to move to the stronger spaces $H^1(m)$ and $H^2_v(m)$. The core idea is to take advantage of the hypodissipativity of the equation by use of the twisted space approach and the Nash-Villani technique \cite{V}. More precisely, choosing two positive constants $\kappa_1 < \kappa_2$ and associated $m_1 = e^{\kappa_1 (\omega-1)}, m_2=e^{\kappa_2(\omega-1)}$, they prove (not without long and tedious technicalities) that for all $\mathcal J$ fixed, there exist constants $K_1,K_2$ and $\delta\in(0,1)$ such that
\begin{equation*}
    \norme{\rho}_{\mathcal H^1}^2 = \norme{\rho}_{L^2(m_2)}^2 + \norme{\partial_w\rho}_{L^2(m_2)}^2 + \pscal{\partial_w\rho}{\partial_v\rho}_{L^2(m_1)} +\norme{\partial_v\rho}_{L^2(m_2)}^2
\end{equation*}
defines a Hilbert norm equivalent to the usual norm on $H^1(m_2)$ and satisfying
\begin{equation*}
\forall \rho\in H^1(m_2),\qquad \pscal{\mathcal Q_{\mathcal J,\varepsilon}\,\rho}{\rho}_{\mathcal H^1} \leqslant K_1 \norme{\rho}_{L^2(\R^2)}^2  - K_2 \norme{\rho}_{\mathcal H^1}^2 .
\end{equation*}
From there, choosing $m = m_2$, Nash's inequality and the equivalence between $\norme{\cdot}_{\mathcal H^1}$ and $\norme{\cdot}_{H^1(m)}$ allow to derive the differential inequality
\[  \dfrac{\diff}{\diff t} \norme{\rho(\cdot,\cdot,t)}^2_{\mathcal H^1} \leqslant \pscal{\mathcal Q_{\mathcal J(t),\varepsilon}\,\rho(\cdot,\cdot,t)}{\rho(\cdot,\cdot,t)} \leqslant K_3 - K_4 \norme{\rho(\cdot,\cdot,t)}^2_{\mathcal H^1},  \]
which proves \textit{a priori} control of $\norme{\rho(\cdot,\cdot,t)}_{\mathcal H^1}$ for any solution of \eqref{eq:FN}.

Then, \cite{MQT} proves, using the partial Fischer information
\[ I_v(\rho) = \int_{-\infty}^{+\infty} \int_{-\infty}^{+\infty} \dfrac{|\partial_v \rho (v,w)|^2}{f(v,w)} \diff v \diff w ,  \]
that solutions to \eqref{eq:FN} with finite initial entropy remain in the space $L^1\log L^1$ of finite entropy functions.

\begin{lemma}
    For all $\rho^0\in L^1(\omega)\cap L^1\log L^1\cap \mathcal P(\R^2)$, for any solution to \eqref{eq:FN} and any time $T>0$, there exists $C(T)>0$ such that
    \[  \sup_{t\in[0,T]}\mathscr H\big(\rho(\cdot,\cdot,t)\big) + \int_0^t I_v\big(\rho(\cdot,\cdot,s)\big)\diff s \leqslant C(T). \]
\end{lemma}
This entropy bound yields the well-posedness bound
\[  \sup_{t\in[0,T]} \norme{\rho_1(\cdot,\cdot,t)-\rho_2(\cdot,\cdot,t)}_{L^1(\omega^\frac12)}\leqslant C_2(T) \norme{\rho_1^0-\rho_2^0}_{L^1(\omega^\frac12)},  \]
which ensures uniqueness of any prospective solution. Although standard results don't apply and the equation has a degenerate structure and unbounded drifts, \cite{MQT} finds its way to the existence of solutions using the a priori bounds above along with a cutoff and a tailored coercive bilinear form. An Aubin-Lions compactness argument yields a solution for any given external input $t\mapsto \mathcal J(t)\in L^\infty(\R_+)$, and then the full solution to \eqref{eq:FN} can be recovered by a fixed point argument on the functional that maps $\mathcal J\in L^\infty([0,T])$ to the internal input $\mathcal J_{\rho_{\mathcal J}}$ associated to the solution $\rho_{\mathcal J}$ obtained from the linear equation with external input $\mathcal J$.

\begin{theorem}[Mischler, Qui\~ninao, Touboul, \cite{MQT}]
    For any $\rho^0\in L^1(\omega)\cap L^1\log L^1 \cap \mathcal P(\R^2)$, there exists a unique global-in-time weak solution $\rho$ to \eqref{eq:FN}. It depends continuously on $\rho^0$ in $L^1(\omega^\frac12)$ and satisfies for some constant $C_0 > 0$,
    \[  \norme{\rho(\cdot,\cdot,t)}_{L^1(\omega)} \leqslant \max\left(C_0,\norme{\rho^0}_{L^1(\omega)}\right).   \]
    If the initial condition enjoys the higher decay rate $\rho^0\in L^1(m)$, then for some constant $C_1 > 0$ depending on $\rho^0$,
    \[  \norme{\rho(\cdot,\cdot,t)}_{L^1(m)} \leqslant C_1.   \]
    There exist norms $\norme{\cdot}_{\mathcal H^1}$ and $\norme{\cdot}_{\mathcal H^2_v}$ equivalent respectively to $\norme{\cdot}_{H^1(m)}$ and $\norme{\cdot}_{\mathcal H^2_v(m)}$ such that if the initial condition enjoys the higher regularity $\rho^0\in H^1(m)$ or $\rho^0\in H^2_v(m)$, then for some constants $C_2,C_3 >0$  depending on $\rho^0$,
    \[  \norme{\rho(\cdot,\cdot,t)}_{H^1(\omega)} \leqslant C_2 ,\qquad \mathrm{or}\qquad  \norme{\rho(\cdot,\cdot,t)}_{H^2_v(\omega)} \leqslant C_3. \]
\end{theorem}

\subsection{Long time behaviour of the mesoscopic system}

Consider the following variant of the Brouwer fixed point theorem:
\begin{theorem}[ \cite{EMR}, \cite{GPV} ]
    Let $\mathcal Z$ a convex and compact subset of a Banach space $\mathcal X$ and $(S_t)_{t\geqslant 0}$ a continuous semigroup on $\mathcal Z$. Assume $\mathcal Z$ is invariant under the action of $(S_t)_{t\geqslant 0}$. Then, there exists $z_0\in\mathcal Z$ such that for any $t\geqslant 0$, $S_t z_0 = z_0$.
\end{theorem}

The article \cite{MQT} uses this result with the space
\[  \mathcal X = H^2_v(m)\cap L^1\log L^1 \cap \mathcal P(\R^2),  \]
an appropriate $\mathcal Z$ and the semigroup generating the solutions of \eqref{eq:FN} to prove existence of at least one stationary solution to the nonlinear system.

\begin{theorem}[Mischler, Qui\~ninao, Touboul, \cite{MQT}]
    For all $\varepsilon\geqslant 0$, there exists at least one stationary state $\rho_\infty\in H^2_v(m)$ to \eqref{eq:FN}. When $\varepsilon\to 0$, it converges strongly in $L^2(m)$ towards the linear stationary state obtained for $\varepsilon=0$. 
\end{theorem}

\paragraph{The linearised problem}

Let $\rho_\infty$ and $\rho$ a stationary state and a solution of \eqref{eq:FN} and consider the variation $ \rho = \rho_\infty + h$. Then the linearised operator for $h$ is
\[  \mathscr L_\varepsilon h = \mathcal Q_{\mathcal J_\infty,\varepsilon}\, h + \varepsilon \mathcal J_h(t) \dfrac{\partial \rho_\infty}{\partial v}, \]
with $\mathcal J_h(t) =\int_{\R^2} vh(v,w,t)\diff v\diff w$. Similar to the method described\footnote{Note that \cite{MW} was published two years later than \cite{MQT} and precedes it only in our exposition.} in Subsection \ref{sec:TE_MW}, the idea is to split the operator in two:
\[   \mathcal A = N \chi_R ,\qquad \mathcal B_\varepsilon = \mathcal L_\varepsilon - \mathcal A,  \]
with $N>0$ a constant and $\chi$ a smooth cutoff function equal to 1 in $B(0,R)$ and supported in $B(0,2R)$. The goal is to take advantage of the regularity of $\mathcal A$ and the dissipativity of $\mathcal B_\varepsilon$ to conclude that for small $\varepsilon$ the spectrum of $\mathcal L_\varepsilon$ is close to the spectrum of $\mathcal B_\varepsilon$; \cite{MQT} emphasises that this procedures corresponds to Weyl's theorem or to a generalised Krein-Rutman \cite{MS2} (see also the survey \cite{FGM}).

Recall the definition of hypodissipativity above (Definition \ref{def:T-E_hyp}).

\begin{lemma}[Mischler, Qui\~ninao, Touboul, \cite{MQT}]
    For $N,R>0$ large enough, $\mathcal B_\varepsilon$ is (-1)-hypodissipative in $H^2_v(m)$.
\end{lemma}
This allows to perform a detailed spectral analysis in the disconnected case and to extend it to $\varepsilon\simeq 0$.

\paragraph{Weakly-nonlinear case}

When $\varepsilon$ is small enough, \cite{MQT} proves that the stationary state is unique. Then, building upon the study of the linear spectrum, they proceed to characterise $\Sigma(\mathcal L_\varepsilon)$.

Coming back to the perturbation $h=\rho-\rho_\infty$, they write the following Duhamel formula:
\[  h  = S_{\mathscr L_\varepsilon}(t) h^0 + \int_0^t S_{\mathscr L_\varepsilon}(t-s)\Big( \varepsilon\mathcal J_h(t) \dfrac{\partial h}{\partial v} \Big)\diff s.  \]
Combining it with the spectral information, it is possible to prove the decay
\[  \norme{h}_{L^2(m)}\leqslant C e^{-\alpha t} \norme{h^0}_{L^2(m)},  \]
for some $C\geqslant 1,\alpha > 0$.

\begin{theorem}[Mischler, Qui\~ninao, Touboul, \cite{MQT}]
    There exists $\varepsilon_0 > 0$ such that for all $\varepsilon\in(0,\varepsilon_0)$, there exists a unique stationary state $\rho_\infty$ to \eqref{eq:FN} and $\eta(\varepsilon)>0$ such that for all $\rho^0\in H^1(m)\cap \mathcal P(\R^2)$, if $\norme{\rho^0-\rho_\infty}_{L^2(m)}\leqslant \eta(\varepsilon)$, there exist $C\geqslant 1$ and $\nu>0$ such that the solution to \eqref{eq:FN} satisfies
    \begin{equation*}
        \forall t\geqslant 0,\qquad \norme{\rho(\cdot,\cdot,t)-\rho_\infty}_{L^2(m)}\leqslant C e^{-\nu t}.
    \end{equation*}
    Moreover, $\lim_{\varepsilon\to 0}\eta(\varepsilon) = + \infty$.
\end{theorem}

Numerical solutions on the particle system $(N=2000)$ in \cite{MQT} confirm that for small $\varepsilon$ the unique stationary state is stable. For large $\varepsilon$, the system can have two stationary states and there is convergence towards one of them. In the intermediate connectivity regime, stable self-sustained oscilations are observed.

\subsection{Periodic solutions arising through noise and interaction}

As we just mentioned, numerical simulations in \cite{MQT} displayed convergence towards a periodic solution in the intermediate conectivity regime. This behaviour was rigorously investigated in \cite{LP2,LP3} using probabilistic methods and noise-driven bifurcations. They study, very similar to the mean-field system \eqref{eq:FN}, the system of stochastic differential equations
\begin{equation}\label{eq:FN_MF_LP}
    \left\{
    \begin{array}{l}
         \diff V_t = \delta\left( N(V_t) - W_t - \dfrac{K}{\delta}(V_t-\mathbb E[V_t])  \right) \diff t + \sqrt2 \sigma \diff B_t, \\
         \diff W_t = \dfrac{\delta}{c}\big(V_t + a - bW_t\big)\diff t, 
    \end{array}
    \right.
\end{equation}
with $b,c,\delta,K$ positive parameters, $a\in\R$ and $N(v)=v-\frac{v^3}{3}$. The main idea is to look at the small $\delta$ range, where the term $V_t-\mathbb E[V_t]$ is close to 0 and thus the density in $v$ is close to a Gaussian. It allows \cite{LP2} to perform a slow-fast perturbative analysis. In a parameter range where the system without noise ($\sigma=0$) nor interactions ($K=0$) would converge to equilibrium, the noise-alone could send the individual neuron dynamics (no interaction) towards a periodic Carnard-like dynamics. Strong interaction ($K\neq0$ and $\delta$ small) then allow these individual oscillations to happen collectively.

The slow-fast dynamics approach proceeds by replacing, for small enough $\delta$, the voltage variable by a random variable of Gaussian distribution $\mathcal N\left(\mathbb E[V_t],\frac{\sigma^2}{K}\right)$. The dynamics of $W_t$ being linear in $V_t$, it will also be Gaussian at first order. The couple $V_t,W_t$ will then behave like $\mathcal N(m_t,\Sigma_\delta)$, and an evolution equation can be written for $m_t$ in the form
\[  \dfrac{\diff }{\diff t}m_t = \delta F_{\frac{\sigma^2}{K}}(m_t),\qquad F_u(v,w) = \left(
\begin{matrix}
    (1-u)v - \frac{v^3}{3} - w\\\frac1c (v+a-bw)
\end{matrix}
\right), \]
which is also of FitzHug-Nagumo type. The covariance matrix $\Gamma_\delta$ can be characterized as
\[  \Gamma_\delta = \dfrac{\sigma^2}{K}{}\left(\begin{matrix}
   1 &\frac{\delta}{Kc+b\delta}\\
   \frac{\delta}{Kc+b\delta} & \frac{\delta}{b(Kc+b\delta)}
\end{matrix}\right).  \]
Making this analysis rigorous, \cite{LP2} proves that the nonlinear mean-field system \eqref{eq:FN_MF_LP} behaves close enough to the system $(v,w)=F_u(v,w)$, which, for a suitable choice of parameters, has a stable fixed point for $u=0$ and a stable periodic solution for some range of the parameter $u$. Because $u=\frac{\sigma^2}{K}$ in the equation for $m_t$, the mean-field oscillations are induced by both noise and interaction.

\begin{theorem}[Luçon, Poquet \cite{LP2}]
    There exist a choice for $a,b,c,K,\sigma$ such that for all $\delta$, the dynamical system $\frac 1\delta (v'(t),w'(t)) = F_{\frac{\sigma^2}{K}}(v(t),w(t))$ admits a stable limit cycle. There exists $\delta^*$ such that for all $\delta\in(0,\delta^*)$, there exists a periodic solution $\rho^{\mathrm{per}}$ to \eqref{eq:FN_MF_LP}. This solution is associated to an invariant limit cycle $C_\delta$ and there exist $C_1,C_2$ independent of $\delta$, $C(\delta),\lambda(\delta)$ such that if $\rho^0 \in \mathcal P_2(\R^2)$, $\int_{\R^2} \norme{(v,w)}^6 \diff\rho^0\leqslant C_1$ and
    $  W_2( \rho^0, C_\delta )\leqslant C_2 \delta$, then the solution $\rho$ to \eqref{eq:FN_MF_LP} satisfies
    \[  W_2(\rho_t, C_\delta)\leqslant C(\delta)e^{-\lambda(\delta)t} W_2(\rho^0,C_\delta). \]
\end{theorem}
They obtain in the process that $C(\delta)\to+\infty$ and $\lambda(\delta)\to\min( \lambda(1),\frac{b}{c} )$ when $\delta\to0$.

\subsection{A spatially extended kinetic FitzHugh-Nagumo model and its macroscopic limit}

From the same base used in the model \eqref{eq:FN}, it is possible to provide a mesoscopic description of a brain area while also taking into account the spatial structure. As we shall see, this allows to rigorously construct macroscopic equations.

The article \cite{BF} considers a system of $n$ interacting FitzHug-Nagumo neurons $((V_t^i,W_t^i,x_i))_{1\leqslant i\leqslant n}$ subject to independent Brownian noises at the voltage level,
\begin{equation}\label{eq:FN_sto_se}
    \left\{\begin{array}{rcl}
     \diff V_t^i &=\ &\displaystyle \Big(N(V_t^i) - W_t^i + \frac{1}{n}\sum_{j=1}^{n} \Phi(x_i,x_j)(V_t^i - V_t^j)  \Big)\diff t + \sigma \diff B_t^i ,\\[0.3cm]
     \diff W_t^i &=\ &\displaystyle ( a V_t^i - b W_t^i+c )\diff t,
\end{array}\right.
\end{equation}
with $c\in\R$, $a,b>0$ and again for simplicity $\sigma^2=2$. Note that the case $\sigma=0$ has also been thoroughly investigated \cite{CFF,C4,C5}. Each neuron is described at time $t$ by a voltage $V_t^i$, an adaptation variable $W_t^i$ and a position $x_i\in  K$ lying in a compact set $K$. The coupling function $\Phi(x,y)$ describes the strength of interactions between points $x$ and $y$. Following \cite{B6,LS4,QT}, the coupling function is composed of strong local connections and weaker long range interactions
\begin{equation}\label{eq:seFN_phi}
      \Phi(x,y)= \dfrac{1}{\varepsilon} \delta_{x=y} + \Psi(x,y).  
\end{equation}
In the formal limit $n\to +\infty$, the system is described by the probability density $\rho(x,v,w,t)$ that we will often write as $\rho(x,u,t)$ with $u=(v,w)$ for clarity.
\begin{equation}\label{eq:FN_se} \tag{se-FhN}
    \left\{\begin{array}{l}
     \displaystyle \dfrac{\partial \rho}{\partial t} + \dfrac{\partial }{\partial w}\Big( (av-bw+c)\rho\Big) + \dfrac{\partial }{\partial v}\Big( \big(N(v)-w-\mathcal J(x,v,t)\big) \rho \Big) - \dfrac{\partial^2 \rho}{\partial v^2} = 0,\\[0.3cm]
     \displaystyle \mathcal J(x,v,t) = \int_K \int_{-\infty}^{+\infty}\int_{-\infty}^{+\infty} \Phi(x,x')(v-v') \rho(x',v',w',t)\diff v'\diff w'\diff x',\\
     \displaystyle \rho(x,v,w,0) = \rho^0(x,v,w)\geqslant 0,\qquad \int_K\int_{-\infty}^{+\infty}\int_{-\infty}^{+\infty} \rho^0(x,v,w)\diff v\diff w\diff x = 1.
\end{array}\right.
\end{equation}
Given the choice \eqref{eq:seFN_phi}, consider the following parameter assumptions.
\begin{hyp}[Regularity of parameters]\label{as:seFhN}
    Assume $N\in C^2(\R)$, $\Psi\in C(K,L^1(K))$ and there exists $p\geqslant 2$ and $r\in(1,+\infty]$ such that
    \[  \limsup_{|v|\to +\infty} \frac{N(v)}{v}=-\infty,\qquad \sup_{|v|\geqslant 1} \frac{|N(v)|}{|v|^p} + \sup_{x'\in K} \int_K \Psi(x,x)\diff x + \sup_{x\in K} \int_K |\Psi(x,x')|^r\diff x' < +\infty. \]
\end{hyp}

The articles \cite{BF,B4,BB} investigate the limit of strong local interactions $\varepsilon\to0$ in \eqref{eq:seFN_phi}, the main idea being that the density should concentrate on a particular time-dependent voltage $\mathcal V(t,x)$ solution to a macroscopic model. In view of this mesoscopic to macroscopic limit, consider the density $\rho^\varepsilon$ obtained as a solution to \eqref{eq:FN_se} from an initial condition $\rho^{0,\varepsilon}$ with the interaction function as in \eqref{eq:seFN_phi}. Then, renormalise $\rho$ by the quantity of neurons $m^\varepsilon(x)$ at location $x$ and define the macroscopic voltage and adaptation variables $\mathcal V^\varepsilon,\mathcal W^\varepsilon$:
\begin{multline*}
    m^\varepsilon(x) = \int_{\R^2}\rho^{0,\varepsilon}(x,u)\diff u,\qquad m^\varepsilon(x)\mu_{t,x}^\varepsilon = \rho(x,\cdot,t),\qquad \bar\mu_{t,x}^\varepsilon = \int_{-\infty}^{+\infty}\mu_{t,x}^\varepsilon(v,\cdot)\diff v\\m^\varepsilon(x)\mathcal V^\varepsilon(t,x)=\int_{\R^2}v\rho^\varepsilon(x,u,t)\diff u,\qquad m^\varepsilon(x)\mathcal W^\varepsilon(t,x)=\int_{\R^2} w\rho^\varepsilon(x,u,t)\diff u. 
\end{multline*}
The system \eqref{eq:FN_se} can then be rewritten as
\begin{equation}\label{eq:seFN_mu}
    \left\{\begin{array}{l}
     \displaystyle \dfrac{\partial \mu^\varepsilon}{\partial t} + \nabla_u \cdot \left[\left(\begin{matrix} &N(v)-w-\mathcal J\\&av-bw+c\end{matrix} \right)\mu^\varepsilon\right] - \dfrac{\partial^2 \mu^\varepsilon}{\partial v^2} = \frac1\varepsilon m^\varepsilon(x)\dfrac{\partial }{\partial v}\Big[ (v-\mathcal V^\varepsilon)\mu^\varepsilon \Big],\\[0.3cm]
     \displaystyle \mathcal J(x,v,t) = \int_K  \Psi(x,x') m^\varepsilon(x)(v-\mathcal V^\varepsilon(t,x))\diff x'.
\end{array}\right.
\end{equation}
Because of the leading order term in $\frac1\varepsilon$, it is natural to expect that in the strong local interaction limit $ (v-\mathcal V^\varepsilon)\mu^\varepsilon\to 0$. To better quantity this concentration, \cite{BF} uses the $W_2$ Wasserstein distance over the space $\mathcal P_2(\R^2)$ of probability measures with finite second moment and they find that the measure $\mu^\varepsilon$ concentrates in the $v$ variable as a Dirac mass centred at the macroscopic variable $\mathcal V(x,t)$ and in the $w$ variable as the continuous marginal $\bar \mu$. More precisely, they obtain an estimate of the form
\[   W_2\left( \,\mu\,,\, \delta_{\mathcal V^\varepsilon } \otimes \bar \mu\,\right) \sim \sqrt\varepsilon.   \]
The dynamics is then asymptotically reduced to a macroscopic model with a voltage $\mathcal V(t,x)$ and an adaptation variable $\mathcal W (t,x) = \int w\diff \bar\mu$ defined through the average of the limit marginal $\bar\mu$. The long-range interaction through the kernel $\Psi$ becomes the macroscopic convolution operator
\[ \mathcal L_m[\mathcal V] = \int_K \Psi(x,x') m(x')(\mathcal V(x)-\mathcal V(x'))\diff x'.  \]
The macroscopic system is then of the form
\begin{equation}\label{eq:seFN_macro} 
    \left\{\begin{array}{l}
     \displaystyle \dfrac{\partial \mathcal V}{\partial t} = N(\mathcal V(t,x)) - \mathcal W(t,x) - \mathcal L_{m}[\mathcal V],\\[0.3cm]
     \displaystyle \dfrac{\partial \bar \mu}{\partial t} + \dfrac{\partial}{\partial w}[ (a\mathcal V(t,x) - b \mathcal W(t,x)  +c )\bar \mu]=0,\\[0.3cm]
     \displaystyle \mathcal W(t,x) = \int_{-\infty}^{+\infty} w \diff \bar\mu_{t,x}(w).
\end{array}\right.
\end{equation}

Consider the following assumptions on initial data.
\begin{hyp}[Initial conditions]\label{as:seFhN_2}
    For all $\varepsilon>0$, assume $x\mapsto \mu_{0,x}^\varepsilon \in C(K,\mathcal P_2(\R^2))$, $m^\varepsilon \in C(K)$ and there exist $p\geqslant 2$, $r\in(1,+\infty]$ of conjugate exponent $r'$ and $C_p,\bar C_p>0$ independent of $\varepsilon$ such that
    \[ m^*\leqslant m^\varepsilon(x)\leqslant \frac1{m^*},\qquad \sup_{x\in K} \int_{\R^2} |u|^{2p}\diff \mu^\varepsilon_{0,x}(u)\leqslant C_p,\qquad \int_{K\times\R^2} |u|^{2p r'}m^\varepsilon(x)\diff \mu^\varepsilon_{0,x}(u) \leqslant \bar C_p, \]
    and that there exist $c_1^\varepsilon,c_2^\varepsilon>0$ such that
    \[  \sup_{x\in K} \int_{\R^2} e^{\frac{|u|^2}{2}}\diff \mu^\varepsilon_{0,x}(u)\leqslant c_1^\varepsilon,\quad \sup_{x\in K} \int_{\R^2}\log(\mu^\varepsilon_{0,x}) \diff \mu^\varepsilon_{0,x}(u)\leqslant c_2^\varepsilon,\quad \sup_{x\in K} \norme{\nabla_u \sqrt{\mu^\varepsilon_{0,x}}}^2_{L^2(\R^2)}\leqslant c_2^\varepsilon. \]
\end{hyp}
Then both the mesoscopic and macroscopic problems are well-posed.
\begin{theorem}[Blaustein, Filbet, \cite{BF}]
    Grant assumptions \ref{as:seFhN} and \ref{as:seFhN_2}. There exists a unique weak solution
    \[ \mu^\varepsilon \in C(\R_+\times K, L^1(\R^2))\cap L^\infty_{\mathrm{loc}}(\R_+\times K, \mathcal P_2(\R^2)) \]
    to \eqref{eq:seFN_mu} and the macroscopic quantities satisfy $\mathcal V^\varepsilon,\mathcal V^\varepsilon\in C(\R_+\times K)$. There exists a unique weak solution $\mathcal V \in C(\R_+\times K)$, $\bar mu \in C(\R_+\times K, L^1(\R))$ to \eqref{eq:seFN_macro}.
\end{theorem}
Using this setting, it is possible to make the asymptotic reduction above rigorous. Rescale the densities $\mu^\varepsilon$ and $\bar \mu^\varepsilon$ as
\begin{equation}
    \mu^\varepsilon_{t,x}(v,w) = \dfrac{1}{\varepsilon^\alpha} \nu^\varepsilon_{t,x}\left(\dfrac{1}{\varepsilon^\alpha}(v-\mathcal V^\varepsilon),w-\mathcal W^\varepsilon\right),\qquad \bar \nu^\varepsilon_{t,x}(w) = \int_{-\infty}^{+\infty} \nu^\varepsilon_{t,x}(v,w)\diff v,
\end{equation}
with $\varepsilon^\alpha$ a concentration profile. Inputing this change of variables in \eqref{eq:seFN_mu} allows to find that the correct profile should be with $\alpha = \frac12$. Then it can be proved that
\[  \nu^\varepsilon = \mathcal M_{m^\varepsilon,0} \otimes \bar \nu^\varepsilon + O(\sqrt{\varepsilon}), \qquad \mathcal M_{m^\varepsilon(x),\mathcal V(x)} (v) = \sqrt{\dfrac{m^\varepsilon(x)}{2\pi}} e^{-\frac{-m^\varepsilon(x) |v-\mathcal V(x)|^2}{2}}.\]
\begin{theorem}[Blaustein, Filbet, \cite{BF}]
    Grant assumptions \ref{as:seFhN} and \ref{as:seFhN_2}. Consider the solutions $\mu^\varepsilon$ to \eqref{eq:seFN_mu} and $(\mathcal V,\bar\mu)$ to \eqref{eq:seFN_macro}. Define the macroscopic and mesoscopic errors
    \[ \mathcal E_{\mathrm{mac}} = \norme{(\mathcal V^0,\mathcal W^0)-(\mathcal V^{0,\varepsilon},\mathcal W^{0,\varepsilon})}_{L^\infty(K)} + \norme{m-m^\varepsilon}_{L^\infty(K)},\qquad \mathcal E_{\mathrm{mes}} = \sup_{x\in K} W_2(\bar\nu_{0,x}^\varepsilon,\bar\nu_{0,x}).   \]
    Then there exists $C,\varepsilon_0 > 0$ such that for all $ 0 < \varepsilon\leqslant \varepsilon_0$ and all $t,x\in \R_+\times K$,
    \begin{equation}\label{bf}    
    W_2(\mu^\varepsilon,\mathcal M_{\frac1\varepsilon m,\mathcal V} \otimes \bar\mu)\leqslant C(\min((\mathcal E_{\mathrm{mac}}+\varepsilon)e^{Ct},1) +\mathcal E_{\mathrm{mes}} e^{-bt} + e^{-\frac{t}{\varepsilon} \rho}).  
    \end{equation} 
\end{theorem}
If $\mathcal E_{\mathrm{mac}} + \mathcal E_{\mathrm{mes}}  = O(\varepsilon)$, the right-hand side in \eqref{bf} can be improved into
$C(\min(\varepsilon e^{Ct},1) + e^{-\frac{t}{\varepsilon} \rho})$ and uniformly on all time compacts,
\[  \lim_{\varepsilon\to 0} \sup_{x\in K} (|\mathcal V^\varepsilon(t,x)-\mathcal V(t,x) + W_2(\bar\mu^\varepsilon_{t,x},\bar\mu_{t,x}) ) = 0,\qquad \lim_{\varepsilon\to 0} \sup_{x\in K}  W_2(\mu^\varepsilon_{t,x},\delta_{\mathcal V(t,x)}\otimes \bar \mu_{t,x})  = 0,  \]
with convergence rates explicitly given in terms of $t,m^*, \varepsilon$ and universal constants not depending on $t,\varepsilon$. Refined estimates on the Gaussian profile during the $\varepsilon\to 0$ blow-up are provided in \cite{B4,BB}, the former proving estimates in a $L^1$ and a weighted $L^2$ framework, the latter using a Hopf-Cole approach to derive $L^\infty$ bounds.

\subsection*{Acknowledgements}
JAC was supported by the Advanced Grant Nonlocal-CPD (Nonlocal PDEs for Complex Particle Dynamics: Phase Transitions, Patterns and Synchronization) of the European Research Council Executive Agency (ERC) under the European Union's Horizon 2020 research and innovation programme (grant agreement No. 883363).
JAC and PR want to thank Westlake University, Hangzhou, China, where a significant part of this review article was written, for its warm hospitality. JAC and PR thank Zhennan Zhou and Alejandro Ramos-Lora for valuable feedback on the manuscript.

{\small
\bibliographystyle{abbrv}

}
\end{document}